\documentclass{article}
\usepackage{amsmath}
\numberwithin{equation}{section}
\usepackage{amssymb}
\usepackage[dvipsnames]{xcolor}
\usepackage{comment}
\usepackage{hyperref}
\usepackage{graphicx}
\usepackage{yhmath}
\usepackage[margin=2.5cm]{geometry}
\usepackage{subcaption}
\newtheorem{theorem}{Theorem}[section]
\newtheorem{proposition}[theorem]{Proposition}
\newtheorem{definition}[theorem]{Definition}
\newtheorem{corollary}[theorem]{Corollary}
\newtheorem{lemma}[theorem]{Lemma}
\DeclareMathOperator{\Op}{Op}

\begin{document}

\title{Feynman Propagators on Exponentially Decaying Perturbations of Stationary Spacetimes with Compact Cauchy Surface}

\author{Selim Amar}

\maketitle
\noindent
Email: selama@stanford.edu \\
Department of Mathematics, Stanford University

\vspace{10pt}

\begin{abstract}

\noindent In this paper, we construct Feynman propagators on exponentially decaying perturbations of stationary spacetime with compact Cauchy surface. The propagators are defined as parametrices of the wave operator on exponentially weighted Sobolev spaces defined by complex contours $\Gamma$. We then apply this theory to the computation of indexes of Dirac operators using a continuous Wick-like rotation to relate it to the Riemannian case. 

\end{abstract}

\tableofcontents

\section{Introduction}

Motivated by applications to Quantum Field Theories on curved spacetimes, the study of Feynman/anti-Feynman propagators on Lorentzian geometries has flourished in the past years. In layman's terms these are inverses of wave operators which propagate future directed frequencies forward/backward in time and past directed frequencies backward/forward in time. More succinctly, it propagates them forward/backward along the lifted geodesic flow on phase space. In practice, ensuring this condition exactly is tricky as one first needs an inverse for the wave operator and then obtain exact control on the support of its kernel. Thus, one typically define Feynman(+)/anti-Feynman(-) parametrices $Q_{\pm F}$ instead, which invert the wave operator up to compact errors. The condition on the propagation of frequencies is then only imposed at high enough frequencies, in which case the propagation is simply given by a choice of direction along the lifted geodesic flow on phase space: propagating along/opposit to the flow gives the Feynman/anti-Feynman condition. This is summarized by the \textit{Hadamard condition} on their wavefront sets \footnote{For an introduction to the mathematical theory of QFT on curved spacetimes, see \cite{bar2009quantum}. Section 4.3.4 discusses the wavefront set of propagators.}
\begin{equation}
    WF(Q_{\pm F}) = \Delta_{T^{\ast}M} \cup \bigcup_{t \geq 0} \varphi_{\pm t}(\Delta_{T^{\ast} M} \cap \pi^{-1} char(p))
\end{equation}
Here $\varphi_{\pm t}$ is the null-geodesic flow on the left factor of $T^{\ast}M \times T^{\ast}M$, $char(p)$ is the characteristic set of $g^{-1}$ i.e. the null covectors and $\pi: T^{\ast}M \times T^{\ast}M$ is the projection on the left factor. Put this in contrast with the retarded and advanced propagators, which propagate \textit{all} frequencies either forward or backward in time. One should think of retard and forward propagators as solving the initial value problems, while Feynman/anti-Feynman propagators solve scattering problems. Because of that, they are mainly relevant in the computation of scattering amplitudes in perturbative QFT. On Minkowski spacetime, they can be explicitly constructed using complex deformation contour techniques, which gives an explicit formula for its kernel
\begin{equation}
    \Box_{\pm F}^{-1} = \lim_{\epsilon \rightarrow 0} \frac{1}{(2\pi)^{n+1}}\int_{\mathbb{R}^{n+1}} \frac{e^{-i(\sigma(t - t') - \xi(x - y))}}{\sigma^2 - p^2 \pm i\epsilon} d\sigma d\xi \quad (\sigma, \xi) \in \mathbb{R}^{n+1}
\end{equation}
where the integral should be computed using an appropriate complex contour. If a time symmetry is available on a spacetime, similar methods roughly corresponding to separation of variables can also yield the Feynman/anti-Feynman propagators (see for example a construction on static spacetime in \cite{Derezi_ski_2018} based on the limiting absorption principle). Unfortunately, this is not always the case, and even when it is, it is often preferable to have a construction of it which works for \textit{perturbations} of the spacetime. This is particularly crucial for nonlinear applications, which physically corresponds to \textit{interacting} QFT. One of the ways this can be achieved is through microlocal methods. The approach we take here is to use \textit{non-elliptic} Fredholm theory to define the Feynman parametrix as a Fredholm inverse of the wave operator on appropriate function spaces. The non-elliptic Fredholm theory we shall use was first proposed by A. Vasy in \cite{vasy2011microlocalanalysisasymptoticallyhyperbolic} in the context of (Riemannian) hyperbolic spaces, and afterwards applied to wave equations in numerous works, including \cite{vasy2014positivitypropagatordifferences} \cite{Gell_Redman_2016} \cite{vasy2018quantumfieldsglobalpropagators} \cite{vasy2017essentialselfadjointnesswaveoperator}. In this work, the spaces will be weighted Sobolev spaces depending on contours $\Gamma$, and can be interpreted as choices of 'boundary conditions' for the wave equation. The purpose of this text is to rigorously define such function spaces, show that exponential perturbations of stationary wave operators are Fredholm on them and finally apply it to the computation of the index of Dirac operators. We start by summarizing the results in subsection \ref{subsec:SummaryResults} and discussing how they relate to previous work in subsection \ref{subsec:PreviousWork}. Afterwards, we quickly recall the Fourier-Laplace transform and construction of Feynman operators for stationary spacetimes in section \ref{sec:FeynmanStationary}, then extend the construction to perturbation of stationary spacetimes in section \ref{sec:PerturbationStationary} and finally apply it to the computation of the index of Dirac operators in section \ref{sec:IndexDirac}. 

\subsection{Summary of Results}
\label{subsec:SummaryResults}
We will being with a short summary of the theory and its main applications. The purpose of this project is to provide a definition of Feynman propagators on perturbations of stationary spacetimes with compact Cauchy surfaces, and then apply it to compute the index of Lorentzian Dirac operator on such spacetimes. Furthermore, we wish for the theory to be robust enough that a rigorous Wick rotation could be applied in this setting. The idea is simple: in the stationary case, given a wave equation $\Box u = f$ on $\mathbb{R}_{t} \times \Sigma$, one can solve it by taking a Fourier transform in time, which gives a one-parameter family $\widehat{\Box}(\sigma)$ of operators on the spatial variables, and then taking an inverse Fourier transform of $\widehat{\Box}(\sigma)^{-1}$ on a contour that does not cross the \textit{resonances}, which are when $\widehat{\Box}(\sigma)$ fails to be invertible:
\begin{equation}
    \Box_{\Gamma}^{-1} = \int_{\Gamma} e^{i\sigma t} \Box(\sigma)^{-1} d\sigma + (\text{finite rank error})
\end{equation}
This can be thought of as the technique of separation of variables under time symmetry. The precise causal nature of $\Box_{\Gamma}^{-1}$ is then determined by $\Gamma$. The first part of the project was rigorously defining Sobolev spaces on which this operation can be done. Furthermore, since we wish for our results to apply to perturbations $\sim O(e^{-t})$ of stationary spacetimes, we wish for them to be acted upon by such perturbations of stationary operators as well. In order for the spaces to be well behaved, a very natural condition is that $\Gamma$ lies in a horizontal strip in $\mathbb{C}$, given as a graph over the real line: $\Gamma = \phi(\mathbb{R})$, $\phi(\sigma) = \sigma + ig(\sigma)$. With that imposed, the first result of the paper takes the following form:
\begin{proposition}(See definition \ref{def:WeightedSobolevSpacesStationary}, \ref{def:b_WeightedSobolevSpace} and proposition \ref{prop:PseudoDiffActContinuously})
    Given two contours $\{ \Gamma_{+}, \Gamma_{-} \} = \Gamma$, there exists Banach spaces $H^s_{\Gamma}$ whose distributions $v \in H^s_{\Gamma}$ have asymptotics at $t = \pm \infty$ given by
    \begin{equation}
        v|_{t = \pm \infty} \simeq \int_{\Gamma_{\pm}} e^{i\sigma t} \widehat{v}(\sigma) d\sigma \quad \widehat{v} \in L^2(\Gamma, \left<\sigma\right>^{2s} d\sigma; H^s_{\Sigma, \sigma}) \quad \text{$H^s_{\Sigma, \sigma} =$ semiclassical Sobolev spaces on $\Sigma$}
    \end{equation}
    Furthermore, a class of pseudo-differential operators $\Psi^{m}_{b, \Gamma}$ acts continuously on $A:H^s_{\Gamma} \rightarrow H^{s-m}_{\Gamma}$, which take the following form near $t = \pm \infty$
    \begin{equation}
        A = A_0 + \sum_{l, k; Im(\eta_l) > 0} e^{i\eta_l |t|} |t|^k A_{lk} +  e^{-(\mu + \delta) |t|} \Psi^{m}_b (M) \quad \mu = \sup_{\sigma \in \Gamma{\pm} - \Gamma{\pm}} |Im(\sigma)|, \delta > 0
    \end{equation}
    Where $A_0$ and $A_{lk}$ are stationary pseudo-differential operators and the $\eta_l$ must satisfy some geometric conditions with respect to $\Gamma$. 
\end{proposition}
The choice of contours we can use will thus be influenced by the decay of our operators to their stationary models. For a perturbation such that the wave operator $\Box$ is in $\Psi^2_{b, \Gamma}$, one can then ask when it is invertible on these spaces. As usual, invertibility is quite delicate and depends on the exact nature of the perturbation. On the other end, invertibility up to finite dimensional obstructions, i.e. $Q \Box = Id - R$, is less senseitive and is amenable to a general theory. Following the (now standard) microlocal approach to non-elliptic Fredholm theory developed in \cite{vasy2011microlocalanalysisasymptoticallyhyperbolic}, the main ingredient we need is an estimate near $t = \pm \infty$. We then use propagation of singularity to cover the whole spacetime and obtain a global estimate. Depending on $\Gamma$, this estimate will take different forms and will in turn dictate the causal properties of the specific parametrix $Q_{\Gamma}$. We present here the estimate that corresponds to the \textit{Feynman} parametrix in a simplified and non-rigorous form. 
\begin{proposition}(Imprecise; see proposition \ref{prop:RadialPointEstimates})
\label{prop:EstimateImprecise}
    Let $\Box$ be a wave operator decaying exponentially fast to stationary models and let $\Gamma$ be contours in $\mathbb{C}$ satisfying assumptions depending only on $\Box$, roughly of the form $Im(\sigma) Re(\sigma) \leq C < 0$ for $\sigma \in \Gamma_{+}$ with $|Re(\sigma)|$ large and some constants $C$ depending on $\Box$ and such that the asymptotics of $\Box$ imply $\Box \in \Psi^2_{b, \Gamma}$. Let $\sigma$ be the frequency dual to the time coordinate $t$. Then for any solutions of $\Box u = f$, we have estimates of the form
    \begin{equation}
    \begin{aligned}
        &\| \Op(\psi(\sigma < 0)) \chi(t >> 1) u \|_{H^s_{\Gamma}} \lesssim \| \Op(\psi(\sigma < 0)) f \|_{H^{s-1}_{\Gamma}} + \| u \|_{H^{-N}_{\Gamma}} \\
        &\| \Op(\psi(\sigma > 0)) \chi(t >> 1) u \|_{H^s_{\Gamma}} \lesssim \| \Op(\psi(\sigma > 0)) f \|_{H^{s-1}_{\Gamma}} + \| h(t)u \|_{H^s} + \| u \|_{H^{-N}_{\Gamma}}
    \end{aligned}
    \end{equation}
    Where $\chi$ is a cutoff near $t = \infty$, $h(t)$ is a compactly supported function of $t$ and $\psi(\pm \sigma > 0)$ is a frequency cutoff in future directed frequencies ($\sigma > 0$) or past directed frequencies ($\sigma < 0$). Similar estimates holds at $t = -\infty$, with the role of $\sigma > 0$ and $\sigma < 0$ reversed and assuming similar properties for $\Gamma_{-}$. 
\end{proposition}
Note that by setting $f = 0$, this corresponds to requiring no high-frequency past directed covectors in the infinite future, and no high-frequency future-directed covectors in the infinite past, which is the high-frequency Feynman condition. Pairing this estimate with the standard propagation of singularity, elliptic estimates and the invertibility theory of $\Box$ in the stationary case then yields a global Fredholm theory for the wave operator on such spacetimes. In particular, it allows us to define Feynman parametrices for it. This is the main result of the paper. A similar theorem also holds for anti-Feynman parametrices, and retarded/advanced propagators. 
\begin{theorem}(Imprecise; see theorem \ref{theo:PropagatorConstruction})
    Let $\Box \in \Psi^2_{b, \Gamma}$ be a wave operator on a globally hyperbolic spacetime $\mathbb{R}_t \times \Sigma$ with compact Cauchy surface and whose stationary model has timelike killing vector fields $\partial_t$. Let $\Gamma$ be contours satisfying the assumptions of proposition \ref{prop:EstimateImprecise}. Assume furthermore that the Fourier-Laplace transform of the stationary model $\widehat{\Box}_{0}(\sigma)$ is invertible on the contours $\Gamma$ and satisfy some estimates. Then there exists a Feynman parametrix for $\Box$ between $\{ u \in H^s_{\Gamma}: \Box u \in H^{s-1}_{\Gamma} \} \rightarrow H^{s-1}_{\Gamma}$, in the sense that $\Box$ is Fredholm on these spaces and a parametrix $E$ satisfies
    \begin{equation}
        WF(E f)|_{t = \pm \infty} \subset \{ \pm \sigma > 0 \} \cup WF(f)|_{t = \pm(\infty)}
    \end{equation}
\end{theorem}
The final part of the paper was applying this theory to Dirac operators. In particular to the computation of their index using a continuous Wick rotation argument, which is possible thanks to the robustness of the Fredholm theory under complex perturbations of the operator. Directly from the Fredholmness of wave operators, we will also obtain the Fredholmness of Dirac operators $D$ on the same spaces. As for the index, Wick rotation is guaranteed to be well behaved as long as the model operator is asymptotically \textit{ultra-static}, in the sense that the metric features no $dt dx$ terms and $\partial_t$ has norm $1$ asymptotically. In that case, we can rotate the $dt$ field in the complex plane to interpolate between our Lorentzian and a Riemannian metric. Since the Fredholm theory is uniform, we then obtain that the index agrees. The Riemannian index is then given by Melrose's result \cite{Melrose1993}. All and all, we obtain the following. The result applies to more general Dirac operators than what is stated here, but these are the main examples.
\begin{theorem}(Imprecise; see theorem \ref{theo:IndexTheorem})
    Suppose the spacetime is even-dimensional and asymptotically static. Let $\Gamma$ be contours corresponding to the Feynman problem such that  they cross $\mathbb{R}$ exactly once at some $\lambda = \pm \epsilon$ small ($+$ for the future and $-$ for the past). Let $D \in \Psi^1_{b, \Gamma}$ be a Dirac operator acting on the spinor bundle or $D = d + \delta$. Then there exists a uniformly Fredholm family of operator $D_{\theta}$ on $H^s_{b, \Gamma}$ such that
    \begin{equation}
        D_0 = D \quad D_{\pi/2} = D_R
    \end{equation}
    where $D_R$ is a Riemannian $b-$Dirac operator. The (graded) index of $D_0$ is thus given by the one of $D_{\theta}$, which shows for the spinor Dirac operator that
    \begin{equation}
        ind_{\lambda}(D) = {}^b\int_{\Bar{M}} \widehat{A}(\nabla) - \frac{(h(A^{+}_{0}) + h(A^{+}_{1}) + \eta(A^{+}_{0}) - \eta(A^{+}_{1}))}{2}
    \end{equation}
    where $\widehat{A}(\nabla)$ is the $\widehat{A}$-form of the connection while $\eta$ is the Eta-invariant and $h$ the dimension of the kernel. $A^{+}_{0}$ denotes the (positive) induced Dirac operator on $\Sigma$ in the past and $A^{+}_{1}$ the (positive) induced Dirac operator on $\Sigma$ in the future. Similarly for the $d+\delta$ operator, on both odd and even dimensional spacetimes which are asymptotically ultra-static, its (graded) index is thus given by minus the total dimension of the deRham cohomology of $\Sigma$
    \begin{equation}
        ind_{\lambda}(d+\delta) = -dim(H_{dR}(\Sigma))
    \end{equation}
\end{theorem}
This formula for the spinor Dirac operator was first found in \cite{B_r_2019} and agrees with the preceding index computations done in the literature, which we will cover in the next subsection. These also made use of Wick rotation in some way. The difference in this treatment is that we do a continuous Wick rotation, which directly allows us to set the Lorentzian and Riemannian indexes equal, rather than a 'one-shot' Wick rotation, which requires further argument to relate the indexes of the two end points of the Wick rotation. Hopefully this setup can also be used as a testbed for the theory of Feynman propagators and Wick rotation, which ideally can then be generalized to the more physical case of non-compact Cauchy surface.
\subsection{Comparison with Previous Work}
\label{subsec:PreviousWork}
Here we give a survey of previous works related to the results established here. The definition of Feynman propagator has attracted much work in the past from the mathematical community. The first rigorous result can be traced back to the work of C. Baer and A. Strohmaier in \cite{B_r_2019}, where they showed that under APS boundary condition, the Dirac operator on a compact globally hyperbolic Lorentzian spacetime with spacelike Cauchy boundary (in particular with finite time interval) is Fredholm, with an index given by the same formula as above (i.e. the same formulas as in the Riemannian case). Subsequent generalizations were done in \cite{B_R_2025} by the same authors, where they relaxed the assumptions and derived a local index theorem, and in \cite{baer2017boundaryvalueproblemslorentzian} by C. Baer and S. Hannes, where different boundary conditions were considered. In the non-compact Cauchy surface setting, M. Braverman gave an index for Dirac operators with growing potential in \cite{braverman2019indexstronglycalliasoperators} and Damaschke gave one in \cite{damaschke2021atiyahsingerdiracoperatorspacetimes} when the Cauchy surfaces are assumed to be Galois coverings. Most relevant to this paper is the result of M. Wrochna and D. Shen in \cite{Shen_Wrochna_2022}, where they generalized for the first time the index theorem and APS boundary conditions to the infinite-time setting with compact Cauchy surface and showed the same index formula as above on polynomially weighted spaces in $t$. As the APS boundary condition has been the main framework for defining and computing indexes of Lorentzian Dirac operators, we give a quick comparison with the one of this paper. The APS boundary condition proposed in \cite{B_r_2019} works by imposing boundary conditions at the initial and final Cauchy surfaces $N_0$, $N_1$ given by
\begin{equation}
    H_{APS} = \{ v \in H([0, 1] \times N) ; \Pi_{A_1 < 0} v|_{t = 1} = 0 \quad \Pi_{A_0 \geq 0} v|_{t=0} = 0 \}
\end{equation}
Where $H$ is some function space and $\Pi$ are spectral projectors of the induced Dirac operators $A_0, A_1$ on $N_0, N_1$. Thus, it imposes that $v$ has no positive/negative frequency modes in the past/future. In the case of an infinite time interval as in \cite{Shen_Wrochna_2022}, this is imposed asymptotically as $t \rightarrow \pm \infty$. To see that this is a Feynman condition, it is useful to think of what solutions to the Dirac equation look like on a ultra-static spacetime $dt^2 - h$. There, the Dirac operator on spinors takes the form
\begin{equation}
    D = cl(\partial_t) (-i\partial_t - A) \quad cl(\partial_t) = \text{Clifford multiplication by $\partial_t$}
\end{equation}
Thus, the solutions $Du = 0$ are given by superpositions of eigenfunctions $\psi_{\lambda}$ of $i A$, with oscillation $e^{i \lambda t}$. The APS condition then translates to the expansion being constrained at the initial and final Cauchy surfaces:
\begin{equation}
    u|_{t = 0} \sim \sum_{\lambda < 0 } \psi_{\lambda} e^{i\lambda t} \quad u|_{t = 1} \sim \sum_{\lambda \geq 0 } \psi_{\lambda} e^{i\lambda t}
\end{equation}
In terms of causality theory, we restrict to only past-directed data at the initial Cauchy surface and future-directed data at the final Cauchy surface, which is exactly the Feynman condition. Now compare it to how this is imposed by the spaces in this paper. Here we take the approach of the boundary condition being prescribed by a contour $\Gamma$. If we expand the solution of our equation into quasi-normal modes $v \sim \sum_{\sigma} e^{i\sigma t} + O(e^{-L |t|})$, then membership in our spaces $H^s_{\Gamma}$ requires 
\begin{equation}
    u \in H^s_{\Gamma} \iff u \in H^s; \quad u|_{t = -\infty} \sim \sum_{\sigma \text{ below } \Gamma_{-}} e^{i\sigma t} + O(e^{Lt}) \text{ and } u|_{t = \infty} \sim \sum_{\sigma \text{ above } \Gamma_{+}} e^{i\sigma t} + O(e^{-Lt})
\end{equation}
In particular in the ultra-static case, where all the modes $\sigma$ are real, the two conditions agree \footnote{Note that in our convention, $\sigma>0$ corresponds to future directed covectors which is the opposite of what physicists typically use.} if we pick the contours $\Gamma$ such that $\Gamma_{-} = \Gamma_{+}$ and $\mp g(\sigma) > 0$ for $\pm \sigma > 0$ (i.e. $\Gamma_{+}$ is below the real axis for positive $\sigma$ and above it otherwise). The main advantage of using the contours $\Gamma$ is that our spaces are directly global in time (of course \cite{Shen_Wrochna_2022} solves that at the expense of introducing limits as $t \rightarrow \pm \infty$ for the APS boundary condition) and do not use the structure of the Dirac equation in any way, hence can also be used to define Feynman parametrices for wave operators that do not come from Dirac operators. It also allows us to handle certain non-static spacetimes. Furthermore, it fits nicely into the $b-$theory developed by Richard Melrose in \cite{Melrose1993} hence the index theorem can be understood from this geometric point of view rather than spectral flow techniques as in \cite{B_r_2019} and the following work.
\newline
Now, parallel to this is the development of non-elliptic Fredholm theory using microlocal tools based on radial set estimates, pioneered by A. Vasy in \cite{vasy2011microlocalanalysisasymptoticallyhyperbolic}. The approach developed there is to define non-elliptic parametrices as Fredholm inverses on weighted Sobolev spaces. The Fredholm estimate itself is shown using two key tools from microlocal analysis: radial sets and propagation of singularities (we refer to Hintz's book for an introduction to microlocal analysis \cite{Hintz2025Microlocal} including what these words mean in general). In the theory of wave equations, this approach has been successfully applied to many problems, nicely summarized in \cite{Vasy2018Minicourse}, and ultimately culminated in applications to nonlinear wave equations in the proof of the stability of the spacetimes Minkowski \cite{Hintz_2020}, Kerr-de Sitter \cite{Hintz_2018} and Kerr \cite{hintz2026nonlinearstabilitysubextremalkerr} \footnote{Of course for these nonlinear equations, substantial additional work has to be done to get the nonlinear estimate. We refer to these results here as an example of applications where a good linear theory is the bedrock of the result.}. It has also been successfully applied to the theory of Anosov flows in dynamics; see the pioneering work of S. Dyatlov and M. Zworski in \cite{dyatlov2016dynamicalzetafunctionsanosov}. The later also made us of the approach to Fredholm theory of \cite{faure2011upper} based on conjugation by exponential weights involving escape functions for the vector field. For the purpose of defining Feynman propagators, the theory was first applied in \cite{Gell_Redman_2016} to construct Feynman propagators on perturbations of Minkowski spacetimes roughly of size $\sim O(\left<(t, r)\right>^{-\delta})$. More precisely, they construct Feynman propagators in the category of \textit{scattering metrics}, which should be thought of in this case as symbolic smooth perturbations of Minkowski that are decaying jointly in $t$ and the radius variable $r$. Other works include a proof of the essential self-adjointness of the wave operator and the limiting absorption principle in \cite{vasy2017essentialselfadjointnesswaveoperator} and a construction of quantum fields from the Feynman propagator in \cite{vasy2018quantumfieldsglobalpropagators}. The theory can also handle spacetimes with gravitational radiation of size $\sim O(\left< r \right>^{-1})$, see for example \cite{molodyk2025feynmanpropagatormassivekleingordon}. As the theory in this paper take this route to define the Feynman propagator, it is important to describe how it relates to these past results. The spaces used for defining Feynman propagators using Vasy's method are typically graph spaces of the following form
\begin{equation}
    H = \{ v \in H^r_{sc}; \Op(\tau^{-l} \rho^{s}) v \in L^2 \} \quad s, l \text{ some functions on $S_{sc}^{\ast}M$}
\end{equation}
In words, they are scattering Sobolev spaces with variable weights in $\tau \sim \left<(t, r)\right>^{-1}$ and $\rho \sim \left<(\sigma, \xi)\right>$. By choosing $s$ and $l$ carefully with respect to the lifted geodesic flow, one can show using radial point estimates that scattering wave operators $\Box$ satisfy Fredholm estimates. The key property to make this definition well behaved is that $\Op(\tau^{-l} \rho^s)$ can be given a precise sense as a pseudo-differential operator in the category of scattering metrics. When working with $b-$metrics, which corresponds to perturbations of size $\sim e^{-|t|}$, this no longer works \footnote{This is related to logarithmic loses in $\tau \sim e^{-t}$ when composing $b-$operators with symbols of the form $\tau^{-l(\sigma, \xi)}$, which we shall encounter in some form as well.}. Nonetheless, the philosophy remains that having a variable weight $\tau^{-l}$ in our space should give an estimate. From that point of view, the work in this paper can be thought of as extending variable order Sobolev spaces to the $b-$setting of exponential decay. We solve these issues by defining the spaces not as a graph space with respect to a chosen quantization $\Op(\tau^{-l} \rho^{s})$, but rather as an exact coordinate-based definition with roots in complex analysis. 
\newline
Finally, one should of course mention that index theorem in the Riemannian setting have been a topic of research for quite a while now. The original index theorem for elliptic operator on closed manifold was proven by M. Atiyah and I. Singer in \cite{Atiyah1968TheIO}, while Dirac operators were discussed in \cite{Atiyah1968TheIO3}. In the case of a compact manifold with boundary, the first result, called the APS index theorem, is also attributed to the same author and V. Patodi in a series of three papers starting with \cite{Atiyah1975SpectralAA}. A nice reference to the topic is the book \cite{BerlineGetzlerVergne1992}. Most relevant to us is the treatment by R. Melrose in \cite{Melrose1993}, where a proof is given for $b-$Dirac operators in terms of the $b-$calculus. The work here should be thought of as an extension of this treatment to the Lorentzian setting. Furthermore, the Wick rotation technique we will be using in this paper to compute the index also has precedent; for example in \cite{Banerjee:2024vcu} the authors discuss the lapse-based Wick rotation of the metric tensor we will be using. Regarding Wick rotation, the novelty in this paper is to construct a robust enough analytical theory so that Wick rotating a Dirac operator can be done while keeping its analytical index constant.

\subsection{Notation}

We here set the notation that will be in used for the remainder of the paper so that our conventions are clear. These will be recalled during the paper, but we divide them in categories here for ease of use. We start by describing the short-hand notation we'll use in our argument:

\begin{itemize}
    \item The symbol $\lesssim$ will be used to denoted $\leq$ up to a constant $C$ that is uniform over the parameters of the inequality. For example $\| u \| \lesssim \| P u \|$ denotes a bound uniform over $u$.
    
    \item For $r \in \mathbb{C}$, the symbol $\left< r \right>$ will denote $\sqrt{1 + |r|^2}$, while if $r$ is a vector it will denote $\sqrt{1 + \| r \|^2}$

    \item Given a multi-index $J = (j_1, j_2, ..., j_k)$ we denote by $\partial^J u$ the derivative $\partial^{j_1} \partial^{j_2} ... \partial^{j_k} u$.
\end{itemize}
We now describe the spaces of smooth functions and the corresponding space of distributions we'll be using throughout the paper. 

\begin{itemize}
    \item On $\mathbb{R}^n$, we denote by $\mathcal{S}(\mathbb{R}^n)$ the Schwartz functions, $\mathcal{S}'(\mathbb{R}^n)$ the tempered distributions, $\mathcal{C}^\infty_c(\mathbb{R}^n)$ the compactly supported smooth functions and $\mathcal{D}'(\mathbb{R}^n)$ the space of distributions on compactly supported functions. Similarly, $\mathcal{C}^k((\mathbb{R}^n)$ denotes function whose first $k$ derivatives exist and are continuous and bounded.
    \item We denote by $\dot{\mathcal{C}}^{\infty}_{exp}(\mathbb{R}^n)$ the super-exponentially decaying smooth functions, i.e. the smooth functions such that for all multi-index $J$ and $N \in \mathbb{R}$ there exists a constant $C$ such that

    \begin{equation}
        |\partial^J u| \leq C_{J, N} \exp(-N \left<x \right>)
    \end{equation}
    The minimal such constants $C_{J, N}$ form seminorms, and a distribution on $\dot{\mathcal{C}}^{\infty}_{exp}(\mathbb{R}^n)$ will be with respect to those.
    \item The space of rapidly decaying smooth functions on a compact manifold with boundary $\Bar{M}$ will be denoted by $\dot{\mathcal{C}}^{\infty}(\Bar{M})$. If we denote by $\tau$ a boundary defining function, then these consists of smooth functions $u$ on $\Bar{M}$ such that for all multi-index $J$ and $N \in \mathbb{R}$ there exists a constant $C$ such that

    \begin{equation}
        | \partial^J u | \leq C_{J, N} \tau^N 
    \end{equation}
    The space of distributions on these test functions with respect to the semi-norms defined by the $C_{J, N}$ will be denoted $\mathcal{C}^{-\infty}(\Bar{M})$. Similarly if $V$ is a smooth vector bundle on it, then $\dot{\mathcal{C}}^{\infty}(\Bar{M}; V)$ denotes rapidly decaying smooth sections of $V$, and $\mathcal{C}^{-\infty}(\Bar{M}; V)$ distributions on them.
    \item Given a space of distribution, say $\mathcal{D}'(\mathbb{R})$, we will denote the subspace of distributions having support in some open set $U$ by $\mathcal{D}'(U)$. 
\end{itemize}

Finally we collect conventions about Sobolev spaces, Fourier-Laplace transform and wave operators.

\begin{itemize}
    \item Let $\sigma$ be a real parameter (pick $Re(\sigma)$ or $|\sigma|$ otherwise). Then $H^s_{\sigma}(\mathbb{R}^n)$ will denote the high-energy Sobolev space of order $s$ with $\sigma$ dependent norm

    \begin{equation}
        \| u \|^2_{H^s_{\sigma}} = \int_{\mathbb{R}^n} \left(1 + \frac{\| \xi \|^2}{\left< \sigma \right>^2} \right)^s |\widehat{u}(\xi)|^2 d\xi
    \end{equation}
    Their usage in this work follows from the fact that the usual Sobolev spaces $H^s$ on $\mathbb{R}^{n+1}$ can be split along one dimension into $H^s_{\sigma}$ valued functions:
    \begin{equation}
        H^s(\mathbb{R}^{n+1}) \simeq \left \{ u(\sigma): \mathbb{R} \rightarrow H^s_{\sigma}; \quad \| u \|^2_{H^s} = \int_{\mathbb{R}} \left< \sigma \right>^{2s} \| u \|^2_{H^s_{\sigma}} d\sigma \right\}
    \end{equation}
    \item We use the following convention for the Fourier-Laplace transforms:

    \begin{equation}
        \widehat{f}(\sigma) = \frac{1}{\sqrt{2\pi}} \int_{\mathbb{R}} e^{-i\sigma t} f(t) dt \quad \sigma \in \mathbb{C}
    \end{equation}
    \item Throughout the paper, $\Gamma$ will denote a complex contour in $\mathbb{C}$ of the form $\Gamma = \phi(\mathbb{R})$

    \begin{equation}
        \phi(\sigma) = \sigma + ig(\sigma) \quad g \in \mathcal{C}^{\infty}(\mathbb{R}) \text{ bounded with bounded derivative $|g'| \leq C < 1$} \quad \mu = sup_{\sigma, \eta \in \mathbb{R}}(g(\sigma) - g(\eta))
    \end{equation}
    As we need them, we'll specify more assumptions on $g$ hence $\Gamma$. 
    \item Given a Lorentzian metric $g$ and inverse metric $g^{-1}$ with sign convention $+, -, -, -, ...$, we define a wave-type operator to be any pseudo-differential operator on vector bundles whose principal symbol $p$ is principally scalar and given by
    \begin{equation}
        p(\zeta) = g^{-1}(\zeta, \zeta) Id \quad \zeta \in T^{\ast}M
    \end{equation}
    In particular, its Hamiltonian flow is given by the geodesic flow on $T^{\ast}M$, and if we write $\eta = (\sigma, \xi)$ with $\sigma$ dual to time, then $\sigma > 0$ corresponds to future-directed frequencies and $\sigma < 0$ the past-directed ones. 'The' wave operator on functions is then given by
    \begin{equation}
        \Box u = - \frac{1}{\sqrt{|det g|}} \partial_{\mu}(\sqrt{|det g|} g^{\mu \nu} \partial_{\nu} u)
    \end{equation}
    In this work we will be mostly interested in wave-type operators between vector bundles, as   these are the ones that arise from Dirac operators. Nonetheless, most of the material will be written for functions, with comments added to explain how to generalize to vector bundles. This is done in order to ease up notation, as in practice going from functions to vector bundles is not difficult. The main distinction to keep in mind between functions and sections of some vector bundle is the existence of a canonical positive definite metric for the first, which makes the wave operator on functions above formally self-adjoint.
    \item Given an operator with asymptotic expansion $P = P_0 + decaying$ on a compact manifold with boundary $\Bar{M}$, $N(P)$ will be defined as the operator $P_0$. It amounts to 'restricting' $P$ to the boundary $\partial \Bar{M}$; hence it serves as a stationary model for the operator $P$. For $b-$pseudo-differential operators, this is a well defined object. Its Mellin transformed family $\widehat{N(P)}(\sigma)$ is called the \textit{indicial family}.
\end{itemize}

\newpage

\section{Feynman Propagators on Stationary Spacetimes}
\label{sec:FeynmanStationary}
We begin by describing the construction of Feynman propagators using the Fourier-Laplace transform of their Kernel. This construction is particularly relevant to us, as it makes inverting non-elliptic operators formally the same as inverting elliptic operators. This method was also one of the motivations for this project. The technique is well known, but finding a specific reference focused on stationary spacetimes turned out to be difficult. We thus refer to \cite{Hormander2003AnalysisI} section 7. for the general theory. Subsections \ref{subsec:FourierLaplaceTransform} and \ref{subsec:ConstructionStationary} recall properties of the Fourier-Laplace transform and its use in construction Green functions, while \ref{subsec:WeightedSpacesStationary} gives the novel definition of spaces on which we can apply the method.

\subsection{The Fourier-Laplace Transform and Associated Distributions}
\label{subsec:FourierLaplaceTransform}
Before going into the construction, it will be useful to collect results about the Fourier-Laplace transform of distribution. Most of the material was taken from Hormander's book \cite{Hormander2003AnalysisI} section 7. Let's begin by defining the Fourier transform of tempered distribution in order to set notation and convention straight.

\begin{definition}
    Let $u \in \mathcal{S}'(\mathbb{R}^n)$ be a tempered distribution. Then its Fourier transform is the tempered distribution $\widehat{u}$ characterized by

    \begin{equation}
        \widehat{u}(\psi) = u(\widehat{\psi}) \quad \psi \in \mathcal{S}(\mathbb{R}^n)
    \end{equation}
    where the Fourier transform of Schwartz functions is defined by
    \begin{equation}
        \widehat{\psi}(\eta) = \frac{1}{\sqrt{(2 \pi)^n}} \int_{\mathbb{R}^n} e^{-i \left< \eta, x \right>} \psi(x) dx
    \end{equation}
\end{definition}
\vspace{0.2in}
The key property requires to take the Fourier transform of tempered distributions is that $\widehat{\psi}$ is also a Schwartz function if $\psi$ is. Unfortunately, if we wish to define it for complex frequencies $\eta$, which amounts to taking the Fourier-Laplace transform, we no longer can define it as easily because there are no natural spaces preserved by the transform. For smooth functions, the most natural class of test functions to take the Fourier-Laplace transform off are the smooth compactly supported function.

\begin{proposition}
\label{prop:PaleyWinerSchawrtz}
(\cite{Hormander2003AnalysisI} Theorem 7.3.1; Paley-Wiener-Schwartz)
    Let $\psi$ be a smooth compactly supported function in $\mathcal{C}_c^{\infty}(\mathbb{R}^n)$. Then its Fourier-Laplace function is the function on $\mathbb{C}^n$ defined by

    \begin{equation}
        \widehat{\psi}(\eta) = \frac{1}{\sqrt{(2 \pi)^n}} \int_{\mathbb{R}^n} e^{-i \left< \eta, x \right>} \psi(x) dx \quad \eta \in \mathbb{C}
    \end{equation}
    It is an entire function of $\eta$ and for any $N \in \mathbb{R}$ it satisfies the estimate

    \begin{equation}
        |\widehat{\psi}(\eta)| \lesssim \left< \eta \right>^{-N} \exp(H_{\psi}(Im(\eta)))
    \end{equation}
    where $H_{\psi}$ is the support function of $\psi$ defined by
    \begin{equation}
        H_{\psi}(\eta) = \sup_{x \in supp(\psi)} \left< \eta, x\right>
    \end{equation}
    Similarly, a tempered distribution $u$ of order $N$ with compact support has an entire Fourier-Laplace transform bounded by $\lesssim \left<\eta\right>^N \exp(H_{u}(Im(\eta)))$. The converse is also true, that is if the Fourier-Laplace transform of $\psi$ satisfies estimates of the above form with $H(\eta) = sup_{x \in K} \left< \eta, x\right>$ for some $K$ compact and convex, then $\psi$ is supported in $K$.
\end{proposition}
Thus, the Fourier-Laplace transform sends smooth compactly supported function to holomorphic functions on $\mathbb{C}^n$ satisfying certain bounds, hence defining the transform by duality is a bit trickier in this case. There are many approaches to this, but the simplest is to consider distributions which are tempered once multiplied by exponential factors. This allows us to define a tempered distribution on the corresponding horizontal line in $\mathbb{C}^n$, and in good situations, this family of distributions will actually come from an holomorphic function. This is the content of the next proposition.

\begin{figure}[htbp]
    \centering
    \includegraphics[width=0.75\textwidth]{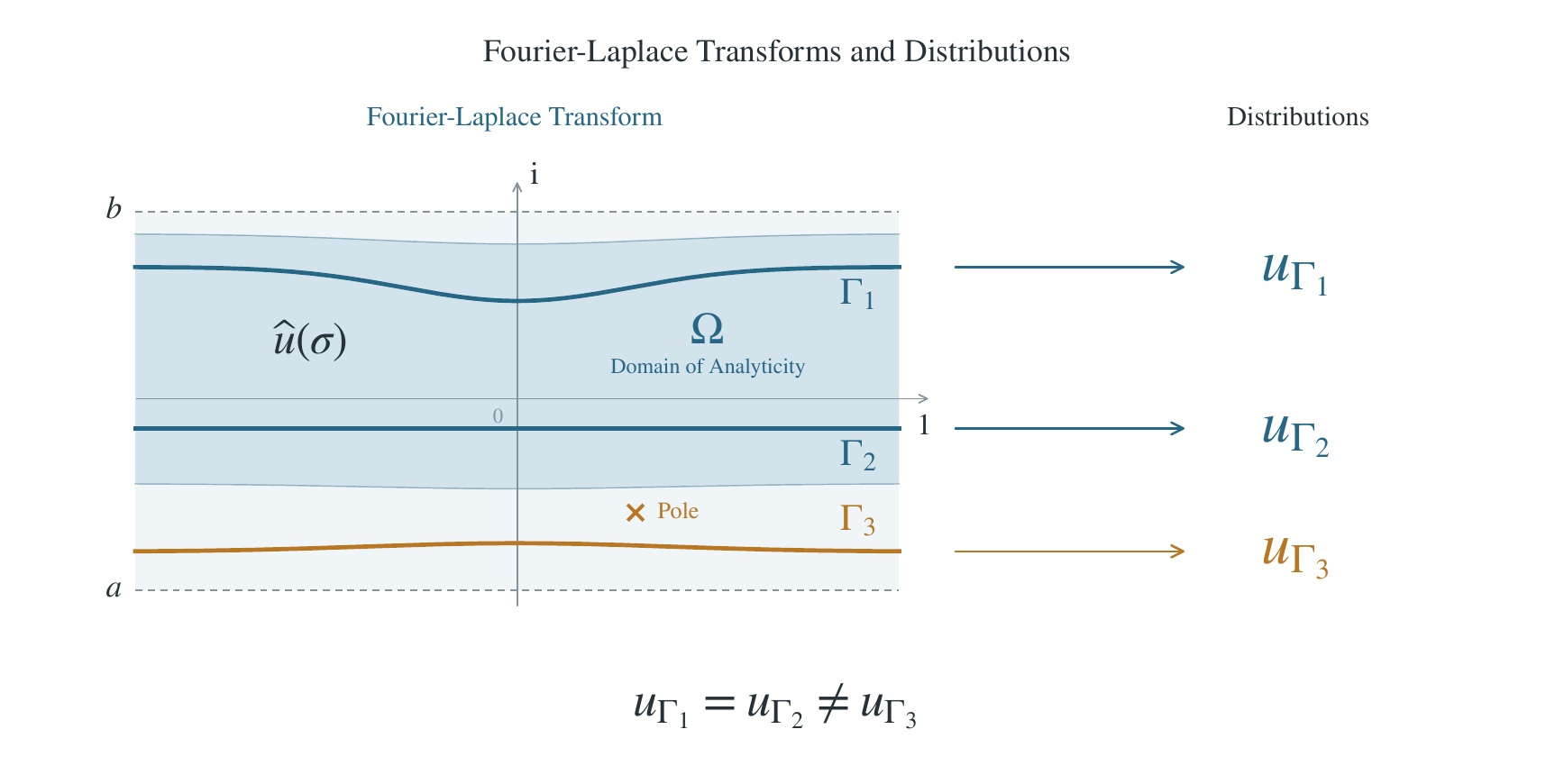}
    \caption{An illustration of the distributional correspondence between Fourier-Laplace transform on contours $\Gamma$ and distributions on $\mathbb{R}$. Analyticity of $\widehat{u}$ makes the correspondence $\Gamma$-independent, while poles contributes a residual term to $u_{\Gamma_2} - u_{\Gamma_3}$. For a general distribution $\widehat{u}|_{\Gamma}$ there might be no way to relate it to a different choice of contour.}
    \label{fig:fourier_laplace_strip}
\end{figure}

\begin{proposition}(\cite{Hormander2003AnalysisI} Theorem 7.4.2)
\label{prop:FLTransform_Domains}
    Let $u \in \mathcal{S}'(\mathbb{R}^n)$ be a tempered distribution and consider the set $\mathcal{F}_u = \{ \eta \in \mathbb{C}: e^{-i\left< \eta, x \right>} u \in \mathcal{S}'(\mathbb{R}^n)\}$. Suppose that $\mathcal{F}$ has a non-empty interior. Then the Fourier-Laplace transform of $u$ exists as an holomorphic function in $\mathcal{F}^{\circ}$ which for some $N$ satisfies the following estimate on horizontal strips

    \begin{equation}
        | \widehat{u} |(\eta) \lesssim \left< \eta \right>^N \quad \eta \in \mathcal{F}^{\circ}, Im(\eta) \text{ varies in a compact set}
    \end{equation}
     Conversely, given an open convex set $\mathcal{I}$ and a holomorphic function $\widehat{u}$ on $\mathbb{R} + i \mathcal{I}$ satisfying an estimate as above, then there exists a (exponentially weighted tempered) distribution $u$ with $\mathcal{F}_u \supset \mathbb{R} + i\mathcal{I}$ and such that the Fourier transform of $e^{-i \left< \eta, x \right>} u$ is $\widehat{u}(\cdot + \eta)$.
\end{proposition}
The canonical example of these distributions are compactly supported distributions, but there are many others. Behind this construction is the fact that one can construct distributions out of various objects that play the role of the 'Fourier-Laplace' transforms on $\mathbb{C}^n$ using duality. One of the key idea in this paper is to work with distributions whose Fourier-Laplace transform in time exists on a well chosen subset of $\mathbb{C}$, which will be a contour $\Gamma$. We will make use of the case $n=1$ hence we focus on that. Take a contour 

\begin{equation}
    \phi(\mathbb{R}) = \Gamma \subset \mathbb{C} \quad \phi(\sigma) = \sigma + ig(\sigma)
\end{equation}
parametrized by a bounded map $g \in \mathcal{C}^{\infty}(\mathbb{R})$ with bounded derivatives. Using duality, we can make sense of distributions whose Fourier-Laplace transform exists on $\Gamma$. In order to define a nice class of test function for $\widehat{u}$ to act on, we work with smooth functions with super-exponential decay:

\begin{equation}
    \dot{\mathcal{C}}^{\infty}_{exp}(\mathbb{R}) = \{ \psi \in \mathcal{C}^{\infty}(\mathbb{R}); \forall k, N \quad |\partial_x^k \psi| \lesssim \exp(-N\left< x \right>) \}
\end{equation}
Proposition \ref{prop:FLTransform_Domains} then tells us that such super-exponentially decaying functions have entire Fourier-Laplace transform, which are Schwartz on any such contours $\Gamma$. Given any two such functions $\psi_1, \psi_2 \in \dot{\mathcal{C}}^{\infty}_{exp}(\mathbb{R})$, by using Fourier inversion on $\Gamma$, we can thus write their distributional pairing as
\begin{equation}
    \int_{\mathbb{R}} \psi_1(t) \psi_2(t) dt = \frac{1}{\sqrt{(2\pi)^{n}}} \int_{\mathbb{R}} \int_{\Gamma} \psi_2(t) \times  e^{i\eta t} \widehat{\psi_1}(\eta) d\eta dt = \int_{\Gamma} \psi_1(\eta) \psi_2(-\eta) d\eta
\end{equation}
We can thus extend this pairing to arbitrary distribution $\widehat{u}$ on $\Gamma$ by defining a distribution $u$ on $\dot{\mathcal{C}}^{\infty}_{exp}(\mathbb{R})$ out of it:
\begin{equation}
    u(\psi) = \widehat{u}(\widehat{\psi}|_{\Gamma}(-\cdot))
\end{equation}
If $\widehat{u}$ comes from an holomorphic function on some strip containing $\Gamma$ as above, then $\widehat{u}(\widehat{\psi}|_{\Gamma}(-\cdot))$ is just integration against $\psi|_{\Gamma}(-\cdot)$ and the resulting value is independent of $\Gamma$ by contour deformation arguments. We will also be interested in the case where $u$ is defined by a meromorphic function $\widehat{u}$ on some horizontal strip containing $\Gamma$. In that case, the choice of contour does matter as different ones will give different value for $u(\psi)$ due to residue terms when deforming the contour. In the general case of an arbitrary distribution, there is no simple way to relate it to a different contour $\Gamma'$. We summarize this discussion in a lemma and figure \ref{fig:fourier_laplace_strip}.

\begin{lemma}
\label{lemma:DistributionalFL_Transform}
    Let $\widehat{u}$ be a tempered distribution on a contour $\Gamma = \phi(\mathbb{R})$. Then it defines a distribution on $\dot{\mathcal{C}}^{\infty}_{exp}(\mathbb{R})$ by 
    \begin{equation}
        u_{\Gamma}(\psi) = \widehat{u}(\widehat{\psi}|_{\Gamma}(-\cdot))
    \end{equation}
    If $\widehat{u}$ is the restriction of a holomorphic function defined on $\Omega$ to $\Gamma$, where $\Omega \subset \mathbb{C}$ is some convex domain containing $\Gamma$, and satisfies $|\widehat{u}(\sigma)| \lesssim \left< \sigma \right>^N$ in $\Omega$, then $\widehat{u}(\widehat{\psi}|_{\Gamma}(-\cdot))$ is independent of $\Gamma \subset \Omega$. In particular if $\Omega = \mathbb{R} + i \mathcal{I}$ as in proposition \ref{prop:FLTransform_Domains}, then $\widehat{u}$ corresponds to the Fourier-Laplace transform of some exponentially weighted tempered distributions, and the two distributions agree for all $\psi \in \dot{\mathcal{C}}^{\infty}_{exp}(\mathbb{R})$. Finally, if $\widehat{u}$ is the restriction of some meromorphic function defined on $\Omega \supset \Gamma$ with finitely many poles and satisfying the same bounds, the value depends on $\Gamma$ but only up to residuals at the poles. We will usually drop the $\Gamma$ in $u_{\Gamma}$ and just write $u$.
\end{lemma}
The most important example of this are Fourier-Laplace functions of exponentials, often called \textit{quasinormal modes}. 

\vspace{0.2in}
\noindent Example: The Fourier-Laplace transform of $e^{i\eta t} t^k$ is 

\begin{equation}
\label{eq:FourierTransformPowers}
    \widehat{u}(\sigma) = \frac{1}{\sqrt{2\pi}} \frac{\partial^k}{\partial \sigma^k} \frac{1}{i(\sigma - \eta)} = \frac{1}{\sqrt{2\pi}}  \frac{k!}{i^{k+1}(\sigma - \eta)^{k+1}}
\end{equation}
in the sense that for any $\psi \in \dot{\mathcal{C}}^{\infty}_{exp}(\mathbb{R})$ supported in $(0, \infty)$ we have

\begin{equation}
    \int_{\mathbb{R}} e^{i\eta t} t^k \psi(t) = \frac{1}{\sqrt{2\pi}}  \int_{\Gamma} \frac{k! \widehat{\psi}(-\sigma)}{i^{k+1}(\sigma - \eta)^{k+1}} d\sigma
\end{equation} 
for any contour $\Gamma$ with $Im(\sigma) < Im(\eta)$ for all $\sigma \in \Gamma$. The proof follows by deforming the contour to $Im(\sigma) = +\infty$. The residue term you pick up is going to be the integral on the left, while $\widehat{\psi}(-\sigma)$ decays exponentially fast as $Im(\sigma) \rightarrow \infty$ due to its support function $H_{\phi}$ satisfying $H_{\psi}(-Im(\sigma)) = -Im(\sigma) \inf_t supp(f)$, hence the integral along $Im(\sigma) = +\infty$ goes to 0. On the other end, the same argument would give 0 if $Im(\eta) > Im(\sigma) \quad \forall \sigma \in \Gamma$, which illustrates the dependence on the contour.
\vspace{0.2in}
\newline
As we will see, distributions with Fourier-Laplace transform on $\Gamma$ as above will be the constituents on which the Fredholm theory of this paper is based. Eventually, we will need to be able to take restrictions of these distributions, hence we quickly discuss how multiplication by a smooth function behaves. Just like in the usual theory of Fourier transforms, multiplication by some function $\chi$ can be expressed as a convolution against the Fourier-Laplace transform $\widehat{\chi}$. The main subtlety is that this is now a convolution in the complex plane, hence apriori requires $\chi$ to have sufficient decay to be well defined. 

\begin{proposition}
    Let $\chi, \psi \in \dot{\mathcal{C}}^{\infty}_{exp}(\mathbb{R})$ be two smooth super-exponentially decaying functions. Then, the Fourier-Laplace transform of $\chi \psi$ is given by complex convolution

    \begin{equation}
        \widehat{\chi \psi}(\sigma) = \frac{1}{\sqrt{2\pi}} \widehat{\chi} \ast \widehat{\psi}(\sigma) = \frac{1}{\sqrt{2\pi}} \int_{\Gamma} \widehat{\chi}(\sigma - \eta) \widehat{\psi}(\eta) d\eta
    \end{equation}
    Where $\phi(\mathbb{R}) = \Gamma$ is any contour as above. In particular, this still holds if one of the two functions, say $\chi$, is in $\exp(-(\mu + \delta) \left< t \right>) \mathcal{C}^{\infty}(\mathbb{R})$ for some $\delta > 0$.
\end{proposition}
Where it gets interesting is multiplication by some bounded smooth function. In that case, $\widehat{\chi}$ is not entire hence one must consider convolution by singular kernels instead. A natural class of functions $\chi$ for which we can still make sense of it are the functions which are eventually constant, in which case $\widehat{\chi}$ has a well behaved meromorphic structure that allows us to define the convolution. 

\begin{proposition}
\label{prop:FourierTransformCutoffs}
    Let $\chi \in C^{\infty}(\mathbb{R})$ be a bounded smooth function approaching a constant $\chi(\infty)$ at a rate $O(e^{-(\mu + \epsilon) \left< t \right>})$ as $t \rightarrow \infty$ with support in $t > 0$. Thus, it is of the form

    \begin{equation}
        \chi \in \mathcal{C}^{\infty}(\mathbb{R}), supp(\chi) \subset \{ t > 0\} \quad \chi' \in \exp(-(\mu + \delta)\left< t \right>) \mathcal{C}^{\infty}
    \end{equation}
    Then, the Fourier-Laplace function of $\chi$ exists as a holomorphic function in $Im(\sigma) < 0$. Its distributional Fourier transform on $\mathbb{R}$ is then given as the boundary value of this holomorphic function:
    \begin{equation}
        \widehat{\chi}(\sigma) = \frac{\widehat{\chi'}}{i(\sigma - i0)} = \sqrt{\frac{\pi}{2}}\chi(\infty) \delta(\sigma) - i p.v. \frac{\chi'(\sigma)}{\sigma}
    \end{equation}
    In particular, it has a distributional meromorphic extension to $Im(\sigma) < \mu + \delta$ in the sense that for any (graphical) contours $\gamma$ such that $\sup_{\gamma} Im(\sigma) < \mu + \epsilon$ and $\{ 0 \} \in \gamma$ we have
    \begin{equation}
        \int_{\mathbb{R}} \chi \psi dt = \int_{\gamma} \widehat{\chi}(\sigma) \widehat{\psi}(-\sigma) d\sigma = \sqrt{\frac{\pi}{2}}\chi(\infty) \widehat{\psi}(0) - i p.v. \int_{\gamma} \frac{\widehat{\chi'}(\sigma) \widehat{\psi}(-\sigma)}{\sigma} \quad \forall \psi \in \dot{\mathcal{C}}^{\infty}_{exp}(\mathbb{R})
    \end{equation}
    Therefore, for any $\psi \in \dot{\mathcal{C}}^{\infty}_{exp}(\mathbb{R})$ and contour $\Gamma = \phi(\mathbb{R})$ such that $sup_{\Gamma - \Gamma} |Im(\sigma)| < \mu + \delta$ the Fourier-Laplace transform of $\chi \psi$ on $\Gamma$ is given by convolution against $\widehat{\chi}$:

    \begin{equation}
        \widehat{\chi \psi}(\sigma)|_{\Gamma} = \frac{1}{\sqrt{2\pi}} \widehat{\chi} \ast \widehat{\psi}(\sigma) = \frac{1}{2} \chi(\infty) \widehat{\psi}(\sigma) - \frac{i}{\sqrt{2 \pi}} p.v. \int_{\Gamma} \frac{\widehat{\chi'}(\sigma - \eta) \widehat{\psi}(\eta)}{\sigma - \eta} d\eta
    \end{equation}
\end{proposition}
\textbf{Proof:} This is a standard computation except the last statement about contours. The claimed holomorphicity of $\widehat{\chi}$ in $Im(\sigma) < 0$ follows from proposition \ref{prop:FLTransform_Domains} along with the fact that $e^{-r t} \chi \in C^{\infty}(\mathbb{R})$ for any $r > 0$. As for $\widehat{\chi}$ on $\mathbb{R}$ being the boundary value of this holomorphic function, this follows from
\begin{equation}
    \int_{\mathbb{R}} \chi \psi dt = \lim_{\epsilon \rightarrow 0} \int_{\mathbb{R}} \exp(-\epsilon t) \chi \psi dt = \lim_{\epsilon \rightarrow 0} \int_{\mathbb{R}} \widehat{\chi}(\sigma - i\epsilon) \widehat{\psi}(-\sigma) d\sigma
\end{equation}
Now, given any (graphical) contour $\gamma$ passing through $\{ 0 \}$ and such that $\sup_{\Gamma} Im(\sigma) < \mu + \delta$, we can use the holomorphicity of $\widehat{\chi'}$ in $Im(\sigma) < \mu + \delta$ to deform the contour $\mathbb{R}$ to $\gamma$ and obtain
\begin{equation}
\begin{aligned}
    \int_{\mathbb{R}} \chi \psi dt &= \lim_{\epsilon \rightarrow 0} \int_{\mathbb{R}} \widehat{\chi}(\sigma - i\epsilon) \widehat{\psi}(-\sigma) d\sigma \\
    &= \lim_{\epsilon \rightarrow 0} \int_{\mathbb{\gamma}} \widehat{\chi}(\sigma - i\epsilon) \widehat{\psi}(-\sigma) d\sigma
\end{aligned}
\end{equation}
which is the desired claim. The one about convolution then follows, since we can deform the convolution integral from $\mathbb{R} - \mathbb{R}$ to $\Gamma - \Gamma$. 
$\Box$
\vspace{0.2in}
\newline
This proposition will end up being crucial later on when we wish to localize distributions near $t = \pm \infty$, as it tells us that one can express the localization of super-exponentially decaying functions $\psi$ to $t = \pm \infty$ purely in terms of their Fourier-Laplace transform on $\Gamma$. Now in practice, the 'time' direction in which we wish to take the Fourier-Laplace transform is only going to be one of the dimensions of our space $\sim \mathbb{R} \times \Sigma$. Thus, we will really be working with \textit{partial} Fourier-Laplace transforms. For smooth functions, this poses no problem but for distributions one needs to be a bit more careful. Roughly speaking, all the results mentioned above extend to this setting, except that now the Fourier-Laplace transform has values in distributions defined on the remaining variables and the estimates given above are the same except that $|\widehat{u}|$ denotes (a sum of finitely many) distributional semi-norms. We summarize this in the following proposition.

\begin{proposition}
\label{prop:rapid_decay_holomorphicity}
(Proposition \ref{prop:FLTransform_Domains} for product distributions)Let $u \in \mathcal{S}'(\mathbb{R}^k \times \mathbb{R}^m)$ be a tempered distribution. Then it's partial Fourier transform in the first $k$ variables is a tempered distribution on $\mathbb{R}^k$ valued in $\mathcal{S}'(\mathbb{R}^m)$. If we denote by $\mathcal{F}_{u} = \{ \sigma \in \mathbb{C}^k: e^{-i\left< \sigma, y \right>}u \in \mathcal{S}'(\mathbb{R}^k \times \mathbb{R}^m)); y \in \mathbb{R}^k \}$, then the Fourier-Laplace transform in the first $k$ variables exists as a (weakly) holomorphic function in $\mathcal{F}^{\circ}$ with values in $\mathcal{S}'(\mathbb{R}^m)$ which for some $N$ satisfies the following estimate on horizontal strips

    \begin{equation}
        | \widehat{u} |(\sigma) \lesssim \left< \sigma \right>^N \quad \sigma \in \mathcal{F}^{\circ}, Im(\sigma) \text{ varies in a compact set}
    \end{equation}
Where $| \widehat{u}(\sigma) | $ denotes a sum of finitely many semi-norms of $\widehat{u}$. Conversely, given an open convex set $\mathcal{I}$ and a (weakly) holomorphic function $\hat{u}$ on $\mathbb{R}^k + i \mathcal{I}$ with values in $\mathcal{S}'(\mathbb{R}^m)$ satisfying an estimate as above, then there exists a (exponentially weighted tempered) distribution $u$ with values in $\mathcal{S}'(\mathbb{R}^m)$ as above with $\mathcal{F}_u \supset \mathbb{R}^k + i \mathcal{I}$ and such that the partial Fourier transform of $e^{-i \left< \sigma, x \right>} u$ is $\widehat{u}(\cdot + \sigma)$.
\end{proposition}
Similarly, if we have a distributional family in $\sigma \in \mathbb{R}^k$ of distributions on $\mathbb{R}^m$, this again defines a distribution on $\mathbb{R}^{k+m}$ by the same construction as lemma \ref{lemma:DistributionalFL_Transform}. As we will be working with class of distributions that are best described by their Fourier-Laplace transform, it will be useful to have a notation for them. 

\begin{definition}
\label{def:Distributions_Gamma}
    Let $\Gamma_1 \times \Gamma_2 \times ... \times \Gamma_k $ be contours in $\mathbb{C}^k$ as above . We will denote by $\mathcal{D}'_{\Pi_i \Gamma_i}$ the set of distributions on super-exponentially decaying smooth functions in the first $k-$variables and tempered in the rest

    \begin{equation}
        \dot{\mathcal{C}}^{\infty}_{exp}(\mathbb{R}^k) \widehat{\otimes} \mathcal{S}(\mathbb{R}^m) = \left \{ \psi \in    \mathcal{C}^{\infty}(\mathbb{R}^{k+m}); \forall k, N, M \quad    |\partial^k \psi| \lesssim \exp(-N \left< y \right>) \left< z \right>^{-M}; y \in \mathbb{R}^k, z \in \mathbb{R}^m \right \}
    \end{equation}
    which are defined by a tempered distribution $\widehat{u}$ on $\Gamma_1 \times \Gamma_2 \times ... \times \Gamma_k$ with values in $\mathcal{S}'(\mathbb{R}^m)$:

    \begin{equation}
        u(\psi(x_1, ..., x_k) \chi(x_{k+1}, ... x_m)) = \left< \left< \widehat{u}(\sigma), \widehat{\psi}(-\sigma) \right>, \chi \right> 
    \end{equation}
    We will denote the distributional family with values in $\mathcal{S}'(\mathbb{R}^m)$ by 

    \begin{equation}
        \left< \widehat{u}(\sigma), \widehat{\psi}(-\sigma) \right> \equiv \int_{\Gamma_1 \times ... \times \Gamma_k} \widehat{u}(\sigma) \widehat{\psi}(-\sigma) d\sigma
    \end{equation}
    If $\widehat{u}$ comes from a holomorphic function on some domain $\Omega$ containing $\Gamma_1 \times \Gamma_2 \times ... \times \Gamma_k $ satisfying estimates of the form $|\widehat{u}(\sigma)|\lesssim \left<\sigma\right>^N$, then we say it is in $\mathcal{H}_{\Omega}$. 
\end{definition}
In particular, as we'll see later, the natural class of operators which acts on the spaces in this paper will have kernels in this set. For that eventual purpose, it is worth mentioning that if a distribution is in $\mathcal{D}'_{\Pi_i \Gamma_i}$, we can also let it act on spaces of test functions with finite exponential growth/decay. Indeed, to justify $u(\psi)$, we only really need that $\widehat{\psi}$ makes sense on $\Pi_i (-\Gamma_i)$, which only requires limited decay. As this will be relevant later on, we make precise the class of test functions on which they act.

\begin{lemma}
\label{lemma:Gamma_Test}
    Let $u \in \mathcal{D}'_{\Pi_i \Gamma_i}$. Then $u$ acts on the space of test functions of the form \footnote{We show later in proposition \ref{prop:ElementaryProperties_Hs} that $\psi \in C^{\infty}_{\Gamma}$ uniquely defines $\widehat{\psi}$ on $\Gamma$ if $|g'| \leq C < 1$, hence this definition is well behaved under our assumptions on $\Gamma$.}

    \begin{equation}
        \mathcal{C}^{\infty}_{\Pi_i (-\Gamma_i)} = \left\{ \psi(y,z) = \frac{1}{(2\pi)^{k/2}}\int_{(-\Gamma_1) \times (-\Gamma_2) \times ... \times (-\Gamma_k)} e^{i \left< \sigma, y \right>} \widehat{\psi}(\sigma, z) d\sigma ; \widehat{\psi} \in \mathcal{S}(\mathbb{R}^k \times \mathbb{R}^m) \right\}
    \end{equation}
    In particular, it acts on test functions decaying with rate $\exp(-(\nu + \delta) \left< y \right>)$ with $\delta > 0$ and $\nu = max_j sup_{\Gamma_j} |Im(\sigma_j)|$, i.e. the ones whose Fourier-Laplace transform is holomorphic in a strip $max_j |Im(\sigma_j)| < \nu + \delta$.
\end{lemma}
We conclude this subsection by an important remark. When working with perturbations of spacetimes, we will be using the framework of the $b-$calculus to reduce the problem to a compact geometry. This is done by applying a compactification $\tau = \exp(-\left< t \right>)$ to the time interval $\mathbb{R}$. This has the effect of converting the Fourier-Laplace transform into the complex Mellin transform. The results in that setting are exactly the same under the transformation $\tau = e^{-t}$: Schwartz functions get replaced by functions decaying faster than any logarithm, while smooth functions decaying faster than any exponential are now rapidly decaying in $\tau$. We will thus use the same notation to denote the Mellin transforms and the corresponding distributions as we did for the Fourier-Laplace transform, with the understanding that all the results stated in this subsection also hold for complex Mellin transform. 
\vspace{0.2in}
\newline
\textbf{Comments for Vector Bundles:} The material in this section extends easily to vector-valued functions, by working component wise. A more invariant way to do it is to consider a distribution $u$ with value in some vector space $V$, i.e. a distribution on functions with values in $V^{\ast}$. Then for any dual element $v \in V^{\ast}$, $\left< u, v \right>$ is a distribution hence we can take its Fourier-Laplace transform, and all the results above apply to it. In particular, for vector-valued functions, saying $u \in \mathcal{D}'_{\Pi_i \Gamma_i}$ will denote $\left< u, v \right> \in \mathcal{D}'_{\Pi_i \Gamma_i}$ for any $v$ and similarly for $\mathcal{H}_{\Omega}$. One could also require that $u(\psi)$ can be expressed in terms of the Fourier-Laplace transform of $\psi$ as above, with $\widehat{\psi}$ valued in $V^{\ast}$ instead and $\widehat{u}$ valued in $V$.

\subsection{Construction of Feynman Propagators using the Fourier-Laplace Transform}
\label{subsec:ConstructionStationary}
We now give a review of the construction of Feynman propagators on stationary spacetimes using complex analysis. I could not find a reference in the literature treating the case of stationary spacetime with compact Cauchy surface in this fashion, but the idea is well known and quite similar to the construction on Minkowski space. Hörmander treats the construction of the Green's function of the wave operator on Minkowski using the Fourier-Laplace transform in section 7.4 of \cite{Hormander2003AnalysisI}. The idea is that due to stationarity, we can consider the Fourier transform in time of the equation, solve for $\widehat{u}(\sigma)$ and then integrate along a contour which avoids its poles. Depending on our choice of contour, this will yield the forward/backward propagators or the Feynman/anti-Feynman propagators. Let's start by recalling the definition of a stationary spacetime with compact Cauchy surface, and then we will describe how the wave operator can be inverted on those.

\begin{definition}
    A stationary spacetime M with compact cauchy surface is a smooth Lorentzian manifold $(M, g)$ of the form

    \begin{equation}
        M \simeq \mathbb{R} \times \Sigma
    \end{equation}
    where $\Sigma$ is a compact Riemannian manifold such that the time translation vector field $\partial_t$ is a timelike Killing vector field. The metric then takes the form

    \begin{equation}
        ds^2 = \lambda (dt - \omega dy_i)^2 - \lambda^{-1} h_{ij} dy^i dy^j
    \end{equation}
    with $h_{ij}$ the Riemannian metric on $\Sigma$, and all components being time-independent. We say it is static if the Cauchy surfaces $\Sigma$ are orthogonal to $\partial_t$, in which case the metric takes the form
    \begin{equation}
        ds^2 = \lambda dt^2 - \lambda^{-1} h_{ij} dy^i dy^j
    \end{equation}
\end{definition}
The Lorentzian metric then defines a time-invariant wave operator $\Box \sim -g^{\ij} \partial_i \partial_j + O(\partial)$ on the spacetime, and our goal is to invert it as a linear operator between function spaces. We do so by considering the Fourier-Laplace transform in time of the wave equation:

\begin{equation}
    \widehat{\Box}(\sigma) \widehat{u}(\sigma) = \widehat{f}(\sigma)
\end{equation}
In order to solve for $\widehat{f}(\sigma)$, we need to invert the family of operators $\widehat{\Box}(\sigma)$. In typical situations, this family will be analytic and Fredholm, hence by the analytic Fredholm theorem will be invertible everywhere besides a discrete set of poles in the complex plane. Thus, as long as $\widehat{f}(\sigma)$ is defined for $\sigma$ off of these poles, one can invert the $\widehat{\Box}(\sigma)$ there and obtain $\widehat{u}(\sigma)$. So suppose $\Gamma$ is a contour in the complex plane of the form $\Gamma = \phi(\mathbb{R})$, $\phi(\sigma) = \sigma + ig(\sigma)$. Assume that $\widehat{\Box}(\sigma)$ is invertible on $\Gamma$. Then, for say $f$ compactly supported in time, we can solve for $\widehat{u}(\sigma) = \widehat{\Box}(\sigma)^{-1} \widehat{f}(\sigma)$ and integrate it on $\Gamma$ with prefactor $e^{i\sigma t}$. By deforming the contour, we can change the integration domain to $\mathbb{R}$ and then obtain $u$ itself:

\begin{equation}
    u = \int_{\Gamma} \frac{e^{i\sigma t} \widehat{f}(\sigma)}{\widehat{\Box}(\sigma)} d\sigma + \sum_k a_k e^{i\sigma_k t}
\end{equation}
where $\sigma_k$ are the poles of $\widehat{\Box}(\sigma)^{-1}$ between $\Gamma$ and the real line $\mathbb{R}$. Thus, up to finitely many 'quasinormal modes', the integral of $e^{i\sigma t} \widehat{\Box}(\sigma)^{-1}$ over a contour gives an inverse for the wave operator. Let's now examine the example of Minkowski in order to understand how the choice of $\Gamma$ impacts which propagator we obtain. Of course Minkowski does not have a compact Cauchy surface, but the example is still instructive as it gives the general idea.

\begin{definition}
    Minkowski spacetime $\mathbb{R}^{n, 1}$ is the Lorentzian metric defined by $ds^2 = dt^2 - \sum_i dx_i^2$ on $\mathbb{R}^{n+1}$.
\end{definition}
The wave operator associated to it is the standard one $\Box = -\partial_t^2 + \sum_i \partial_{x_i}^2$ and can be analyzed using the Fourier-Laplace transform on $\mathbb{R}^{n+1}$. In order to apply the approach illustrated above, we will also take the Fourier transform in the spatial variables so as to avoid dealing with $(\Delta - \sigma^2)^{-1}$ on $\mathbb{R}^n$. We then need to study

\begin{equation}
    \int_{\Gamma} \frac{\widehat{f}(\sigma, \xi)}{\sigma^2 - \xi^2} d\sigma
\end{equation}
for a fixed frequency $\xi \neq 0$. The poles of the integrand are at $\sigma = \pm |\xi|$, and up to homotopy not crossing the poles there are four different contours one can consider, yielding the four propagators on Minkowski, i.e. the possible inverses of the wave operator \footnote{We would still have to treat $\xi = 0$ to get an actual inverse rather than just a parametrix, but this still gives the idea hence we won't try to be rigorous.}. See figure \ref{fig:propagator_contours} for an illustration of each contours.

\begin{definition}
\label{def:MinkowskiPropagator}
    The four different propagators, named the retarded, advanced, Feynman and Anti-Feynman propagators, on the Minkowski spacetime are given by the following integrals along contours $\Gamma$:

    \begin{equation}
    \begin{aligned}
        &\Box^{-1}_{\Gamma} = \int_{\mathbb{R}^n} \int_{\Gamma} \frac{e^{i\sigma t +  \left<x, \xi \right>} \widehat{f}(\sigma, \xi)}{\sigma^2 - \xi^2} d\sigma d\xi \\
        & \Gamma_{retarded} = \{ \sigma: Im(\sigma) = -\epsilon \}\\
        & \Gamma_{advanced} = \{ \sigma: Im(\sigma) = \epsilon \}\\
        & \Gamma_{Feynman} = \phi(\mathbb{R}) \quad \phi(\sigma) = \sigma - i \tanh(\sigma) \\
        & \Gamma_{Anti-Feynman} = \phi(\mathbb{R}) \quad \phi(\sigma) = \sigma + i \tanh(\sigma)\\
    \end{aligned}
    \end{equation}
    For $\widehat{f}(\sigma)$ holomorphic, the choice of contour doesn't matter as long as the placement of the poles $\sigma = \pm |\xi|$ relative to $\Gamma$ (i.e. below or above) is preserved. Otherwise, one often writes a limit as the contour approaches $\mathbb{R}$. 
\end{definition}
\begin{figure}[htbp]
    \centering
    \includegraphics[width=0.5\textwidth]{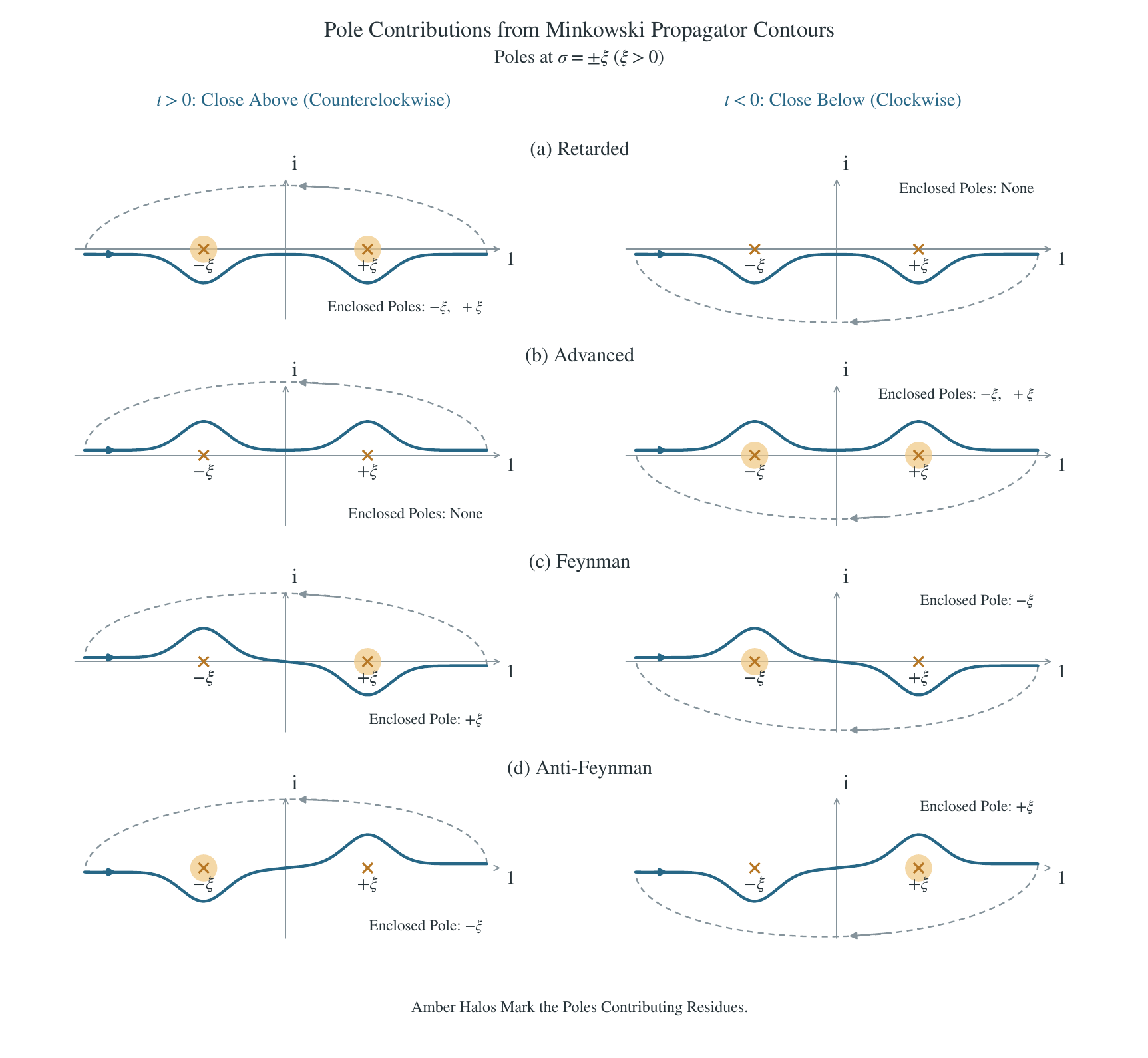}
    \caption{The Minkowski propagators in terms of complex contours in $\mathbb{C}$. The contour impacts which frequency $\sigma$ are included in $t > 0$ or $t < 0$ since one closes the contour in two different ways to compute the integral.}
    \label{fig:propagator_contours}
\end{figure}
We can use the residue theorem to compute these integrals, and must close the loop in the upper half plane for $t > 0$ and the lower half plane for $t < 0$. Thus, we pick up different poles for each of the four contours, which correspond to the different boundary/causality conditions defining each propagators:

\begin{equation}
\begin{aligned}
    &\Box^{-1}_{ret}: \{ f \in \mathcal{C}_c^{\infty}: f|_{t \leq 0} = 0 \} \rightarrow \{ u \in \mathcal{C}^{\infty}: u|_{t \leq 0} = 0 \} \\
    &\Box^{-1}_{adv}: \{ f \in \mathcal{C}_c^{\infty}: f|_{t \geq 0} = 0 \} \rightarrow \{ u \in \mathcal{C}^{\infty}: u|_{t \geq 0} = 0 \} \\
    &\Box^{-1}_{Feyn}: \{ f \in \mathcal{C}_c^{\infty}: \widehat{f}(\sigma > 0)|_{t < 0} = 0, \widehat{f}(\sigma < 0)|_{t > 0} = 0 \} \rightarrow \{ u \in \mathcal{C}^{\infty}: \widehat{u}(\sigma > 0)|_{t < 0} = 0, \widehat{u}(\sigma < 0)|_{t > 0} = 0 \} \\
    &\Box^{-1}_{Ant-Feyn}: \{ f \in \mathcal{C}_c^{\infty}: \widehat{f}(\sigma < 0)|_{t < 0} = 0, \widehat{f}(\sigma > 0)|_{t > 0} = 0 \} \rightarrow \{ u \in \mathcal{C}^{\infty}: \widehat{u}(\sigma < 0)|_{t < 0} = 0, \widehat{u}(\sigma > 0)|_{t > 0} = 0 \} \\
\end{aligned}
\end{equation}
This gives the general idea that we'll try to apply: different contours define propagators on function spaces with boundary conditions on $\widehat{u}$ specified by the contour. In order to make this rigorous, we're left to define weighted Sobolev spaces defined by these contours $\Gamma$ and study when $\widehat{\Box}(\sigma)^{-1}$ acts on them. As the definition of the inverses above involves the Fourier-Laplace transform of the source term only along $\Gamma$, it is natural to work on function spaces defined directly by their Fourier-Laplace transform, such as the distributions defined in \ref{def:Distributions_Gamma}. Furthermore, in the case of compact Cauchy surface, we no longer need to take the Fourier transform in space because $\widehat{\Box}(\sigma)$ will typically be invertible everywhere except discrete $\sigma$, hence we can pick $\Gamma$ to avoid these poles. For convenience, we give the analogue of definition \ref{def:MinkowskiPropagator} in the compact Cauchy surface setting:

\begin{definition}
    Let $M \simeq \mathbb{R} \times  \Sigma$ be a stationary spacetime with compact Cauchy surface and $\widehat{\Box}(\sigma)$ be the Fourier-transformed wave operator. Suppose that $\widehat{\Box}(\sigma)$ is a meromorphic family of Fredholm operators \footnote{We don't specify here the spaces on which they act, but it will be on $L^2$ based spaces} on $\Sigma$ satisfying an estimate of the form

    \begin{equation}
        \| \widehat{\Box}(\sigma)^{-1} \| \lesssim \left<\sigma\right>^{N} \quad \text{for some $N$ and } Im(\sigma) \text{ varying in some bounded strip in $\mathbb{C}$}
    \end{equation}
    Let $\Gamma$ be a contour $\phi(\mathbb{R})$ on which $\widehat{\Box}(\sigma)^{-1}$ exists, i.e. avoiding the set $\mathcal{P}$ of poles of $\widehat{\Box}(\sigma)^{-1}$. Then the wave operator $\Box$ has the explicit inverse

    \begin{equation}
        \Box^{-1}_{\Gamma}(f) = \int_{\Gamma} e^{i\sigma t} \times \widehat{\Box}(\sigma)^{-1} \widehat{f}(\sigma) d\sigma
    \end{equation}
    on smooth functions with decay prescribed by $\Gamma$, as defined in lemma \ref{lemma:Gamma_Test}:

    \begin{equation}
        \Box^{-1}_{\Gamma}(f): \mathcal{C}^{\infty}_{\Gamma} \rightarrow \mathcal{C}^{\infty}_{\Gamma}
    \end{equation}
    Suppose $\mathcal{P}$ is a subset of the real line $\mathbb{R}$. We then define the retard, advanced, Feynman and anti-Feynman operators as the   inverses arising from the following contours:
    \begin{equation}
    \begin{aligned}
        & \Gamma_{retarded} = \{ \sigma: Im(\sigma) = -1 \}\\
        & \Gamma_{advanced} = \{ \sigma: Im(\sigma) = 1 \}\\
        & \Gamma_{Feynman} = \phi(\mathbb{R}) \quad \phi(\sigma) = \sigma - i \tanh(\sigma) + \delta(\mathcal{P}) \\
        & \Gamma_{Anti-Feynman} = \phi(\mathbb{R}) \quad \phi(\sigma) = \sigma + i \tanh(\sigma) + \delta(\mathcal{P}) \\
    \end{aligned}
    \end{equation}
    where $\delta(\mathcal{P})$ is equal to $\frac{1}{2} \inf_{\mathcal{P} \backslash \{ 0\} } |\sigma|$ if $0 \in \mathcal{P}$ and 0 otherwise \footnote{This is just to ensure our contours stay disjoint from $\mathcal{P}$. Whether the 0 frequency is above or below $\Gamma$ is a choice to be made, similar to whether or not 0 should be included in the past or future when imposing the $APS$ boundary condition.}.
\end{definition}
\begin{figure}[h]
    \centering
    \includegraphics[width=0.75\textwidth]{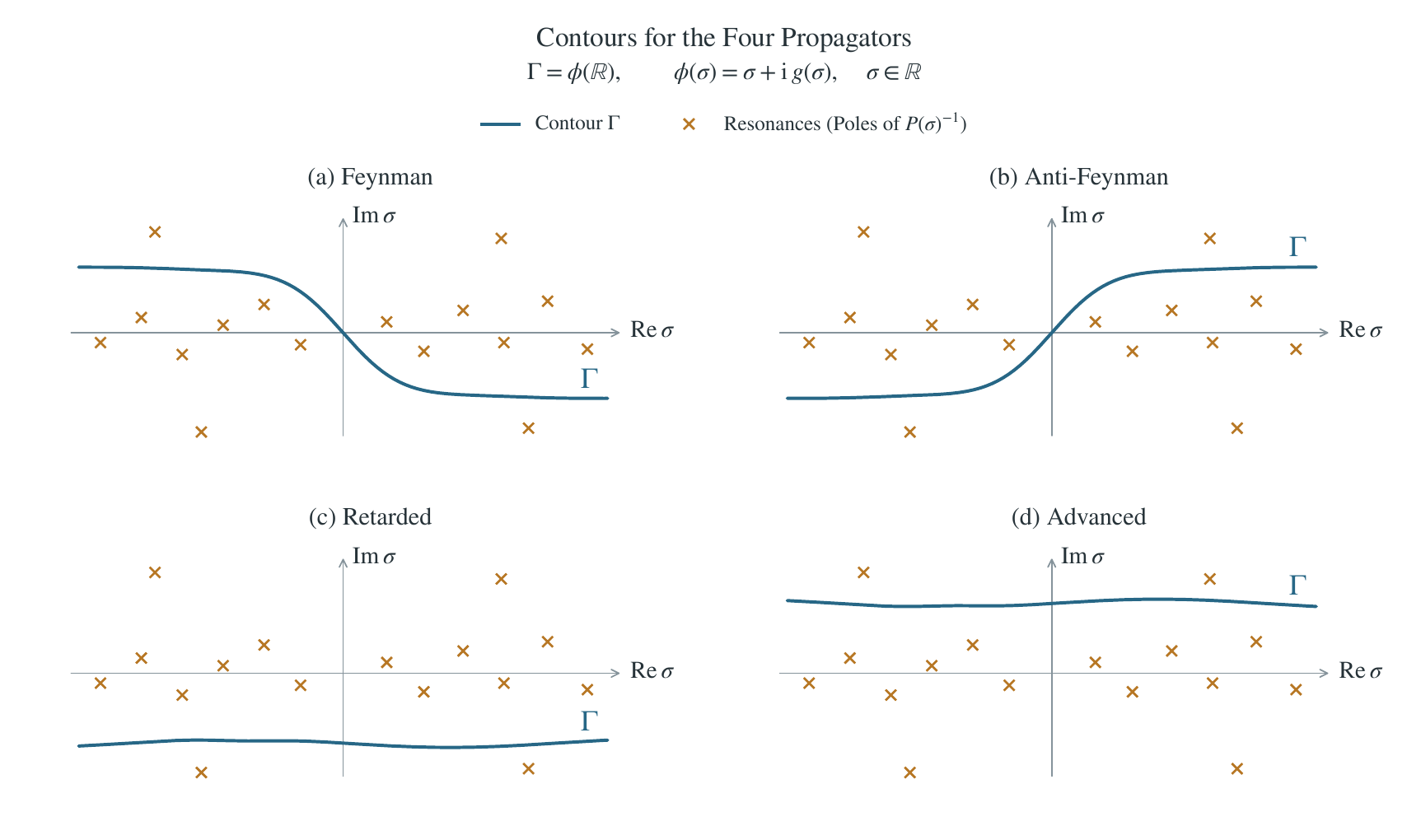}
    \caption{An illustration of different contours and how they correspond to their respective propagators. One can of course have different contours for the same propagators, depending on a choice of finitely many resonances to include or exclude from solutions. What determines the causality of the propagator associated to $\Gamma$ is the sign of $Re(\sigma) \times Im(\sigma)$ on $\Gamma$. Furthermore, if one wishes for exact support conditions on the propagators, rather than just decay/growth, then one should assume that there is no resonances above/below the contour for the retard/advanced propagators, and similarly for Feynman/anti-Feynman in the appropriate quadrants. }
    \label{fig:propagator_contours}
\end{figure}
Other inverses could be considered depending on the relative positions of the poles of $\widehat{\Box}(\sigma)^{-1}$ to $\Gamma$, but they don't have a clean interpretation in terms of causality. The different inverses can thus roughly be considered in bijection with the set of contours $\Gamma$ up to homotopy not crossing the poles of $\widehat{\Box}(\sigma)^{-1}$. The goal of this paper is to obtain a similar theory for \textit{perturbations} of stationary spacetimes. While we won't be able to obtain an explicit inverse in this case due to the problem not being time-translation invariant, we will be able to show that the wave operator is almost invertible (Fredholm) on function spaces with boundary conditions at $t = \pm \infty$ defined by contours $\Gamma$. In the next subsection, we thus define the model function spaces on which our theory will take place and study their basic properties.
\vspace{0.2in}
\newline
\textbf{Comments for Vector Bundles:} Again, for a wave-type operator on vector valued functions, the inversion can be done component-wise, all along the same contour $\Gamma$ avoiding the poles. 

\subsection{Exponentially Weighted in Time Spaces Defined by Contours and Basic Properties}
\label{subsec:WeightedSpacesStationary}
As we've discussed in the previous subsection, to solve a stationary wave equation with time independent coefficients one can integrate against a holomorphic function of the frequency on a contour which avoids poles of $\widehat{\Box}(\sigma)^{-1}$. We now define function spaces on $\mathbb{R} \times \Sigma$ for which this operation can be done in a continuous manner. They will consists of variants of weighted Sobolev spaces defined using a choice of contour $\Gamma \subset \mathbb{C}$ playing the role of where we will integrate $\widehat{\Box}^{-1}(\sigma)$. We assume that it is the image of the real line $\mathbb{R}$ under a map $\phi(\sigma) = \sigma + ig(\sigma)$ where $g \in \mathcal{C}^{\infty}(\mathbb{R})$ has first derivative bounded by $|g'| \leq M < 1$ \footnote{We will need stronger assumptions for the theory to work out. We will state them as needed.}:

\begin{equation}
    \Gamma = \phi(\mathbb{\mathbb{R}}) \quad \phi(\sigma) = \sigma + ig(\sigma) \text{ with $g \in \mathcal{C}^{\infty}(\mathbb{R})$ real bounded with bounded derivatives}; \quad |g'(\sigma)| \leq M < 1
\end{equation}
The example to keep in mind is $g(\sigma) = \tanh(\sigma)$. Roughly speaking, the idea is to consider $H^s$ superpositions of modes $e^{i \phi(\sigma)t}$ of the form

\begin{equation}
    u = \frac{1}{\sqrt{2 \pi}} \int_{\Gamma} e^{i\sigma t} \widehat{u}(\sigma) d\sigma
\end{equation}
One can think of this as exponentially weighted Sobolev spaces where the weight depends on the frequency $\sigma$. This is naturally isomorphic to $H^s$ functions on $\Gamma$ with values in functions on $\Sigma$. More formally, these integrals can be interpreted as distributions acting on super-exponentially decaying functions by

\begin{equation}
    u(v) = \frac{1}{\sqrt{2\pi}} \int_{\mathbb{R}} \int_{\Gamma} e^{i \sigma \cdot t } \left< v(t, x), \widehat{u}(\sigma) \right> dt = \int_{\Gamma} \left< \widehat{v}(-\sigma), \widehat{u}(\sigma) \right> d\sigma
\end{equation}
 They are thus distributions in $\mathcal{D}'_{\Gamma}$ as defined in definition \ref{def:Distributions_Gamma} whose Fourier-Laplace transform on $\Gamma$ is in $H^s$. We next give the proper definition.

\begin{definition}
\label{def:WeightedSobolevSpacesStationary}
    We define the space of exponentially weighted Sobolev spaces  by $\Gamma$ to be the following subspace of distributions $\mathcal{D}'_{\Gamma}$ on smooth, super-exponentially decaying in time, functions:
    \begin{equation}
    H^s_{\Gamma} = \left \{ \frac{1}{\sqrt{2\pi}} \int_{\Gamma} e^{i \sigma \cdot t } \widehat{u}(\sigma) \quad \widehat{u}: \Gamma \rightarrow H^s_{\Sigma, \sigma} \text{ weakly measurable }, \| \widehat{u}(\phi(\sigma)) \|_{H^s_{\Sigma, \sigma}} \in L^2(\mathbb{R}, \left< \sigma \right>^{2s} d\sigma) \right \}
\end{equation}
    where $H^s_{\Sigma, \sigma}$ are the high-energy spaces on $\Sigma$ with norm $\| h(\sigma) \|^2_{H^s_{\Sigma, \sigma}} = \int_{\mathbb{R}}\int_{\Sigma} (1 + \frac{1}{\left< \sigma\right>^2}\Delta_{\Sigma})^s h(\sigma) d\sigma dV$ and $\Delta_{\Sigma}$ is the Laplacian associated to a Riemannian metric on $\Sigma$. The norm on $H^s_{\Gamma}$ is then
    \begin{equation}
        \| u \|^2_{H^s_{\Gamma}} = \int_{\mathbb{R}} \| \widehat{u}(\phi(\sigma)) \|^2_{H^s_{\Sigma, \sigma}} \left< \sigma \right>^{2s} d\sigma
    \end{equation}
    which makes it into a Banach space \footnote{Proposition \ref{prop:ElementaryProperties_Hs} will show that this norm is well defined.}.
\end{definition}
 The case of $\Gamma = \mathbb{R}$ corresponds to the standard Sobolev space $H^s$ while $\Gamma = \{ \sigma: Im(\sigma) = l \}$ corresponds to the exponentially weighted Sobolev space $e^{-lt} H^s$. Furthermore, it will follow from a density argument that $\widehat{u}(\sigma)$ is uniquely specified as a function in $L^2(\Gamma, \left<\sigma\right>^{2s})$, so that these spaces are naturally isomorphic to $L^2(\Gamma, \left<\sigma\right>^{2s})$ functions with value in $H^s_{\Sigma, \sigma}$. For this to be reasonable spaces on which to work, we'd like some basic results about it:

\begin{itemize}
    \item A characterization of membership in $H^s_{\Gamma}$ in terms of $\widehat{u}$ for distributions $u$.
    \item A dense subset which is easy to work with, say the compactly supported smooth functions.
\end{itemize}
These properties all follow from simple properties of the Fourier-Laplace transform, hence we summarize them in the next proposition.

\begin{proposition}(Elementary Properties of $H^s_{\Gamma}$)
\label{prop:ElementaryProperties_Hs}
The weighted spaces $H^s_{\Gamma}$ satisfy the following properties:

\begin{itemize}
    \item $H^s(\mathbb{R} \times \Sigma)$ functions whose partial Fourier transform $\widehat{u}(\sigma)$ in time extends holomorphically from $\mathbb{R}$ to a holomorphic function on some bounded strip $\Omega = \{\sigma: a < Im(\sigma) < b \}$ containing $\Gamma$ and satisfying uniform bounds on its $H^s$ norms along horizontal lines:

    \begin{equation}
        \int_{\mathbb{R}} \left< \sigma \right>^{2s} \| \widehat{u}(\sigma + i \eta) \|^2_{H^s_{\Sigma, \sigma}} d\sigma \lesssim 1 \quad \forall  \mathbb{R}+ i\eta \in \Omega
    \end{equation}
    
    are members of this space. In particular, if $\nu > \sup_{\Gamma} |Im(\sigma)|$, $\exp(-\nu \left< t \right>)H^s(\mathbb{R} \times \Sigma)$ functions are members of $H^s_{\Gamma}$. 

    \item The smooth compactly supported functions are dense in these spaces. In particular, the action of $\int_{\Gamma} e^{i\sigma t} \widehat{u}(\sigma) d\sigma \in H^s_{\Gamma}$ as a distribution on super-exponentially decaying in time functions uniquely determines $\widehat{u}(\sigma)$ hence $H^s_{\Gamma} \simeq L^2(\Gamma, \left< \sigma \right>^{2s})$.

    \item The $L^2$ pairing on super-exponentially decaying functions

    \begin{equation}
        \left< \psi_1, \psi_2 \right> = \int_{\mathbb{R} \times \Sigma} \psi_1(t, x) \overline{\psi_2(t, x)} dt   dx = \int_{\Gamma} \left< \widehat{\psi_1}(\sigma), \overline{\widehat{\psi_2}(\Bar{\sigma})} \right> d\sigma
    \end{equation}
    extends to a non-degenerate continuous pairing between $H^s_{\Gamma}$ and $H^{-s}_{\Bar{\Gamma}}$
    \begin{equation}
        \left< u, v \right> = \int_{\Gamma} \left< \widehat{u}(\sigma), \overline{\widehat{v}(\overline{\sigma})} \right> d\sigma
    \end{equation}
    which makes them dual to each other.

\end{itemize}

\end{proposition}
\textbf{Proof:} These properties are essentially 1d statements in $t$ hence we will suppress the tangential variable in $\Sigma$ for the proof. Furthermore, by adding $\left<\sigma\right>^s$ factors, one can obtain the $H^s$ statements from the $L^2$ one, so we focus on $s = 0$. We prove each of the statement step by step:
\vspace{0.1in}
\newline
1) The first property follows by contour deformation, i.e. by the holomorphic assumption and the uniform estimate on $L^2$ norms we have

\begin{equation}
    u = \int_{\mathbb{R}} e^{i \sigma t} \widehat{u}(\sigma) d\sigma = \int_{\Gamma} e^{i \sigma t} \widehat{u}(\sigma) d\sigma                 
\end{equation}
in the sense of distributions. The $\exp(-\nu \left< t \right>)H^s$ functions have holomorphic Fourier-Laplace transform in $\Omega = \{ |Im(\sigma)|\leq \nu \} \supset \Gamma$ satisfying estimates of this form, hence are members of this space.
\vspace{0.1in}
\newline
2) As for the smooth compactly supported functions being dense, the statement is equivalent to $h \in L^2(\mathbb{R})$ satisfying

\begin{equation}
    \int_\mathbb{R} h(\sigma) \widehat{u}(\phi(\sigma)) d\sigma = 0 \quad \forall \text{ $u$ smooth compactly supported}
\end{equation}
implies $h = 0$. One can show this by an explicit computation: pick for some $\eta \in \mathbb{R}$ the function $\widehat{u}(\sigma) = \exp(i \eta \cdot \sigma) \widehat{\chi}(\sigma)$, where $\widehat{\chi}$ is the Fourier Laplace transform of a smooth compactly supported function $\chi$. $\widehat{u}(\sigma)$ is the transform of a translate of $\chi$ hence is also the Fourier-Laplace transform of a compactly supported function. Define the function

\begin{equation}
    \Tilde{h}(\sigma) = h(\sigma) \times \widehat{\chi}(\phi(\sigma))
\end{equation}
This is a rapidly decaying function of $\sigma$ and the property for $h$ implies

\begin{equation}
    \int_{\mathbb{R}} e^{i \eta \phi(\sigma)} \Tilde{h}(\sigma) d\sigma = 0 \quad \forall \eta \in \mathbb{R}
\end{equation}
This integral is absolutely convergent and defines a smooth function. Integrating it against super-exponentially decaying function of $\eta$ gives

\begin{equation}
    \int_{\eta \in \mathbb{R}} \psi(\eta) \int_{\mathbb{R}} e^{i \eta \phi(\sigma)} \Tilde{h}(\sigma) d\sigma d\eta = \int_{\mathbb{R}} \widehat{\psi}(-\phi(\sigma)) \Tilde{h}(\sigma) d\sigma =  0
\end{equation}
i.e. this function of $\eta$ is 0 as a distribution on super-exponentially decaying function. Thus, it suffices to show that such functions have dense Fourier-Laplace transform in $L^2(\Gamma)$, since that would imply $\Tilde{h}(\sigma) = 0$, which implies $h = 0$ since the Fourier-Laplace transform of compactly supported functions separates points. We show this by a convolution argument. Let $f \in \mathcal{C}^{\infty}_c$ be a compactly supported smooth function on $\Gamma$ and consider the functions

\begin{equation}
    f_{\epsilon}(\phi(\sigma)) = \frac{1}{\sqrt{4\pi \epsilon}} \int_{\mathbb{R}} \exp(-(\phi(\sigma) - \phi(\eta))^2 / 4\epsilon) f(\phi(\eta)) \phi'(\eta) d\eta
\end{equation}
i.e. $f_{\epsilon}$ is the convolution of $f$ with the usual Gaussian kernel but now computed on $\Gamma \times \Gamma$. These correspond to the Fourier-Laplace transform of smooth and super-exponentially decaying functions since they are entire with uniform Schwartz estimates in horizontal strips. We claim that 

\begin{equation}
    lim_{\epsilon \rightarrow 0} f_{\epsilon} = f \quad \text{in }L^2(\Gamma)
\end{equation}
Indeed, first note that we can change variables to write the convolution as
\begin{equation}
    f_{\epsilon}(\phi(\sigma)) = \frac{1}{\sqrt{4\pi}} \int_{\mathbb{R}} \exp(-(\phi(\sigma) - \phi(\sigma - \sqrt{\epsilon} \xi))^2 / 4\epsilon) f(\phi(\sigma - \sqrt{\epsilon} \xi)) \phi'(\sigma - \sqrt{\epsilon} \xi) d\xi
\end{equation}
The limit of the integrant pointwise as $\epsilon \rightarrow 0$ is

\begin{equation}
    \frac{1}{\sqrt{4 \pi}} \exp(-\xi^2 \phi'(\sigma)^2 / 4) f(\phi(\sigma)) \phi'(\sigma)
\end{equation}
We now wish to apply the dominated convergence theorem. Note that the exponential above is uniformly bounded over $\epsilon$ by

\begin{equation}
    |\exp(-(\phi(\sigma) - \phi(\sigma - \sqrt{\epsilon} \xi))^2 / 4\epsilon)| \leq \exp(- \alpha \xi^2 / 4)
\end{equation}
where $\alpha = inf Re(\phi'(\sigma)^2) > 0$ by our assumption that $|g'(\sigma)| < 1$. Thus we can apply DCT hence we obtain in $L^2$
\begin{equation}
\begin{aligned}
    lim_{\epsilon \rightarrow 0} f_{\epsilon}(\phi(\sigma)) &= \frac{1}{\sqrt{4\pi}} \int_{\mathbb{R}} \exp(-\xi^2 \phi'(\sigma)^2/ 4) f(\phi(\sigma)) \phi'(\sigma) d\xi \\
    &= f(\phi(\sigma))
\end{aligned}
\end{equation}
Where in the last step we used the formula for integral of complex Gaussians. To get the statement for general $f \in L^2(\Gamma)$, use density of $\mathcal{C}^{\infty}_c(\Gamma)$ in $L^2(\Gamma)$. 
\vspace{0.1in}
\newline
3) The duality statement follows from $L^2(\Gamma) \simeq L^2(\overline{\Gamma})$ under the map $\sigma \rightarrow \overline{\sigma}$, then use that $L^2$ is self-dual. As for the equality with the usual pairings of functions on $\mathbb{R} \times \Sigma$, proceed by changing the order of integration:

\begin{equation}
\begin{aligned}
    \int_{\mathbb{R} \times \Sigma} \psi_1(t, x) \overline{\psi_2(t, x)} dt   dx &= \int_{\Sigma} \int_{\mathbb{R}} \int_{\Gamma} e^{i\sigma t} \widehat{\psi_1}(\sigma, x) d\sigma \overline{\psi_2(t, x)} dt dx \\
    &= \int_{\Gamma} \int_{\Sigma} \widehat{\psi_1}(\sigma, x) \int_{\mathbb{R}} e^{i\sigma t} \overline{\psi_2(t, x)} dt dx d\sigma \\
    &= \int_{\Gamma}  \left< \widehat{\psi_1}(\sigma), \overline{\widehat{\psi_2}(\overline{\sigma})} \right>
\end{aligned}    
\end{equation}
$\Box$
\vspace{0.2in}
\newline
By the last proposition, these spaces behave very similarly to the usual weighted Sobolev spaces $e^{-lt} H^s$ hence it is reasonable to try to invert the wave operator $\Box$ on them through the Fourier-Laplace transformed operator $\widehat{\Box}^{-1}(\sigma)$. Assuming that $\widehat{\Box}(\sigma)^{-1}$ has no poles on $\Gamma$ and satisfies \textit{high-energy} estimates

\begin{equation}
    \| \widehat{\Box}(\sigma)^{-1} \|_{H^s_{\Sigma, \sigma} \rightarrow H^{s+k}_{\Sigma, \sigma}} \lesssim \left< \sigma \right>^{-k} \quad \sigma \in \Gamma
\end{equation}
the operator $\int_{\Gamma} e^{i\sigma t} \widehat{\Box}(\sigma)^{-1} d\sigma$ acts continuously on this space as a Fourier-multiplier operator

\begin{equation}
   \widehat{\Box}(\sigma)^{-1}: H^s_{\Gamma} \rightarrow H^{s+k}_{\Gamma}
\end{equation}
and this provides an inverse for the wave operator on these spaces. More generally, any Fourier multiplier acts on these spaces, which we think off as 'stationary' operators.

\begin{lemma}
\label{lemma:StationaryOperatorContinuous}
    Let $A(\sigma): \Gamma \rightarrow \mathcal{L}(H_{\sigma}^s(\Sigma), H^{s-k}_{\sigma}(\Sigma))$ be a one-parameter family of bounded linear maps on $H^s_{\sigma}(\Sigma)$ satisfying

    \begin{equation}
        \| A(\sigma) \|_{H^s_{\Sigma, \sigma} \rightarrow H^{s-k}_{\Sigma, \sigma}} \lesssim \left< \sigma \right>^{k}
    \end{equation}
    Then $A$ defines a continuous map $A: H^s_{\Gamma} \rightarrow H^{s-k}_{\Gamma}$. In particular, if $\widehat{\Box}(\sigma)^{-1}$ has no pole on $\Gamma$ and $\widehat{\Box}(\sigma)^{-1}$ satisfies high-energy estimates, then $\Box$ is invertible between $\{ u \in H^s_{\Gamma}: \Box u \in H^{s-k}_{\Gamma}\} \rightarrow H^{s-k}_{\Gamma}$.
\end{lemma}
\textbf{Proof:} By proposition \ref{prop:ElementaryProperties_Hs} elements of $H^s_{\Gamma}$ are determined by their Fourier-Laplace transform in $H^s(\Gamma, H^{s}_{\Sigma, \sigma})$ hence it suffices to show that the action is continuous on this space. So given $\widehat{v}(\sigma) \in H^s(\Gamma, H^{s}_{\Sigma, \sigma})$, we have

\begin{equation}
    \|A(\phi(\sigma)) \widehat{v}(\phi(\sigma)) \|_{H^{s-k}_\Gamma} = \int_{\mathbb{R}} \| A(\phi(\sigma)) \widehat{v}(\phi(\sigma)) \|^2_{H^{s-k}_{\Sigma, \sigma}} \left< \sigma \right>^{2(s-k)} d\sigma  \lesssim \int_{\mathbb{R}} \| \widehat{v}(\phi(\sigma)) \|^2_{H^{s}_{\Sigma, \sigma}} \left< \sigma \right>^{2s} d\sigma = \| \widehat{v}(\phi(\sigma)) \|_{H^s_{\Gamma}}
\end{equation} 
Thus, it is continuous. The statement about the wave operator follows from directly defining the inverse as a bounded operator.
$\Box$
\vspace{0.2in}
\newline
Thus, we've defined weighted Sobolev spaces on which Feynman propagators in the stationary case naturally act on. For this to be a well-behaved theory, we'd like that it generalizes well to \textit{perturbations} of stationary spacetimes. This will involve studying the action of operators that are stationary up to decaying perturbative terms. We need two key tools for that:

\begin{itemize}
    \item That decaying powers of $e^{-t}$ act on them continuously when the support is restricted to $[0, \infty)$.
    \item That these spaces are isomorphic to $H^s(\mathbb{R} \times \Sigma)$ under some reasonable operator which would allow us to compute the norms using pairings of distributions $\left< u , v \right>$. This is required so that we can use standard arguments in microlocal analysis based on pairings, which is what will eventually give us the Fredholm estimates.
\end{itemize}
The first property ammounts to a containment statement between these different spaces according to their rate of decay, since powers of $e^{-t}$ naturally translate the contour $\Gamma$ in the complex plane. We expect that if $\Gamma$ implies more decay for $L^2_{\Gamma}$ than $\Gamma'$, than $L^2_{\Gamma} \subset L^2_{\Gamma'}$. The next proposition confirms this for reasonable $\Gamma$ and $\Gamma'$, and is illustrated in figure \ref{fig:IntermediateContour}:
\begin{figure}[h]
    \centering
    \includegraphics[width=0.75\linewidth]{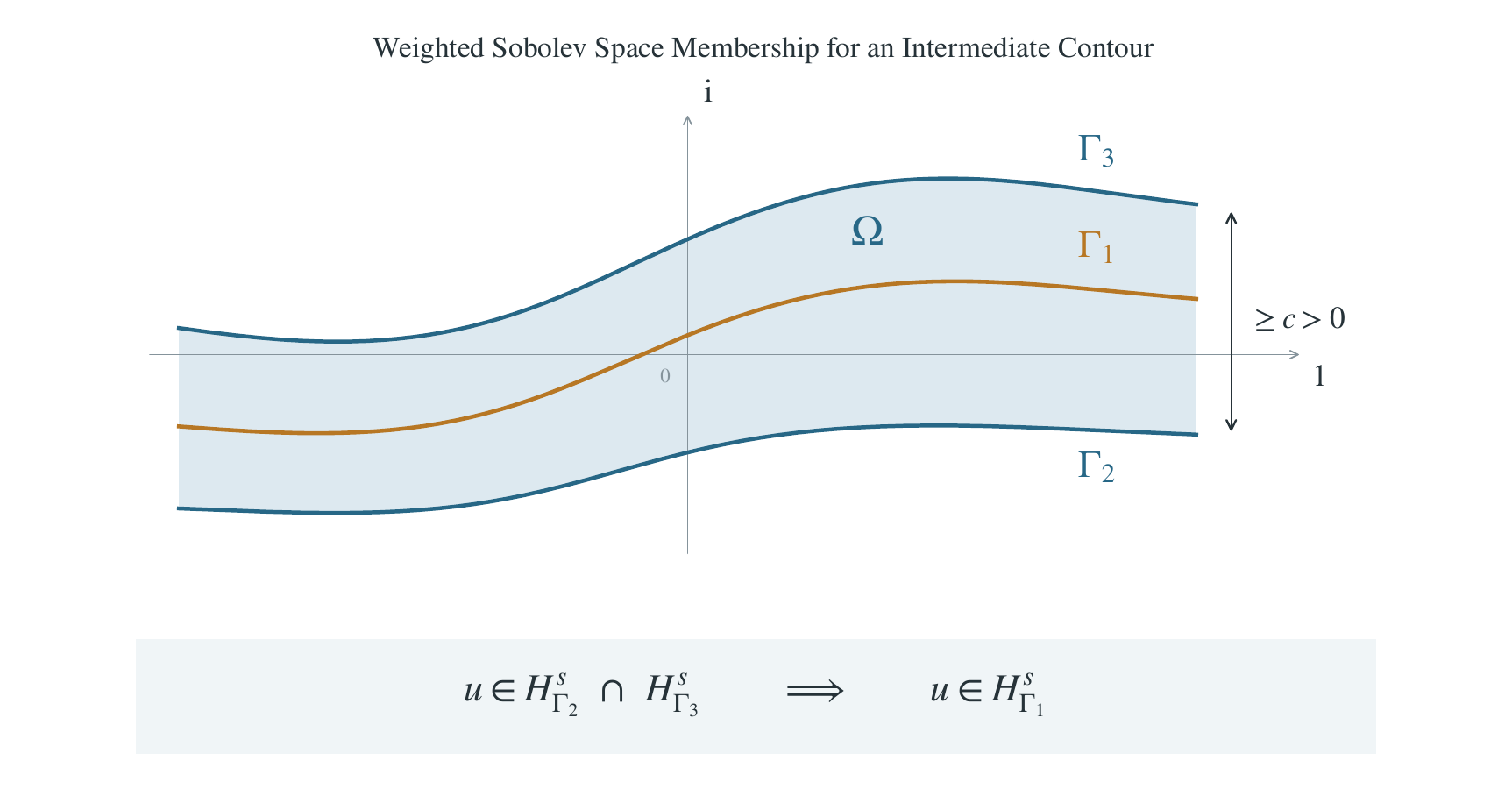}
    \caption{An illustration of how different $H^s_{\Gamma}$ are related. Membership in the boundary contours $\Gamma_2$ and $\Gamma_3$ guarantees membership in all the intermediate contours $\Gamma_1$. This can be understood from the fact that the boundary values of $\widehat{u}$ defines a unique holomorphic extension to the domain $\Omega$ between the contours, given by the Cauchy transform.}
    \label{fig:IntermediateContour}
\end{figure}

\begin{proposition}
\label{prop:InclusionStationarySpaces}
Suppose $\psi$ is a super-exponentially decaying function in time and let $\Gamma_i = \phi_i(\mathbb{R})$, $i= 1,2 ,3$, be three contours such that $Im(\phi_3(\sigma)) > Im(\phi_1(\sigma)) > Im(\phi_2(\sigma))$ for all $\sigma \in \mathbb{R}$ such that

\begin{equation}
    Im(\phi_3(\sigma)) - Im(\phi_2(\sigma)) \geq c > 0 \quad \text{ uniformly over $\sigma$}
\end{equation}
and similarly for $\Gamma_1$, i.e. it is uniformly vertically separated from $\Gamma_2$ and $\Gamma_3$. Then there exists a constant independent of $\psi$ such that

    \begin{equation}
    \label{eq:Bounds_Different_Gamma}
        \| \widehat{\psi}(\sigma)|_{\Gamma_1} \|_{L^2(\Gamma_1, \left< \sigma \right>^{2s} d\sigma)} \leq C \left( \| \widehat{\psi}(\sigma)|_{\Gamma_2} \|_{L^2(\Gamma_2, \left< \sigma \right>^{2s} d\sigma)} + \| \widehat{\psi}(\sigma)|_{\Gamma_3} \|_{L^2(\Gamma_3, \left< \sigma \right>^{2s} d\sigma)}\right)
    \end{equation}
    and we have the continuous inclusion:
    \begin{equation}
        H^s_{\Gamma_2} \cap H^s_{\Gamma_3}  \hookrightarrow H^s_{\Gamma_1} 
    \end{equation}
    given by the weighted Cauchy transform on $\Omega = \{ \sigma: Im(\phi_2(Re(\sigma))) < Im(\sigma) < Im(\phi_3(Re(\sigma))) \}$ with weight $w(\eta) = (C + \sqrt{\Delta} + i\eta)^{s}$ for $C >\sup_{\mathbb{R}}|g(\sigma)| $ and $\Delta$ some positive Laplacian on $\Sigma$:
    \begin{equation}
        \widehat{u}|_{\Gamma_1}(\sigma) = \frac{1}{w(\sigma)} \left \{ \frac{1}{2\pi i} \int_{\Gamma_2} \frac{w(\eta) \widehat{u}|_{\Gamma_2}(\eta)}{\eta - \sigma} d\eta - \frac{1}{2\pi i} \int_{\Gamma_3} \frac{w(\eta) \widehat{u}|_{\Gamma_3}(\eta)}{\eta - \sigma} d\eta \right \}
    \end{equation}
    In particular, $u$ is then in $\mathcal{H}_{\Omega}$ in the sense that is defined by a holomorphic Fourier-Laplace transform in $\Omega$. Furthermore, the inclusion $H^s_{\Gamma_2} \cap H^s_{\Gamma_3} \hookrightarrow H^r_{\Gamma_1}$ is compact for $s > r$. Finally, if one restricts the support to positive time $t$ in $(0, \infty) \times \Sigma$ we furthermore have a continuous inclusion

    \begin{equation}
        H^s_{\Gamma_1} \cap \mathcal{D}'_{\Gamma_1}((0, \infty) \times \Sigma)  \hookrightarrow H^s_{\Gamma_2} \cap \mathcal{D}'_{\Gamma_2}((0, \infty) \times \Sigma) 
    \end{equation} 
    given also by the (one-sided, weighted) Cauchy transform:
    \begin{equation}
    \label{eq:OneSidedCauchy}
        \widehat{u}|_{\Gamma_2}(\sigma) = \frac{-1}{2\pi i w(\sigma)} \int_{\Gamma_1} \frac{w(\eta) \widehat{u}|_{\Gamma_1}(\eta)}{\eta - \sigma} d\eta
    \end{equation}
    in which case $\widehat{u}$ is holomorphic below $\Gamma_1$, i.e. in $\mathcal{H}_{\{ Im(\sigma) < g_1(Re(\sigma)) \}}$. Again, $H^s_{\Gamma_1} \cap \mathcal{D}'_{\Gamma_1}((0, \infty) \times \Sigma)  \hookrightarrow H^r_{\Gamma_2} \cap \mathcal{D}'_{\Gamma_2}((0, \infty) \times \Sigma)$ is compact for $s>r$. In particular, one can think of elements of $H^s_{\Gamma} \cap \mathcal{D}'_{\Gamma}((0, \infty))$ as being given by holomorphic Fourier-Laplace transforms in $\{ Im(\sigma) < g(Re(\sigma)) \}$ with $L^2(\Gamma, \left<\sigma\right>^{2s})$ data at the boundary $\Gamma$.
\end{proposition}
\textbf{Proof: } Again this is a $1d$ statement hence we suppress the tangential variables. Given any super-exponentially decaying function $\psi$, its Fourier-Laplace transform $\widehat{\psi}(\sigma)$ is holomorphic in $\sigma$ with rapid decay on horizontal lines, hence its value in the region $\Omega$ bounded by $\Gamma_2$ and $\Gamma_3$ is given by the Cauchy transform of its boundary values \footnote{This follows directly from the residue theorem applied to increasingly large contours.}:
\begin{equation}
    \widehat{\psi}(\sigma) = \frac{1}{2\pi i} \int_{\Gamma_2} \frac{\widehat{\psi}|_{\Gamma_2}(\eta)}{\eta - \sigma} d\eta - \frac{1}{2\pi i} \int_{\Gamma_3} \frac{\widehat{\psi}|_{\Gamma_3}(\eta)}{\eta - \sigma} d\eta \quad \sigma \in \Omega
\end{equation}
In particular, taking the restriction to $\Gamma_1$ gives a continuous map from $L^2(\Gamma_2, d\sigma) \bigoplus L^2(\Gamma_3, d\sigma)$ to $L^2(\Gamma_1, d\sigma)$ since by our uniform separation assumption on the contour, the kernel (of say the first term) of this integral transform on $\mathbb{R} \times \mathbb{R} \simeq \Gamma_2 \times \Gamma_3$ can be split as
\begin{equation}
\begin{aligned}
    K(\eta, \sigma) &= \frac{1}{2\pi i} \frac{1 + ig'_2(\eta)}{Re(\eta) - Re(\sigma) + i(g_2(\eta) - g_1(\sigma))} \\
    & = \frac{1}{2\pi i} \left( \frac{1 + ig'_2(\eta)}{Re(\eta) - Re(\sigma) - i} - \frac{i(1 + ig'_2(\eta))(g_2(\eta) - g_1(\sigma) + 1)}{(Re(\eta) - Re(\sigma) - i)(Re(\eta) - Re(\sigma) + i(g_2(\eta) - g_1(\sigma)))} \right)
\end{aligned}
\end{equation}
The first term is the kernel for flat contours (times a bounded factor) hence is $L^2$ bounded, while the second term is $L^2$ bounded by Schur's test. The weighted statement similarly follow by applying the result to $w(\eta) \widehat{\psi}$ and using that $w$ is holomorphic. Now, to turn this into a containment statement, we show that given $u \in H^s_{\Gamma_2} \cap H^s_{\Gamma_3}$, one can approximate $\widehat{u}$ by a holomorphic function on $\Omega$ in both norms simultaneously. Indeed, by the construction from above, we can approximate it in $H^s_{\Gamma_2}$ by

\begin{equation}
    \widehat{u}_{\epsilon}(\phi_2(\sigma)) = \frac{1}{\sqrt{4\pi \epsilon}} \int_{\mathbb{R}} \exp(-(\phi_2(\sigma) - \phi_2(\eta))^2 / 4\epsilon) \widehat{u}(\phi_2(\eta)) \phi_2'(\eta) d\eta
\end{equation}
Now using that $u$ in $H^s_{\Gamma_3}$ defines the same distribution on super-exponentially decaying functions, this integral can also be written as

\begin{equation}
    \widehat{u}_{\epsilon}(\phi_2(\sigma)) = \frac{1}{\sqrt{4\pi \epsilon}} \int_{\mathbb{R}} \exp(-(\phi_2(\sigma) - \phi_3(\eta))^2 / 4\epsilon) \widehat{u}(\phi_3(\eta)) \phi_3'(\eta) d\eta
\end{equation}
Furthermore, by a contour deformation argument applied to the holomorphic functions $u_{\epsilon}$, the following also defines the same distribution

\begin{equation}
    \widehat{u}_{\epsilon}(\phi_3(\sigma)) = \frac{1}{\sqrt{4\pi \epsilon}} \int_{\mathbb{R}} \exp(-(\phi_3(\sigma) - \phi_3(\eta))^2 / 4\epsilon) \widehat{u}(\phi_3(\eta)) \phi_3'(\eta) d\eta
\end{equation}
So $u_{\epsilon}$ tend to $u$ in both of $H^s_{\Gamma_{2, 3}}$ hence by inequality \ref{eq:Bounds_Different_Gamma} (which applies to any holomorphic function in $\Omega$ with the required decay by what we showed) it defines a distribution in $H^s_{\Gamma_1}$ that is equal to $u$, which shows the containment. The Cauchy formula then also holds for $u$ by the same argument. This also explains why the inclusion is compact for $s > r$, since then the kernel between $L^2$ spaces are of the form (say for the $\Gamma_2$ to $\Gamma_1$ part)
\begin{equation}
    K(\eta, \sigma) = \frac{1}{2\pi i (C + i\phi_1(\sigma))^{s-r}} \frac{1 + ig'_2(\eta)}{Re(\eta) - Re(\sigma) + i(g_2(\eta) - g_1(\sigma))}
\end{equation}
Taking a cutoff in $|Re(\sigma)| < C$ then give a compact operator, and as $C \rightarrow \infty$ these operators tend in operator norm to $K(\eta, \sigma)$ thanks to the extra decay from $s > r$. Finally, if the support is restrained to $(0, \infty)$, we then have by what we showed

\begin{equation}
    H^s_{\Gamma_1} \cap \mathcal{D}'_{\Gamma_1}((0, \infty) \subset H^s_{\Gamma_1} \cap H^s_{Im(\sigma) = -C} \subset H^s_{\Gamma_2}
\end{equation}
where $C > 0$ is chosen so that $-C < inf_{\Gamma_2} Im(\sigma) < inf_{\Gamma_1} Im(\sigma)$. The first containment follows from $u \in H^s_{\Gamma_1} \cap \mathcal{D}'_{\Gamma_1}((0, \infty))$ implying that $u$ is a distributions on functions of the form $\psi \in e^{-Ct} H^{-s}$ since for some cutoff $\chi$ equal to 1 on $(-1, \infty)$ we have

\begin{equation}
    u(\psi) = u(\chi \psi) = \int_{\Gamma_1} \widehat{u}(\sigma) \widehat{\chi \psi}(-\sigma) d\sigma
\end{equation}
where $\widehat{\chi \psi}$ exists on $\Gamma_1$ since $| \chi \psi | \lesssim e^{-C\left< t \right>}$ implies that it is has holomorphic Fourier-Laplace transform in a strip $|Im(\sigma)| < C$ which contains $\Gamma_1$. Thus by duality $u$ is in $e^{Ct}H^s = H^s_{Im(\sigma) = -C}$ with norm bounded by the $L^2(\Gamma_1, \left<\sigma\right>^{2s} d\sigma)$ norm of $\widehat{u}_{\Gamma_1}$ \footnote{Technically this uses that $\widehat{\chi \psi}(\sigma)|_{\Gamma_1}$ and $\widehat{\chi \psi}(\sigma)|_{Im(\sigma) = -C}$ are comparable, which again follows from the Cauchy formula.}. Since the support of a distribution is independent of the subspace of distributions in which it is defined, we conclude that we have a continuous inclusion
\begin{equation}
       H^s_{\Gamma_1} \cap \mathcal{D}_{\Gamma_1}'((0, \infty)) \hookrightarrow H^s_{\Gamma_2} \cap \mathcal{D}'_{\Gamma_2}((0, \infty))
\end{equation}
as desired. It is given by the one-sided Cauchy formula \ref{eq:OneSidedCauchy} since we can apply it on contours $\Gamma_1$ and $Im(\sigma) = -C$, then take $C \rightarrow \infty$. The corresponding integral vanishes in the limit since $\widehat{u}(\sigma)|_{Im(\sigma) = -C}$ decays exponentially fast in $C$, by the fact that $\widehat{u}(\sigma)|_{Im(\sigma) = -C}$ defines the distribution $e^{-Ct} u$ which tends in norm to $0$ as $C \rightarrow \infty$. Once again, compactness follow from the Cauchy integral formula. $\Box$
\vspace{0.2in}
\newline
Now the last property we required of our spaces was that they should be naturally isomorphic to the usual $H^s(\mathbb{R} \times \Sigma)$ spaces so that we can apply pairing arguments. Directly from the definition, the space is trivially isomorphic to $L^2(\Gamma, \left<\sigma\right>^{2s})$ functions with value in distributions on $\Sigma$ under the Fourier-Laplace transform, which is itself isomorphic to $H^s(\mathbb{R} \times \Sigma)$ under pull-back by $\phi$. Nonetheless, it is still worth summarizing that as a proposition since it is the key property that will allow us to obtain non-trivial estimates for wave operators on these spaces.

\begin{proposition}
\label{prop:IsoContourSobolev}
    The Fourier-Laplace transform in time restricted to $\Gamma \subset \mathbb{C}$ 

    \begin{equation}
        \mathcal{F}_{\Gamma}(\psi) = \frac{1}{\sqrt{2\pi}} \int_{\mathbb{R}} e^{-i\sigma t} \psi(t, x) dt \quad \sigma \in \Gamma \quad x \in \Sigma
    \end{equation}
    acting on smooth super-exponentially decaying functions in time extends to an isometry between $H^s_{\Gamma}$ and $L^2(\Gamma, \left< \sigma \right>^{2s} d\sigma)$ functions with values in $H^s_{\Sigma, \sigma}$. Thus, composing this with the pullback map by $\phi$ gives an isometry $T$ with $H^s(\mathbb{R} \times \Sigma)$:

    \begin{equation}
        T(u)(t, x) = \frac{1}{2\pi} \int_{\mathbb{R}} \int_{\mathbb{R}} e^{i \sigma t} e^{-i \phi(\sigma) t'} u(t', x)  dt' d\sigma
    \end{equation}
    The inverse map is given by taking the Fourier transform in time and then taking the pullback map for the projection $\Gamma \rightarrow \mathbb{R}$ to obtain a member of $L^2(\Gamma, \left< \sigma \right>^{2s} d\sigma)$, which can be mapped to $H^s_{\Gamma}$:

    \begin{equation}
        T^{-1}(v)(t, x) = \frac{1}{2\pi} \int_{\Gamma} \int_{\mathbb{R}} e^{i\sigma t} e^{-i \phi^{-1}(\sigma) t'} v(t', x) dt' d\sigma
    \end{equation}
     Furthermore, $T$ maps super-exponentially decaying smooth functions to Schwartz functions:
     \begin{equation}
        T: \dot{\mathcal{C}}^{\infty}_{exp}(\mathbb{R} \times \Sigma) \rightarrow \mathcal{S}(\mathbb{R} \times \Sigma)
     \end{equation}
      And if $g$ extends to a bounded holomorphic function in a strip $|Im(\sigma)| \leq C$, then $T$ maps super-exponentially decaying smooth functions to functions of decay $e^{-C\left< t \right>}$:
      \begin{equation}
       T: \dot{\mathcal{C}}^{\infty}_{exp}(\mathbb{R} \times \Sigma) \rightarrow e^{-C\left< t \right>} \mathcal{C}^{\infty}(\mathbb{R} \times \Sigma)
      \end{equation}
\end{proposition}
\textbf{Proof:} This essentially follows from the definition of our space. By density of super-exponentially decaying functions, defining the map on them suffices. It is an isometry since $\| u \|^2_{H^s_{\Gamma}} = \int_{\mathbb{R}} \left< \sigma \right>^{2s} \| \widehat{u}(\phi(\sigma)) \|^2_{H_{\sigma}^s(\Sigma)} d\sigma$ on those which is exactly the $H^s$ norm of $T(u)$. It is left to check that $T^{-1}$ as above is indeed the inverse of $T$. Note that we have

\begin{equation}
\begin{aligned}
    T^{-1} \circ T &= \mathcal{F}_{\Gamma}^{-1} \circ (\phi^{-1})^{\ast} \circ \mathcal{F}^{-1} \circ \mathcal{F} \circ \phi^{\ast} \circ \mathcal{F}_{\Gamma}(v) \\
    &= \mathcal{F}_{\Gamma}^{-1} \circ \mathcal{F}_{\Gamma}(v) \\
    &= Id
\end{aligned}
\end{equation}
hence the two maps are inverse of each other on $H^s_{\Gamma}$. As for the mapping properties, they follow from rapid decay in $\sigma$ being preserved under the pullback maps $\phi^{\ast}$ while if $\phi$ is holomorphic we can deform the contour of integration to obtain decay.
$\Box$
\vspace{0.2in}
\newline
\noindent Note that this takes a form similar to a pseudo-differential operators, except that the 'symbol' $e^{-g(\sigma) t}$ isn't symbolic in $\sigma$ in general since when differentiating it in $\sigma$ one gets losses in powers of $t$ and cannot obtain improvement in the joint frequency $\left< (\sigma, \xi) \right>$. Nonetheless, the theory will behave very similarly to the usual theory of weighted Sobolev spaces. Now, as explained above, in practice one wants to solve the wave equations when the coefficient are dependent on time in a control way. For that, we will need to use these spaces on perturbations of stationary spacetimes, and develop a more robust theory than writing down explicit inverses to show that the wave operator $\Box$ is invertible. In the next section, we construct such a theory using standard microlocal tools: radial set estimates and the propagation of singularity theorem.
\vspace{0.2in}
\newline
\textbf{Comments for Vector Bundles:} For vector-valued distributions, the theory of these spaces proceeds in the exact same fashion by working component-wise and requiring $\Gamma$ to be the same for all components. Better yet, one chooses a norm on the vector space $V$ and uses that to define the $H^s_{\Gamma}$ spaces, with $\widehat{u}(\sigma)$ now being valued in the vector space $V$.
\newpage

\section{Perturbations of Stationary Spacetimes}
\label{sec:PerturbationStationary}
We now move on to constructing the theory for perturbations of stationary spacetimes. In order to work with them, we need a setup which naturally incorporates the perturbations we will consider. One such setting is provided by the \textit{b-calculus}, which was developed by Melrose to solve differential equations on compact manifold with boundary, hence we will be working within that framework. An exposition can be found in his book \cite{Melrose1993} chapter 2, 4 and 5; we will recall the basics as we need them. For the stationary spacetimes mentionned above, this will mean we will be working with \textit{exponentially decaying} perturbation of stationary spacetimes. This is a rather natural assumption, since the spaces defined above in the stationary case are weighted by exponential in time factors. Setting ourselves on a compact geometry makes the analysis much easier, hence we compactify the spacetime $\mathbb{R} \times \Sigma$ by compactifying the time interval $\mathbb{R}$. We thus define the compactified spacetime $\Bar{M}$ by

\begin{equation}
    \Bar{M} = [0, 1] \times \Sigma \quad M \hookrightarrow \Bar{M}, (t, x) \rightarrow (\rho(t), x) \in [0, 1] \times \Sigma
\end{equation}
where $\rho(t)$ is any bijective smooth monotonic function equal to $1 - e^{-t}$ for $t >> 0$ large and $e^t$ for $t << 0$ very negative. Under the change of variable $t \rightarrow \tau = e^{-t}$, the Fourier transform becomes the Mellin transform 

\begin{equation}
    \widehat{u}(\sigma) = \frac{1}{\sqrt{2\pi}} \int_{\mathbb{R}} \tau^{i\sigma} u(\tau) \frac{(-d\tau)}{\tau}
\end{equation}
and the standard calculus on $\mathbb{R} \times \Sigma$ corresponds to the $b-calculus$ on $\Bar{M}$, which is defined in terms of the Mellin transform and can be thought of as allowing differentiation in $\partial_t = -\tau \partial_{\tau}$ only rather than $\partial_{\tau}$. In subsection \ref{subs:Contour_Weighted_Sobolev}, we define the weighted Sobolev spaces we wish to work on directly on the compactification $\Bar{M}$, while in the remaining subsections we prove the properties it satisfies, and study when a pseudo-differential operator will statisfy Fredholm estimates on it.

\subsection{Contour Weighted b-Sobolev Spaces $H^s_{b, \Gamma}$}
\label{subs:Contour_Weighted_Sobolev}
We are interested in defining weighted Sobolev space on the compactified spacetime $\Bar{M}$ on which we could invert the wave operator. As in the stationary case \ref{subsec:WeightedSpacesStationary}, we define them through a choice of contours. As this part does not require an exact product structure, we will work in the general case even though the spaces we will apply the theory too will be of the form $[0, 1] \times \Sigma$. Let $\Bar{M}$ be an arbitrary compact manifold with boundary, $\partial \Bar{M}_i$ be the connected components of it's boundary and $M$ be its interior $int(\Bar{M}) = M$. Pick normal neighborhoods $N_{i}$ with product structure for each of them

\begin{equation}
    N_i \simeq [0, b)_{\tau_i} \times \partial M_i
\end{equation}
Here, $\tau_i$ is a defining function for the boundary component $\partial \Bar{M}_i$. On this normal neighborhood and with respect to this product structure, we can define the spaces analogously as we did in the stationary case in definition \ref{def:WeightedSobolevSpacesStationary}, but now using the Mellin transform and a contour $\Gamma_i$ depending on $N_i$:

\begin{equation}
     H^s_{N_{i}, \Gamma_i} = \left \{ \frac{1}{\sqrt{2\pi}} \int_{\Gamma_i} \tau^{-i\sigma t} \widehat{u}(\sigma) d\sigma ; \widehat{u}: \Gamma_i \rightarrow H^s_{\partial M_i, \sigma}, \| \widehat{u}(\phi_i(\sigma)) \|_{H^s_{\partial M_i, \sigma}} \in L^2(\mathbb{R}, \left< \sigma \right>^{2s} d\sigma) \right \}
\end{equation}
They satisfy the same property as in the stationary case by the correspondence $\log(\tau) = -t$ between the Mellin and the Fourier transform. In particular, these are distributions on smooth rapidly decaying functions in $\tau$ (i.e. super-exponentially decaying in $t$) in the class $\mathcal{D}'_{\Gamma_i}$. We would then like to define the global weighted space as the distributions on $\Bar{M}$ which in $N_i$ are in $H^s_{N_i, \Gamma_i}$. For this to be reasonable, we need to be able to localize the elements of $H^s_{N_i, \Gamma_i}$ near $\tau_i = 0$. The following lemma states exactly that. 

\begin{lemma}
\label{lemma:CutoffContinuous}
    Let $\chi(\tau) \in C^{\infty}([0, \infty))$ be a smooth bounded function with bounded derivatives supported in $[0, b)$ and equal to a constant $c$ on $[0, \epsilon)$, for $0 < \epsilon < b$. Then multiplication by $\chi$ is a continuous map on $H^s_{\Gamma}([0, \infty) \times \Sigma)$.
\end{lemma}
\textbf{Proof:} By proposition \ref{prop:FourierTransformCutoffs} given any super-exponentially decaying function $\psi$ the Mellin transform in time of $\chi \psi$ on $\Gamma$ is given by

\begin{equation}
    \widehat{\chi \psi}(\sigma) = \frac{c}{2} \widehat{\psi}(\sigma) - \frac{i}{\sqrt{2 \pi}} p.v. \int_{\Gamma} \frac{\widehat{-\tau \partial_{\tau}(\chi)}(\sigma - \eta) \widehat{\psi}(\eta)}{\sigma - \eta} d\eta
\end{equation}
So it remains to show that this is a bounded map on $L^2(\Gamma, \left< \sigma \right>^{2s} d\sigma)$. The first term is obviously bounded so we focus on the principal value ($p.v.$) integral. For Lipschitz contours $\Gamma$, in \cite{CoifmanMcIntoshMeyer1982} the following operator

\begin{equation}
\label{eq:p.v.Contour}
    p.v. \int_{\Gamma} \frac{\widehat{\psi}(\eta)}{\sigma - \eta} d\eta
\end{equation}
was proven to be continuous on $L^2(\Gamma)$. Note that for smooth $\Gamma$, i.e. when $\phi$ is smooth, this is a pseudo-differential operator on $\Gamma$ hence it's mapping property on $L^2$ are a special case of the pseudo-differential calculus. To get the same statement with a $\widehat{\tau \partial_{\tau}(\chi)}(\sigma - \eta)$ factor in the integral, use a Taylor series for $\widehat{\tau \partial_{\tau}(\chi)}(\sigma - \eta)$ around $\eta = \sigma$ to split the integral into a singular term (a $p.v.$ integral as in eq. \ref{eq:p.v.Contour}) plus a convolution by a smooth rapidly decaying function. Both terms are then bounded on $L^2(\Gamma)$. Finally to get the continuity on weighted spaces $L^2(\Gamma, \left< \sigma \right>^{2s} d\sigma)$, apply the unweighted result to $\left< \eta \right>^{s} \widehat{\psi}(\eta)$ by absorbing the extra $\left< \eta \right>^{-s}$ factor in $\widehat{\tau \partial_{\tau}(\chi)}(\sigma - \eta)$. The statement for general distribution $u \in H^s_{\Gamma}$ then follows by density of smooth super-exponentially decaying functions. 
$\Box$
\vspace{0.2in}
\newline
With this lemma, we can 'glue' together the asymptotics defined by $H^s_{N_i, \Gamma_i}$ to form a global distribution on $\Bar{M}$. Thus, the following definition of our weighted spaces on compact manifolds with boundary is now well behaved. 
\begin{figure}[h]
    \centering
    \includegraphics[width=0.75\linewidth]{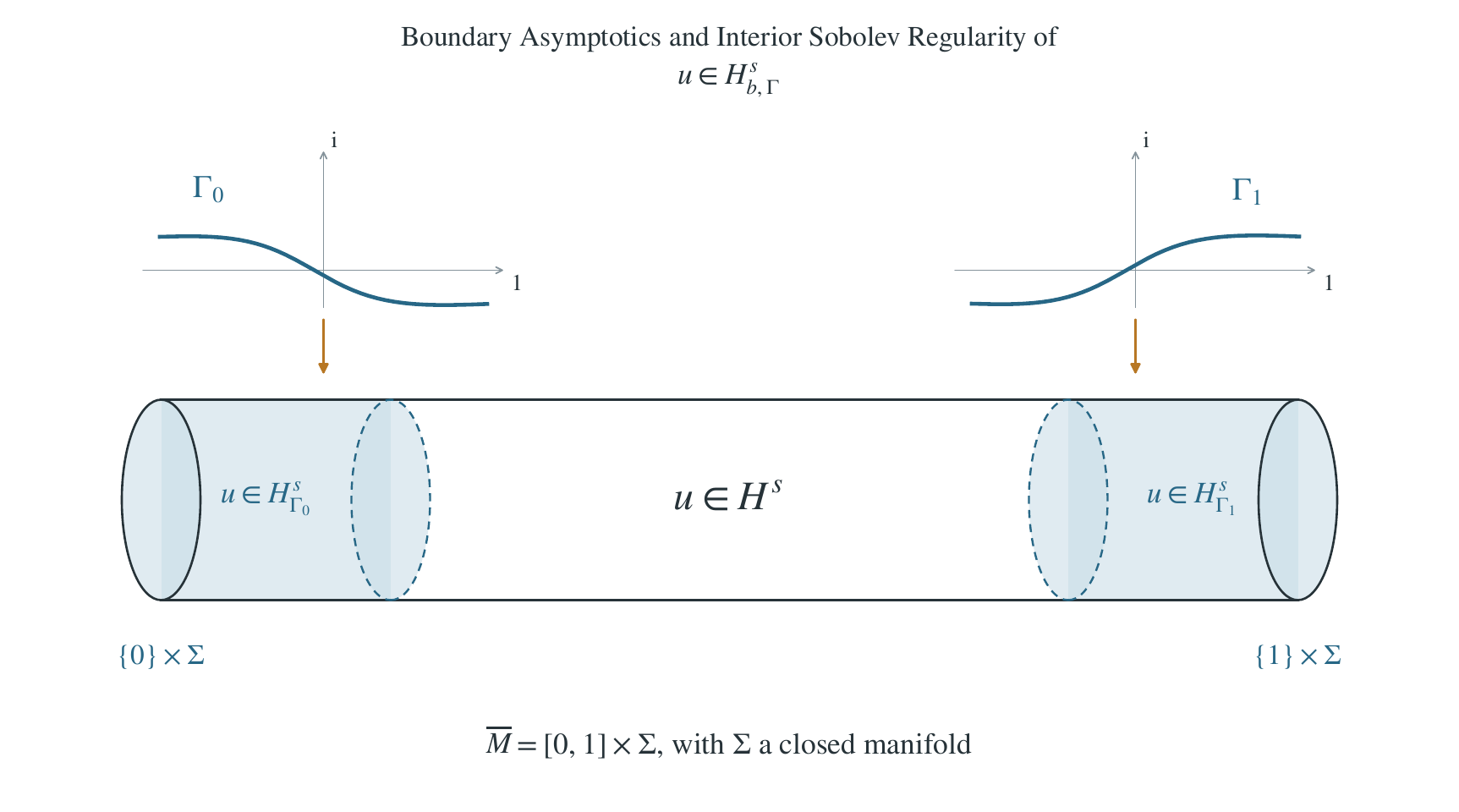}
    \caption{An illustration of how we define members of $H^s_{b, \Gamma}$ in the case of $\Bar{M} = [0, 1] \times \mathbb{S}^1$. In the interior, these distribution look like the usual $H^s$ distributions, while at the boundary components they have asymptotics given by the contours $\Gamma_i$. }
    \label{fig:boundary_asymptotics}
\end{figure}

\begin{definition} 
\label{def:b_WeightedSobolevSpace}
Let $\Bar{M}$ be a compact manifold with boundary and $\partial \Bar{M}_i$ the connected components of its boundary equipped with collar neighborhoods $N_i \simeq [0, b)_{\tau_i} \times \partial \Bar{M}_i$. Pick a $b-$density $\sim \frac{d\tau}{\tau} dy$. Let $V$ be a smooth vector bundle equipped with a positive-definite metric and trivializations of the normal direction \footnote{That is, an identification of sections $u$ with functions $\tau \rightarrow u(\tau, \cdot)$ valued in sections on $\partial \Bar{M}_i$. This comes directly from the product structure $N_i \simeq [0, b_i) \times \partial \Bar{M}_i$.} in each collar neighborhoods $N_i$. We define the $b-$Sobolev space weighted by the contours $\Gamma = \{\Gamma_i\}_i$ to be the following subspace of distributions on rapidly decaying smooth sections of $V$:

\begin{equation}
    H^s_{b, \Gamma}(\Bar{M}; V) = \left\{ u \in \mathcal{C}^{-\infty}(\Bar{M}; V); u|_{N_i} \in H^s_{N_i, \Gamma_i} \text{ and } u \in H^s_{loc}(M) \right\} \subset \mathcal{C}^{-\infty}(\Bar{M}; V)
\end{equation}
Given a finite partition of unity $\{ \chi_i \} \cup \chi$ subordinate to a cover $\{N_i\} \cup U$ satisfying $N_i \cap N_j = \empty$ for $i \neq j$, $U$ some some open set such that $\Bar{U}$ is disjoint from $\partial M$ and $\chi_j$ equal to 1 on a neighborhood of $\partial M_i$, we define the norm on this space as

\begin{equation}
    \| u \|_{H^s_{b, \Gamma}} = \sum_i \| \chi_i u \|_{H^s_{N_i, \Gamma_i}} + \| \chi u \|_{H^s}
\end{equation}
Clearly, by proposition $\ref{prop:ElementaryProperties_Hs}$, it contains rapidly decaying smooth sections $\dot{\mathcal{C}}^{\infty}(\Bar{M}; V)$ as a dense subset. Note that in a trivialization for $V$ near the boundary components $\partial \Bar{M}_i$, elements of $H^s_{b, \Gamma}(\Bar{M}; V)$ have components in the scalar spaces $H^s_{b, \Gamma}(\Bar{M}; \mathbb{R})$.
\end{definition}
\textbf{Remarks:} Note that the definition of this space is not a typical graphical space $\Op(b)^{-1} H^s_b(\Bar{M})$. Indeed, while the operators $T_i$ defined in \ref{prop:IsoContourSobolev} act as local isomorphisms from $H^s_{N_i, \Gamma_i}$ to $H^s_b(N_i)$, they cannot be combined together to yield a global isomorphism with $H^s_b(\Bar{M})$ since they are \textit{global} in time, even up to smoothing terms. Furthermore, since the underlying bundle $V$ will always be implicit, we denote these spaces by $H^s_{b, \Gamma}$ in what follows and write as if dealing with functions $V = \mathbb{R}$. All results apply to any vector bundle $V$. The main distinction from scalar functions to vector bundle $V$ will actually be the choice of positive-definite metric, which in the Lorentzian setting are not canonical.
\vspace{0.2in}
\newline
By $u|_{N_i}$ we of course really mean the pullback to $[0, b_i) \times \partial \Bar{M}_i$ as a distribution on rapidly decaying functions. By lemma \ref{lemma:CutoffContinuous} $\chi_i u \in H^s_{N_i, \Gamma_i}$ since $\chi_i$ equals $1$ near $\tau_i = 0$, hence the norm is well defined. In words, this space is $H^s$ in the interior while near the $i^{th}$ boundary component it is $H^s_{\Gamma_i}$ with decay prescribed by the contour $\Gamma_i$. For this definition to be reasonable, we would first like to say that on our previous setting of $\mathbb{R} \times \Sigma$, this space agrees with the global definition of $H^s_{\Gamma}$ we gave in \ref{def:WeightedSobolevSpacesStationary}. This is indeed the case, precisely because we can localize members of $H^s_{\Gamma}$.

\begin{proposition}
    Let $\Bar{M} = [0, 1] \times \Sigma$ be the exponential compactification of $\mathbb{R} \times \Sigma$ with boundary defining functions $\tau_0 = e^{t}$ at $t = -\infty$ and $\tau_1 = e^{-t}$ at $t = \infty$. Let $-\Gamma_0 = \Gamma_1 = \Gamma$, where $\Gamma$ is some contour in $\mathbb{C}$. Then $H^s_{b, \{ \Gamma_0, \Gamma_1 \}}(\Bar{M})$ as defined above in definition \ref{def:b_WeightedSobolevSpace} is isomorphic to the global $H^s_{\Gamma}(\mathbb{R} \times \Sigma)$ defined in \ref{def:WeightedSobolevSpacesStationary}.
\end{proposition}
\textbf{Proof:} This follows mostly follows from lemma \ref{lemma:CutoffContinuous} and proposition \ref{prop:ElementaryProperties_Hs}. Indeed, cover the real line $\mathbb{R}$ as
\begin{equation}
    \mathbb{R} = (-\infty, -1) \cup (-2, 2) \cup (1, \infty)
\end{equation}
and pick a partition of unity $\chi_{-\infty}, \chi, \chi_{\infty}$ subordinate to this cover, with $\chi_{\pm\infty}(t)$ equal to 1 for $\pm t$ large. Then by lemma \ref{lemma:CutoffContinuous} and proposition \ref{prop:ElementaryProperties_Hs} for any $u \in H^s_{\Gamma}(\mathbb{R} \times \Sigma)$ we have

\begin{equation}
    \chi_{-\infty}  u \in H^s_{\Gamma}(\mathbb{R} \times \Sigma) \cap \mathcal{D}'_{\Gamma}((-\infty, -1)) \quad \chi_{+\infty}  u \in H^s_{\Gamma}(\mathbb{R} \times \Sigma) \cap \mathcal{D}'_{\Gamma}((1, \infty)) \quad \chi u \in H^s((-2, 2) \times \Sigma)
\end{equation}
and conversely if each terms is in its corresponding space than $u \in H^s_{\Gamma}$. In terms of the compactified spacetime $\Bar{M}$ with boundary defining functions $e^{\pm t}$, this cover corresponds to
\begin{equation}
    [0, 1] = [0, e^{-1}) \cup U \cup (1-e^{-1}, 1] \quad U \subset (0, 1)
\end{equation}
and under the change of variable $\tau_0 = e^{t}$, $\tau_1 = e^{-t}$ in the normal neighborhood $[0, e^{-1})$, $(1 - e^{-1}, 1]$ we have
\begin{equation}
    \chi_{-\infty} u \in H^s_{[0, e^{-1}), -\Gamma} \quad \chi_{+\infty} u \in H^s_{(1 - e^{-1}, 1], \Gamma} \quad \chi u \in H^s(U \times \Sigma)
\end{equation}
Thus, the pullback map under $\mathbb{R} \times \Sigma \hookrightarrow \Bar{M}$ gives an isomorphism between the two spaces.
$\Box$
\vspace{0.2in}
\newline
For general contours $\Gamma_i$, these spaces are quite easy to construct. Given a partition of unity $\chi_i, \chi$ subordinate to a cover $\{ N_i \} \cup U$ of normal neighborhoods as above, one can simply define distributions $u$ by  elements of 

\begin{equation}
    H^s_{b, \Gamma} = \sum_i \chi_i H^s_{N_i, \Gamma_i} + \chi H^s_{loc}(M) \subset C^{-\infty}(\Bar{M})
\end{equation}
As before, we wish to show that this space satisfy some reasonable properties. In this manifold setting, we'd like to check the following:

\begin{itemize}
    \item $H^s_{N_i, \Gamma_i}$ does not depend on the choices we've made too heavily, i.e. on the choice of collar neighborhoods $N_i$, their product structures and on the choice of partition of unity.
    \item They obey natural inclusions between each other based on their regularity order $s$ and decay $\Gamma$.
    \item A class of pseudo-differential operators act on them, which we will use to solve time-dependent PDEs.
    \item We have analogs of elliptic estimates, propagation of singularity and subprincipal-type estimates on them.
\end{itemize}
We conclude this subsection by showing the independence of $H^s_{b, \Gamma}$ on our choices and the basic inclusion they satisfy between each other. The first statement is rather simple; the spaces do not depend on the choices we made.

\begin{proposition}
    Suppose $|g_i'(\sigma)| \lesssim \left<\sigma\right>^{-1}$ for all contours $\Gamma_i$ in $\Gamma$. Then the choice of collar neighborhoods $N_i$ with product structure $[0, b_i)_{\tau_i} \times \partial \Bar{M}_i$ does not matter in the definition of $H^s_{b, \Gamma}$ as long as we preserve the curves of $-\tau \partial_{\tau} \sim \partial_t$. More precisely, for any other choices of collar neighborhoods $\Tilde{N}_i \simeq [0, \Tilde{b}_{\Tilde{\tau}_i}) \times \partial \Bar{M}_i$ related by a coordinate transformations of the form
    \begin{equation}
        \varphi(\tau, x) = (a(\tau, x) \tau, y(x)) \quad a > 0, y \text{ some diffeomorphism of $\Sigma$}
    \end{equation}
    we have that $u|_{\Tilde{N}_i} \in H^s_{\Tilde{N}_i, \Gamma_i}$ with comparable norms. Similarly, the norm on $H^s_{b, \Gamma}$ is independent of the choice of partition of unity, in that the norms are equivalent for some other choice of partition of unity $\Tilde{\chi}_j$ subordinate to $\{ \Tilde{N}_i \} \cup \Tilde{U}$ and equal to 1 in a neighborhood of $\Tilde{N}_i$:
    \begin{equation}
        \sum_i \| \Tilde{\chi}_i u \|_{H^s_{\Tilde{N}_i, \Gamma_i}} + \| \Tilde{\chi} u \|_{H^s} \simeq \sum_i \| \chi_i u \|_{H^s_{N_i, \Gamma_i}} + \| \chi u \|_{H^s}
    \end{equation}
\end{proposition}
\textbf{Proof:} Let us start by reducing to the case where $N_i = \Tilde{N_i}$ so that we only have to worry about the choice of product structure. Given two neighborhoods $N_i$ and $\Tilde{N}_i$, we show that $u|_{N_i} \in H^s_{N_i, \Gamma_i}$ iff $u|_{N_i \cap \Tilde{N}_i} \in H^s_{N_i \cap \Tilde{N}_i, \Gamma_i}$ where $N_i \cap \Tilde{N}_i$ is thought of as an open set in $[0, b_i) \times \partial M_i$, and similarly we show $u|_{\Tilde{N}_i} \in H^s_{\Tilde{N}_i, \Gamma_i}$ iff $u|_{N_i \cap \Tilde{N}_i} \in H^s_{N_i \cap \Tilde{N}_i, \Gamma_i}$ when $N_i \cap \Tilde{N}_i$ is thought of as an open set in $[0, \Tilde{b}_i) \times \partial M_i$. This comes down to the following statement:

\begin{lemma}
    Let $u\in \mathcal{C}^{-\infty}$ be a distribution on smooth rapidly decaying functions in time with support in $[0, b) \times \Sigma$ and let $U \subset [0, b) \times \Sigma$ be an open subset containing $\{0\} \times \Sigma$. Suppose that $u \in H^s_{loc}((0, \infty) \times \Sigma)$. Then $u \in H^s_{\Gamma}$ iff $u|_{U} \in H^s_{\Gamma}$.
\end{lemma}
\textbf{Proof}: Given $\chi$ a smooth cutoff function supported in $U$ and equal to $1$ on a neighborhood of the boundary $\{ 0 \} \times \Sigma$, split $u$ as $\chi u + (1 - \chi)u$. Because $(1 - \chi) u$ is supported away from the boundary and $u \in H^s_{loc}((0, \infty) \times \Sigma)$, $(1 - \chi) u$ is in $H^s$ and compactly supported hence is in $H^s_{\Gamma}$ by proposition \ref{prop:ElementaryProperties_Hs}. It follows that $u \in H^s_{\Gamma}$ iff $\chi u \in H^s_{\Gamma}$. That last statement is true iff $u|_{U} \in H^s_{\Gamma}$, since the two differ by a distribution supported away from the boundary hence in $H^s_{\Gamma}$ too.
$\Box$
\vspace{0.2in}
\newline
Thus, it suffices to consider the case $N = \Tilde{N}$ with different product structures. Let $\varphi: N \rightarrow \Tilde{N}$ be the transition map between the two collar neighborhoods. It is of the form

\begin{equation}
    \varphi(\tau, x) = (a(\tau, x) \tau, y(x))
\end{equation}
with both $a$ and $y$ smooth. We must show that pullback under this map gives a continuous map $\varphi^{\ast}: H^s_{\Tilde{N}, \Gamma} \rightarrow H^s_{N, \Gamma}$. Clearly, the pullback under a diffeomorphism of the form $(a, \Phi)$, $(\tau, x) \rightarrow (a(x)\tau, \Phi(x))$, is a continuous map on our spaces since the Fourier-Laplace transform changes as $\widehat{\psi}(\sigma) \rightarrow a^{-i\sigma} \Phi^{\ast} \widehat{\psi}(\sigma)$, which is a continuous map on high-energy spaces $H^s_{\Sigma, \sigma}$. We can thus compose by the inverse of $(a(0, x), y(x))$ and reduce to studying diffeomorphisms of the form $\varphi(\tau, x) = (a(\tau, x) \tau, x)$ with $a(\tau, x) = 1 + O(\tau)$. Let $\psi$ be a rapidly decaying function supported near $\tau = 0$. Then $\varphi^{\ast} \psi(\tau, x) = \psi(a(\tau, x)\tau, x)$ and we need to estimate its Fourier-Laplace transform on $\Gamma$:
\begin{equation}
    \widehat{\varphi^{\ast} \psi}(\phi(\sigma), x) = \frac{1}{\sqrt{2\pi}} \int \psi(a(\tau, x) \tau, x) \tau^{i \phi(\sigma)}\frac{d\tau}{\tau}
\end{equation}
Consider the family $\psi_{\lambda}$ given by $\psi(a(\tau,x)^{\lambda} \tau, x)$. We have
\begin{equation}
\begin{aligned}
    \partial_{\lambda} \widehat{\psi_{\lambda}} &= \frac{1}{\sqrt{2\pi}} \int q_{\lambda} \tau \partial_{\tau}(\psi(a(\tau, x)^{\lambda} \tau, x)) \tau^{i \phi(\sigma)} \frac{d\tau}{\tau} \quad q_{\lambda} = \frac{\log(a)}{1 + \lambda \tau \partial_{\tau} \log(a)} \\
        &=-\frac{1}{\sqrt{2\pi}} \int \psi(a(\tau, x)^{\lambda} \tau, x) \tau \partial_{\tau}(q_{\lambda} \tau^{i\phi(\sigma)})  \frac{d\tau}{\tau}
\end{aligned}
\end{equation}
We can write this as convolution along $\Gamma - \Gamma$:
\begin{equation}
    \partial_{\lambda} \widehat{\psi_{\lambda}} = - \frac{i}{\sqrt{2\pi}} \int_{\mathbb{R}} \phi(\eta) \widehat{q_{\lambda}}(\phi(\sigma) - \phi(\eta)) \widehat{\psi_{\lambda}}(\eta, x) \phi'(\eta) d\eta
\end{equation}
 We have that $\log(a)$ with $a = 1 + O(\tau)$ is a smooth $O(\tau)$ function, hence $\widehat{q_{\lambda}}$ exists on $\Gamma - \Gamma$ by considering its finite Taylor series in $\tau$ up to a rapidly decaying term (the powers $\tau^k$ then have meromorphic Mellin transforms with poles located above $\Gamma - \Gamma$ hence existing on $\Gamma - \Gamma$; a similar argument will be explained in more details in \ref{prop:PseudoDiffActContinuously}). Write the kernel as
 \begin{equation}
 \begin{aligned}
    & \frac{-i}{\sqrt{2\pi}} \phi(\eta) \phi'(\eta) \widehat{q_{\lambda}}(\phi(\sigma) - \phi(\eta)) = \frac{-i}{\sqrt{2\pi}} \eta \widehat{q_{\lambda}}(\sigma - \eta) + R(\sigma, \eta)\\
    & R(\sigma, \eta) = \frac{-i}{\sqrt{2\pi}} \{ \eta(\widehat{q_{\lambda}}(\phi(\sigma) - \phi(\eta)) - \widehat{q_{\lambda}}(\sigma - \eta)) + ig(\eta)(1+ig'(\eta))\widehat{q_{\lambda}}(\phi(\sigma) - \phi(\eta)) + i \eta g'(\eta) \widehat{q_{\lambda}}(\phi(\sigma) - \phi(\eta)) \}
\end{aligned}
 \end{equation}
Using that $|g'(\eta)| \lesssim \left<\eta \right>^{-1}$, one can bound $|R(\sigma, \eta)|$ by $\left< \sigma - \eta \right>^{-2}$, hence it defines a bounded map on $L^2(\Gamma)$ by Schur's test. We thus have
\begin{equation}
\begin{aligned}
    & \partial_{\lambda} \| \widehat{\psi_{\lambda}} \|_{L^2(\Gamma)}^2 = \left<(R + R^{\dagger})  \widehat{\psi_{\lambda}},  \widehat{\psi_{\lambda}} \right> + \left< (A + A^{\dagger})\widehat{\psi_{\lambda}},  \widehat{\psi_{\lambda}} \right> \\
    & (A + A^{\dagger}) \widehat{\psi_{\lambda}} = \int_{\mathbb{R}} \frac{i}{\sqrt{2\pi}} (\sigma - \eta) \widehat{q_{\lambda}}(\sigma - \eta) \widehat{\psi_{\lambda}}(\eta, x) d\eta
\end{aligned}
\end{equation}
$A + A^{\dagger}$ is thus also $L^2$ bounded by Schur's test, hence we conclude that
\begin{equation}
    \partial_{\lambda} \| \widehat{\psi_{\lambda}} \|_{L^2(\Gamma)}^2 \lesssim  \| \widehat{\psi_{\lambda}} \|_{L^2(\Gamma)}^2 \implies  \| \widehat{\psi_{\lambda = 1}} \|_{L^2(\Gamma)}^2 \lesssim  \| \widehat{\psi} \|_{L^2(\Gamma)}^2
\end{equation}
Since $\lambda = 1$ corresponds to the pullback of $\psi$, we are done. The statement for $L^2(\Gamma, \left< \sigma \right>^{2s} d\sigma)$ is proved similarly by inserting weights in these estimates (or equivalently doing the same proof for derivatives of $\psi$, which gives the result for integer $s$, and applying interpolation + duality). Finally, for the equivalence of the norms defined by different partitions of unity, given two partitions of unity $\chi_j$ and $\Tilde{\chi}_j$ subordinate to $\{ N_i \} \cup U$, $\{ \Tilde{N}_i \} \cup \Tilde{U}$ and equal to $1$ near $\partial M_i$, by what we showed the norms on $N_i \cap \Tilde{N}_i$ are equivalent (assuming the product structure is related by $\varphi$ as above). By combining this with the fact that multiplication by smooth bounded functions is continuous on $H^s$, multiplication by a smooth compactly supported smooth functions is continuous on $H^s_{\Gamma}$ (again from proposition \ref{prop:PseudoDiffActContinuously}) and the fact that the $H^s$ norm and $H^s_{\Gamma}$ norms on compactly supported functions are equivalent we get the desired bound:

\begin{equation}
\begin{aligned}
    \sum_i \| \chi_i u \|_{H^s_{N_i, \Gamma_i}} + \| \chi u \|_{H^s} &\lesssim \sum_i \sum_j \| \Tilde{\chi}_j \chi_i u \|_{H^s_{N_i, \Gamma_i}} + \| \Tilde{\chi} \chi_i u \|_{H^s_{N_i, \Gamma_i}} + \sum_j \| \Tilde{\chi}_j \chi u \|_{H^s} + \| \Tilde{\chi} \chi u \|_{H^s} \\
    &\lesssim \sum_{i, j} \| \Tilde{\chi}_j \chi_i u \|_{H^s_{\Tilde{N}_j, \Gamma_j}} + \| \Tilde{\chi} \chi_i u \|_{H^s} + \sum_j \| \Tilde{\chi}_j \chi u \|_{H^s_{\Tilde{N}_j, \Gamma_j}} + \| \Tilde{\chi} \chi u \|_{H^s}  \\
    &\lesssim \sum_j \| \Tilde{\chi}_j u \|_{H^s_{\Tilde{N}_j, \Gamma_j}} + \sum_{i, j}\| (1 - \chi_i) \Tilde{\chi}_j u \|_{H^s_{\Tilde{N}_j, \Gamma_j}} + \| (1 - \chi) \Tilde{\chi}_j u \|_{H^s_{\Tilde{N}_j, \Gamma_j}} \\
    &+ \| \chi_i \Tilde{\chi} u \|_{H^s} + \| \chi \Tilde{\chi} u \|_{H^s} \\
    & \lesssim \sum_j \| \Tilde{\chi}_j u \|_{H^s_{\Tilde{N}_j, \Gamma_j}} + \| \Tilde{\chi} u \|_{H^s}
\end{aligned}
\end{equation}
Thus, the two defined norms are equivalent. We note here that what was important was that $\chi_j$ and $\Tilde{\chi}_j$ both had the same asymptotic at $\partial M_j$, in this case that they were both equal to $1$ nearby. 
 $\Box$
\vspace{0.2in}
\newline
Following this proposition, when speaking of $H^s_{b, \Gamma}$, we will always implicitly fix a collar neighborhood structure on the normal neighborhoods $N_i$ up to equivalence under change of coordinates of the form $(a(\tau, x), y(x))$. In the wave setting, this is naturally given by the integral curves of the vector field $\partial_t$. Another crucial property for Fredholm theory is the inclusion of these spaces with each other and the duality between them. For the standard weighted Sobolev spaces $H^{s, l}_b = \tau^l H^s_b$, they are included within each other according to their regularity order $s$ and rate of decay $l$, and the dual of $H^{s, l}_b$ is simply $H^{-s, -l}_b$. Furthermore, and crucially for the Fredholm theory, the inclusion $H^{s, l}_b \hookrightarrow H^{r, k}_b$ is compact for $s>r$, $l > k$. A similar statement holds for the spaces $H^s_{b, \Gamma}$ we defined here, which is going to allow us to obtain Fredholm estimates on these spaces:

\begin{proposition}
\label{prop:InclusionBWeightedSobolev}
    Let $\Gamma_1$ and $\Gamma_2$ be two collections of contours such that $Im(\phi_1(\sigma)) - Im(\phi_2(\sigma)) \geq c > 0$ for all $\sigma \in \mathbb{R}$, some $c \in \mathbb{R}$ and all $\phi_i(\mathbb{R}) \in \Gamma_i$. Then we have a continuous inclusion

    \begin{equation}
        H^s_{b, \Gamma_1} \hookrightarrow H^r_{b, \Gamma_2} \quad \forall s \geq r
    \end{equation}
    Furthermore, if $s > r$, then the inclusion is compact. Finally, $H^{s}_{b, \Gamma}$ is dual to $H^{-s}_{b, \Bar{\Gamma}}$ under the $L^2_b$ pairing on $\Bar{M}$.
\end{proposition}
\textbf{Proof:} This essentially follows from proposition \ref{prop:InclusionStationarySpaces} where we showed this containment for distributions supported in the half-line $[0, \infty)$ in terms of the Fourier-Laplace transform in $t = -\ln(\tau)$. Indeed, fixing collar neighborhoods $N_i \simeq [0, b_i) \times \partial \Bar{M}_i$ and a partition of unity $\chi_i$, $\chi$ subordinate to a cover $\{ N_i \} \cup U$ of $\Bar{M}$, we have by proposition \ref{prop:InclusionStationarySpaces} that for $u \in H^s_{b, \Gamma_1}$
\begin{equation}
    \chi_i u \hookrightarrow H^r_{N_i, (\Gamma_2)_i} \quad \forall i
\end{equation}
which gives a continuous inclusion $H^s_{b, \Gamma_1} \hookrightarrow H^r_{b, \Gamma_2}$. As for its compactness when $s > r$, each inclusion $H^s_{N_i, (\Gamma_1)_i} \hookrightarrow H^r_{N_i, (\Gamma_2)_i}$, $H^s(U) \hookrightarrow H^r(U)$ is then compact by proposition \ref{prop:InclusionStationarySpaces}, hence $H^s_{b, \Gamma_1} \hookrightarrow H^r_{b, \Gamma_2}$ is then also compact. The duality statement follows from the corresponding statement in proposition \ref{prop:ElementaryProperties_Hs} that said that $H^s_{N_i, \Gamma_i}$ is dual to $H^{-s}_{N_i, \Bar{\Gamma_i}}$. $\Box$
\vspace{0.2in}

\subsection{Mapping Properties of Operators on $H^s_{b, \Gamma}$}
We now wish to study which operators $A: \dot{\mathcal{C}}^{\infty}(\Bar{M}; V) \rightarrow \mathcal{C}^{-\infty}(\Bar{M}; W)$ act continuously on these spaces, where $V$ and $W$ are smooth vector bundles on $\Bar{M}$ equipped with positive-definite metrics. Again, we discuss the case $V = W = \mathbb{R}$, with the extension to a general vector bundles implicit. This will be easiest to understand from the properties of their kernel
\begin{equation}
    K_A \in \mathcal{C}^{-\infty}(\Bar{M} \times \Bar{M})
\end{equation}
thought of as a distribution on rapidly decaying smooth functions $\dot{\mathcal{C}}^{\infty}(\Bar{M} \times \Bar{M})$. It turns out that the mapping properties are best stated in terms of the Mellin transform in the complex plane (its 'Fourier-Laplace' transform) of $K_A$. To apply the transform to our study of the mapping properties of operators, we will need to work in a normal neighborhood $N$ of the boundary. So pick normal neighborhoods $N_i \simeq [0, b_i) \times \partial \Bar{M}_i$ for the connected components of the boundary and a cover $\{ N_i \} \cup U$ with a partition of unity $\chi_i^2, \chi^2$ subordinate to it. We split up the operator $A$ as

\begin{equation}
\label{eq:ExpansionA}
    A = \sum_{i, j} \chi_i^2 A \chi_j^2 + \chi^2 A \chi_j^2 + \chi_i^2 A \chi^2
\end{equation}
$A:H^s_{b, \Gamma^1} \rightarrow H^{r}_{b, \Gamma^2}$ will then be bounded if each terms are bounded on their respective spaces: 

\begin{equation}
    \chi_i A \chi_j: H^s_{N_j, \Gamma_j^1} \rightarrow H^{r}_{N_i, \Gamma_i^2}, \quad \chi A \chi_j: H^s_{N_j, \Gamma_j^1} \rightarrow H^{r}_{loc}, \quad \chi_i A \chi: H^s_{loc} \rightarrow H^{r}_{N_i, \Gamma_i^2}
\end{equation}
In practice, the second and third term will be harmless, as they are respectively supported away from the boundary or involve the action of $A$ on a function $\chi v$ supported away from the boundary, which will be rapidly decaying. Thus, we focus on the terms of the form $\chi_i A \chi_j$, which we can study in normal coordinates for $N_i$ and $N_j$. Moving forward in this subsection, we will assume that the kernel $K_A$ is supported in $N_i \times N_j$. We'll understand its mapping property by studying its behavior under conjugation with the isomorphisms $T_i: H^s_{N_i, \Gamma_i} \rightarrow H^s_b$ from proposition \ref{prop:IsoContourSobolev}. By showing that the conjugated operator $T_i \circ \chi_i A \chi_j \circ T_{j}^{-1}$ is a continuous map $H^s_{b} \rightarrow H^{s'}_b$ we'd show that $\chi_i A \chi_j$ is continuous as a map $H^s_{N_j, \Gamma_j} \rightarrow H^{s'}_{N_i, \Gamma_i}$. This can be done by demonstrating that the conjugated kernel has enough decay.

\begin{figure}[h]
    \centering
    \includegraphics[width=0.75\linewidth]{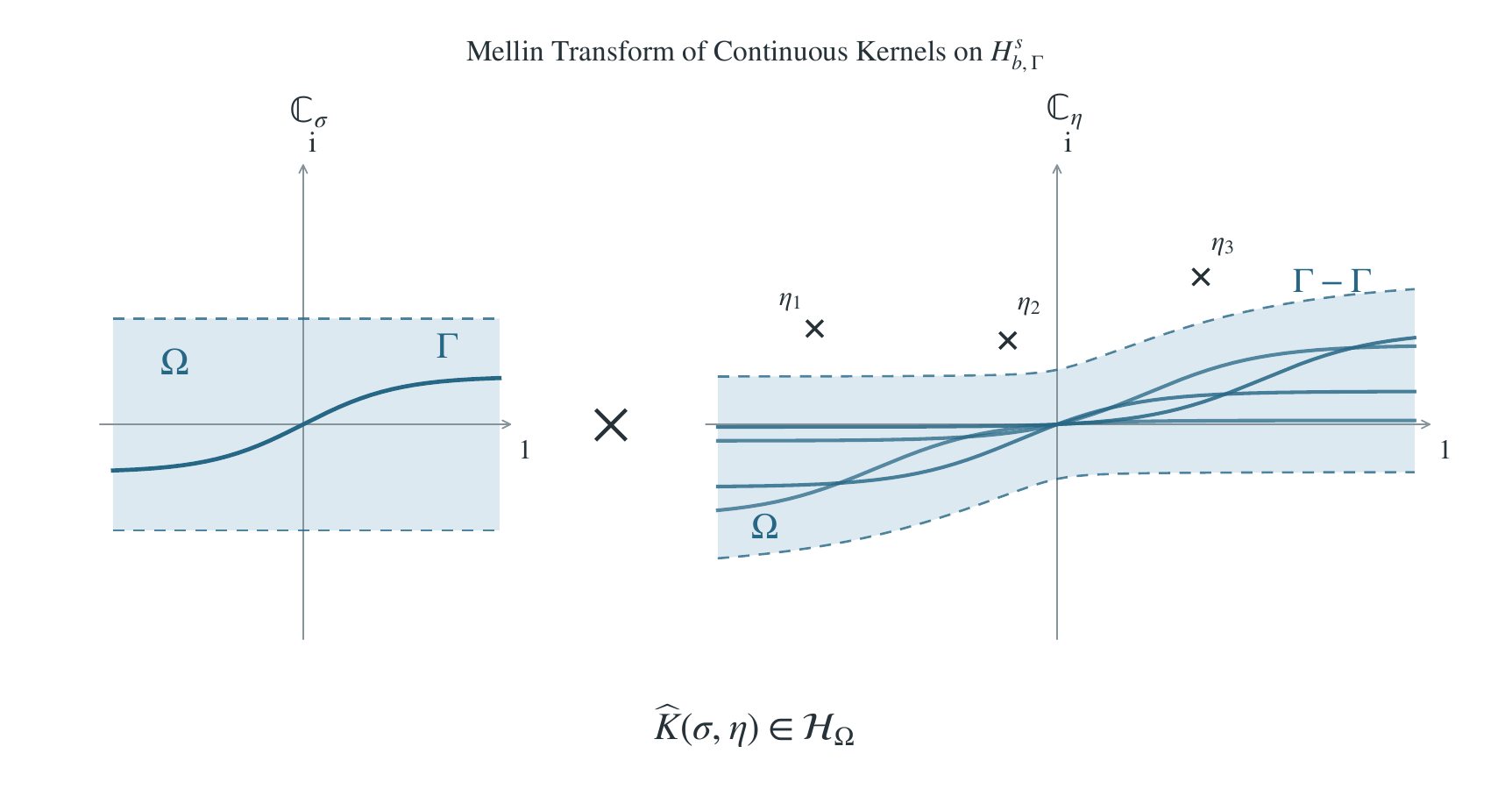}
    \caption{An illustration of the Mellin transform $\widehat{K}(\sigma, \eta)$ in $\kappa = \tau / \rho$ and $\rho$ for a typical continuous kernel acting on $H^s_{b, \Gamma}$. We need for it to make sense on $\Gamma \times (\Gamma - \Gamma)$, which typically means that it is given as a holomorphic function in some region $\Omega$ around it.}
    \label{fig:difference_contour_domains}
\end{figure}

\begin{proposition}
\label{prop:ConjugationFormula}
    Let $K_A$ be a kernel supported in a normal neighborhood $N_i \times N_j$ of $\partial \Bar{M}_i \times \partial \Bar{M}_j$. Suppose that $K_A \in \mathcal{D}'_{\Gamma_i \times -\Gamma_j}$ as a distribution on $(\tau, \tau', x, y) \in N_i \times N_j$ as defined in \ref{def:Distributions_Gamma} and let $T_i$ be the isomorphism $T_i: H^s_{b, \Gamma_i} \rightarrow H^s_b$ from \ref{prop:IsoContourSobolev}. Then the Mellin-transform of the conjugated kernel in $\tau, \tau'$ is given by
    
    \begin{equation}
        \widehat{K}_{T_i \circ A \circ T_j^{-1}}(\sigma, \eta) = \phi_i^{\ast} \times (-\phi_j(-\cdot))^{\ast} \widehat{K}_A \times \phi_j'(\eta)
    \end{equation}
    which we denote by
    
    \begin{equation}
       K_{T_i \circ A \circ T_j^{-1}}(\tau, \rho) \equiv \frac{1}{2\pi} \int_{\mathbb{R}} \int_{\mathbb{R}} \tau^{-i\sigma} \rho^{i\eta} \widehat{K}(\phi_i(\sigma), - \phi_j(\eta)) \phi_j'(\eta) d\sigma d\eta
    \end{equation}
    In words, the complex Mellin transform of the conjugated kernel equals the pullback of $\widehat{K}$ under the map $\phi_i \times (-\phi_j(-\cdot)): \mathbb{R} \times \mathbb{R} \rightarrow \Gamma_i \times (-\Gamma_j)$ times a derivative factor $\phi_j'(\eta)$.
\end{proposition}
\textbf{Proof:}
As usual we will suppress the behavior of tangential variable as the conjugation does not affect them. Let $\chi(\tau), \psi(\rho)$ be rapidly decaying smooth functions. We compute the pairing:

\begin{equation}
\begin{aligned}
    \left< K_{T_i \circ A \circ T_j^{-1}}, \chi(\tau) \psi(\rho) \right> &= \left< K_{A}, T_i^{\dagger} \chi(\tau) T_j^{-1} \psi(\rho) \right> \\
    &= K_A\left(\frac{1}{(2\pi)^2} \int_{\mathbb{R}_{\sigma}} \int_{\tau'} \tau^{i\phi_i(\sigma)} \tau'^{- i \sigma} \chi(\tau') \frac{d\tau'}{\tau'} \int_{\Gamma_{\eta}} \int_{\rho'} \rho^{-i\eta} \rho'^{i \phi_j^{-1}(\eta)} \psi(\rho') \frac{d\rho'}{\rho'} \right) \\
    &= \int_{\mathbb{R}_{\sigma}} \int_{\Gamma_{\eta}} \left< \widehat{K_A}(\phi_i(\sigma),- \eta), \widehat{\chi}(-\sigma) \widehat{\psi}(\phi_j^{-1}(\eta)) \right>
\end{aligned}
\end{equation}
where in the last line we've used the definition of the action of $K_A$ on $\mathcal{C}^{\infty}_{\Gamma_i \times -\Gamma_j}$ in terms of $\widehat{K}_A$. This is the desired expression.
$\Box$
\vspace{0.2in}
\newline
Immediately from this, we see that $\chi_i A \chi_j$ defines a bounded map between $H^s_{N_j, \Gamma^1_j}$ and $H^r_{N_i, \Gamma^2_i}$ when the Mellin transform of $K_A$ in the variables $\tau, \rho$ exists on $\Gamma^1_i \times -\Gamma^2_j$ and has some decay in $\sigma, \eta$. This has to do with the decay (and smoothness) property of the kernel in these variables. Now, in practice, the kernels we will be working with will be concentrated (usually singular) at the diagonal; hence, directly deriving it's properties in terms of $\tau, \rho$ can be awkward. To capture this, we instead rewrite the Mellin transform in $\tau, \rho$ in terms of the Mellin transform in the variable $\kappa = \tau / \rho$ which measures the distance off the diagonal. In that case, we find that

\begin{equation}
    \widehat{K_A}(\phi_i(\sigma), -\phi_j(-\eta)) \times \phi_j'(\eta)= \widehat{K}(\phi_i(\sigma), \phi_i(\sigma) - \phi_j(-\eta)) \phi_j'(\eta) \quad K(\kappa, \rho) = K_A(\kappa \rho, \rho) 
\end{equation}
More precisely, in terms of distributions

\begin{equation}
    \phi_i^{\ast} \times (-\phi_j)^{\ast} \widehat{K}_A \times \phi_j'(\eta) =   \phi_i^{\ast} \times (\phi_i-\phi_j(-\cdot))^{\ast} \widehat{K} \times \phi_j'(\eta)
\end{equation}
i.e. the pullback under the map $\phi_i \times (\phi_i - \phi_j(-\cdot)): \mathbb{R} \times \mathbb{R} \rightarrow \Gamma_i \times (\Gamma_i - \Gamma_j)$. So we need that $\widehat{K}$ is a distribution on $\Gamma_i \times (\Gamma_i - \Gamma_j)$ to make sense of the expression above. One could define general objects which have a Mellin transform transform on such a set, but in practice it would be difficult to check membership in such classes hence we will rather go with the simple approach of requiring that the Mellin transform of $\widehat{K}$ is holomorphic on an open set $\Omega \subset \mathbb{C}^2$ containing $\Gamma_i \times (\Gamma_i - \Gamma_j)$, i.e. $K \in \mathcal{H}_{\Omega}$ as defined in definition \ref{def:Distributions_Gamma}:

\begin{equation}
    \left< K, \psi(s) \chi(\rho) \right> = \int_{\Gamma_i \times \Gamma'} \left< \widehat{K}(\sigma, \eta), \widehat{\psi}(-\sigma) \widehat{\chi}(-\eta) \right> \quad \Gamma_i \times (\Gamma_i - \Gamma_j) \subset \Omega \quad \text{ and $\Gamma' \subset \Gamma_i - \Gamma_j$}
\end{equation}
In terms of decay, this is guaranteed if $K$ is decaying rapidly off the diagonal uniformly over $\rho$, which by proposition \ref{prop:rapid_decay_holomorphicity} distributionally means that:

\begin{equation}
    \kappa^{-i\sigma} K(\kappa, \rho) \text{ is a distribution on smooth functions rapidly decaying in $|\log(\kappa)|$ for $|Im(\sigma)| < C$ }
\end{equation}
where $ C > \sup_{\Gamma_i} |Im(\sigma)|$. The Mellin transform of $K(\kappa, \rho)$ in $\kappa$ will then be a holomorphic function over a half-plane $|Im(\sigma)| < C$ which for $C$ large enough will contain $\Gamma_i$. As for the behavior in the $\rho$ variable, to guarantee existence on $\Gamma_i - \Gamma_j$ we can assume that it has a polyhomogeneous expansion in decaying terms in $\rho$ up to a remainder term of size $O(\rho^{\mu + \epsilon})$ for some small number $\delta>0$ and $\mu_{ij} = \sup_{\sigma \in \Gamma_i - \Gamma_j}|Im(\sigma)|$:

\begin{equation}
    K(\kappa, \rho) = \sum_{l, k; Im(\eta_l) > 0} k_{j, k}(\kappa) \rho^{-i\eta_l} \log(\rho)^k + O_{\kappa}(\rho^{\mu_{ij} + \delta})
\end{equation}
In that case, the Mellin transform in $\rho$ of $K(\kappa, \rho)$ is a meromorphic function in a half plane $Im(\eta) < \mu + \delta$, with poles of order $k+1$ at the exponents $\eta_l$. For this to be holomorphic on $\Gamma_i - \Gamma_j$, we must assume that the exponents $\eta_l$ are above $\Gamma_i - \Gamma_j$ in the complex plane:

\begin{equation}
    Im(\eta_l) + g_j(\sigma) - g_i(\sigma + Re(\eta_l)) > \epsilon > 0 \quad \forall \sigma \in \mathbb{R}
\end{equation}
That assumption may seem strange at first, but makes sense geometrically. The spaces $H^s_{b, \Gamma}$ specify the rate of decay of each real frequencies $Re(\sigma)$ individually, so if one were to multiply by $\tau^{-i\eta_j}$ such that $(\Gamma + \eta_j) \cap \Gamma \neq \emptyset$, which amounts to translating $\Gamma$ by it, we'd have a contradiction between the two rates of decay imposed on $\eta \in (\Gamma + \eta_j)$: one by $\Gamma$ and the other by $(\Gamma + \eta_j)$. Thus, $H^s_{\Gamma + \eta_j}$ may not lie in $H^s_{\Gamma}$. Heuristically, this assumption says that $\Gamma + \eta_l$ is 'a bit more decaying' than $\Gamma$ for all possible real frequencies $Re(\sigma)$, which agrees with our containment statement in \ref{prop:InclusionStationarySpaces}. These properties of the Mellin transform of $K$ in the complex plane is what we will use to derive the continuity of these operators; see figure \ref{fig:difference_contour_domains} for an illustration of them. We can have different variant of the following proposition, depending on which integral test we'd like to use, but they all roughly say the same thing: given strong enough estimates for $\widehat{K}(\sigma, \eta)$ on $\Gamma_i \times (\Gamma_i - \Gamma_j)$ we obtain continuity.
\begin{proposition}
\label{prop:Continuous_Map_Conjugated}
    Let $A$ be an operator $A:\dot{\mathcal{C}}^{\infty}(\Bar{M}; V) \rightarrow \mathcal{C}^{-\infty}(\Bar{M}; W)$ with kernel $K \in \mathcal{H}_{\Omega}$ having holomorphic Melin transform near $\Gamma_i \times (\Gamma_i - \Gamma_j)$ and supported in a normal neighborhood $N_i \times N_j$ of $\partial \Bar{M}_i \times \partial \Bar{M}_j$. Suppose that the Mellin transform of $K$ in $(\kappa, \rho)$ $\widehat{K}$ satisfies for some $N>1 + |s|$

    \begin{equation}
    \label{eq:BoundsKernel}
        \|\widehat{K}(\sigma, \eta)\|_{H^s_{\partial M_j, \sigma} \rightarrow H^r_{\partial M_i, \sigma}} \lesssim \frac{\left< \sigma \right>^{s - r}}{\left< \eta \right>^{N}} \quad \forall (\sigma, \eta) \in \Omega 
    \end{equation}
    In words, for each $\sigma, \eta$, the Mellin transform of the kernel $\widehat{K}(\sigma, \eta)$ defines a bounded map on semiclassical spaces $H^s_{\partial M_j, \sigma} \rightarrow H^r_{\partial M_i, \sigma}$, with norms satisfying the bound \ref{eq:BoundsKernel}. Then $A$ acts continuously on 
    \begin{equation}
        A: H^s_{b, \Gamma} \rightarrow H^{r}_{b, \Gamma}
    \end{equation}
    Note that this is the classical bound for $H^s_b \rightarrow H^r_b$ continuity, except that now we require it near $\Gamma_i \times (\Gamma_i - \Gamma_j)$.
\end{proposition}
\textbf{Proof:} We reduced to showing continuity of some operator $H^s_b([0, \infty) \times \partial M_j) \rightarrow H^r_b([0, \infty) \times \partial M_j)$ by conjugating by the isomorphism $T: H^s_{N_i, \Gamma_i} \rightarrow H^s_b([0, \infty) \times \partial M_i)$. The key point is to show that this conjugation preserves the decay properties of the kernel so that the resulting distribution defines a continuous pairing for $H^r_b \times H^s_b$ functions. From proposition \ref{prop:ConjugationFormula}, the conjugated kernel involves the Mellin transform in $\kappa$ on $\Gamma_i$ and the transform in $\rho$ on $\Gamma_i - \Gamma_j$. Let $\chi$ and $\psi$ be smooth compactly supported functions. The pairing $\left< K_{T_i \circ A \circ T_j^{-1}}, \chi(\tau) \psi(\rho) \right>$ is

\begin{equation}
    \int_{\mathbb{R}_{\sigma}} \int_{\mathbb{R}_{\eta}} \left< \widehat{K}(\phi_i(\sigma), \phi_i(\sigma) - \phi_j(\eta)), \widehat{\chi}(-\sigma) \widehat{\psi}(\eta) \right> \phi_j'(\eta) d\sigma d\eta
\end{equation}
By assumption $\widehat{K}$ satisfies the estimate:

\begin{equation}
\begin{aligned}
    \left| \left< \widehat{K}(\phi_i(\sigma), \phi_i(\sigma) - \phi_j(\eta)), \widehat{\chi}(-\sigma) \widehat{\psi}(\eta) \right> \right| &\lesssim \frac{ \left< \sigma \right>^{s-r} \| \widehat{\chi}(-\sigma) \|_{H^{-r}_{\partial M_i, \sigma}} \| \widehat{\psi}(\eta) \|_{H^s_{\partial M_j}} }{\left< \phi_i(\sigma) - \phi_j(\eta) \right>^{N}}  \\
    &= \left< \sigma \right>^{-r} \| \widehat{\chi}(-\sigma) \|_{H^{-r}_{\partial M_i, \sigma}} \left< \eta \right>^{s} \| \widehat{\psi}(\eta) \|_{H^s_{\partial M_j}} \times \frac{\left<\sigma\right>^s}{\left<\eta\right>^s \left< \phi_i(\sigma) - \phi_j(\eta) \right>^{N}}
\end{aligned}
\end{equation}
Note that for $N > 1 + |s|$, the multiplicative factor is uniformly integrable in each variables, hence by the Schur integral test it defines a $L^2$ pairing between $\left< \sigma \right>^{-r} \| \widehat{\chi}(-\sigma) \|_{H^{-r}_{\partial M_i, \sigma}}$ and $\left< \eta \right>^{s} \| \widehat{\psi}(\eta) \|_{H^s_{\partial M_j}}$. Thus, after integrating we obtain the bound

\begin{equation}
    \left| \left< K_{T_i \circ A \circ T_j^{-1}}, \chi(\tau) \psi(\tau') \right> \right| \lesssim \| \psi \|_{H^s_b} \| \chi \|_{H^{-r}_b}  
\end{equation}
which is the desired continuity statement. 
$\Box$
\vspace{0.2in}
\newline
By applying this theorem to the $\chi_i A \chi_j$ terms in the expansion \ref{eq:ExpansionA}, we can show that a very large class of operators are bounded on $H^s_{b, \Gamma} \rightarrow H^r_{b, \Gamma}$. In particular, all the following operators will define continuous maps:

\begin{itemize}
    \item Multiplication by a smooth function $h$ which is equal to a constant in a neighborhood of each of the boundary components $\partial M_i$. In that case, we can take a partition of unity $\chi_i^2, \chi^2$ which is subordinate to normal neighborhoods $N_i$ on which $h$ is constant. $\chi_i h \chi_j$ are then trivially bounded, while in the interior $h$ defines a bounded map on Sobolev spaces. 
    \item Multiplication by powers $\tau^{-i\eta}$, $\eta \in \mathbb{C}$, defines a bounded map $H^s_{b, \Gamma} \rightarrow H^s_{b, \Gamma + \eta}$ and if $\Gamma_i + \eta$ lies above $\Gamma_i$, i.e. $Im(\phi_i(\sigma) + \eta) - Im(\phi_i(\sigma + Re(\eta))) > \epsilon > 0$ for all $\sigma \in \mathbb{R}$ and $i$, then we have a continuous inclusion $H^s_{\Gamma + \eta} \rightarrow H^s_{\Gamma}$ by proposition \ref{prop:InclusionBWeightedSobolev}. Thus, for such $\eta$, $\tau^{-i\eta}$ defines a bounded map on $H^s_{b, \Gamma} \rightarrow H^s_{b, \Gamma}$.
    \item A class of $b-$pseudodifferential operator will also act on them, as described in section \ref{sec:PSDO_Section}.
\end{itemize}
Note that multiplication by an arbitrary smooth function will also be continuous, as we can expand it into powers of $\tau$, but multiplication by a \textit{conormal} function, say for example smooth in $t = -\log(\tau)$, won't be continuous in general. We now would like to use the above mapping theory to understand the action of \textit{pseudo-differential} operators rather than general kernels. This is essential for us to apply microlocal tools to the Fredholm theory of wave operators. Thus, in the next subsection \ref{sec:PSDO_Section}, we will recall $b-$pseudo-differential operators and then discuss the subclass which acts on our spaces $H^s_{b, \Gamma}$.

\subsection{The Action of b-Pseudo-Differential Operators}
\label{sec:PSDO_Section}
In order to solve PDEs with variable coefficients, we wish to allow the action of $b-$pseudo-differential operators on these spaces. We start with a very short introduction of what they are than discuss the subclass which will act on the $H^s_{b, \Gamma}$ spaces. Within the interior, these are regular pseudo-differential operators but they differ when acting near the boundary $\partial \Bar{M}$. We start by recalling the definition in charts associated to normal neighborhoods $N_i \simeq [0, b_i) \times \partial \Bar{M}_i$ of the boundary and then give the global coordinate free description in the manner of Melrose \cite{Melrose1993} Chapter 4.  

\begin{definition}(Local Definition of b-Pseudo-Differential Operators) Suppose $a(\tau, \sigma, x,\xi) \in S^m$ is a symbol of order $m$ on some chart $[0, \infty) \times \mathbb{R} \times \mathbb{R}^{2n}$, meaning that it satisfies the following bounds:

\begin{equation}
    |(\tau \partial_{\tau})^i \partial_{x}^j \partial_{\sigma}^k \partial_{\xi}^l a(\tau, \sigma, x, \xi)| \lesssim \left< (\sigma, \xi) \right>^{m - k - l} 
\end{equation}
Then in a chart $[0, \infty) \times \mathbb{R}^n$ the action of $A = \Op(a)$ on rapidly decaying function in $\tau$ is given by

\begin{equation}
    Av = \frac{1}{(\sqrt{2\pi})^{(n+1)}} \int_{\mathbb{R}} \int_{\mathbb{R}^n} \tau^{-i\sigma} e^{i x \cdot \xi} a(\tau, \sigma, x, \xi) \widehat{v}(\sigma, \xi) d\xi d\sigma + Rv
\end{equation}
Where $\widehat{v}(\sigma, \xi)$ denotes the Mellin transform in $\tau$ and Fourier transform in $x$ while $R$ is some smoothing operator. We will say $A$ is \textit{stationary} if $a$ is independent of $\tau$ in some product neighborhood of the boundary $[0, \infty) \times \Sigma$, i.e. $a(\tau, x, \sigma, \xi) = a(0, x, \sigma, \xi)$. We will denote this class of operators by $A \in \Psi^m$, and the stationary ones by $\Psi^m_0$. Under $-\log(\tau) = t$, it corresponds to the usual Hormander class of pseudo-differential operators on $\mathbb{R} \times \mathbb{R}^{n-1}$. Similarly, considering matrix-valued symbols $a$ leads to pseudo-differential operators $\Psi^m$ acting on vector bundles.
\end{definition}
Equivalently, after fixing a product structure, one could think of them as 1-parameter family of pseudo-differential operators on $\Sigma \simeq \partial \Bar{M}$ which we might do so when convenient, i.e.
\begin{equation}
    Av = \frac{1}{\sqrt{2\pi}}\int_{\mathbb{R}} \tau^{-i\sigma} a(\tau, \sigma) \widehat{v}(\sigma) d\sigma
\end{equation}
where $a(\tau, \sigma)$ takes values in (semiclassical) pseudo-differential operators on $\Sigma$ with symbolic estimates on its operator norm \footnote{Semiclassical PSDOs $\Psi^m_{\sigma}(\Sigma)$ are simply the usual PSDOs but defined with respect to the semiclassical Fourier transform instead, i.e. scaling the frequency by $Re(\sigma)^{-1}$: $e^{-i\left< \xi, x\right>} \rightarrow e^{-i\left< \xi, x \right> / Re(\sigma)}$.} :
\begin{equation}
    a(\tau, \sigma) \in \Psi^m_{\sigma}(\Sigma)
\end{equation}
In order to extend it to $H^s_{\Gamma}$ one needs to show estimates on the integral kernel of the operator, which is the distribution

\begin{equation}
    K_A(\tau, \rho, x, y) = \frac{1}{(2\pi)^{n+1}}\int_{\mathbb{R}} \int_{\mathbb{R}} \left( \frac{\rho}{\tau} \right)^{i\sigma} e^{i (x-y) \cdot \xi} a(\tau, \sigma, x, \xi) d\xi d\sigma
\end{equation}
In practice, one wishes for these operators to act on functions with finite decay $O(\tau^l)$. As it stands, this kernel does not due to only have $\log(\tau)$ decay off the diagonal $\tau / \rho = 1, x = y'$. In order to impose rapid decay of the kernel off the diagonal, we can manually multiply it by a cutoff $h(\tau/\rho)$ and instead define the action by the kernel:

\begin{equation}
    K_A(\tau, \rho, x, y) = \frac{1}{(2\pi)^{n+1}} \int_{\mathbb{R}} \int_{\mathbb{R}} \left( \frac{\rho}{\tau} \right)^{i\sigma} e^{i(x - y) \cdot \xi} a(\tau, \sigma, x, \xi) h(\tau / \rho) d\xi d\sigma
\end{equation}
The 'true' symbol of this operator, which is the Mellin transform in the $\tau/ \rho$ variable, is thus the convolution of $a$ by the Mellin transform of $h$:
\begin{equation}
    \Tilde{a}(\tau, \sigma, x, \xi) = \int_{\mathbb{R}} \widehat{h}(\sigma - \eta) a(\tau, \eta, x, \xi) d\eta
\end{equation}
It has an holomorphic extension to $\sigma \in \mathbb{C}$ hence makes $\Tilde{a}$ into an entire function. It is convenient to point out for later on that this microlocally preserves the order of $a$ in finite horizontal strip $Im(\sigma) \leq L$. Indeed, if $a$ is globally of order $M$ but of order $m$ in some open cone
\begin{equation}
    C = \left \{ (\tau, \sigma, x, \xi); (\tau, x, (\sigma, \xi) / \left< (\sigma, \xi) \right>) \in V \subset [0, \infty) \times \mathbb{R} \times \mathbb{S}^{n+1} \right\}
\end{equation}
then so is $\Tilde{a}$ for $(\tau, Re(\sigma), x, \xi)$ in a slightly smaller cone $\Tilde{C}$. Indeed, for small enough $\delta>0$ we can split the integral as

\begin{equation}
    \Tilde{a}(\tau, \sigma, x, \xi) = \int_{|Re(\sigma) - \eta| \geq \delta \left<(\sigma, \xi)\right>} \widehat{h}(\sigma - \eta) a(\tau, \eta, x, \xi) dr + \int_{|Re(\sigma) - \eta| \leq \delta \left<(\sigma, \xi)\right>} \widehat{h}(\sigma - \eta) a(\tau, \eta, x, \xi) dr
\end{equation}
The first term is $O(\left<(Re(\sigma), \xi)\right>^{-N})$ for any $N$ by using the rapid decay of $\widehat{h}$ and the global bound on $a$, while the second term is $O(\left<(Re(\sigma), \xi)\right>^{m})$ in a smaller cone $\Tilde{C}$ by using the estimate of $a$ within $C$. The proof for the derivatives follow in the exact same fashion \footnote{Even more is true: this convolution preserves the symbol up to smoothing terms, in the sense that $\Tilde{a}(\tau, \sigma, x, \xi) - a(\tau, \sigma, x, \xi) = O(\left<(\sigma, \xi)\right>^{-N})$ for all $N$ on $\sigma \in \mathbb{R}$, hence is harmless for all microlocal purposes. The proof is quite similar.}. Thus, up to smoothing operators, one can always assume our symbols are holomorphic with microlocal estimates uniform in horizontal strip $Im(\sigma) \leq L$, and that the resulting operator $A$ acts on weighted spaces. Because of this cutoff away from the diagonal, the global definition of $b-$pseudo differential operators on manifold is better stated in terms of the regularity properties of the kernel $K_A$, as exposed in Melrose \cite{Melrose1993} Chapter 4 Section 4.8. We recall here the rigorous definition rapidly and summarize afterwards its important consequences.

\begin{definition}(Global definition of b-Pseudo-differential operators)
    A pseudo-differential operator $A \in \Psi^m_b(\Bar{M})$ of order $m$ (in the small calculus) is an operator $A: \dot{\mathcal{C}}^{\infty}(\Bar{M}) \rightarrow \mathcal{C}^{-\infty}(\Bar{M})$ whose kernel $K_A \in \mathcal{C}^{-\infty}(\Bar{M} \times \Bar{M})$ is a distribution conormal to the $b-diagonal$ of order $m$ and vanishing to infinite order at the right and front faces $\tau = 0 \backslash \{ \tau = \rho \}$, $\rho = 0 \backslash \{ \tau = \rho \}$ :

    \begin{equation}
        \Psi^m_b(\Bar{M}) = \left \{ A: \dot{\mathcal{C}}^{\infty}(\Bar{M}) \rightarrow \mathcal{C}^{-\infty}(\Bar{M}); K_A \in I^m(\Bar{M}_b^2; \Delta^2_b), K_A|_{lf, rf} = 0 \right \}
    \end{equation}
    We denote by $\Psi^m_{b, 0}$ the stationary operators, i.e. the ones defined in some collar neighborhood $N_i \times N_j$ and stationary there. Similarly, for general vector bundles $V$ and $W$, an operator $A: \dot{\mathcal{C}}^{\infty}(\Bar{M}; V) \rightarrow \mathcal{C}^{-\infty}(\Bar{M};W)$ will be said to be in $A \in \Psi^m_b(\Bar{M})$ if its kernel is conormal as above, which is equivalent to asking all matrix entries in a trivialization to be scalar b-pseudo-differential operators of order $m$.
\end{definition}
For a reference on what all of these object means, we refer again to Melrose \cite{Melrose1993} Chapter 4. One can also just think of them as a compatible family of local $b-$pseudodifferential operators in charts as we defined above, with the rapid decay off the diagonal being a global condition. For now, the important properties that this implies are that $A$ acts on weighted b-Sobolev spaces:
\begin{equation}
    A: H^{s, l}_{b}(M) = \tau^l H^s_{b}(\Bar{M}) \rightarrow H^{s-m, l}_b(\Bar{M})
\end{equation}
This can be understood from the fact that rapid decay off the diagonal guarantees the Mellin transform of the kernel $\widehat{K_A}(\sigma, \eta)$ in $\kappa = \tau / \rho$ and $\rho$ to exist in $\mathbb{C} \times \mathbb{R}$, i.e. is entire in $\sigma$, hence is a member of
\begin{equation}
    \widehat{K_A}(\sigma, \eta) \in \mathcal{D}'_{\mathbb{C} \times \mathbb{R}}
\end{equation}
By a strengthening of proposition \ref{prop:Continuous_Map_Conjugated} it then acts on weighted Sobolev space \footnote{This is because in proposition \ref{prop:Continuous_Map_Conjugated} we assumed holomorphicity around the second factor too with some decay in $\eta$, while here we are dealing with the borderline case of having only a distribution on $\mathbb{R}_{\eta}$.}. We are interested in these operators acting on the spaces $H^s_{b, \Gamma}$ we defined. A motivation of why they should do so is that by holomorphicity of the symbol, for any contour $\Gamma$ we can write the action on super-exponentially decaying function as
\begin{equation}
\label{eq:ActionOnGamma}
    Av = \frac{1}{\sqrt{2\pi}} \int_{\Gamma} \tau^{-i\sigma} a(\tau, \sigma) \widehat{v}(\sigma) d\sigma
\end{equation}
This only involves the Mellin transform of $v$ on $\Gamma$, hence should be extendible to other Mellin transform on $\Gamma$ coming from elements of $H^s_{b, \Gamma}$ \footnote{Note that this only need $a(\tau, \sigma)$ to make sense for $\sigma \in \Gamma$, hence a more general class of pseudo-differential operators would also act on these spaces, where we consider symbols $a(\tau, \sigma)$ on $[0, \infty) \times \Gamma$. We won't need this here.}. If $A$ was stationary, then this would just be a Mellin multiplier hence trivially acts on the space. Now if $a(\tau, \sigma)$ depends on $\tau$, then one must check that the result equation \ref{eq:ActionOnGamma} can still be expressed as an integral over $\Gamma$. The easiest way to do so is to conjugate $A$ by the local isomorphism $T_i:H^s_{N_i, \Gamma_i} \rightarrow H^s_b([0, \infty) \times \partial \Bar{M}_i)$ and use proposition \ref{prop:Continuous_Map_Conjugated} to conclude continuity.

\begin{proposition}
\label{prop:PseudoDiffActContinuously}
    Let $A \in \Psi^m_{b}(\Bar{M})$ be a b-pseudo-differential operator of order $m$ which near every components $\partial \Bar{M}_i$ of the boundary has a collar neighborhood $N_i$ of $\partial \Bar{M}_i$ in which it has an expression of the form

    \begin{equation}
    \begin{aligned}
        &\chi_i A \chi_i = \chi_i A_0 \chi_i + \sum_{l, k; Im(\eta_l) > 0} \chi_i \tau^{-i\eta_l} (-\log(\tau))^k A_{lk} \chi_i +  \chi_i \tau^{\mu + \delta} \Psi^{m}_b (\Bar{M}) \chi_i \\
        &\mu = \sup_{\sigma \in \Gamma_i - \Gamma_i, i} |Im(\sigma)|, \quad \delta > 0, \quad A_0, A_{lk} \in \Psi^m_{b, 0}(\Bar{M})
    \end{aligned}
    \end{equation}
    where $\tau$ is a defining function for $\partial \Bar{M}$, $\chi_i(\tau)$ is a smooth cutoff near $\partial \Bar{M}_i$ and the exponents $\eta_l$ are such that $\eta_l$ sit above $\Gamma_i - \Gamma_i$ i.e. satisfy for some $\epsilon$

    \begin{equation}
        Im(\eta_l) + g_i(\sigma) - g_i(Re(\eta_l) + \sigma) > \epsilon > 0 \quad \forall \sigma \in \mathbb{R}
    \end{equation}
    Note that the exponents $\eta_l$ can be different for different components $\partial M_i$. In words, $A$ near the boundary is a sum of stationary pseudo-differential operators plus perturbations that are constrained powers of $\tau$ plus a perturbation of order $O(\tau^{\mu + \delta})$. Then $A$ acts continuously on the spaces $H^s_{b, \Gamma}(\Bar{M})$ and its adjoint $A^{\dagger}$ on $H^s_{b, \Bar{\Gamma}}(\Bar{M})$:

    \begin{equation}
        A: H^s_{b, \Gamma} \rightarrow H^{s-m}_{b, \Gamma} \quad A^{\dagger}: H^s_{b, \Bar{\Gamma}} \rightarrow H^{s-m}_{b, \Bar{\Gamma}}
    \end{equation}
    We will denote by $\Psi^m_{b, \Gamma}(\Bar{M})$ the class of such b-pseudo-differential operators.
\end{proposition}
\textbf{Proof:} Pick a partition of unity $\chi_i^2, \chi^2$ subordinate to collar neighborhoods $N_i \simeq [0, b_i) \times \partial \Bar{M}_i$ and an interior set $U$. We must show that every term in the exansion \ref{eq:ExpansionA} is continuous:
\begin{equation}
    \chi_i A \chi_j: H^s_{N_j, \Gamma_j} \rightarrow H^{s-m}_{N_i, \Gamma_i}, \quad \chi A \chi_j: H^s_{N_j, \Gamma_j} \rightarrow H^{s-m}_{loc}, \quad \chi_i A \chi: H^s_{loc} \rightarrow H^{s-m}_{N_i, \Gamma_i}
\end{equation}
Let's begin with the last two terms. Any $v \in H^s_{N_j, \Gamma_j}$ satisfies $\chi_j v \in H^s_{N_j, Im(\sigma) = -l}$ for a large enough growing weight $l$, and the inclusion is continuous. Then, $A$ maps $H^s_{N_j, Im(\sigma) = -l}$ to $H^{s-m}_{N_j, Im(\sigma) = -l}$, and finally $\chi$, being supported away from the boundary, maps it to $H^{s-m}_{loc}$. Thus $\chi A \chi_j$ is bounded. As for $\chi_i A \chi$, use again that for $v \in H^s_{loc}$, $\chi v$ is compactly supported hence in $H^s_{N_i, Im(\sigma) = l}$ for any $l$ with all norms equivalent. Then, $A$ maps $H^s_{N_i, Im(\sigma) = l}$ to $H^{s-m}_{N_i, Im(\sigma) = l}$ continuously, and finally the continuous inclusion in $H^{s-m}_{N_i, \Gamma_i}$ for large enough $l$ gives the remaining map. So we focus now on $\chi_i A \chi_j$. If $i \neq j$, $\chi_i A \chi_j$ has kernel supported away from the diagonal hence is rapidly decaying at $\tau_i = 0$ and $\tau_j = 0$ since $K_A$ is smooth and rapidly decaying off the diagonal. In particular, $\chi_i A \chi_j$ maps $H^s_{N_j, Im(\sigma) = -l} \rightarrow H^{s-m}_{N_j, Im(\sigma) = l}$ for any $l$, hence $H^s_{N_j, \Gamma_j} \rightarrow H^{s-m}_{N_i, \Gamma_i}$ continuously. Finally, we are left with $\chi_i A \chi_i$, which is where we use our assumptions on the expansion of $A$. By proposition \ref{prop:Continuous_Map_Conjugated}, $\chi_i A \chi_i$ will act continuously if its kernel $K$ has Mellin transform in $(\kappa, \rho)$ defined on $\Gamma_i \times (\Gamma_i - \Gamma_i)$ with appropriate estimates. In our case the Mellin transform of the kernel is given in terms of the Mellin transform of the symbol in $\tau$ as
\begin{equation}
    \widehat{K}(\sigma, \eta) = \widehat{\chi_i a \chi_i}(\eta, \sigma) \quad (\eta, \sigma) \in \mathbb{C}
\end{equation}
This with proposition \ref{prop:Continuous_Map_Conjugated} aldready shows the boudness of the fast decaying term
\begin{equation}
    \chi_i \tau^{\mu + \delta} \Psi^{m}_b (M) \chi_i 
\end{equation}
since the Mellin transform in $\kappa, \tau$ then exists in a strip $\mathbb{C} \times \{Im(\eta) < \mu + \delta \}$, which contains $\Gamma_i \times (\Gamma_i - \Gamma_i)$, and is bounded by $\lesssim \left<\sigma\right>^m \times \left< \eta \right>^{-N}$ for any $N$ in the operator norm of maps $H^s_{\partial \Bar{M}_i, \sigma} \rightarrow H^{s-m}_{\partial \Bar{M}_i, \sigma}$. Indeed, for each fixed $\sigma, \tau$ we have a (semiclassical) pseudo-differential operator
\begin{equation}
    \Op(\widehat{\chi_i a \chi_i}(\eta, \sigma)) \in \Psi^m_{\sigma}(\partial M_i)
\end{equation}
with uniform control on its operator norm in horizontal strips in $\mathbb{C} \times \{Im(\eta) < \mu + \delta \}$ (this follows from the operator norm being controlled by symbolic norms of the total symbol, which are uniform in strips). As for the expansion in powers of $\tau$, the Mellin transform $\widehat{K}(\sigma, \eta)$ is then decaying too slowly to apply $\ref{prop:Continuous_Map_Conjugated}$ directly, since
\begin{equation}
    \widehat{\tau^{-i\eta_l} (-\log(\tau))^k} = \frac{1}{\sqrt{2\pi}}  \frac{k!}{i^{k+1}(\eta - \eta_l)^{k+1}}
\end{equation}
which prevents fast decay in $\eta$. One could do the continuity proof directly by computing the Kernel, this is essentially what we will do when computing the symbol of the conjugated operator in proposition \ref{prop:ConjugatedSymbol}, but here instead we proceed with a simpler argument. First the stationary part $A_{lk}$ is clearly a continuous operator by lemma \ref{lemma:StationaryOperatorContinuous} since they are just Mellin multiplier (and again the uniformity of the operator norms $H^s_{\partial \Bar{M}_i, \sigma} \rightarrow H^{s-m}_{\partial \Bar{M}_i, \sigma}$ follows from the symbolic estimates on the total symbol). So we're left to show that multiplication by $\tau^{-i\eta_l} (-\log(\tau)^k)$ is continuous when the support is restricted to $\tau < 1$. That was shown for $k=0$ in proposition \ref{prop:InclusionStationarySpaces}, while for $k>1$ then we have enough decay to apply proposition \ref{prop:Continuous_Map_Conjugated} in the $s = 0$ case (since we get $\left< \eta \right>^{-2}$ decay), and we can handle any $s$ after by scaling $A$ appropriately. Finally, the statement about the adjoint $A^{\dagger}$ follows from $A \in \Psi^m_{b, \Gamma}$ implying $A^{\dagger} \in \Psi^m_{b, \Bar{\Gamma}}$, due to the condition on $\eta_l$ being symmetric under $\eta_l \rightarrow -\Bar{\eta}_l$.
$\Box$
\vspace{0.2in}
\newline
We thus have quite a large class of pseudo-differential operators acting on our spaces $H^s_{b, \Gamma}$. In practice, one would check membership in $\Psi^m_{b, \Gamma}$ by doing a partial expansion at the boundary $\partial \Bar{M}$ of the total symbol of $A$ up to some term decaying fast enough. Now, in order for this subclass of $b$-pseudo-differential operators to be useful for solving PDEs, it is very desirable that they form an algebra under composition. This is the content of the next lemma.

\begin{lemma}
\label{lemma:AlgebraGamma}
    $\Psi^{\circ}_{b, \Gamma}(\Bar{M})$ forms a graded algebra, that is
    \begin{equation}
        \Psi^m_{b, \Gamma}(\Bar{M}) \circ \Psi^r_{b, \Gamma}(\Bar{M}) \subset \Psi^{m+r}_{b, \Gamma}(\Bar{M})
     \end{equation}
     Furthermore, it is closed under asymptotic summation, that is if we have a collection of operators $A_j \in \Psi^{m_j}_{b, \Gamma}$ sharing the same exponents $\eta_l$ and with $m_j \rightarrow -\infty$ monotonically, then there exists an operator $A \in \Psi^m_{b, \Gamma}$ such that
     \begin{equation}
         A - \sum_{j = 1}^k A_j \in \Psi^{m_k}_{b, \Gamma} \quad \forall k
     \end{equation}
\end{lemma}
\textbf{Proof:} We have to check two things: that the algebra is well behaved under commutators and that the condition on the exponents $\eta_l$ is additive. Let's start with the exponents: if $\eta_l$ and $\eta_k$ satisfy it, then we find $\forall \eta \in \mathbb{R}$
\begin{equation}
    Im(\eta_l) + Im(\eta_k) > g(Re(\eta_l + Re(\eta_k) + \eta) - g(\eta + Re(\eta_k)) + g(Re(\eta_k) + \eta) - g(\eta)
\end{equation}
which after canceling terms is the desired property. As for the commutators, note that for a stationary operator $A_0 \in \Psi^{m}_{b, 0}$ and power $\tau^{-i\eta_l}$ we have
\begin{equation}
    \left[A_0, \tau^{-i\eta_l} \right] = \tau^{-i\eta_l} B_0
\end{equation}
for some stationary operator $B_0 \in \Psi^{m}_{b, 0}$ with symbol $a(0, \sigma + \eta_l) - a(0, \sigma)$. Thus when composing two elements of $\Psi^{\circ}_{b, \Gamma}(M)$ one can simply commute the powers of $\tau$ and the stationary part of the operators to bring it to the required form. As for the closedness under asymptotic summation, it follows from each term in the sum being closed under it. $\Box$
\vspace{0.2in}
\newline
We now have a natural algebra of operators acting on our spaces and it remains to show that under the usual hypothesis (ellipticity or dynamical) we can obtain estimates for them. This will be the topic of the next subsection \ref{subs:EllipticNonEllipticEstimates}. In preparation for that, we need to introduce one last class of pseudo-differential operators. These operators will come up naturally in the proof hence it is useful to gather results about them. They are the \textit{product-type} pseudo-differential operators and are modeled on product operators of the form
\begin{equation}
    K_A = \delta(x - y) \times A_x(\tau, \sigma) \quad A_x(\tau, \sigma) \in \Psi^m_b([0, \infty))
\end{equation}
i.e. a family of pseudo-differential operators on $[0, \infty)$ parametrized by $\Sigma$. They were first discussed in \cite{MelroseRochon2006} for the purpose of defining an algebra of operators adapted to a fibration, which in our case is simply the projection $[0, \infty) \times \Sigma \rightarrow \Sigma$. The symbolic estimate can be found in [Appendix A] while the relevant boundness statement in section 10. Another reference is Melrose's unpublished notes \cite{Melrose2006Lectures} section 18, which we refer for the precise definition of the calculus. We give here for the reader a simplified definition of the calculus:

\begin{definition}
\label{def:ProductTypeOperators}
    We will denote by $\Psi^{(m, m)}([0, \infty) \times \Sigma)$ the calculus of (properly supported) pseudo-differential operators $A$ whose kernel $K_A = K_{A_1} + K_{A_2}$ (modulo schwartz kernel) near the diagonal $\Delta = \{ \tau = \rho, x = y \}$ is given locally by the oscillatory integral

    \begin{equation}
        K_{A_1}(\tau, \rho, x, y) = \frac{1}{(2\pi)^{dim(\Sigma) + 1}} \int \left( \frac{\tau}{\rho} \right)^{-i\sigma} e^{i(x-y)\cdot\xi} a_1(\tau, \rho, \sigma, x, y, \xi) d\sigma d\xi
    \end{equation}
    Where the symbol $a_1$ is assumed to be conormal on the phase space obtained by blowing-up $\sigma = 0$ in frequency space. In particular it satisfies the estimate
    \begin{equation}
    \label{eq:ProductTypeEstimates}
        |(\tau \partial_{\tau})^i (\rho \partial_{\rho})^{i'}  \partial_x^j \partial_{y}^{j'} \partial_{\sigma}^k \partial_{\xi}^l a_1| \lesssim \left< (\sigma, \xi) \right>^{m} \left< \sigma \right>^{-k} \left< (\sigma, \xi) \right>^{-l}
    \end{equation}
    Furthermore, it can admit an additional term near the partial diagonal $\Delta_{\Sigma} = \{ x = y\}$ and away from the full diagonal $\tau = \rho$ given locally by
    \begin{equation}
        K_{A_2}(\tau, \rho, x, y) = \frac{1}{(2\pi)^{dim(\Sigma)}} \int e^{i(x-y)\cdot\xi} a_2(\tau, \rho, x, y, \xi) d\xi
    \end{equation}
    where the symbol $a_2$ is $O(|\log(\tau)|^{-\infty}, |\log(\rho)|^{-\infty})$, conormal in $\tau, \rho$ but symbolic of order $m$ in $(x, y, \xi)$. 
\end{definition}
An equivalent definition can be given more precisely in terms of the singularities of the kernels. While ordinary pseudo-differential are singular only at the diagonal $\tau = \rho, x = y$, such product type pseudo-differential operators are allowed to have singularities on the partial diagonal $\Delta_{\Sigma} = \{ x = y \}$. They can thus be defined as operators whose kernel is conormal in a precise sense to both the diagonal $\Delta$ and the partial one $\Delta_{\Sigma}$:
\begin{equation}
    K_A \in I^{(m, m)}(([0, \infty) \times \Sigma)^2, \Delta, \Delta_{\Sigma})
\end{equation}
On a higher abstraction still, they can be understood as \textit{paired Lagrangian distributions}, in that their kernel $K_A$ is singular at both Lagrangians $N^{\ast} \Delta$ and $N^{\ast}_{\Sigma}$ with a compatibility condition at the intersection. An introduction to those can be found in \cite{dehoop2012diffractionconormalsingularities} section 5. The estimates there correspond exactly to the class of symbols in equation \ref{eq:ProductTypeEstimates}. Anyhow, for our purposes, these arise naturally since after conjugation by the isomorphisms $T_i:H^s_{N_i, \Gamma_i} \rightarrow H^s_b$, we naturally obtain operators in this class (see proposition \ref{prop:ConjugatedSymbol}):
\begin{equation}
    T_i \circ \chi_i \Psi^m_{b, \Gamma} \chi_i \circ T_i^{-1} \subset \Psi^{(m, m)}([0, \infty) \times \Sigma)
\end{equation}
They fail to be standard pseudo-differential operator since $\partial_{\sigma}$ derivatives do not improve in the joint frequency $\left< (\sigma, \xi) \right>$, and in particular should be thought of as global in 'time' $[0, \infty)$. Nonetheless, they behave quite reasonably. The two results we collect about them are that they act continuously on the $b$-Sobolev spaces $H^s_b$ and that after localizing to the elliptic set of $\sigma$, i.e. $\sigma \neq 0$, they reduce to the usual pseudo-differential operators $\Psi^m$.
\begin{lemma}
\label{lemma:MicrolocalizationAgree}
    Let $A \in \Psi^{(m, m)}([0, \infty) \times \Sigma)$. Then $A: H^{s}_b([0, \infty) \times \Sigma) \rightarrow H^{s-m}_b([0, \infty) \times \Sigma)$ continuously. Furthermore, suppose $B \in \Psi^0$ has wavefront set in $\{ \sigma \neq 0 \}$, i.e. $\sigma$ is elliptic on it. Then
    \begin{equation}
        AB, BA \in \Psi^{m}([0, \infty) \times \Sigma)
    \end{equation}
    with wavefront set in $WF(B)$.
\end{lemma}
\textbf{Proof:} The continuity was discussed by R. Melrose and F. Rochon in \cite{MelroseRochon2006} section 10, see equation (10.11) (technically this is in the compact setting with classical symbols but the arguments are the same; see also their Appendix A). As for the wavefront set statement, note that by our assumptions on $WF(B)$, for some $\psi$ we have
\begin{equation}
    AB =  A\Op(\psi(\left< \sigma \right> / \left< (\sigma, \xi) \right>))B + R \quad \text{$R$ smoothing}
\end{equation}
where $\psi$ localizes on $\left< \sigma \right> / \left< (\sigma, \xi) \right> \geq \epsilon$. Then the total symbol of $A \Op(\psi(\left< \sigma \right> / \left< (\sigma, \xi) \right>))$ is just $\psi(\left< \sigma \right> / \left< (\sigma, \xi) \right>) a$, which is a (standard) symbol since within its support $\left<\sigma\right>$ and $\left<(\sigma, \xi)\right>$ are comparable. Furthermore, it trivially has wavefront set contained in $WF(B)$. Similarly for $BA$. 
$\Box$
\vspace{0.2in}
\newline
This last lemma tells us that within $\sigma \neq 0$, the two algebras agree. In particular, $\Psi^{(m, m)}$ becomes \textit{microlocal}, in the sense that it preserves wavefront sets and other microlocal properties. This will be essential to make the commutator estimators go through, which is the topic of the next section \ref{subs:EllipticNonEllipticEstimates}.

\subsection{Elliptic and Non-Elliptic Estimates on $H^s_{b, \Gamma}$}
\label{subs:EllipticNonEllipticEstimates}
In order to obtain Fredholm estimate, we will use a combination of propagation arguments and apriori control. The first step will be to derive microlocal elliptic estimates so that we control what happens in the elliptic set of $P$. In our case, this is rather simple since it amounts to show that elliptic parametrix of operators act on our spaces in which case proposition \ref{prop:PseudoDiffActContinuously} can be used. Thus, the statement end ups very similar to the usual one.  
\begin{proposition}(Elliptic Estimates)
\label{prop:EllipticEstimates}
    Let $P \in \Psi^m_{b, \Gamma}(\Bar{M})$ be a b-pseudo-differential operator acting on vector bundles $V$ and $W$ with expansion as in proposition \ref{prop:PseudoDiffActContinuously} and let $B \in \Psi^0_{b, \Gamma}(\Bar{M})$ be microlocally supported in the elliptic set of $P$, i.e. $WF(B) \subset ell(P)$. Then there exists a microlocalizer $Q \in \Psi^{-m}_{b, \Gamma}(\Bar{M})$ such that $B = QP + R$ with $R \in \Psi^{-\infty}_{b, \Gamma}(\Bar{M})$. Thus, there exists $G \in \Psi^{0}_{b, \Gamma}(\Bar{M})$ such that for any $N$ we have the following elliptic estimate

    \begin{equation}
        \| Bu \|_{H^s_{b, \Gamma}} \lesssim \| G Pu \|_{H^{s-m}_{b, \Gamma}} + \| u \|_{H^{-N}_{b, \Gamma}} 
    \end{equation}
    In particular, $\Psi^{\circ}_{b, \Gamma}(\Bar{M})$ is closed under inversion (up to smoothing) of elliptic operators.
\end{proposition}
\textbf{Proof:} By $WF(B) \subset ell(P)$ we can construct a local parametrix $Q \in \Psi^{-m}_b(\Bar{M})$ for $P$ satisfying
\begin{equation}
    B = QP + R
\end{equation}
where $R \in \Psi^{-\infty}_b(\Bar{M})$. So we only have to check that $Q$ can be taken to lie in $\Psi^{-m}_{b, \Gamma}(\Bar{M})$ and $R$ in $\Psi^{-\infty}_{b, \Gamma}(\Bar{M})$ so that they act on the $H^s_{b, \Gamma}$ spaces. We recall the construction of $Q$ and check that every steps preserve the decay in $\tau$ we need. Let $p$ be the principal symbol of $P$. By assumption, it is of the form
\begin{equation}
    p = p_0 + \sum_{l,k; Im(\eta_l) > 0} \tau^{-i\eta_l} (-\log(\tau))^k p_{lk} + \tau^{\mu + \delta} S^m_b
\end{equation}
We define $Q_1$ by a quantization of $\frac{b}{p}$, where $b$ is the principal symbol of $B$. This is well defined since $b$ vanishes outside the elliptic set of $p$. Note that $\frac{b}{p}$ is then also of this form since the space form an algebra by lemma \ref{lemma:AlgebraGamma} and is closed under inversion of symbols (by a Taylor series argument). For the quantization $\Op(\frac{b}{p})$ to land in $\Psi^{-m}_{b, \Gamma}(\Bar{M})$, one can thus simply work in a chart. Then, 
\begin{equation}
    B = Q_1 P + R_1
\end{equation}
Where $R_1 \in \Psi^{-1}_{b, \Gamma}(\Bar{M})$. As usual to improve this to a smoothing operator, we improve the residual iteratively and then use a Neumann serie to sum them all up. We do the same process with $B - Q_1 P$ to find $Q_2 \in \Psi^{-m-1}_{b, \Gamma}(\Bar{M})$ such that $B - Q_1P = Q_2P + R_2$ with $R_2 \in \Psi^{-2}_{b, \Gamma}(\Bar{M})$. Doing this to all order we obtain a collection $Q_i$, which by proposition \ref{lemma:AlgebraGamma} we can asymptotically sum to obtain $Q$ such that $B = QP + R$, $R \in \Psi^{-\infty}_{b, \Gamma}(\Bar{M})$. $\Box$
\vspace{0.2in}
\newline
This already shows that elliptic operators will be Fredholm on our spaces assuming that their indicial family $\widehat{N(P)}(\sigma)$ is invertible on $\Gamma$ \footnote{This is just the Mellin transformed family of $P_0$, the stationary part of $P$.}. In order to obtain estimates within the characteristic set of $P$, we will use a propagation of singularity statement paired with an estimate close to the boundary. The usual propagation of singularity theorem can be used away from the boundary hence it thus remain to control $u$ by $Pu$ near the boundary ('timelike infinity'), which we'll handle using a sort of radial point-type estimate. We start by giving an heuristic on why we should be able to obtain such an estimate for the wave operator on these spaces, then move on to give the proof. A wave operator has principal symbol $\sim \sigma^2 - \xi^2$. Now roughly speaking, if we choose our contour $\Gamma$ so that $Im(\phi(\sigma)) = g(\sigma)$ is negative for $\sigma \geq 0$ and positive otherwise, i.e. we choose our space to be growing for future directed $\sigma > 0$ frequencies and decaying for past directed frequencies $\sigma < 0$, then the Hamiltonian derivative of the weight $\tau^{-g(\sigma)}$

\begin{equation}
    H_p(\tau^{-g(\sigma)}) = -g(\sigma) \tau^{g(\sigma) - 1}H_p(\tau) - \log (\tau)\tau^{-g(-\sigma)} H_p(g(\sigma)) \sim -2 g(\sigma) \sigma \tau^{-g(\sigma)}
\end{equation}
is a positive weight times $\tau^{-g(\sigma)}$. It thus follows that the commutator

\begin{equation}
    \frac{1}{i} [P, \Op(\tau^{-g(\sigma)})] \sim \Op(-H_p(\tau^{-g(\sigma)})) + R
\end{equation}
has a negative weight in its principal symbol, which can be used in a commutator estimate to control $u$ in terms of $Pu$. Now the reason that this argument can be tricky to make rigorous is that the function $\tau^{g(\sigma)}$ is not symbolic and does not define a standard pseudo-differential operator. It has two separate issues: every $\partial_{\sigma}$ derivative of its symbol introduces a $\log(\tau)$ factor, which is growing at the boundary, and secondly even assuming that $\phi(\sigma)$ is symbolic, i.e. every $\partial_{\sigma}$ improves by $\left< \sigma \right>^{-1}$, the symbol won't be jointly symbolic in $(\sigma, \xi)$. Thus, we will instead work with the conjugated operator. These issues will come up in some way nonetheless, and are respectively handled by the assumed expansion of $P$ and by assuming that $\sigma$ is elliptic in the characteristic set of $P$ \footnote{I am quite confident that a commutator estimate could still be made, but one would need to be quite careful about the errors in $[P, T]$ factoring through $T$.}.

\begin{figure}
    \centering
    \includegraphics[width=0.7\linewidth]{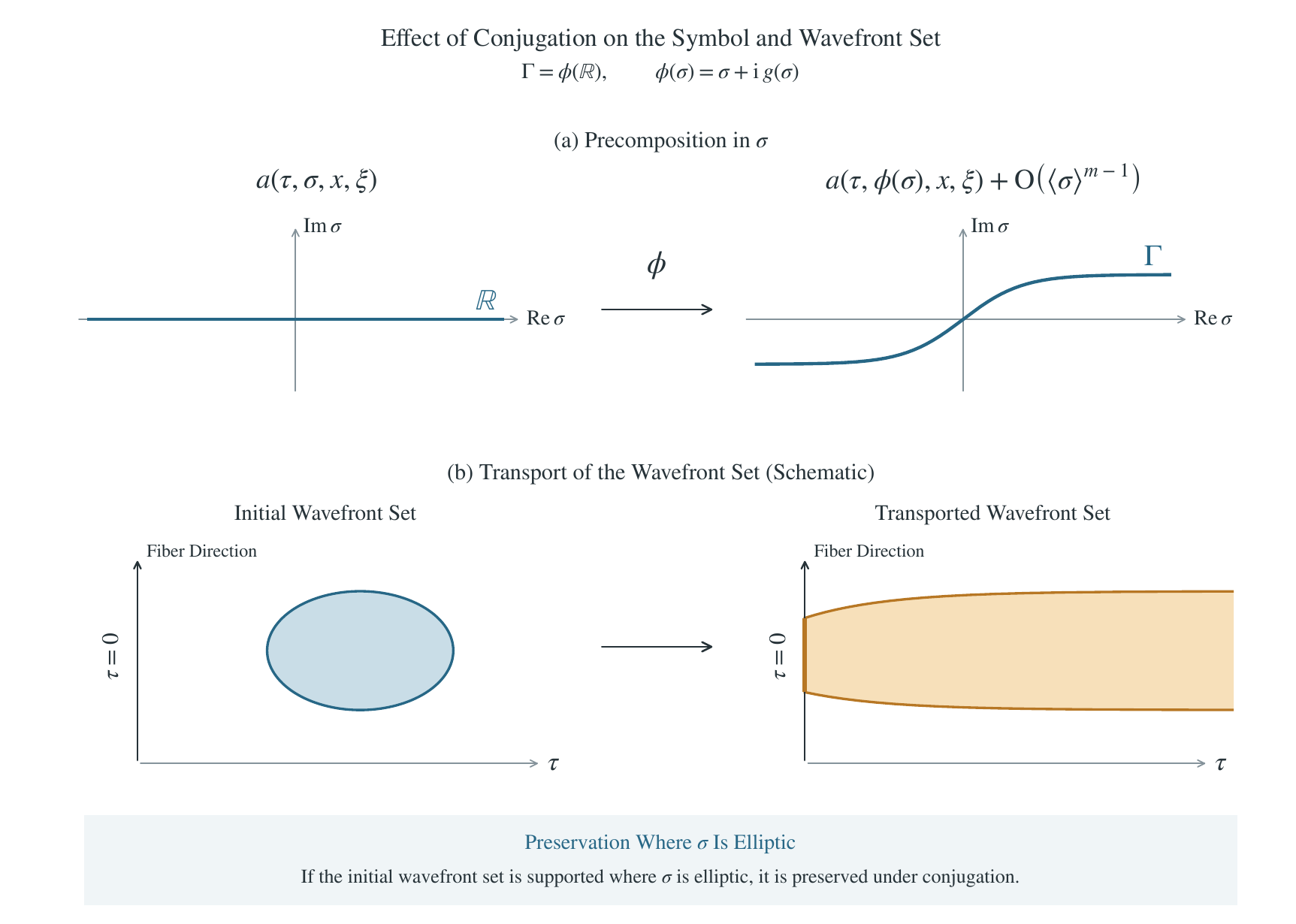}
    \caption{An illustration of the effect of conjugation on pseudo-differential operators. Roughly speaking, it can be expressed as a pull-back by $\phi$ on the frequency $\sigma$ and a deformation in $\tau$ of the wavefront set. Lower terms also arise, which are the expected one from the commutator formula.}
    \label{fig:symbol_conjugation}
\end{figure}

\begin{proposition}(Formula for the Conjugated Operator)
\label{prop:ConjugatedSymbol}
    Let $A \in \Psi^m_{b, \Gamma}(\Bar{M})$ be a pseudo-differential operator acting on vector bundles $V$ and $W$ with expansion as in proposition \ref{prop:PseudoDiffActContinuously}

    \begin{equation}
    \label{eq:AssumedExpansionForConjugation}
        \chi_i A \chi_i = \chi_i A_0 \chi_i + \sum_{l, k; Im(\eta_l) > 0} \chi_i \tau^{-i\eta_l} (-\log(\tau))^k A_{lk} \chi_i +  \chi_i \tau^{\mu + \delta} \Psi^{m}_b (\Bar{M}) \chi_i \quad \mu = \sup_{\Gamma - \Gamma} |Im(\sigma)|
    \end{equation}
    Assume that $g_i(\sigma) = \frac{1}{i} (\phi_i(\sigma) - \sigma)$ are symbols of order $0$ on $\mathbb{R}$. Let $\chi_i^2, \chi^2$ be a partition of unity subordinate to a cover of normal neighborhoods $N_i \simeq [0, b_i) \times \partial \Bar{M}_i$ and $U$, and set $\phi =  \phi_i$. Then the conjugated operator $T_i \circ \chi_i A \chi_i \circ T_i^{-1}$ under the isomorphism $T_i: H^s_{N_i, \Gamma_i} \rightarrow H^s_b$ acting on $H^s_{b}([0, \infty) \times \partial \Bar{M}_i)$ is a pseudo-differential in the class $\Psi^{m, m}$ as defined in \ref{def:ProductTypeOperators} with (right) total symbol in $N_i$ given by

    \begin{equation}
    \begin{aligned}
    \label{eq:conjugated_symbol}
        &\sigma(T_i \circ \chi_i A \chi_i \circ T_i^{-1})(\sigma, \rho) = \frac{1}{\sqrt{2\pi}} \int_{\mathbb{R}_{\vartheta}} \rho^{-i\vartheta} \widehat{\chi_i a \chi_i}(\vartheta + i\vartheta q(\sigma, \vartheta), \phi(\sigma)) \phi'(\sigma - \vartheta) d\vartheta \\
        &=\chi_i a \chi_i(\rho, \phi(\sigma)) -ig'(\sigma) \log(\rho) \rho \partial_{\rho}(\chi_i a \chi_i)(\rho, \phi(\sigma)) + O(\left<\sigma\right>^{m-2}, |\log(\rho)^{-\infty}|)  \\
        &= \chi_i a \chi_i(\rho, \sigma) - i \frac{H_{\chi_i a \chi_i}(\rho)}{\rho} g(\sigma) + O(\left<\sigma\right>^{m-1}, |\log(\rho)|^{-\infty}) + O(\left< \sigma \right>^{m-2})\\
        \newline \\
        &q(\sigma, \vartheta) = \frac{g(\sigma) - g(\sigma - \vartheta)}{\vartheta} = \int_{0}^{1} g'(\sigma - r \vartheta) dr
    \end{aligned}
    \end{equation}
    where $\chi_i a \chi_i$ is the (right) total symbol of $\chi_i A \chi_i$ in $N_i$ and $ O(\left<\sigma\right>^{r}, |\log(\rho)^{-\infty}|)$ denotes symbol of order $r$ (in $\sigma$) that are $O(|\log(\rho)|^{-\infty})$. The conjugated operator is thus of the form
    \begin{equation}
        T_i \circ \chi_i A \chi_i \circ T_i^{-1} = \chi_i A \chi_i + i H + \Op(O(\left< \sigma \right>^{m-2})) + O(|\log(\tau)^{-\infty}|)
    \end{equation}
    where $H$ has right total symbol $\frac{-H_{\chi_i a \chi_i}(\tau)}{\tau} g(\sigma)$.
    Finally, if $\sigma$ is elliptic in $WF(\chi_i A \chi_i)$, then it is as well in $WF(T_i \circ \chi_i A \chi_i \circ T_i^{-1})$. In that case, the conjugated operator is a regular pseudo-differential operator $\Psi^m$ on $[0, \infty) \times \partial \Bar{M}_i$ and $WF(T_i \circ \chi_i A \chi_i \circ T_i^{-1}) = WF(\chi_i A \chi_i)$.
\end{proposition}
\textbf{Remark:} While the conjugated operator is not a standard pseudo-differential operator it behaves quite similarly to the usual ones except one should think of the action to be global (even up to smoothing terms) in the $[0, \infty)$ factor when $\sigma = 0$. Note that the main effect of conjugation is to pullback the symbol $a$ by the map $\phi: \mathbb{R} \rightarrow \Gamma$, with some error terms that are well behaved for our class of operators. For stationary operators $A_0$, conjugation is \textit{exactly} precomposition by $\phi$.  
\vspace{0.2in}
\newline
\textbf{Proof:} In proposition \ref{prop:ConjugationFormula} we computed the Mellin transform of the kernel of the conjugated operator to be the Mellin transform of $K(\kappa, \rho)$, the kernel as a function of $\kappa = \tau / \rho$ and $\rho$, in $\Gamma_i \times (\Gamma_i - \Gamma_i)$:

\begin{equation}
    K_{T_i \circ \chi_i A \chi_i \circ T_i^{-1}} = \frac{1}{2\pi} \int_{\mathbb{R}_{\sigma}} \int_{\mathbb{R}_{\eta}} \tau^{-i\sigma} \rho^{i\eta} \widehat{K_{\chi_i A \chi_i}}(\phi_i(\sigma), \phi_i(\sigma) - \phi_i(\eta)) \phi_i'(\eta) d\sigma d\eta
\end{equation}
For pseudo-differential operators, the Mellin transform in $\kappa$ of the Kernel is simply the (right) total symbol of the operator, so $\widehat{K}$ is just

\begin{equation}
    \widehat{\chi_i a \chi_i}(\phi_i(\sigma) - \phi_i(\eta), \phi_i(\sigma), x, \xi) \phi_i'(\eta) 
\end{equation}
where $\widehat{\chi_i a \chi_i}$ denotes the Mellin transform of $\chi_i a \chi_i$ in the $\rho$ variable. We then take the inverse Mellin transform over $\Gamma_{\eta}$ to get the symbol, or equivalently over $\vartheta = \sigma - \eta$:

\begin{equation}
\label{eq:total_conjugation}
    \sigma(T_i \circ \chi_i A \chi_i \circ T_i^{-1})(\rho, \sigma) = \frac{1}{\sqrt{2\pi}} \int_{\mathbb{R}_{\vartheta}} \rho^{-i \vartheta} \widehat{\chi_i a \chi_i}(\vartheta + i\vartheta q(\sigma, \vartheta), \phi_i(\sigma), x, \xi) \phi_i'(\sigma - \vartheta)  d\vartheta
\end{equation}
We will show the expansion by considering terms of the expansion \ref{eq:AssumedExpansionForConjugation} separately. Set $\phi_i = \phi$ for the remainder of the proof.
\vspace{0.1in}
\newline
1) Let's start with the highly decaying term $\tau^{\mu + \epsilon} \Psi^m_b(\Bar{M})$ as it is the easiest to understand. In that case, $\widehat{\chi_i a \chi_i}(\eta, \phi(\sigma), x, \xi)$ exists and is holomorphic in $Im(\eta) < \mu + \epsilon$ with rapid decay along horizontal lines in $\eta$. Note that this integral is a $0,0$ order symbol as defined in \ref{def:ProductTypeOperators}, i.e. well behaved under differentiation by $\rho \partial_{\rho}$, $\left< \sigma \right> \partial_{\sigma}$, $\partial_{x}$ and $\left<(\sigma, \xi)\right> \partial_{\xi}$. Indeed, differentiation by $\rho \partial_{\rho}$ creates powers of $\vartheta$, which can be absorbed by the fast decay of $\widehat{\chi_i a \chi_i}(\vartheta  + i\vartheta q(\sigma, \vartheta), \phi(\sigma), x, \xi)$ in $\left< \vartheta + i\vartheta q(\sigma, \vartheta) \right> \sim \left< \vartheta \right>$, while differentiation in $\partial_{\sigma}$ gives

\begin{equation}
\begin{aligned}
    &\sqrt{2\pi} \partial_{\sigma} \sigma(T_i \circ \chi_i A \chi_i \circ T_i^{-1})(\rho, \sigma) = \int_{\mathbb{R}_{\vartheta}} \rho^{-i \vartheta} (\widehat{\chi_i a \chi_i}(\vartheta + i \vartheta q(\sigma, \vartheta), \phi(\sigma), x, \xi) ig''(\sigma - \vartheta) \\
    &+ \widehat{\partial_{\sigma} \chi_i a \chi_i} (\vartheta + i\vartheta q(\sigma, \vartheta), \phi(\sigma), x, \xi) \phi'(\sigma) \phi'(\sigma - \vartheta) +  \partial_{\eta} \widehat{\chi_i a \chi_i}(\vartheta + i\vartheta q(\sigma, \vartheta), \phi(\sigma), x, \xi) \partial_{\sigma}(i\vartheta q(\sigma, \vartheta)) \phi'(\sigma - \vartheta)) d\vartheta
\end{aligned}
\end{equation}
We claim each terms is uniformly bounded by $\left<(\sigma, \xi) \right>^{m} \left< \sigma \right>^{-1}$. Indeed, for the first term $|g''(\sigma - \vartheta)| \lesssim \left< \sigma - \vartheta \right>^{-2}$ which after integration against a rapidly decaying function of $\vartheta$ gives the weight $\left< \sigma \right>^{-1}$. For the second term $\partial_{\sigma} \chi_i a \chi_i$ is directly $\left< (\sigma, \xi) \right>^{m-1}$ while the other factors are bounded. Finally, for the third term we have that $\partial_{\sigma} \vartheta q(\sigma, \vartheta) = \partial_{\sigma}(g(\sigma) - g(\sigma - \vartheta)) = \vartheta \int_{0}^{1} g''(\sigma - r \vartheta) dr$ is bounded by 

\begin{equation}
    |\vartheta \partial_{\sigma} q| \lesssim |\vartheta| \int_0^1 \left< \sigma - r \vartheta \right>^{-2} dr \lesssim |\vartheta| \left< \vartheta \right>^2 \left< \sigma \right>^{-2}
\end{equation}
Integrating this factor against the rapidly decaying in $\vartheta$ function $\partial_{\eta} \widehat{\chi_ia  \chi_i}$ then gives the weight $\left< \sigma \right>^{-1}$. Similarly, derivatives in $\partial_{\sigma}^N$ gives symbolic bound $\lesssim \left< (\sigma, \xi) \right>^m \left< \sigma \right>^{-N}$ by the same argument. Indeed, the $N^{th}$ derivative with respect to $\sigma$ is of the form (overall factors of $\phi'$ are absorbed in big $O$ since they are bounded; the full formula would also involve mix derivative terms $\partial^j_{\eta} \partial^k_{\sigma}$ of $\widehat{\chi_i a \chi_i}$ which are handled in the exact same way)
\begin{equation}
\begin{aligned}
    &\sqrt{2 \pi}\partial_{\sigma}^{N} \sigma(T_i \circ \chi_i A \chi_i \circ T_i^{-1})(\rho, \sigma) \\
    &\simeq \int_{\mathbb{R}_{\vartheta}} \rho^{-i \vartheta} \sum_{j, k; k+j = N} \Bigg\{ \widehat{\partial_{\sigma}^k \chi_i a \chi_i}(\vartheta + i\vartheta q(\sigma, \vartheta), \phi(\sigma), x, \xi) \times \sum_{n, l_i; n + \sum l_i = j} O\left( \partial_{\sigma}^n \phi'(\sigma - \vartheta) \Pi_{i=1}^k \partial^{l_i}_{\sigma} \phi'(\sigma) \right) \\
    &+ \partial^{k}_{\eta} \widehat{\chi_i a \chi_i}(\vartheta + i\vartheta q(\sigma, \vartheta), \phi(\sigma), x, \xi) \times  \sum_{n, l_i; n + \sum l_i = j}O\left( \partial^{n}_{\sigma}\phi'(\sigma - \vartheta) \Pi_{i=1}^k \partial^{l_i}_{\sigma}(\vartheta q'(\sigma, \vartheta)) \right) + \text{mixed $\partial$} \Bigg\} d\vartheta
\end{aligned}
\end{equation}
The first term satisfies the bound since $\partial^k_{\sigma} \chi_i a \chi_i$ bounded by $\left< (\sigma, \xi) \right>^{m-k}$ while $\partial_{\sigma}^n \phi'(\sigma - \vartheta) \Pi_{i=1}^k \partial^{l_i}_{\sigma} \phi'(\sigma)$ is bounded by $\left< \sigma - \vartheta \right>^{-n} \left<\sigma\right>^{-\sum l_i}$, which after integrations against a rapidly decaying function of $\vartheta$ gives a weight $\left< \sigma \right>^{-j}$. As for the second term we have

\begin{equation}
    \Pi_{i=1}^k |\partial^{l_i}_{\sigma} (\vartheta q')| \leq |\vartheta|^k \Pi_{i=1}^k \int_{0}^{1} |g^{(l_i+2)}(\sigma - r \vartheta)| dr \lesssim |\vartheta|^k \Pi_{i = 1}^k \left< \vartheta \right>^{l_i + 2} \left< \sigma \right>^{-l_i-2}
\end{equation}
Taking the product with $\partial_{\sigma}^{n} \phi'(\sigma -\vartheta) \lesssim \left< \sigma - \vartheta \right>^{-n}$ and with $\partial_{\eta}^k \widehat{\chi_i a \chi_i} \lesssim \left< (\sigma, \xi) \right>^{m}$ then gives a product bounded by $\lesssim \left< (\sigma, \xi) \right>^{m} \left< \sigma \right>^{-n-\sum l_i -k}$. Integrating against a rapidly decaying function of $\vartheta$ then gives the desired bound. Finally, it is of order $O(|\log(\tau)|^{-\infty})$ uniformly over $\sigma$ because we are taking the inverse Mellin transform of a rapidly decaying function of $\vartheta$. It remains to show the expansion in powers of $\left< \sigma \right>$. This is exactly the same way we obtain the asymptotic expansion for composition of pseudo-differential operator using Taylor series, hence we will be concise. Expanding the integrant by doing a Taylor series in $\vartheta$ gives
\begin{equation}
\begin{aligned}
    \widehat{\chi_i a \chi_i}(\vartheta + i\vartheta q(\sigma, \vartheta), \phi(\sigma), x, \xi) \phi'(\sigma - \vartheta) = \widehat{\chi_i a \chi_i}(\vartheta, \phi(\sigma), x, \xi)\phi'(\sigma) \\
    + i \vartheta g'(\sigma) \partial_{\eta} \widehat{\chi_i a \chi_i}(\vartheta, \phi(\sigma), x, \xi) + O(\left< \sigma - \vartheta \right>^{m-2})
\end{aligned}
\end{equation}
Taking the inverse Mellin transform then gives the first expansion in \ref{eq:conjugated_symbol}. As for the second, simply use a Taylor series in $\sigma$ now, and use that $-H_{\chi_i a \chi_i}(\tau)/\tau = \partial_{\sigma} \chi_i a \chi_i$ while the $\partial_{\rho}$ derivatives are rapidly decaying in $|\log(\rho)|^{-\infty}$.
\vspace{0.1in}
\newline
2) We now show that the conjugation of the stationary part $\chi_i A_0 \chi_i$ of the operator is also a pseudodifferential operator of class $\Psi^{m,m}$ with the expansion \ref{eq:conjugated_symbol}. As in proposition \ref{prop:FourierTransformCutoffs} the Mellin transform of a one sided-cutoff $\chi_i$ is the boundary value of $\widehat{-\rho \partial_{\rho} \chi}_i(\eta)/ i\eta$ as you approach from the lower half plane. Similarly here, since $\partial_{\rho} (\chi_i a_0 \chi_i)$ is a compactly supported function, we have that
\begin{equation}
    \widehat{\sigma(\chi_i A_0 \chi_i)}_i(\eta, \phi(\sigma)) = -\frac{\widehat{\rho \partial_{\rho}(\chi_i a_0 \chi_i)}(\eta, \phi(\sigma))}{i(\eta - i0)} \equiv \frac{\widehat{h}(\eta, \phi(\sigma))}{i(\eta - i0)}
\end{equation}
The conjugated operator then has kernel with Mellin transform in $(\kappa, \rho)$ given by the distribution
\begin{equation}
    \widehat{K}(\sigma, \vartheta) = \frac{\widehat{h}(\vartheta + i\vartheta q(\sigma, \vartheta), \phi(\sigma))}{i(\vartheta  + i\vartheta q(\sigma, \vartheta) - i0)} \phi'(\sigma - \vartheta) = \frac{\widehat{h}(\vartheta + i\vartheta q(\sigma, \vartheta), \phi(\sigma))}{i(1 + iq(\sigma, \vartheta))(\vartheta  - i0)} (1 + ig'(\sigma - \vartheta))
 \end{equation}
Under the (distributional) inverse Mellin transform in $\vartheta$, this corresponds to a $0,0$ order symbol satisfying the expansion \ref{eq:conjugated_symbol}. Indeed, as before we wish to Taylor expand in $\vartheta$. Before doing so, we must remove the singularity at $\vartheta = 0$, which we do so by subtracting the Mellin transform of $\chi_i A_0 \chi_i$, resulting in a smooth function $H(\sigma, \vartheta)$ (since the denominator is smooth in $\vartheta$ and both numerator and denominator vanishes at $\vartheta = 0$):
\begin{equation}
\begin{aligned}
    H(\sigma, \vartheta) &= \widehat{K}(\sigma, \vartheta) - \frac{\widehat{h}(\vartheta, \phi(\sigma))}{i(\vartheta - i0)} \\
    &= \frac{\widehat{h}(\vartheta + i\vartheta q(\sigma, \vartheta), \phi(\sigma))}{i(1 + iq(\sigma, \vartheta))(\vartheta -i0)} (1 + ig'(\sigma - \vartheta)) - \frac{\widehat{h}(\vartheta, \phi(\sigma))}{i(\vartheta - i0)} \\
    &= \frac{\widehat{h}(\vartheta(1 + iq(\sigma, \vartheta)), \phi(\sigma))(1 + ig'(\sigma - \vartheta)) - \widehat{h}(\vartheta, \phi(\sigma))(1 + iq(\sigma, \vartheta))}{(i(1 + iq(\sigma, \vartheta))(\vartheta - i0)} \\
\end{aligned} 
\end{equation}
We can thus apply the Taylor series argument to $H(\sigma, \vartheta)$. To see that it gives expansion \ref{eq:conjugated_symbol}, note that we could do the computation with $\frac{\widehat{h}(\eta, \phi(\sigma))}{i(\eta - i\epsilon)}$ for $\epsilon>0$ small (i.e. the Mellin transform in the lower half-plane) and then take the limit $\epsilon \rightarrow 0$ to obtain the real Mellin transform. Since the two computations must agree, and in the lower-half plan it gives expansion \ref{eq:conjugated_symbol}, we are done. 
\vspace{0.1in}
\newline
4) Next we compute the effect of conjugation on the powers $\tau^{-i\eta_l}(-\log(\tau))^k$ similarly to how we did for the stationary case. The Mellin transform acting on functions supported on $\tau < 1$ is given by equation \ref{eq:FourierTransformPowers}
\begin{equation}
    \widehat{\tau^{-i\eta_l}(-\log(\tau))^k} = \frac{1}{\sqrt{2\pi}}  \frac{k!}{i^{k+1}(\eta - \eta_l)^{k+1}} = \frac{1}{\sqrt{2\pi}} \frac{\partial^k}{\partial \eta^k} \frac{1}{i(\eta - i0)}(\eta - \eta_l)
\end{equation}
And similarly, the Mellin transform of $\chi_i \tau^{-i\eta_l} (-\log(\tau))^k a_{lk} \chi_i$ is given by
\begin{equation}
    \widehat{\sigma(\chi_i \tau^{-i\eta_l} (-\log(\tau))^k A_{lk} \chi_i)}(\eta, \phi(\sigma)) = i^k \partial_{\eta}^k \left(\frac{\widehat{\tau \partial_{\tau}(\chi_i a_0 \chi_i)}}{i(\eta - i0)}\right)(\eta - \eta_l, \phi(\sigma))
\end{equation}
Which can be understood to be defined on $Im(\sigma) < Im(\eta_l)$ with a boundary value defined on the line $Im(\sigma) = Im(\sigma_l)$ as you approach it from the lower plane. The Kernel thus has Mellin transform on $\mathbb{R}$ given by
\begin{equation}
    \widehat{K}(\sigma, \vartheta) = i^k \partial_{\eta}^k \left(\frac{\widehat{\tau \partial_{\tau}(\chi_i a_0 \chi_i)}}{(i(\eta - i0))}\right)(\vartheta + i\vartheta q(\sigma, \vartheta) - \eta_l, \phi(\sigma)) \times \phi'(\sigma - \vartheta)
\end{equation}
Because of our assumption on $\eta_l$, this expression makes sense since the contour $\vartheta \rightarrow \vartheta + i\vartheta q(\sigma, \vartheta)$ sits below $\eta_l$ for all $\sigma$ (we don't even need to take a boundary value). Then to obtain the expansion \ref{eq:conjugated_symbol}, we just proceed exactly as before and do a Taylor series in $\vartheta$.
\vspace{0.1in}
\newline
5) Finally we study how the conjugation affects wave front sets. Suppose that $(\sigma_0, x_0, \xi_0) \notin (WF(\chi_i A \chi_i))$, so the total symbol is rapidly decaying around it uniformly over $\tau$:
\begin{equation}
    \forall N \quad |a(\tau, \sigma, x, \xi)| \lesssim \left<(\sigma, \xi) \right>^{-N} \quad \text{for $(Re(\sigma), x, \xi) \in C$, where $C$ is some open cone centered at $(\sigma_0, x_0, \xi_0)$}
\end{equation}
The conjugated total (right) symbol is then
\begin{equation}
    \sigma(T_i \circ \chi_i A \chi_i \circ T_i^{-1})(\rho, \sigma) = \frac{1}{\sqrt{2\pi}} \int_{\mathbb{R}_{\vartheta}} \rho^{-i \vartheta} \widehat{\chi_i a \chi_i}(\vartheta + i\vartheta q(\sigma, \vartheta), \phi(\sigma), x, \xi) \phi_i'(\sigma - \vartheta)  d\vartheta
\end{equation}
By the uniformity of the estimate over $\tau$, we get the same one for the Mellin transform of $a$ in $\tau$, hence for $(Re(\sigma), x, \xi) \in C$, the integral above is bounded by
\begin{equation}
    |\sigma(T_i \circ \chi_i A \chi_i \circ T_i^{-1})(\rho, \sigma)| \lesssim \int_{\mathbb{R}_{\vartheta}} \left<\vartheta\right>^{-2} \left<(\sigma,  \xi)\right>^{-N} \lesssim \left<(\theta, \xi)\right>^{-N}
\end{equation}
The proof is the exact same for derivatives. Thus, $(\sigma_0, x_0, \xi_0) \notin WF(T_i \circ \chi_i A \chi_i \circ T_i^{-1})$ uniformly over $\tau$ hence the conjugation by $T_i$ preserves $\tau$ independent wavefront sets. In particular, if $\sigma$ is elliptic on $WF(\chi_i A \chi_i)$, then it is as well on $WF(T_i \circ \chi_i A \chi_i \circ T_i^{-1})$, which makes the conjugated operator a regular pseudo-differential operator since we can bound $\left<\sigma\right>^{-1} \lesssim \left<(\sigma, \xi)\right>^{-1}$ within its wavefront set as shown in lemma \ref{lemma:MicrolocalizationAgree}. In that case, the wavefront sets are the same, since the expansion $\ref{eq:conjugated_symbol}$ can be done to all orders with a remainder $\left<\sigma\right>^{-N}$. If $\sigma$ is elliptic, this remainder is smoothing hence the wavefront set of the conjugated symbol lies in the wavefront sets of the individual terms in the expansion, which are local in the symbol of $\chi_i A \chi_i$. 
$\Box$
\vspace{0.2in}
\newline
It now remains to show that $\phi(\sigma) = \sigma + ig(\sigma)$ can be chosen so that the conjugated operator satisfies an estimate in a neighborhood of the boundary. Heuristically, this is quite simple to do since the symbols we'll be working with will be of the form $\sigma^2 - \xi^2$, hence the conjugated symbol is

\begin{equation}
    \sigma^2 + i 2 \sigma g(\sigma) - g(\sigma)^2 - \xi^2
\end{equation}
Which has a lower order term $i2 \sigma g(\sigma)$ that has a sign as long as $g$ has a definite sign in $\pm \sigma > 0$. Thus, in general, we can expect to get an estimate on these operators for operator $P$ whose stationary part $P_0$ has a subprincipal symbol with definite sign in $\pm \sigma > 0$ after composition by $\phi_i(\sigma)$ at each boundaries $\partial M_i$. This is summarized by the next proposition, which we call 'radial point estimate' because the existence of such estimates was motivated by a positive commutator heuristic.

\begin{figure}[htbp]
    \centering
    \includegraphics[width=0.75\textwidth]{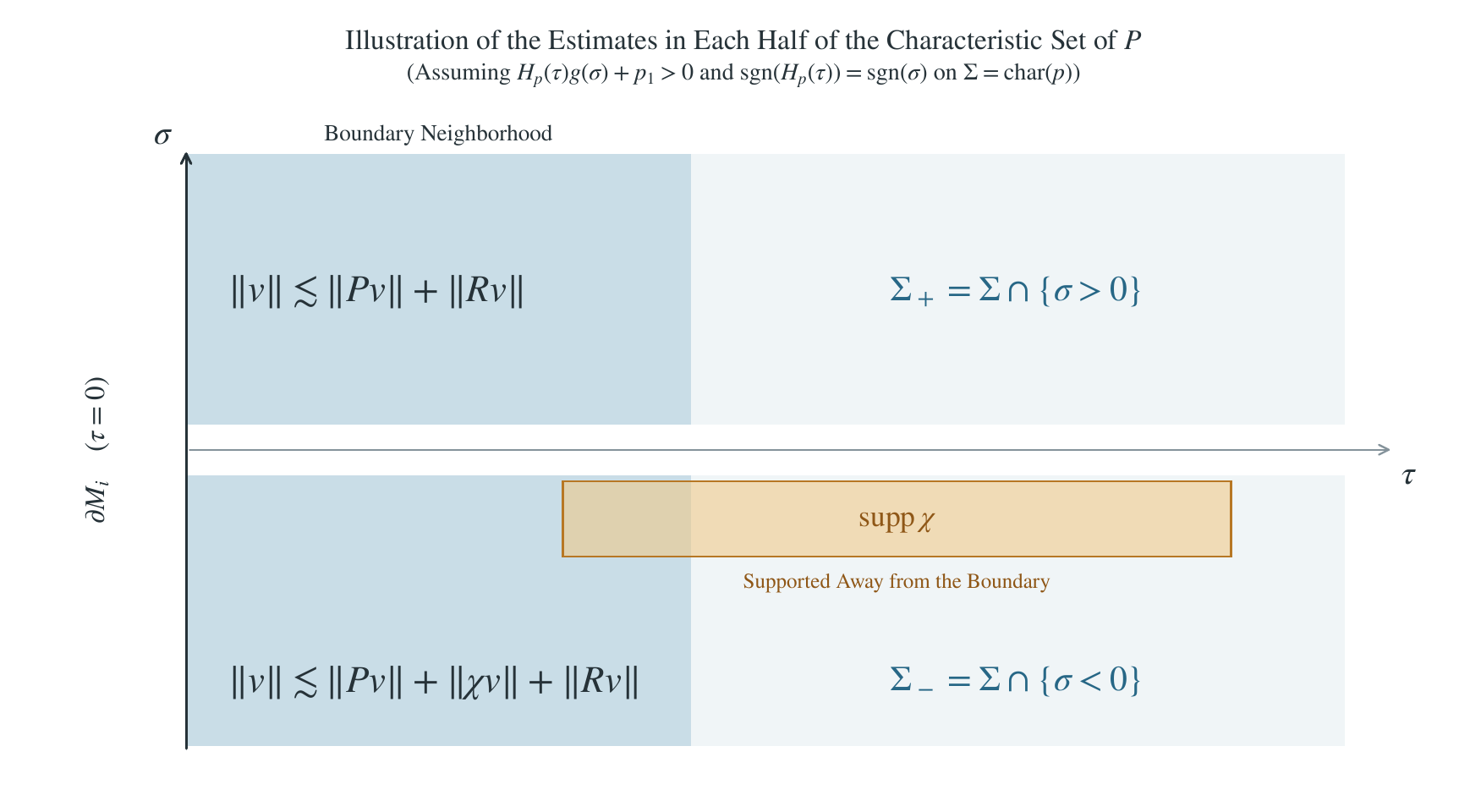}
    \caption{A sketch of the estimate in the Anti-Feynman case, using $\varrho = \sigma / \left< (\sigma, \xi) \right>$ for simplicity. Each half has a separate estimate. In one we have a 'free' estimate, while in the other we require control away from the boundary $\tau = 0$.}
    \label{fig:boundary_estimate_sketch}
\end{figure}

\begin{proposition}
\label{prop:RadialPointEstimates}
(Radial Point Estimates) Let $P, Q \in \Psi^m_{b, \Gamma}(M)$ be classical pseudo-differential operators acting on some vector bundle $V$ equipped with a positive-definite metric. Suppose $g_i(\sigma)$ is a symbol of order $0$ for the normal neighborhood $N_i$ of a connected component $\partial \Bar{M}_i$ of the boundary. Let $p$ be the real scalar principal symbol of $P$ and $p_1$ its subprincipal symbol, i.e. the principal symbol of $\frac{1}{2i}(P - P^{\dagger})$, where $P^{\dagger}$ is the adjoint with respect to the chosen positive definite metric on $V$ and $b-$density on $\Bar{M}$. Suppose that $Q$ has Hermitian principal symbol and that there exists $\epsilon > 0$ such that 
    \begin{equation}
        \sigma \text{, } -\frac{H_p(\tau)}{\tau}g(\sigma) Id + p_1 \text{ are elliptic on } \{ \tau < \epsilon, \frac{|p|}{\left<(\sigma, \xi)\right>^m} < \epsilon \} \cap N_i \text{ and $sgn(-\frac{H_p(\tau)}{\tau} g(\sigma)Id + p_1) = sgn(Q)$ there}
    \end{equation}
where $sgn(Q)$ denotes the signature of the operator $Q$ up to terms of order $m - 2$. More precisely, $Q$ is of the form \footnote{For $Q$ principally scalar, which it usually will be, a sufficient condition for that property is that the total symbol of $Q$ in charts is non-negative (non-positive) in $-(H_p(\tau)/\tau)g(\sigma) Id + p_1>0$ ($<0$) up to terms of order $m-2$. This is the Fefferman-Phong inequality \cite{FeffermanPhong1978}. One could also just require the principal symbol of $Q$ to satisfy the sign condition, in which case its subprincipal term would enter the threshold.}
\begin{equation}
    Q = Q_{+} - Q_{-} + Q_r \quad WF(Q_{\pm}) \subset \left\{ \pm(-\frac{H_p(\tau)}{\tau}g(\sigma) Id + p_1) > 0 \right\}, Q_{\pm} \in \Psi^{m} \text{ and $Re\left<Q_{\pm}v, v\right> \geq 0$}, Q_r \in \Psi^{m-2}_{b, \Gamma}
\end{equation}
Furthermore, assume that $-H_p(\tau)$ shares the same sign as $-\frac{H_p(\tau)}{\tau}g(\sigma) Id + p_1$ when $\varrho > 0$, where $\varrho$ is some smooth function on the b-cosphere bundle $S^{\ast}_b \Bar{M}$ elliptic in the characteristic set $char(p)$ \footnote{For this condition, if $H_p(\tau)$ vanishes at a point, we say it shares the same sign as $-H_p(\tau)/\tau Id + p_1$ automatically.}:
    \begin{equation}
        -sgn(H_p(\tau)) = sgn\left( -\frac{H_p(\tau)}{\tau}g(\sigma) Id + p_1 \right) \quad \text{in } \{ \tau < \epsilon, \frac{|p|}{\left<(\sigma, \xi)\right>^m} < \epsilon, \varrho > 0 \} \cap N_i = W_{\epsilon}^{+}
    \end{equation}
Then, near $\partial M_i$ we have a local estimate in the sense that for any $N$ there exists $B \in \Psi^0_{b, \Gamma}$ elliptic in $W_{\delta}^{+}$ (for some $\delta$) and $R, G \in \Psi^{0}_{b, \Gamma}$, all with wavefront set in $W_{\epsilon}^{+}$, such that

\begin{equation}
    \| B v \|_{H^s_{\Gamma}} \lesssim \| G(P+iQ) v \|_{H^{s-m+1}_{\Gamma}} +  \| R v \|_{H^{-N}_{\Gamma}} \quad \forall v \in H^s_{b, \Gamma}
\end{equation}
On the other end, if $-H_p(\tau)$ does not share the same sign as $\frac{H_p(\tau)}{\tau}g(\sigma) Id + p_1$, then there exists for any $N$ operators $B \in \Psi^0_{b, \Gamma}$ elliptic in $W_{\delta}^{+}$, $E \in \Psi^0_b(M)$ and $R, G \in \Psi^{0}_{b, \Gamma}$, all with wavefront set in $W_{\epsilon}^{+}$, such that
    \begin{equation}
        \| B v \|_{H^s_{\Gamma}} \lesssim \| G(P+iQ) v \|_{H^{s-m+1}_{\Gamma}} + \| \chi Ev \|_{H^{s}_{b}} +  \| R v \|_{H^{-N}_{\Gamma}} \quad \forall v \in H^s_{b, \Gamma}
    \end{equation}
where $\chi(\tau)$ is a cutoff supported away from the boundary $\partial M_i$. Similarly, the same results holds for $W_{\epsilon}^{-} = \{ \tau < \epsilon, \frac{|p|}{\left<(\sigma, \xi)\right>^m} < \epsilon, \varrho < 0 \} \cap N_i$. That is, if the signs agree we get an estimate on $W_{\delta}^{-}$ of $v$ by $Pv$, while if they don't we must also control the solution $v$ away from the boundary. Furthermore, these estimates all hold in the strong sense in that if the right hand-side is finite and $v \in \cup_r H^r_{\Gamma}$ then $Bv \in H^s_{\Gamma}$ and the inequality holds. Finally, these estimates are uniform over $P$ and $Q$, in that if $P_{\theta} + iQ_{\theta}$ is a bounded family of operators (in the sense that all the symbolic norms on their total symbol in a family of charts are uniformly bounded) satisfying our assumptions uniformly (in the sense that $Q_{\theta}$ still satisfy the sign conditions and $P_{\theta}$ has a uniform lower bound for $|-\frac{H_p(\tau)}{\tau}g(\sigma) Id + p_1|$ in the appropriate regions), then the constants in the estimates above are uniform over $\theta$. 
\end{proposition}
\textbf{Remark:} The assumptions are quite simple for wave-type operators with principal symbol $\sim \sigma^2 - \xi^2$. In that case, $-H_p(\tau)/\tau = \sigma$ and we take $\varrho = \sigma / \left< (\sigma, \xi) \right>$, hence the assumptions are on the sign of $\sigma g(\sigma) + p_1$ in the two halves of the characteristic sets $\pm \sigma >0$. Sometimes it can also be convenient to choose $\varrho = -H_p(\tau)/\tau$ when it is elliptic, which is the case when $p_1 = 0$. As for the uniformity over the imaginary part $Q$, this will be useful when we do Wick rotation. Finally, it will follow from the proof that the estimates can be further microlocalized in any subregion defined by $|\varsigma| < \epsilon$, where $\varsigma$ is some smooth function on the $b-$cosphere bundle $S^{\ast}_b \Bar{M}$ such that $-sgn(H_p(\varsigma)) = sgn(\varsigma) \times sgn(-H_p(\tau)/\tau g(\sigma) Id + p_1)$.
\vspace{0.2in}
\newline
\textbf{Proof:} We will show the case of functions, i.e. $V = \mathbb{R}$; generalization to vector bundles is simple. In order to estimate the solution within the characteristic set, we will work in the conjugated picture under the (local) isomorphism $T_i: H^s_{N_i, \Gamma_i} \rightarrow H^s_b([0, \infty) \times \partial \Bar{M}_i)$ defined in proposition \ref{prop:IsoContourSobolev} together with a standard positive commutator estimate \footnote{Technically here we are using the subprincipal term of the conjugated operator to obtain an estimate, but this could be easily written as a commutator estimate where $[\chi_i P \chi_i, T] = (\chi_i P \chi_i -  T_i \circ \chi_i P \chi_i \circ T_i^{-1}) \circ T$ provides the positive commutator.}. Let $\chi_i$ be a cutoff near a normal neighborhood $N_i$ of $\partial \Bar{M}_i$ and define the conjugated operator
\begin{equation}
    \Tilde{P} + i\Tilde{Q} = T_i \circ \chi_i (P +iQ) \chi_i \circ T_i^{-1}
\end{equation}
 By proposition \ref{prop:ConjugatedSymbol} this is a $m^{th}$ order pseudo-differential operator on $H^s_{b}([0, \infty) \times \partial \Bar{M}_i)$ of class $\Psi^{(m,m)}$, in the sense that $\partial_{\sigma}$ derivatives of its symbol only improves the order by $\left<\sigma\right>^{-1}$ rather than $\left<(\sigma, \xi)\right>^{-1}$. We also computed its symbol; recall for the following that up to $O(\left<\sigma\right>^{m-2})$ terms it was given by 
 \begin{equation}
 \begin{aligned}
     \widetilde{p}_{tot} &= \chi_i^2 p_{tot} - i \frac{H_{\chi_i^2 p}(\tau)}{\tau} g(\sigma) + O(|\log(\tau)|^{-\infty}) + O(\left<\sigma\right>^{m-2}) \\
     \widetilde{q}_{tot} &= \chi_i^2 q_{tot} - i \frac{H_{\chi_i^2 q}(\tau)}{\tau}g(\sigma) + O(|\log(\tau)|^{-\infty}) + O(\left<\sigma\right>^{m-2})
\end{aligned}
 \end{equation}
 where $p_{tot}$ is the total symbol of $P$. For an operator $A$, we have the following identity:
\begin{equation}
\begin{aligned}
    Im \left< A (\Tilde{P} + i\widetilde{Q}) v, A v \right> &=  Re \left(\left< \frac{1}{i}[A, \Tilde{P}] v , Av \right>\right) + \frac{1}{2i} \left<(\Tilde{P}-\Tilde{P}^{\dagger})Av, Av \right> \\
    & + \frac{1}{2} \left<(\widetilde{Q} + \widetilde{Q}^{\dagger})Av, Av \right> + \frac{1}{2} \left<Mv, v\right>
\end{aligned}
\end{equation}
Where $M = [A,[A,\widetilde{Q}]] + (A^{\dagger}-A)[A, \widetilde{Q}] + [A-A^{\dagger}, \widetilde{Q}^{\dagger}]A + [A, \widetilde{Q} - \widetilde{Q}^{\dagger}]A$. We will work with $\varrho>0$; the $\varrho<0$ case proceeds in the exact fashion. Pick a cutoff $w(s)$ equal to $1$ for $s>\epsilon>0$ and supported in $[0, \infty)$ and a cutoff $\psi(s)$ equal to $1$ in $(-\epsilon/2, \epsilon/2)$ with support in $(-\epsilon, \epsilon)$. We will choose $A$ as follows:

\begin{equation}
    A = h(\tau) \Op(\left< \sigma \right>^{s-(m-1)/2}  \times w(\varrho) \times \psi(p|_{\tau=0}/\left<(\sigma, \xi)\right>^m))
\end{equation}
Where the quantization is chosen so that $WF(A) \subset W_{\epsilon}^{+}$. Note that we can write $A = AB + R = BA + R$, where $B = \Op(\Psi(p|_{\tau=0}/\left<(\sigma, \xi)\right>^m))$ for $\Psi$ a cutoff equal to $1$ on $supp(\psi)$, $B = B^{\dagger}$ up to smoothing operator, $WF(B) \subset W_{\epsilon}^{+}$ and $R$ is smoothing (hence can and will be be discarded from the argument). The expression up to smoothing terms can then be written as
\begin{equation}
\begin{aligned}
    Im \left< A (\Tilde{P} + i\widetilde{Q}) v, A v \right> &= Re \left(\left< \frac{1}{i} [A, B\Tilde{P}B] v , Av \right> \right) + \frac{1}{2i} \left<B(\Tilde{P}-\Tilde{P}^{\dagger})BAv, Av \right> \\
    &+ \frac{1}{2} \left<B(\widetilde{Q} + \widetilde{Q}^{\dagger})BAv, Av \right> + \frac{1}{2} \left<Mv, v\right>
\end{aligned}
\end{equation}
Now $B\Tilde{P}B$ and $B \widetilde{Q}B$ are regular pseudo-differential operator of class $\Psi^{m}$ since $\sigma$ is elliptic within the wavefront set of $B$ by proposition \ref{prop:ConjugatedSymbol}. Thus, the remainder of the proof can proceed with a standard commutator estimate. The principal symbols of $B\Tilde{P}B$ is then (where $b$ is the principal symbol of $B$)
\begin{equation}
    \sigma_{m}(B\Tilde{P}B)= \chi_i(\tau)^2 \times b^2 \times p(\tau = 0, \sigma, x, \xi) + O(|\log(\tau)|^{-\infty})
\end{equation}
Thus, the principal symbol of the commutators involving $P$ is of the form
\begin{equation}
    b^2 \times h'(\tau) \times H_p(\tau) \times w(\varrho) \times \psi(p/\left<(\sigma, \xi)\right>^m)) + e(p+iq) + O(|\log(\tau)|^{-\infty}) a
\end{equation}
Where $e$ is some symbol coming from the derivative term $\psi'$ and from the $H_b$ part of the Hamiltonian vector field. As for the subprincipal term, its principal symbol is
\begin{equation}
    \frac{1}{2i} \sigma_{m-1}(B(\widetilde{P}-\widetilde{P}^{\dagger})B) = \chi_i(\tau)^2 \times b^2 \times \left(-\frac{H_p(\tau)}{\tau} g(\sigma) + p_1\right) + O(|\log(\tau)|^{-\infty})
\end{equation}
Thus, on the wavefront set of $A$, the subprincipal term of $B(\widetilde{P}-\widetilde{P}^{\dagger})B$ has an elliptic principal symbol by our assumptions. Suppose without loss of generality that it is positive. Then Garding's inequality gives
\begin{equation}
    \|Av\|^2_{H^{(m-1)/2}_b} \lesssim Re \left( \frac{1}{2i} \left<B(\Tilde{P}-\Tilde{P}^{\dagger})BAv, Av \right> \right) + \|Av\|^2_{H^{(m-2)/2}_b}
\end{equation}
As for the commutator term, there are two cases to separate whether $H_p(\tau) h'(\tau)$ is positive or negative, which depends on $-sgn(H_p(\tau))$ since $h'$ is negative definite. In the first case, the principal symbol of the commutator has a positive sign up to the error terms $e(p+iq) + O(|\log(\tau)^{-\infty}|)a$. Thus, again by Garding's inequality we have
\begin{equation}
    0 \lesssim  Re \left( \left< \frac{1}{i} [A, B\Tilde{P}B]v , Av \right>\right)  + \left|\left<G(\widetilde{P}+i\widetilde{Q})v, Av\right> \right| + \delta \left<Av, Av\right> + \|Rv\|^2_{H^{s-1/2}_b}
\end{equation}
Where $\delta$ is a number that can be taken arbitrary and $G \in \Psi_b^{s-(m+1)/2}$, $R \in \Psi^0_b$ are some microlocalizers. In particular, $\delta \left<Av, Av\right>$ can be absorbed using the positivity of $B(\widetilde{P} - \widetilde{P}^{\dagger})B$. In the other case where $H_p(\tau) h'(\tau)$ is negative, we cannot use Garding's inequality hence must instead bound the commutator. Nonetheless, because $h'(\tau)$ is supported away from the boundary, we can still bound it as follows:
\begin{equation}
     \left|Re\left(\left< \frac{1}{i}[A, B\Tilde{P}B]v , Av \right> \right)\right|  \lesssim |\left<\chi^2 Ev, Av \right>| + |\left<G(\widetilde{P}+i\widetilde{Q})v, Av\right>| + \delta \left< Av, Av \right> + \| Rv \|^2_{H_b^{s-1/2}}
\end{equation}
Where $\chi(\tau)$ is a cutoff supported away from the boundary and $E \in \Psi^s_b$ a microlocalizer. So the main difference between the two cases is whether or not the term $\left< \chi^2 Ev, Av \right>$ is included. For the remainder of the proof, we will restrict ourselves to the case where it is, with the understanding that when $H_p(\tau) h'(\tau)$ is the right sign (positive here) this term disappears. Finally, for the term involving $Q$, we can use its assumed positivity. Indeed, $M$ is a pseudo-differential operator of $2$ lower orders hence can be bounded by $\| Rv \|^2_{H^{s-1/2}_b}$ while our positivity assumption on $Q$ gives a similar one for $\widetilde{Q}$:
\begin{equation}
    B(\widetilde{Q} + \widetilde{Q}^{\dagger})B = B\chi_i(Q_{+} + Q_{+}^{\dagger} - Q_{-} - Q_{-}^{\dagger})\chi_iB + \widetilde{Q}_r
\end{equation}
For some $\widetilde{Q}'_r \in \Psi^{m-2}$. Note that we are using that the subprincipal term induced by the conjugating cancel when taking the real part of the operator. Hence we have positivity up to two orders lower:
\begin{equation}
    0 \lesssim Re \left( \left<B(\widetilde{Q} + \widetilde{Q}^{\dagger})B Av, Av\right>\right) + \| Av \|^2_{H^{(m-2)/2}}
\end{equation}
Thus, combining all these inequalities and applying Cauchy-Schwartz to the pairings yield:
\begin{equation}
\begin{aligned}
    \| Av \|^2_{H^{(m-1)/2}_b} &\lesssim \left| Im \left< A (\Tilde{P} + i\widetilde{Q}) v, A v \right> \right| + \left|\left<\chi^2 Ev, Av \right>\right| + \left| \left<G(\widetilde{P}+i\widetilde{Q})v, Av\right>\right| + \delta \|Av\|^2_{H^0_b} + \|Rv\|^2_{H^{s-1/2}_b} \\
    &\lesssim \frac{1}{\delta} \| G(\widetilde{P} + i\widetilde{Q})v \|^2_{H^{s-m+1}_b} + \| \chi^2 Ev \|^2_{H^s_b} + \delta \|Av\|^2_{H^{(m-1)/2}_b} + \| Rv \|^2_{H^{s-1/2}_b} 
\end{aligned}
\end{equation}
Taking $\delta$ small enough yields
\begin{equation}
\label{eq:EstimateCommutatorConjugated}
\begin{aligned}
    \| h(\tau) \Op(w(\varrho)) \times \psi(p|_{\tau=0}/\left<(\sigma, \xi)\right>^m)) v \|_{H^s_b} &\lesssim \| Av \|_{H^{(m-1)/2}_b} \\
    &\lesssim \| G(\widetilde{P} + i\widetilde{Q})v \|_{H^{s-m+1}_b} + \| \chi^2 Ev \|_{H^s_b} + \| Rv \|_{H^{s-1/2}_b} 
\end{aligned}
\end{equation}
Now to get the estimate on our original space $H^s_{b, \Gamma}$, we simply undo the conjugation. For $u \in H^s_{b, \Gamma}$, pick $v = T_i \Bar{\chi}_i u \in H^s_b([0, \infty) \times \partial M_i)$ where $\Bar{\chi}_i$ is a cutoff equal to $1$ in $supp(\chi_i)$. Then the estimate becomes
\begin{equation}
\begin{aligned}
    &\| T_i^{-1} \circ h(\tau) \Op(w(\varrho) \times \psi(p|_{\tau=0}/\left<(\sigma, \xi)\right>^m)) \circ T_i \Bar{\chi}_i u \|_{H^s_{b, \Gamma}} \lesssim \\
    &\| T_i^{-1} \circ G \circ T_i (\chi_i P \chi_i + i \chi_i Q \chi_i) \Bar{\chi}_i u\|_{H^s_{b, \Gamma}} + \| T_i^{-1} \circ \chi^2 E \circ T_i \Bar{\chi}_i u \|_{H^s_{b, \Gamma}} + \| T_i^{-1} \circ R \circ T_i \Bar{\chi}_i u \|_{H^{s-1/2}_{b, \Gamma}}
\end{aligned}
\end{equation}
To get rid of the conjugation, note that once we have the estimate \ref{eq:EstimateCommutatorConjugated} above, we can replace the operators involved by any operators in $\Psi^m$ elliptic in the corresponding regions. In particular, we can replace the $h(\tau) \Op(w(\varrho) \times \psi(p|_{\tau=0}/\left<(\sigma, \xi)\right>^m))$, $G$, $\chi^2 E$ and $R$ operators by operators $T_i \circ B \circ T_i^{-1}$, $T_i \circ G' \circ T_i^{-1}$, $T_i \circ \chi E' \circ T_i^{-1}$ and $T_i \circ R' \circ T_i^{-1}$, where $B$, $G'$, $\chi E$ and $R'$ are operators in $\Psi^{\circ}_{b, \Gamma}$ that are elliptic in the appropriate regions. This effectively uses that microlocalization before and after conjugation agree as long as $\sigma$ is elliptic, as shown in proposition \ref{prop:ConjugatedSymbol}. Combining everything together then yields the desired result
\begin{equation}
    \| Bv \|_{H^s_{b, \Gamma}} \lesssim \| G'\chi_i (P+iQ) \Bar{\chi}_i u \|_{H^s_{b, \Gamma}} + \| \chi E' u \|_{H^s_b} + \| R' u \|_{H^{s-1/2}_{b, \Gamma}}
\end{equation}
Finally, we can remove the $\Bar{\chi}_i$ term in front of $u$ because $\chi_i (P + iQ) (1 - \Bar{\chi}_i)$ is smoothing due to $\chi_i$ and $(1-\Bar{\chi}_i)$ having disjoint support. Iterating over the remainder term in $H^{s-1/2}_{b, \Gamma}$ then yields the desired result. Finally, for the claim regarding the inequality being strong, in that we only need to assume $v \in H^r_{b, \Gamma}$ for some $r$, add a regularizer in the commutant $A$ of the form (for some large enough $L$ to justify the pairings)
\begin{equation}
    \Lambda_{\delta} = \left<\delta \sigma\right>^{-L}
\end{equation}
This doesn't change the commutators with $P$ to leading term since $H_p(\sigma)|_{\tau = 0} = 0$ while the rest of the argument remains unchanged. One then takes the limit $\delta \rightarrow 0$; as this is rather standard, we omit the details and refer to Hintz's notes \cite{Hintz2025Microlocal} chapter 8 for the details in a similar argument. Finally, for the uniformness with respect to $P$ and $Q$, note that the only time the size of $Q$ enters the argument is when bounding the double commutator term $M$ and the residual $\widetilde{Q}_r$. If $Q$ has uniform symbolic norm on its symbol, then these bounds are uniform as well. As for $P$, if the ellipticity of $-\frac{H_p(\tau)}{\tau}g(\sigma) Id + p_1$ is uniform and we control its total symbol, then the bounds we used in the positive commutator estimate are all uniform as well. Thus, the overall constant can be taken independent of $\theta$.
$\Box$
\vspace{0.2in}
\newline
Combining this estimate with elliptic estimates away from $char(p)$ then yields control of the solutions $v$ near the boundary component $\partial \Bar{M}_i$. To obtain a global estimate, one then use the propagation of singularity theorem to propagates this control everywhere else in the manifold, assuming that the Hamiltonian flow of $p$ has good properties. We will make the result precise in the next subsection \ref{subsection:FredholmTheory} where we states the assumption necessary for the operators to be Fredholm. 

\subsection{Fredholm Theory on $H^s_{b, \Gamma}$}
\label{subsection:FredholmTheory}
\begin{figure}[h]
    \centering
    \includegraphics[width=0.6\linewidth]{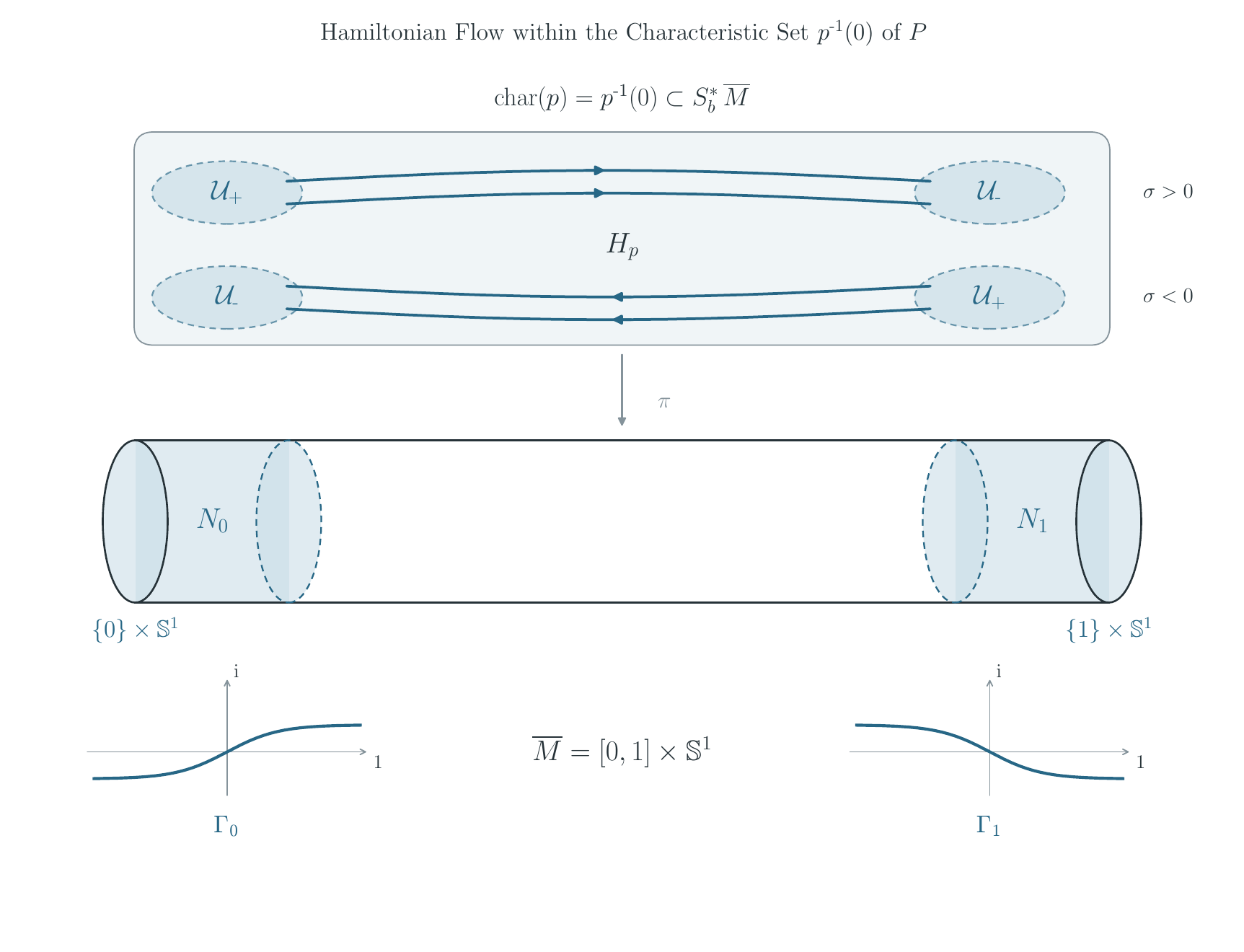}
    \caption{An illustration of a typical global Hamiltonian flow that guarantees a wave operator $P \sim -\partial_t^{2} + \partial_{\theta}^2$ to be Fredholm on the spacetime $\Bar{M} = [0, 1] \times \mathbb{S}^1$. The top figure illustrates the Hamiltonian flow in each halves $\pm \sigma > 0$ of the characteristic set $\sigma^2 = \xi^2$, while the bottom figure illustrates contours defining spaces $H^s_{b, \Gamma}$ on which $P$ is Fredholm. This choice of contour would correspond to a Feynman propagator.}
    \label{fig:HamiltonianFlowFredholm}
\end{figure}
Now that we have all the microlocal estimates required, we can state the Fredholm properties of various operators on our spaces $H^s_{b, \Gamma}$ on a $b-$manifold $\Bar{M}$. As usual, the theory splits between elliptic and non-elliptic operators. Under the usual assumption that the normal family is invertible on $\Gamma$, the elliptic operators will be automatically Fredholm (theorem \ref{theo:EllipticOperatorsAreFredholm}). As for the non-elliptic operators, we must also require some properties regarding the contour and the Hamiltonian flow of the principal symbol $p$ of the operator. Under such assumptions, non-elliptic operators will also be Fredholm on our spaces (see theorem \ref{theo:NonEllipticOperatorsAreFredholm}). Let's begin with the easy case of elliptic operators. As before, $V$ will denote a vector bundle $V$ equipped with a positive-definite metric, while $\Bar{M}$ denotes a compact manifold with boundary equipped with a $b-$density $\sim \frac{d\tau}{\tau} dy$.
\begin{theorem}
\label{theo:EllipticOperatorsAreFredholm}
    Let $P \in \Psi^m_{b, \Gamma}$ be an elliptic pseudo-differential operator on a $b-$manifold $\Bar{M}$ acting on some vector bundle $V$. Suppose that its indicial family on $\partial \Bar{M}_i$ $\widehat{N(P)}(\sigma):H^s_{\partial \Bar{M}_i, \sigma} \rightarrow H^{s-m}_{\partial \Bar{M}_i, \sigma}$ is invertible on $\Gamma_i$ for all $i$, where $\Gamma = \{ \Gamma_i \}_i$. Then $P:H^s_{b, \Gamma} \rightarrow H^{s-m}_{b, \Gamma}$ is a Fredholm operator. 
\end{theorem}
\textbf{Proof:} By proposition \ref{prop:EllipticEstimates} global ellipticity yields the estimates
\begin{equation}
\begin{aligned}
    &\| v \|_{H^s_{b, \Gamma}} \lesssim \| Pv \|_{H^{s-m}_{b, \Gamma}} + \| v \|_{H^{-N}_{b, \Gamma}} \\
    & \| v \|_{H^{-s+m}_{b, \Bar{\Gamma}}} \lesssim \| P^{\dagger}v \|_{H^{-s}_{b, \Bar{\Gamma}}} + \| v \|_{H^{-N}_{b, \Bar{\Gamma}}}
\end{aligned}
\end{equation}
To show a Fredholm estimate, we need the remainder to be compact, which here amounts to showing that we can estimate the remainders $\| v \|_{H^{-N}_{b, \Gamma}}$ on a contour with less decay. This is where we use our assumption on the invertibility of $N(P)(\sigma)$ on $\Gamma$. Pick a partition of unity $\chi_i, \chi$ subordinate to a cover of normal neighborhood $N_i$ and an interior open set $U$. Then
\begin{equation}
    \| v \|_{H^{-N}_{b, \Gamma}} \simeq \| \chi v \|_{H^{-N}_{b}} + \sum_i \| \chi_i v \|_{H^{-N}_{N_i, \Gamma_i}} 
\end{equation}
We need only to improve the decay on $\| \chi_i v \|_{H^{-N}_{\Gamma_i}}$. For that, use the assumed invertibility of $\widehat{N(P)}(\sigma)$ on $\Gamma$. By the continuity criteria for stationary operators on $H^s_{\Gamma}$ in lemma \ref{lemma:StationaryOperatorContinuous}, it gives that $N(P):H^s_{\Gamma} \rightarrow H^{s-m}_{\Gamma}$ is invertible (the necessary control $\| \widehat{N(P)}(\sigma)^{-1} \| \lesssim \left<\sigma\right>^{-m}$ follows from ellipticity).  Then, to improve the norm from a contour $\Gamma_i$ to $\Gamma_i - \delta i$ for $\delta>0$ small enough, we can use $N(P)$:
\begin{equation}
\begin{aligned}
    \| \chi_i v \|_{H^{-N}_{N_i, \Gamma_i}} \lesssim \| N(P) v \|_{H^{-N-m}_{N_i, \Gamma_i}} &\lesssim \| (P-N(P)) v \|_{H^{-N-m}_{N_i, \Gamma_i}} + \| P v \|_{H^{-N-m}_{N_i, \Gamma_i}} \\
    &\lesssim \| v \|_{H^{-N}_{N_i, \Gamma_i - \delta i}} + \| P v \|_{H^{-N-m}_{N_i, \Gamma_i}}
\end{aligned}
\end{equation}
Similarly for $P^{\dagger}$. This gives the global estimate
\begin{equation}
\begin{aligned}
    &\| v \|_{H^s_{b, \Gamma}} \lesssim \| Pv \|_{H^{s-m}_{b, \Gamma}} + \| v \|_{H^{-N}_{b, \Gamma - \delta i}} \\
    & \| v \|_{H^{-s+m}_{b, \Bar{\Gamma}}} \lesssim \| P^{\dagger}v \|_{H^{-s}_{b, \Bar{\Gamma}}} + \| v \|_{H^{-N}_{b, \Bar{\Gamma} - \delta i}}
\end{aligned}
\end{equation}
Proposition \ref{prop:InclusionBWeightedSobolev} implies that $H^{r}_{b, \Gamma} \rightarrow H^{-N}_{b, \Gamma - i\delta}$ is compact for $-N < r$ and $\delta > 0$ and that these two estimates take place on dual spaces, hence the estimate above is a Fredholm estimate for $P$ which shows that it is Fredholm on $H^s_{b, \Gamma}$. 
$\Box$
\vspace{0.2in}
\newline
The analogous theorem for non-elliptic operator has the added dynamical assumptions on the Hamiltonian flow of $p$, but otherwise is similar. Readers familiar with the microlocal non-elliptic Fredholm theory will recognize the assumptions. The ones on the contours are exactly the assumptions needed to apply our estimate near the boundary in proposition \ref{prop:RadialPointEstimates}. This is one of those results where an image is worth a thousand words: the assumptions on the flow are thus illustrated in figure \ref{fig:HamiltonianFlowFredholm}.
\begin{theorem}
\label{theo:NonEllipticOperatorsAreFredholm}
    Let $P + iQ \in \Psi^m_{b, \Gamma}$ be a classical pseudo-differential operator on a b-manifold $\Bar{M}$ acting on some vector bundle $V$. Suppose furthermore that $P+iQ$ and $\Gamma_i$ satisfies the hypothesis of proposition \ref{prop:RadialPointEstimates} for all $i$, i.e. roughly that $g_i$ is symbolic, $-H_p(\tau)/\tau g(\sigma) + p_1$ has a definite sign in each of the halves $\mp H_p(\tau) / \tau > 0$ within the characteristic set of $P$ near $\partial \Bar{M}$, and that $\sigma$ is elliptic there. Suppose that its indicial family on $\partial \Bar{M}_i$ $\widehat{N(P + iQ)}(\sigma):H^s_{\partial \Bar{M}_i, \sigma} \rightarrow H^{s-m + k}_{\partial \Bar{M}_i, \sigma}$ is invertible on $\Gamma_i$ for all $i$, where $\Gamma = \{ \Gamma_i \}_i$, and that it satisfies high-energy estimates there of the form
    \begin{equation}
        \| u \|_{H^s_{\partial \Bar{M}_i, \sigma}} \lesssim \left< \sigma \right>^{k-m} \| \widehat{N(P + iQ)}(\sigma) \|_{H^{s-m+k}_{\partial \Bar{M}_i, \sigma}} \quad \sigma \in \Gamma_i \quad k \geq 1
    \end{equation}
    Finally, suppose that within the region where the principal symbol of $Q$ is non-positive $q \leq 0$ ($q \geq 0$) all the integral curves of the rescaled Hamiltonian vector field $H_p$ on the $b-$cosphere bundle $S^{\ast}_b \Bar{M}$\footnote{Think of this as the usual cosphere bundle within the interior but also including the singular form $d\tau / \tau$ at the boundary.} within the characterstic set $char(p) = p^{-1}(0)$ of $P$ tend backward (forward) to the boundary $S^{\ast}_b \Bar{M} |_{\partial \Bar{M}}$ all while staying in some open neighborhood of $q \leq 0$ ($q \geq 0$) and landing within the region
    \begin{equation}
        \mathcal{U}_{+} = \{ (\sigma, \tau, x, \xi); p(\sigma, \tau, x, \xi) = 0 \quad sgn(-H_p(\tau) / \tau) = sgn(-(H_p(\tau)/\tau) g(\sigma) + p_1) \} \subset char(p) \subset S^{\ast}_b \Bar{M}
    \end{equation}
    and tend forward (backward) to the boundary $S^{\ast}_b \Bar{M} |_{\partial \Bar{M}}$ staying in some open neighborhood of $q \leq 0$ ($q \geq 0$) and landing within the region
    \begin{equation}
        \mathcal{U}_{-} = \{ (\sigma, \tau, x, \xi); p(\sigma, \tau, x, \xi) = 0 \quad -sgn(-H_p(\tau) / \tau) = sgn(-(H_p(\tau)/\tau) g(\sigma) + p_1) \} \subset char(p) \subset S^{\ast}_b \Bar{M}
    \end{equation}
    Then $P +i Q$ is a Fredholm operator between the following spaces:
    \begin{equation}
    \begin{aligned}
        &P + iQ: \mathcal{X}^s_{\Gamma} \rightarrow \mathcal{Y}^{s-m+k}_{\Gamma} \\
        &\mathcal{X}^s_{\Gamma} = \{ u \in H^s_{b, \Gamma}; \quad (P+iQ) u \in H^{s-m+k}_{b, \Gamma} \} \quad \mathcal{Y}^{s-m+k}_{\Gamma} =H^{s-m+k}_{b, \Gamma}
    \end{aligned}
    \end{equation}
    If $Q = 0$, one can instead only assume that integrals curves of $H_p$ tend to $\mathcal{U}_{+}$ along any direction, i.e. either backward (in which case it must tend forward to $\mathcal{U}_{-}$) or forward (in which case it must tend backward to $\mathcal{U}_{-}$). Then $P$ is Fredholm on the same spaces. Furthermore, the Fredholm estimate for $P+iQ$ is uniform over $P$ and $Q$, in the sense that if $P_{\theta}+iQ_{\theta}$ is a bounded family of operators satisfying our assumptions uniformly \footnote{In the sense that both the upper and lower bounds we require are uniform; see theorem \ref{theo:NonEllipticOperatorsAreFredholm} for more details. Note that we also require the high-energy estimates to be uniform over $\theta$, hence in particular that resonances of $P_{\theta} + iQ_{\theta}$ never cross $\Gamma$ as $\theta$ varies.}, then the Fredholm estimates for $P_{\theta}+iQ_{\theta}$ are uniform over $\theta$.
\end{theorem}
\textbf{Remark:} For the wave-type operator on $\Bar{M} = [0, 1] \times \Sigma$ with $P = \Box$ hence $p \sim \sigma^2 - \xi^2$, $-H_p(\tau) / \tau = 2\sigma$ and thus the global condition on the flow of $H_p$ is mostly automatic, since future directed covectors $\sigma > 0$ are gonna to tend forward along $H_p$ to the future $\{ 1 \} \times \Sigma$, and past-directed covectors $\sigma < 0$ are going to tend forward along $H_p$ towards the past $\{ 0 \} \times \Sigma$. So as long as $\Gamma$ is chosen to satisfy proposition \ref{prop:RadialPointEstimates}, we have these assumptions satisfied. The high-energy estimate on the other hand is highly nontrivial. In our typical applications, the indicial family $\widehat{N(\Box)}(\sigma)$ will resemble the spectral family of a self adjoint operator on $\partial \Bar{M}$, in which case we have estimates on any strips avoiding the real line $\sigma \in \mathbb{R}$, which is sufficient for our applications. This is discussed in more details in section \ref{sec:PropagatorsWave}.   
\vspace{0.2in}
\newline
\textbf{Proof:} The proof proceeds the same way as all previous instances of non-elliptic Fredholm theorems following the general technique in \cite{vasy2011microlocalanalysisasymptoticallyhyperbolic}, hence we will be concise. Our assumptions on $P+iQ$ and $\Gamma$ means that we have a 'free' estimate for $u$ in terms of $(P+iQ)u$ in some neighborhood of the boundary $N$ and within $\mathcal{U}_{+}$:
\begin{equation}
\label{eq:FreeEstimate}
    \| Bu \|_{H^s_{b, \Gamma}} \lesssim \| (P+iQ)u \|_{H^{s-m+1}_{b, \Gamma}} + \| u \|_{H^{-N}_{b, \Gamma}} \quad B \text{ elliptic on } \mathcal{U}_{+} \cap S^{\ast}_b M|_{N} 
\end{equation}
In $\mathcal{U}_{-}$, we instead have an estimate of the form
\begin{equation}
    \| Au \|_{H^s_{b, \Gamma}} \lesssim \| \chi E u \|_{H^{s}_{b}} + \| (P+iQ)u \|_{H^{s-m+1}_{b, \Gamma}} + \| u \|_{H^{-N}_{b, \Gamma}} \quad A \text{ elliptic on } \mathcal{U}_{-} \cap S^{\ast}_b M|_{N}
\end{equation}
where $\chi$ is supported away from the boundary $\partial \Bar{M}$. In particular, all integral curves of $H_p$ within the characteristic set $char(p)$ starting from its support $supp(\chi)$ eventually tend backward (in $q \leq 0$)/forward (in $q \geq 0$) to $\mathcal{U}_{+} \cap S^{\ast}_b M|_{N}$ in finite time, where we have the 'free' estimate \ref{eq:FreeEstimate}. We can thus apply the propagation of singularity theorem with imaginary part $Q$ (note that our sign assumptions on $Q$ allow us to apply it in the appropriate direction) to bound $\| \chi E u \|_{H^{s}_{b}}$ by $\| Bu \|_{H^s_{b, \Gamma}} $ from above. Similarly, within the interior of $\Bar{M}$ and away from the boundary, all integral curves of $H_p$ within $char(P)$ tend backward/forward to $\mathcal{U}_{+} \cap S^{\ast}_b M|_{N}$, and again by propagation of singularity and equation \ref{eq:FreeEstimate} this bounds $\| A u \|$ by $\| P u \|$ within the characteristic set. By using the elliptic estimates \ref{prop:EllipticEstimates} on the elliptic set of $P$ and a partition of unity argument, we obtain the global estimate:
\begin{equation}
    \| u \|_{H^s_{b, \Gamma}} \lesssim \| (P+iQ)u \|_{H^{s-m+1}_{b, \Gamma}} + \| u \|_{H^{-N}_{b, \Gamma}}
\end{equation}
Doing the same for the adjoint $P^{\dagger} - iQ^{\dagger}$ on the dual space $H^{-s+m-1}_{b, \Bar{\Gamma}}$, but instead propagating along $H_p$ in the opposite direction (which we can by our assumptions on the reverse flow; note that if $P+iQ$ and $\Gamma_i$ satisfies the hypothesis of proposition \ref{prop:RadialPointEstimates}, then $P^{\dagger} - iQ^{\dagger}$ also satisfies them with the dual contour $\Bar{\Gamma}_i$), gives the estimate
\begin{equation}
    \| u \|_{H^{-s+m-1}_{b, \Bar{\Gamma}}} \lesssim \| (P^{\dagger} - iQ^{\dagger})u \|_{H^{-s}_{b, \Bar{\Gamma}}} + \| u \|_{H^{-N}_{b, \Bar{\Gamma}}}
\end{equation}
The rest of the proof then proceeds as in the elliptic case in proposition \ref{theo:EllipticOperatorsAreFredholm}, in that we use the invertibility of the normal operator $N(P+iQ)$ on $\Gamma$ to improve the remainder term, showing a Fredholm estimate on $H^s_{b, \Gamma}$ with dual space $H^{-s}_{b, \Bar{\Gamma}}$. $\Box$
\vspace{0.2in}
\newline
The main class of operators we wish to apply this theorem to are the wave-type operators $\Box$ on a Lorentzian manifold $[0, 1] \times \Sigma$. We will do so in section \ref{sec:PropagatorsWave}, which will also show that Dirac operators are Fredholm on our spaces. 
\newpage
 
\section{Feynman Parametrices and Indexes of Dirac Operators}
\label{sec:IndexDirac}
Now that we have a working Fredholm theory for non-elliptic operators on (compactified) spacetimes of the form $[0, 1] \times \Sigma$, we can use this to define the index of Dirac operators on them and compute it. We begin by applying theorem \ref{theo:NonEllipticOperatorsAreFredholm} to wave-type operators on such spacetimes in section \ref{sec:PropagatorsWave}, then recall the definition of Dirac operators and their index in section \ref{sec:DiracOperators} and finally use a Wick rotation argument to relate their index to the elliptic one in section \ref{sec:WickRotation}.

\subsection{Propagators of Wave Operators}
\label{sec:PropagatorsWave}
\begin{figure}[h!]
    \centering
    \includegraphics[width=0.75\linewidth]{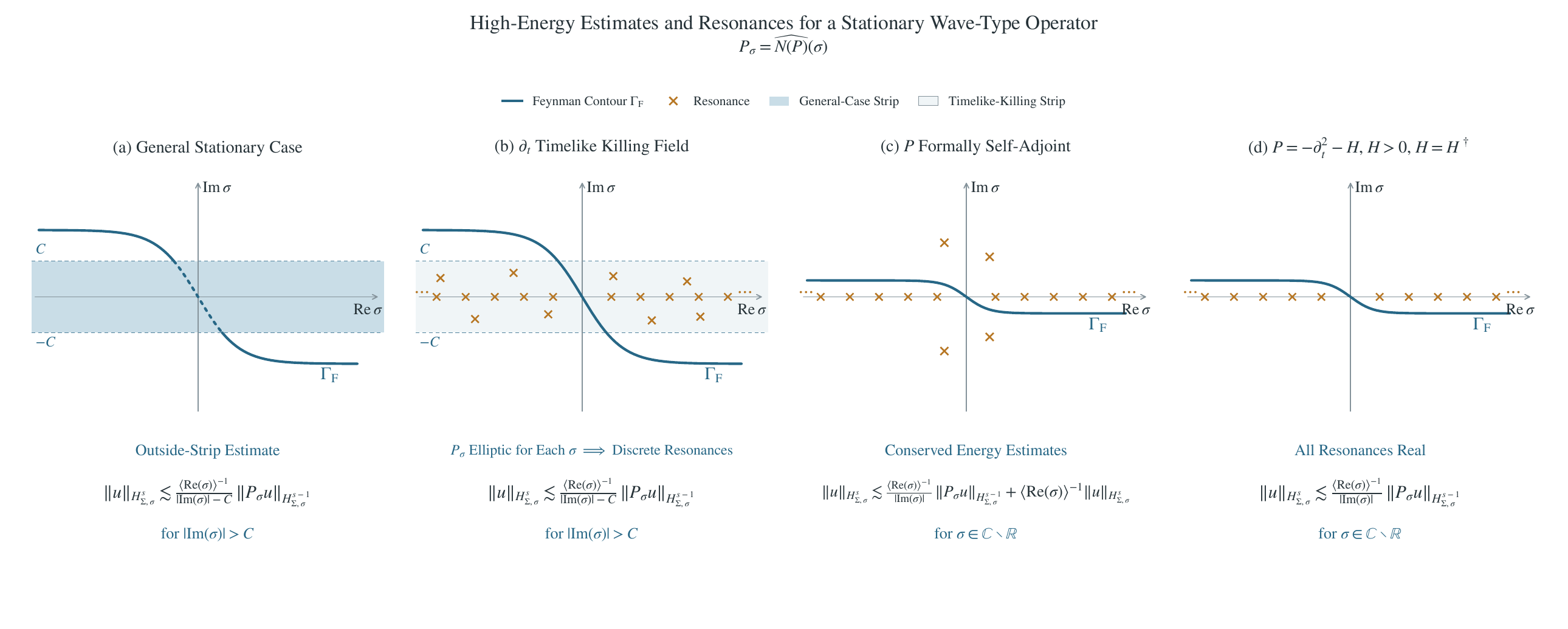}
    \caption{An illustration of the regions of $\mathbb{C}$ where we can obtain the high-energy estimates required to apply theorem \ref{theo:PropagatorConstruction}. In the general stationary case, we are missing esimates within a strip $|Im(\sigma)| \leq C$ in order to 'connect' the tails of $\Gamma$ to each other; furthermore, $\sigma$ is not elliptic within the null covectors. If $\partial_t$ is timelike, we are guaranteed to have the required control, while if $P$ is formally self-adjoint with respect to some positive metric, we can take $\Gamma$ as close as we want to the real line. Furthermore, the poles then come in conjugate pairs. Finally, the best case scenario is if the wave operator is generated by a self-adjoint elliptic positive operator $H$, in which case the resonances are identified with its spectrum and on the real line. Note that here we illustrated the case where the resonances are concentrated on $\mathbb{R}$, but in the presence of a subprincipal term $p_1$ they could be shifted to another line $Im(\sigma) = l$.}
    \label{fig:HighEnergyEstimates}
\end{figure}

We now discuss how the spaces $H^s_{b, \Gamma}$ allow us to define propagators, more specifically Fredholm inverses of wave-type operators, for exponentially decaying perturbations of stationary spacetimes. We start by quickly explaining how globally hyperbolic spacetimes fit into the compactified picture of section \ref{sec:PerturbationStationary}, then move on to apply theorem \ref{theo:NonEllipticOperatorsAreFredholm} to that case.
\begin{definition}
    A globally hyperbolic spacetime with compact Cauchy surface is a Lorentzian manifold $(M, g)$ of the form
    \begin{equation}
        M = \mathbb{R}_t \times \Sigma \quad dt \text{ timelike on $M$ and $(\Sigma, h)$ a closed Riemannian manifold}
    \end{equation}
    In this decomposition, the metric takes the form
    \begin{equation}
        g = N^2 dt^2 - h_{ij}(dx^i + \beta^i dt)(dx^j + \beta^j dt) \quad \text{$N$ is called the lapse function and $\beta$ the shift}
    \end{equation}
    It is stationary if the vector field $\partial_t$ is Killing, i.e. none of the variables above depend on time $t$, and static if $\beta = 0$. The exponential compactification of $M$ is the compact manifold with boundary $\Bar{M}$ obtained by compactifying the time interval $\mathbb{R}_t$ as follows \footnote{One could of course consider any power $e^{\pm \alpha t}$ as valid, which would change the smooth structure on $\Bar{M}$. In practice, say if one wishes for $g$ to be smooth on $\Bar{M}$, one can choose the exponent $\alpha$ from the expansion of the metric $g$ at $t = \pm \infty$.}:
    \begin{equation}
        \mathbb{R}_t \hookrightarrow [0, 1] \quad t \rightarrow \rho(t) \quad  \rho(t)|_{(-\infty, -1)} = e^{t} \text{ and } \rho(t)|_{(1, \infty)} = 1 - e^{-t}
    \end{equation}
    In particular, in terms of $t$, the defining functions in $\Bar{M}$ for the past $\{ 0 \} \times \Sigma$ and future $\{ 1 \} \times \Sigma$ are respectively
    \begin{equation}
        \tau_0 = e^{t} \quad \tau_1 = e^{-t}
    \end{equation}
    hence $dt = d \tau_0 / \tau_0$ at the past and $dt = -d \tau_1 / \tau_1$ at the future.
\end{definition}
Of course, after compactifying, the metric $g$ on $M$ need not extend to the compactification $\Bar{M}$ in any nice way. We will thus restrict our attention to metrics $g$ on $M$ whose asymptotics at $t = \pm \infty$ makes it well behaved on $\Bar{M}$. Before talking about what that imposes on $g$, let us make precise what we mean by 'wave operators' associated to $g$.
\begin{definition}
\label{def:wavetype}
    A wave-type operator $P$ on some vector bundle $V$ of $M$ is a $2^{nd}$ order differential operator whose principal symbol is given by the inverse metric:
    \begin{equation}
        p(\zeta) = g^{-1}(\zeta, \zeta) Id
    \end{equation}
    In particular, its characteristic set $p^{-1}(0)$ is given by the null covectors $g^{-1}(\zeta, \zeta) = 0$. The timelike condition on $dt$ guarantees that any non-zero null covector $g^{-1}(\zeta, \zeta) = 0$ has non-zero pairing with $dt$. We say it is future directed if $g^{-1}(\zeta, dt) > 0$, and past directed if $g^{-1}(\zeta, dt) < 0$. This globally splits the characteristic set minus the $0$ section $p^{-1}(0)\backslash \{ \zeta = 0 \}$ into two connected components $\pm g^{-1}(\zeta, dt) > 0$, which agrees up to a sign with the splitting $\mp H_p(\tau)/\tau > 0$ we used in proposition \ref{theo:NonEllipticOperatorsAreFredholm}:
    \begin{equation}
        -\frac{H_p(\tau)}{\tau} = \pm 2 g^{-1}(\zeta, dt) \quad \text{ near $t = \pm \infty$, i.e. at the past and future}
    \end{equation}
    Finally, note that if the Killing vector field $\partial_t$ is timelike, then its dual coordinate $\sigma$ is elliptic within the characteristic set of null covectors.
\end{definition}
  We note here a change of convention: in the general $b$-manifold theory of section \ref{sec:PerturbationStationary}, $\sigma$ near a boundary component refers to the dual coordinate of $-\tau \partial_{\tau}$, which in the past is $-\partial_t$ rather than $\partial_t$. We switch now as it is more convenient to work with the global $\partial_t$; the two conventions only differ by a sign. Thus, the contours $\Gamma = \{ \Gamma_0, \Gamma_1 \}$ will now denote spaces on the exponential compactification $\Bar{M}$ defined by asymptotics
  \begin{equation}
  \begin{aligned}
       u|_{\{ 0 \} \times \Sigma} = \frac{1}{\sqrt{2\pi}} \int_{\Gamma_{0}} \tau_0^{i\sigma} \widehat{u}(\sigma) \frac{d\tau}{\tau} \quad u|_{\{ 1 \} \times \Sigma} = \frac{1}{\sqrt{2\pi}} \int_{\Gamma_{1}} \tau_1^{-i\sigma} \widehat{u}(\sigma) \frac{(-d\tau)}{\tau}
  \end{aligned}
  \end{equation}
  which in terms of $t$ is just the usual Fourier-Laplace transform. Thus the two conventions are just related by $\sigma \rightarrow -\sigma$ in the past, or in terms of contours $\phi_0 \rightarrow \sigma - ig_0(-\sigma)$. This will also allow us to take $\Gamma_0 = \Gamma_1$ in typical problems, with positive imaginary part corresponding to decay in the future and to growth in the past. To analyze wave operators with the theory of this paper, we need them to be in the class $\Psi^2_{b, \Gamma}$ on $\Bar{M}$ so that we can apply to them the estimates of section \ref{sec:PerturbationStationary}. More precisely, given a choice of contours $\Gamma$ determining the asymptotics at $t = \pm \infty$, we wish for the wave-type operator to be in $\Psi^m_{b, \Gamma}$. This already imposes some constraints on the metric $g$. It terms of the time function $t$, requiring $P \in \Psi^2_{b, \Gamma}$ at the top-order term gives the following asymptotic expansion for $g^{-1}$ in a neighborhood of $t = \pm \infty$:
\begin{equation}
\label{eq:MetricExpansion}
\begin{aligned}
    &g^{-1}|_{\{ \pm \infty \} \times \Sigma} = g_0^{-1} + Re\left( \sum_{l, k; Im(\eta_l) > 0} e^{ i\eta_l |t|} (|t|)^k h_{lk} + e^{-(\mu + \delta) |t|} h_{\delta} \right) \quad h_{lk}, h_{\delta} \in Sym^2(T_b \Bar{M}) \\
    &g^{-1}_0, h_{lk} \text{ stationary}
\end{aligned}
\end{equation}
where $g_0$ is a stationary Lorentzian metric independent of time, $Sym^2(T_b \Bar{M})$ denotes symmetric 2-cotensor and the exponents $\eta_l$ and $- \Bar{\eta}_l$ are above $\Gamma_i - \Gamma_i$ in the complex plane. Note that the stationary model $g_{0}^{-1}$ and the exponents $\eta_l$ can in general differ at $t = \pm \infty$; this is particularly relevant for spacetimes which evolve from one stationary state to another. Lower-order terms must still admit an appropriate expansion at the past and future, and because those vary quite a bit with the operator and bundle $V$ \footnote{They can for example include curvature terms or arbitrary potentials.}, we instead simply fix a wave-type operator $P \in \Psi^{2}_{b, \Gamma}$ for the remainder of the discussion. In that case, a simple application of our Fredholm theorem \ref{theo:NonEllipticOperatorsAreFredholm} allow us to construct Fredholm inverses for $P$, or $P+iQ$ if $Q$ is some real $2^{nd}$ order operator satisfying some sign conditions. The spaces on which they are Fredholm will be as defined in theorem \ref{theo:NonEllipticOperatorsAreFredholm}, i.e.
\begin{equation}
    \begin{aligned}
        &P + iQ: \mathcal{X}^s_{\Gamma} \rightarrow \mathcal{Y}^{s-2+k}_{\Gamma} \quad k \geq 1 \\
        &\mathcal{X}^s_{\Gamma} = \{ u \in H^s_{b, \Gamma}; \quad (P+iQ) u \in H^{s-2+k}_{b, \Gamma} \} \quad \mathcal{Y}^{s-2+k}_{\Gamma} =H^{s-2+k}_{b, \Gamma}
\end{aligned}
\end{equation}
Depending on which contours $\Gamma$ we pick, we then obtain Fredholm inverses that behave qualitatively differently.
\begin{theorem}
\label{theo:PropagatorConstruction}
    Let $\Bar{M} = [0, 1] \times \Sigma$ be the exponential compactification of a globally hyperbolic spacetime with compact Cauchy surface, and let $P$ be a wave-type operator on some vector bundle $V$ and $Q$ a $2^{nd}$ order b-pseudo-differential operator. Let $\Gamma = \{ \Gamma_0, \Gamma_1 \}$ be choices of contours for which $g_i$ is symbolic of order $0$, $g^{-1}$ has the expansion \ref{eq:MetricExpansion} and for which $P + i Q \in \Psi^2_{b, \Gamma}(\Bar{M})$. Suppose the indicial family of the stationary model $\widehat{N(P + iQ)}(\sigma)$ is invertible on $\Gamma_0, \Gamma_1$ at the past $\{ 0\} \times \Sigma$ and future $\{ 1 \} \times \Sigma$, and that we have high energy estimates on them of the form
    \begin{equation}
    \label{eq:HighEnergyEstimateWave}
        \| u \|_{H^s_{\partial \Bar{M}_i, \sigma}} \lesssim \left< \sigma \right>^{k-2} \| \widehat{N(P + iQ)}(\sigma) \|_{H^{s-m+k}_{\partial \Bar{M}_i, \sigma}} \quad \sigma \in \Gamma_i \quad k \geq 1
    \end{equation}
    Suppose furthermore that $\sigma$ is elliptic within the characteristic sets of the stationary models $N(P)$ at both the past an future, i.e. the stationary models for $g$ at the past and future have their Killing vector $\partial_t$ timelike. Then, depending on the relative signs of $g^{-1}(\zeta, dt) g_i(\sigma)$ and $Q$, we obtain different Fredholm inverses for $P + iQ$:
    \begin{itemize}
        \item Near the past $\{ 0 \} \times \Sigma$, suppose $2g^{-1}(\zeta, dt) g_0(\sigma)Id + p_1$ is negative definite within $g^{-1}(\zeta, dt) > 0$, positive definite in $g^{-1}(\zeta, dt) < 0$ and near the future $\{ 1 \} \times \Sigma$, $2g^{-1}(\zeta, dt) g_1(\sigma)Id + p_1$ is negative definite within $g^{-1}(\zeta, dt) > 0$ and positive definite in the other half. Furthermore, we assume that $\pm Q$ is globally non-positive in a neighborhood of $p^{-1}(0)$ up to terms of order $0$ \footnote{We really only need the subprincipal symbol of $Q$ to vanish near the boundary so that it doesn't affect the radial point estimate; in the interior we can impose the sign only on the principal symbol of $Q$.} in $\pm g^{-1}(\zeta, dt) > 0$, as defined in proposition \ref{prop:RadialPointEstimates}. Then $P+iQ: \mathcal{X}_{\Gamma}^s \rightarrow \mathcal{Y}_{\Gamma}^{s-2+k}$ is Fredholm with parametrix $E: \mathcal{Y}_{\Gamma}^{s-2+k} \rightarrow \mathcal{X}_{\Gamma}^s$ satisfying \footnote{$WF_{\Gamma}$ denotes the $\Gamma-$wavefront set, i.e. some covector $\zeta \notin WF_{\Gamma}(u)$ iff $Au \in H^{\infty}_{b, \Gamma}$ for some $b-$operator $A \in \Psi^0_{b, \Gamma}$ elliptic at $\zeta$. }
        \begin{equation}
            WF_{\Gamma}(Ef)|_{{0} \times \Sigma}  \subset WF_{\Gamma}(f)|_{{0} \times \Sigma} \quad  WF_{\Gamma}(Ef)|_{{1} \times \Sigma}  \subset WF_{\Gamma}(f)|_{{1} \times \Sigma} \cup p^{-1}(0)|_{{1} \times \Sigma}
        \end{equation}
        This corresponds to a 'retarded' parametrix since the solution is 'trivial', in the sense of decay and its high-frequency content, in the past. Similarly, reversing the role of the past $\{ 0 \} \times \Sigma$ and the future $\{ 1 \} \times \Sigma$ in these assumptions give also a Fredholm problem, corresponding to the 'advanced' parametrix. Note that the retarded and advanced parametrix are dual to each other.
        \item Near the past $\{ 0 \} \times \Sigma$ and future $\{ 1 \} \times \Sigma$, $2g^{-1}(\zeta, dt) g_0(\sigma)Id + p_1$ and $2g^{-1}(\zeta, dt) g_1(\sigma)Id + p_1$ are negative definite in the whole characteristic set $p^{-1}(0)$. Furthermore, $Q$ is globally non-positive in a neighborhood of $p^{-1}(0)$ up to terms of order $0$, as defined in proposition \ref{prop:RadialPointEstimates}. Then $P+iQ: \mathcal{X}_{\Gamma}^s \rightarrow \mathcal{Y}_{\Gamma}^{s-2+k}$ is Fredholm with parametrix $E: \mathcal{Y}_{\Gamma}^{s-2+k} \rightarrow \mathcal{X}_{\Gamma}^s$ satisfying
        \begin{equation}
            WF_{\Gamma}(Ef)|_{{0} \times \Sigma}  \subset WF_{\Gamma}(f)|_{{0} \times \Sigma} \cup \{g^{-1}(\zeta, dt) < 0 \} \quad  WF_{\Gamma}(Ef)|_{{1} \times \Sigma}  \subset WF_{\Gamma}(f)|_{{1} \times \Sigma} \cup \{g^{-1}(\zeta, dt) > 0 \} 
        \end{equation}
        This corresponds to a 'Feynman' parametrix since the solution is 'trivial' in two different halves of the characteristic sets at the past and future. Similarly, reversing the role of the past $\{ 0 \} \times \Sigma$ and the future $\{ 1 \} \times \Sigma$ in these assumptions give also a Fredholm problem (which ammounts to reversing the signs here), corresponding to the 'anti-Feynman' parametrix. Note that the Feynman and anti-Feynman parametrix are dual to each other.
    \end{itemize}
    Again, if $P$ and $Q$ vary boundely over a one-parameter family of operators $P_{\theta} + iQ_{\theta}$ satisfying the assumptions uniformly, then the corresponding Fredholm estimates are uniform over $\theta$ too. 
\end{theorem}
\textbf{Proof:} This is essentially translating the assumptions of theorem \ref{theo:NonEllipticOperatorsAreFredholm} to the case of wave equations, and using the remark in definition \ref{def:wavetype} about their characteristic set and hamiltonian flow. The main thing to check is the global behavior of the Hamiltonian flow, which here is just the lifted geodesic flow on the cosphere bundle $S^{\ast}_b \Bar{M}$. This is where global hyperbolicity comes in: $g^{-1}(\zeta, dt)$ is \textit{globally} elliptic within the characteristic set, hence in particular any null geodesics eventually converge to the boundary $\partial \Bar{M}$, with future directed covectors $g^{-1}(\zeta, dt) > 0$ converging to the future $\{ 1 \} \times \Sigma$ in the forward direction, and in the past backward. Similarly for past directed covectors $g^{-1}(\zeta, dt)$. Furthermore, our assumptions on the sign of $2 g^{-1}(\zeta, dt) g_1(\sigma)Id + p_1$ guarantees it ends up in the appropriate regions $\mathcal{U}_{\pm}$ from proposition \ref{theo:NonEllipticOperatorsAreFredholm}, and similarly our assumptions on the sign of $Q$ guarantee we can use propagation of singularity in the right direction. This show Fredholmness. As for the wavefront set statements, it follows from our estimate in proposition \ref{prop:RadialPointEstimates} applied with $\varrho = g^{-1}(\sigma, dt)$. Note that $\pm g^{-1}(\sigma, dt) > 0$ is the exact same condition as $\pm \sigma > 0$ in our case since $\sigma$ is assumed elliptic. $\Box$
\vspace{0.2in}
\newline
We are now left with the task of showing the existence of nontrivial examples of our theory. As alluded previously in theorem \ref{theo:NonEllipticOperatorsAreFredholm}, the main assumption to verify is the existence of high energy estimates for $\widehat{N(P+iQ)}(\sigma)$ on contours $\Gamma_0$, $\Gamma_1$. Thankfully, this has been shown in the past to exist for a wide class of $\widehat{N(P+iQ)}(\sigma)$; these are often called \textit{resolvent estimates} in the literature since in the ultra-static case with $N = 1$ and $\beta = 0$, they correspond to estimating the resolvent $(\sigma^2 - \Delta)^{-1}$ within the complex plane of some Laplace operator. In general, the indicial family is given by
\begin{equation}
    \widehat{N(P+iQ)}(\sigma) = \frac{(\sigma + i\beta^j \partial_{x_j})^2}{N^2} - \Delta_h + \text{ lower-order terms } + i\widehat{Q}(\sigma)
\end{equation}
where $\Delta_h$ is the Laplace operator associated to the Riemannian metric $h$. The following are the main examples where it is known that we have high-energy estimates. We focus on the case $Q=0$ for now. The estimates are summarized in figure \ref{fig:HighEnergyEstimates}.
\begin{itemize}
    \item Perturbations of General Stationary Spacetimes: In the general stationary case, high-energy estimates are known to exist outside a large enough strip $|Im(\sigma)| \leq C$ depending on the operator $P$. This corresponds to the well-known fact that we always have exponential bounds on the growth of the energy of solutions to wave equations; see for example \cite{B_r_2015} Theorem 8. In terms of high-energy estimates, it gives the uniform estimate
    \begin{equation}
    \label{eq:StationaryHighEnergy}
        \| u \|_{H^s_{\partial \Bar{M}_i, \sigma}} \lesssim \frac{1}{\left< Re(\sigma) \right> (|Im(\sigma)| - C)} \| \widehat{N(P)}(\sigma) u \|_{H^{s-1}_{\partial \Bar{M}_i, \sigma}} \quad |Im(\sigma)| > C >> 1
    \end{equation}
    In particular, $\widehat{N(P)}(\sigma)$ is invertible outside the strip. As for what happens within the strip $|Im(\sigma)| \leq C$, it heavily depends on the operator $P$. Even in the stationary case, without additional information, we cannot obtain the high-energy estimates required for the Feynman parametrices since the corresponding contours $\Gamma$ would need to pass through the strip $|Im(\sigma)| \leq C$ \footnote{Note that this isn't an issue for the retarded and advanced parametrices, since the corresponding contours can lie above or below the strip.} hence we need some estimates within it. Any estimate along some vertical line $L$ crossing the strip would work for example. On the other end, if we knew that the set of poles of $\widehat{N}(\sigma)^{-1}$ was discrete, then we could pick $\Gamma$ avoiding these poles and which for large enough real frequency $|Re(\sigma)|$ lies in $|Im(\sigma)|>C$, in which case we could then use estimate \ref{eq:StationaryHighEnergy} to obtain the global high-energy estimates. Showing discreteness of the poles is typically done using the analytic Fredholm theorem. In general, $\widehat{N(P)}(\sigma)$ is a non-elliptic family of operator, hence showing that they are Fredholm can be attempted with the general technique of \cite{vasy2011microlocalanalysisasymptoticallyhyperbolic} paired with dynamical assumptions. The best case scenario is when $\widehat{N(P)}(\sigma)$ are elliptic operators, which is our next example. It also happens to coincide with our assumption that $\sigma$ should be elliptic within the characteristic set of $P$.
    \item Perturbations of Stationary Spacetimes with Timelike Killing Vector Field: If the Killing vector field $\partial_t$ is \textit{timelike}, i.e. $g(\partial_t, \partial_t) > 0$ everywhere which requires $N^2 > h_{ij} \beta^i \beta^j$, then the operators $\widehat{N(P)}(\sigma)$ are elliptic hence automatically Fredholm between semiclassical Sobolev spaces. This shows that the poles within the strip $|Im(\sigma)|\leq C$ are discrete, in which case any contour $\Gamma$ that avoid them and which lies outside the strip for sufficiently large $|Re(\sigma)|$ is going to satisfy the high energy estimate \ref{eq:HighEnergyEstimateWave} with $k = 1$. On the other end, it doesn't give by itself estimates within the strip. This shows that for such spacetimes, we can construct Feynman parametrices using theorem \ref{theo:PropagatorConstruction}. Regions where $g(\partial_t, \partial_t)$ changes from timelike to spacelike are usually called \textit{ergoregions}, so our assumptions exclude those. The next best thing is if the all the resonances except finitely many are on the real line, which involves some symmetry assumption on $P$.
    \item Perturbations of Formally Self-Adjoint Wave Operators: Still assuming that $\partial_t$ is timelike, now also suppose that $P$ is formally self-adjoint with respect to some \textit{positive definite} metric on the bundle $V$, i.e. $P = P^{\dagger}$. Then, $\widehat{N(P)}(\sigma)$ satisfies the following estimate in the complex plane $\mathbb{C}$ minus the real line:
    \begin{equation}
        \| u \|_{H^s_{\partial \Bar{M}_i, \sigma}} \lesssim \frac{1}{\left< Re(\sigma) \right> |Im(\sigma)|} \| \widehat{N(P)}(\sigma) u \|_{H^{s-1}_{\partial \Bar{M}_i, \sigma}} + \left< Re(\sigma) \right>^{-1} \| u \|_{H^{-N}_{\partial \Bar{M}_i, \sigma}} \quad \sigma \in \mathbb{C} \backslash \mathbb{R}
    \end{equation}
    In particular, it only has finitely many resonances outside the real line. The proof is quite simple and well-known. It corresponds to an energy estimate again. One uses that for such operators, $\widehat{N(P)}(\sigma)^{\dagger} = \widehat{N(P)}(\Bar{\sigma})$ hence we have the following identity (which follows from the second derivative $\partial_{\sigma}^2\widehat{N((P)}$ being real)
    \begin{equation}
    \begin{aligned}
        &Im \left< \widehat{N(P)}(\sigma) u, u \right> = Im(\sigma) \left< \partial_{\sigma} \widehat{N((P)}(Re(\sigma)) u , u \right> \\
        & \partial_{\sigma} \widehat{N((P)}(Re(\sigma)) = 2\frac{(Re(\sigma) + i \beta^j \partial_{x_j})}{N^2} + \text{ lower order terms}
    \end{aligned}
    \end{equation}
    The principal symbol of $\partial_{\sigma} \widehat{N((P)}(Re(\sigma))$ is thus $2g^{-1}(\zeta, dt)$ hence has a sign in each halves of the (semiclassical) characteristic set of $\widehat{N(P)}(\sigma)$, which paired with (semiclassical) elliptic estimates allows us to control $u$ using the identity above. Equivalently, one could also compute $Im \left< \Bar{\sigma} \widehat{N(P)}(\sigma) u, u \right> = Re \left< \widehat{N(P)}(\sigma) u, i\sigma u \right>$ which amounts to using the usual energy pairing $\left<Pu, \partial_t u \right>$. See \cite{Vasy2013Microlocal} section 3.2 and 3.3 for a similar computation in the semi-classical setting. In that case, one can take our contours $\Gamma$ arbitrarily close to the real line, as long as it stays at least a fixed distance $\epsilon$ away and avoids the resonances. A typical example is the wave-operator on scalar functions, which is formally self-adjoint with respect to the volume form induced by $g$. Furthermore, note that in that case the subprincipal term $p_1$ of $P$ can be chosen to vanish at the boundary $\tau = 0$ by taking the positive definite pairing with respect to which $P$ is self-adjoint.
    \item Perturbations of Wave Operators with a Self-Adjoint Generator: Finally, the simplest case to analyze is when the stationary model metrics are static and the wave operator $P$ is generated by a self-adjoint positive elliptic operator $H$ of order 2 on $\Sigma$ (with respect to a positive-definite pairing on the pullback bundle on $\Sigma$):
    \begin{equation}
        N(P) = (-\partial_t^2 - H) \quad H^{\dagger} = H \text{ positive }
    \end{equation}
    In that case, understanding $N(P)$ reduces to the spectral theory of $H$. In particular, by spectral theory we automatically get high-energy estimates away from the real line and $0$:
    \begin{equation}
         \| u \|_{H^s_{\partial \Bar{M}_i, \sigma}} \lesssim \frac{1}{|Im(\sigma)| \left<Re(\sigma)\right>} \| \widehat{N(P)}(\sigma) u \|_{H^{s-1}_{\partial \Bar{M}_i, \sigma}} \quad \sigma \in \mathbb{C} \backslash \mathbb{R}, |\sigma| \text{ large}
    \end{equation}
    All the resonances lie on $\mathbb{R}$ and are identified with $\sqrt{spec(H)} \cup -\sqrt{spec(H)}$.
\end{itemize}
The estimates when the complex part $Q$ is non-zero will be discussed in section \ref{sec:WickRotation} when we attempt to Wick rotate our operators, but roughly speaking these estimates can be preserved in certain regions of $\mathbb{C}$ if $Q$ satisfy some sign conditions. We summarize the examples in the following corollary, which follows from what we just discussed.
\begin{corollary}
    Let M be a globally hyperbolic spacetime whose stationary models $g_0$ at the past and future has Killing vector field $\partial_t$ timelike. Let $P$ be some wave-type operator in $\Psi^2_{b, \Gamma}$. Then there exists choices of contours $\Gamma$ for which $P:\mathcal{X}^s_{\Gamma} \rightarrow \mathcal{Y}^{s-1}_{\Gamma}$ admits Feynman and anti-Feynman parametrices in the sense of theorem \ref{theo:PropagatorConstruction}. If $P = P^{\dagger}$ with respect to some positive definite metric on $V$, then these contours can be taken to any $\epsilon$ distance away from the real line $\mathbb{R}$. 
\end{corollary}
We conclude this section by explaining under which conditions the Feynman parametrix $E$ of theorem \ref{theo:PropagatorConstruction} satisfies the stronger Feynman condition
\begin{equation}
\label{eq:StrongFeynman}
    \widehat{Ef}(\sigma > 0)|_{t = -\infty} = 0 \quad \widehat{Ef}(\sigma < 0)|_{t = +\infty} = 0
\end{equation}
and similarly for the retarded and advanced parametrix. In other words, when can we convert our weak decay and wavefront set condition of theorem \ref{theo:PropagatorConstruction} to an exact support condition. Rather than requiring only estimates for $\widehat{N(P)}^{-1}(\sigma)$ on the contours $\Gamma$, this corresponds to requiring $\widehat{N(P)}^{-1}(\sigma)$ to have estimates in sub-planes of $\mathbb{C}$ so that we can deform the contour to $Im(\sigma) = \pm \infty$ (or in the retarded/advanced propagator case, simply apply the Paley-Wiener-Schwartz theorem \ref{prop:PaleyWinerSchawrtz} to infer support properties of $u$ from exponential estimates on $\widehat{u}(\sigma)$). Roughly speaking, if we have global control of the meromorphic family $\widehat{N(P)}^{-1}(\sigma)$ and there is no resonances in the appropriate quadrants of $\mathbb{C}$ relative to $\Gamma$, the exact support conditions as in equation \ref{eq:StrongFeynman} hold interpreted in terms of mode expansions. Once again, this expansion in powers of $\tau$ is well-known on standard weighted Sobolev spaces $\tau^l H^s_b$; the resulting distributions are called polyhomogeneous and are developped in \cite{Melrose1993} section 5.10 while the argument for waves can also be found in \cite{vasy2011microlocalanalysisasymptoticallyhyperbolic} proposition 3.5.
\begin{proposition}
\label{prop:ModeExpansion}
    Let $P$ be a wave-type differential operator. Suppose that $\widehat{N(P)}^{-1}(\sigma)$ is a meromorphic family satisfying the high-energy estimate
    \begin{equation}
        \| u \|_{H^s_{\partial \Bar{M}_i, \sigma}} \lesssim \frac{1}{\left< Re(\sigma) \right> (|Im(\sigma)| - C)} \| \widehat{N(P)}(\sigma) u \|_{H^{s-1}_{\partial \Bar{M}_i, \sigma}} \quad |Im(\sigma)| > C > 0
    \end{equation}
    and $\widehat{N(P)}^{-1}(\sigma)$ exists on the contours $\Gamma_i$ for $i = 0, 1$. Suppose furthermore that the contours $\Gamma_i$ are eventually a positive distance away from $|Im(\sigma)| \leq C$, i.e. $|g_i(\sigma)| > C+\epsilon$ for $|\sigma| \in \mathbb{R}$ large enough. Finally, suppose that at both the past and future we have some estimates within the strip $|Im(\sigma)| \leq C$, in the sense that there are vertical lines $L_k = \{ Re(\sigma) = l_k \} \cap \{ |Im(\sigma)| \leq C \}$ in the strip with $l_{2k} \rightarrow \infty$, $l_{2k+1} \rightarrow -\infty$ on which $\widehat{N(P)}^{-1}(\sigma)$ satisfies some polynomial bound\footnote{This is only required for the mode expansion. In the Feynman and anti-Feynman case, we only need it for vertical lines in the appropriate quadrant of $\mathbb{C}$ corresponding to which parametrix we are considering. Also, any polynomial bounds between Sobolev spaces of any order would serve for this purpose, at the expense of a weaker regularity for the remainder.}:
    \begin{equation}
        \| \widehat{N(P)}^{-1}(\sigma) \|_{H^{s-1}_{\partial \Bar{M}_i, \sigma} \rightarrow H^{s}_{\partial \Bar{M}_i, \sigma}} \lesssim \left< l_k \right>^N \quad \sigma \in L_k
    \end{equation}
    Let $u \in H^s_{b, \Gamma}$ be such that $Pu = f \in H^s_{b, \Gamma}$ has compact support. Then the parametrices of theorem \ref{theo:PropagatorConstruction} satisfy the stronger support conditions:
    \begin{itemize}
        \item In the retarded/advanced case, i.e. $\Gamma_0, \Gamma_1$ below/above the strip, the support of the solution $u = E^{-1} f$ is in the future/past. More precisely, there are cutoffs $\chi_0$ and $\chi_1$ supported near $\{ 0 \} \times \Sigma$ and $\{ 1 \} \times \Sigma$ for which
        \begin{equation}
            \chi_0 u = 0 \text{ (retarded case) } \quad \chi_1 u = 0 \text{ (advanced case) }
        \end{equation}
        At the future (retarded) or past (advanced) and assuming the estimate along vertical lines, the solution also admits an expansion into powers of $\tau$ and $-\log(\tau)$, as we'll discuss for the Feynman case.
        \item In the Feynman case, i.e. $\Gamma_i$ eventually below $Im(\sigma) = -C$ for future-directed frequencies $Re(\sigma)>0$ and eventually above $Im(\sigma) = C$ for past-directed frequencies $Re(\sigma) < 0$, suppose there are no resonances below $\Gamma_0$ with $Re(\sigma) > 0$ and similarly no resonances above $\Gamma_1$ with $Re(\sigma) < 0$. Suppose furthermore that $P - N(P)$ admits an expansion into powers of $\tau$ at the past and future up to order $L_i > \mu_i = \sup |Im(\Gamma_i - \Gamma_i)|$ (by definition of $\Psi^m_{b, \Gamma}$ we have this up to order $\mu + \delta$ already):
        \begin{equation}
            P - N(P) = \sum_{l, k} \tau^{-i \eta_l}(-\log(\tau))^k P_{lk} + O(\tau^{L_i + \delta}) \quad \text{$\eta_l$ above $\Gamma_i - \Gamma_i$} \quad \text{near $\partial \Bar{M}_i$}
        \end{equation}
        Then there exists cutoffs $\chi_0$ and $\chi_1$ supported near $\{ 0 \} \times \Sigma$ and $\{ 1 \} \times \Sigma$ for which $\chi_i u$ admits expansions up to decay order $\Gamma_i \mp iL$ into past/future directed quasinormal modes corresponding to the resonances of $\widehat{N(P)}(\sigma)^{-1}$ and the powers $\eta_l$ \footnote{The infinite sum must be taken as a limit of partial sums over the finitely many resonances in the 'box' defined by the contours $\Gamma$, $\Gamma + iL$ and the vertical lines $L_k$, $L_{2k+1}$. The convergence is with respect to the original contour $\Gamma$ but in a weak regularity norm: $H^{-M}_{b, \Gamma}$.} near $\tau = 0$:
        \begin{equation}
        \begin{aligned}
            &\chi_0 u = \sum_{\sigma \text{ below } \Gamma_0, Re(\sigma) \leq 0} \sum_{\alpha, k} e^{i(\sigma - \alpha \cdot \eta) t}(-t)^{k} u^0_{\sigma,k,\alpha} + r_0 \\
            &\chi_1 u = \sum_{\sigma \text{ above } \Gamma_1, Re(\sigma) \geq 0} \sum_{\alpha, k} e^{i(\sigma + \alpha \cdot \eta)t}t^{k} u^1_{\sigma,k,\alpha} + r_1
        \end{aligned}
        \end{equation}
        where we used $\tau_0 = e^{t}$ and $\tau_1 = e^{-t}$ and $r_0$, $r_1$ are in some weak Sobolev space $H^{-M}_{b, \Gamma_i \mp iL}$ . The spatial distributions $u_{\sigma,k,\alpha}$ correspond to the resonant states of $\widehat{N(P)}(\sigma)^{-1}$ with finitely many of the operators $P_{lk}$ applied to them. Furthermore, $\alpha$ denotes a multi-index in $\mathbb{N}^M$, with $\alpha \cdot \eta = \sum \alpha_l \eta_l$, while the sum over $k$ depends on the resonance $\sigma$ and multi-index $\alpha$. Note that by our assumption on the exponents $\eta_l$, $\alpha \cdot \eta_l + \sigma$ is guaranteed to stay above $\Gamma_1$ (and similarly $\sigma - \alpha \cdot \eta$ stays below $\Gamma_0$), while if the $\eta_l$ satisfy $Re(\eta_l) > 0$, then the corresponding exponents $\sigma \pm \alpha \cdot \eta$ at the past/ future also have real part negative/positive. A similar result hold for the anti-Feynman case with the role of the past and future reversed.
    \end{itemize}
\end{proposition}
\textbf{Remark:} The resonance expansion also holds for any contour $\Gamma$ leading to a Fredholm problem as long as we have the appropriate estimate on $\widehat{N(P)}(\sigma)^{-1}$, and in general the resonances that appear in the expansion will be the ones above $\Gamma_1$ and below $\Gamma_0$. There is thus a 'finite-dimensional' choice of resonances to include or exclude at the past or future, corresponding to where exactly the contours $\Gamma$ cross the real axis $\mathbb{R}$. Furthermore, while apriori the remainders $r_0$, $r_1$ are in some weak Sobolev spaces, typically we can apply the strong form of our radial set estimates \ref{prop:RadialPointEstimates} plus propagation of singularity to show that they are in a more regular space.
\vspace{0.1in}
\newline
\textbf{Proof:} The techniques for obtaining these properties are well-known and works the exact same way on our space $H^s_{b, \Gamma}$, hence we rapidly go over them. Pick normal neighborhoods $N_0$ and $N_1$ of $\{ 0 \} \times \Sigma$ and $\{ 1 \} \times \Sigma$ such that $\chi_0 u \in H^s_{b, \Gamma_0}$ and $\chi_1 u \in H^s_{b, \Gamma_1}$. Then we have
\begin{equation}
    P \chi_i u = [P, \chi_i]u + \chi_i Pu
\end{equation}
Furthermore, $P = N(P) + O(\tau^{\nu})$ while $[P, \chi_i] u$ is compactly supported hence so is $\Tilde{f} = \chi_if + [P, \chi_i]u$. Thus we have
\begin{equation}
    N(P) \chi_i u = \Tilde{f} + O(\tau^{\nu}) \chi_i u
\end{equation}
Taking the Mellin transform gives
\begin{equation}
\label{eq:SolvingForMellinTransform}
    \widehat{\chi_i u}(\sigma) = \widehat{N(P)}^{-1}(\sigma)\left(  \widehat{\Tilde{f}} + \widehat{ O(\tau^{\nu}) \chi_i u} \right)(\sigma)
\end{equation}
We know that $\widehat{\Tilde{f}}$ is entire in $\mathbb{C}$ while by proposition \ref{prop:InclusionStationarySpaces} $\widehat{ O(\tau^{\nu}) 
\chi_i u}$ is holomorphic in the region under/above $\Gamma_i \pm i \nu$ at the future/past. In the retarded and advanced propagator case, by our assumptions on the poles of $\widehat{N(P)}^{-1}(\sigma)$, in either the past (retarded, $i = 0$) or future (advanced, $i = 1$), it is holomorphic between $\Gamma_i$ and $\Gamma_i \pm i\nu$ with estimates along vertical lines. Say we work in the retarded case, hence the past. We can thus deform the contour from $\Gamma_0$ to $\Gamma_0 - i\nu$ which shows that $\chi_i u \in H^s_{\Gamma_0 - i \nu}$ there. Doing this argument repeatedly gives that $\widehat{\chi_0 u}$ exists within the whole complex plane $\mathbb{C}$, hence is rapidly decaying. We can then apply energy estimates to conclude that $\chi_0 u = 0$; one uses multiplier of the form $\chi \tau^{-l} \partial_t$ for large $l$ to obtain an exact positive commutator estimate. Equivalently, as can be found in \cite{B_r_2015}, the usual energy estimates on a globally hyperbolic spacetimes give exponential bounds on the $H^s$ norms of solutions in the spatial slice $\Sigma$, hence if $u$ vanishes super-exponentially (in $t$) in a neighborhood of the past, then $\| u|_{t_1} \| \lesssim e^{C|t_1 - t_2|} \| u|_{t_2} \| \lesssim e^{|t_1 - t_2|} e^{-N|t_2|}$ gives $u$ vanishing by picking $N$ large and $t_2$ going to $-\infty$. As for the Feynman case, working at the future, our assumptions on the norm of $\widehat{N(P)}(\sigma)^{-1}$ along vertical lines now translate to deforming the contour from $\Gamma_1$ to $\Gamma_1 + i \nu$ by taking a limit over 'squared' closed contours
\begin{equation}
    L_{2k} \cup L_{2k+1} \cup (\Gamma_1 \cap \{\sigma : l_{2k+1} \leq Re(\sigma) \leq l_{2k} \}) \cup (\Gamma_1 + i\nu \cap \{\sigma : l_{2k+1} \leq Re(\sigma) \leq l_{2k} \})
\end{equation}
In the sense of distributions (pairing with holomorphic Mellin transform $\widehat{\psi}$ decaying rapidly in $Re(\sigma)$), the integrals along the vertical segments go to $0$, which shows that (say in the future for the Feynman case)
\begin{equation}
    \chi_1 u = \sum_{\sigma \text{ above } \Gamma_1 \text{ below } \Gamma_1 + i\nu, Re(\sigma) > 0} \tau_1^{-i\sigma}(-\log(\tau))^k v^1_{\sigma, k} + r_{1, \nu}
\end{equation}
with $r_{1, \nu} \in H^{s - N-2}_{\Gamma_1 + i \nu}$, and the sum over resonances converges in a weak Sobolev norm when ordered by the lines $L_{k}$. We now reiterate this argument with $\chi_1 u$ of this form. Because we've assumed $P - N(P)$ to admit an expansion as well, the remainder $O(\tau^{\nu}) \chi_i u$ now admits an expansion in powers $\tau^{-i\sigma - \eta_l}$ up to a (different) decaying term $\Tilde{r}_{1, \nu} \in H^{s - N - 4}_{\Gamma_1 + i 2\nu}$:
\begin{equation}
    (P - N(P)) \chi_1 u = \sum_{\sigma \text{ above } \Gamma_1 \text{ below } \Gamma_1 + i\nu, Re(\sigma) > 0} \sum_{l, k} \tau_1^{-i\sigma - i\eta_l}(-\log(\tau))^k \Tilde{v}^1_{\sigma, k} + \Tilde{r}_{1, \nu}
\end{equation}
We then solve for $\chi_1 u$ using $\widehat{N(P)}^{-1}(\sigma)$ as above, but this time the right hand-side is meromorphic below $\Gamma + 2i\nu$ with poles at the exponents $\sigma + \eta_l$. This can again be done, with possibly increasing the multiplicity of poles when $\sigma + \eta_l$ coincides with a pole of $\widehat{N(P)}^{-1}(\sigma)$ and adding new resonances when we deform the contour to $\Gamma_1 + 2i\nu$. This gives an expansion with remainder in $H^{s-2N-4}_{\Gamma_0 + 2i\nu}$, and with resonances between $\Gamma_1$ and $\Gamma_1 + 2i \nu$ (along with exponents corresponding to the $\sigma + \eta_l$ terms). Repeating this over and over again gives the full expansion up to any finite decay less than $\Gamma_1 + i L$ (at which we must stop because the expansion is no longer preserved by $P - N(P)$ hence meromorphicity at each iteration is lost).
$\Box$
\vspace{0.2in}
\newline
The remaining part of this paper will be an application of the theory to the computation of indexes of Lorentzian Dirac operators. As explained in the introduction, this has been previously studied using a generalization of the APS boundary condition to Lorentzian geometry together with a spectral flow argument to compute the index. Here, we give an argument based on a continuous Wick rotation of the operator. We quickly recall the definition of Dirac operators in \ref{sec:DiracOperators} and do the Wick rotation in \ref{sec:WickRotation}.
\vspace{0.2in}
\newline
\textbf{Remarks:} We focused in this section on the case $Q = 0$. In the theory of wave equations, there are also interesting situations where $Q$ plays an essential role as a \textit{complex absorbing potential}. The name comes from the interpretation of $Q$ as an absorbing region for the waves. This can be used for example to study the Kerr-de Sitter family of black holes on a compact manifold, where $Q$ is applied outside the event horizons. The theory of this paper still applies to this case and gives Fredholm operators associated to contours, but the interpretation as Feynman propagators is no longer clear.
\subsection{Lorentzian Dirac Operators and their Index}
\label{sec:DiracOperators}
Dirac Operators can be thought of as 'square roots' for wave-type operators since they satisfy the property
\begin{equation}
    D^2 = P \quad \text{ $P$ a wave-type operator }
\end{equation}
Their index have physical meaning in physics since they describe \textit{chiral anomalies}; a mathematical introduction to the topic can be found in \cite{B_r_2016}. Within the Riemannian community, they also attracted a lot of attention since they form a large class of Fredholm operator whose indexes is computable through geometric quantities. Theorems relating their index to geometry are known as \textit{index theorems} in the literature. For this paper, the most relevant is the one R. Melrose derived in \cite{Melrose1993} for Riemannian $b-$Dirac operators. We thus begin by defining Lorentzian $b-$Dirac operators in analogy with Melrose's \cite{Melrose1993} definition 3.7. This is also essentially the same definition as Baer-Strohmaier \cite{B_R_2025} give in section 1 except done in the global $b-$geometry  \footnote{Note that here we are using Melrose's convention for the Dirac operator rather than the one of Baer-Strohmaier; the two differ by a $-i$ factor.}. 
\begin{definition}
    Let $\Bar{M}$ be the exponential compactification of an even-dimensional globally hyperbolic spacetime with compact Cauchy surface. Suppose $V = V^{+} \oplus V^{-}$ is a graded bundle equipped with
    \begin{itemize}
        \item An indefinite Hermitian non-degenerate pairing $\left( \cdot, \cdot \right)$.
        \item A graded b-clifford action compatible with the pairing: $cl(\zeta) cl(\gamma) + cl(\gamma) cl(\zeta) = 2 g^{-1}(\zeta, \gamma) Id \quad cl(\zeta)^{\dagger} = cl(\zeta) \quad G cl(\zeta) + cl(\zeta) G = 0$ where $G$ is the grading.
        \item A grading-preserving $b-$connection $\nabla$ conserving the pairing $( \cdot, \cdot)$ and which is compatible with the Clifford action and Levi-Civita connection $\nabla^L$: $\nabla(cl(\zeta) u) = cl(\nabla^L \zeta) u + cl(\zeta) \nabla u$.
    \end{itemize}
    If we let $C$ denotes the contraction between $T^{\ast}_b M$ and $V$, i.e. $C(\zeta \otimes v) = cl(\zeta) v$, then the associated Dirac operator is defined by:
    \begin{equation}
        Du = -i C \cdot \nabla u
    \end{equation}
    It is automatically a first-order $b-$differential operator $D \in \Psi^1_b(\Bar{M})$ and its principal symbol is given by the Clifford action:
    \begin{equation}
        \sigma_1(D)(\zeta) = cl(\zeta) \quad \sigma_2(D^2)(\zeta) = cl(\zeta)^2 = g^{-1}(\zeta, \zeta) Id
    \end{equation}
\end{definition}
\vspace{0.2in}
The canonical example of a Dirac operator is the familiar $d + \delta$ construction on the vector bundle of $k-forms$, whose data corresponds to
\begin{itemize}
    \item $V = \Lambda^{\circ}(T^{\ast}_b \Bar{M}) \otimes \mathbb{C} = \Lambda^{+}(T^{\ast}_b \Bar{M}) \otimes \mathbb{C} \oplus \Lambda^{-}(T^{\ast}_b \Bar{M}) \otimes \mathbb{C}$ graded into even and odd $k-$forms. \\
    \item The pairing $( \cdot, \cdot )$ is the induced Lorentzian metric on $k-forms$.
    \item The clifford action is $cl(\zeta) = i(\zeta \wedge \cdot - \iota_{\zeta})$ where $\iota_{\zeta}$ denotes contraction with $\zeta$. 
    \item The connection $\nabla$ is the induced Levi-Civita connection on $k-$forms. 
\end{itemize}
Another large class of Dirac operators comes from choosing a $b$-spin structure on $\Bar{M}$ which allows us to define a Dirac operator on the spinor bundle $\mathcal{S}$; we refer to Melrose's paper \cite{Melrose1993} chapter 3 on how this is done in the Riemannian case, while the Lorentzian case proceeds exactly the same way. Immediately from our work in the previous section \ref{sec:PropagatorsWave}, Dirac operators will be Fredholm on our spaces $H^s_{b, \Gamma}$ as long as $D \in \Psi^1_{b, \Gamma}$ and $\Gamma$ satisfies the usual assumptions of our non-elliptic Fredholm theorem \ref{theo:NonEllipticOperatorsAreFredholm}. 
\begin{proposition}
    Let $M$ be a globally hyperbolic spacetime whose stationary models $g_0$ at the past and future has Killing vector field $\partial_t$ timelike. Let $\Gamma_0$ and $\Gamma_1$ be choices of contour in the past $\{ 0 \} \times \Sigma$ and future $\{ 1 \} \times \Sigma$. Suppose that a Dirac operator as above is in $\Psi^1_{b, \Gamma}$. Suppose furthermore that its associated wave operator $P = D^2$ satisfies the assumptions of theorem \ref{theo:PropagatorConstruction}, i.e. that $2 g^{-1}(\zeta, dt) g(\sigma)Id + p_1$ has a definite sign in each halves of the characteristic set at the past and future, and that $\widehat{N(P)}(\sigma)$ satisfies high energy estimates on $\Gamma$ for some $k \geq 1$. Then
    \begin{equation}
        D: \{ u \in H^s_{b, \Gamma}: Du \in H^{s-1+k}_{b, \Gamma} \} \rightarrow H^{s-1+k}_{b, \Gamma}
    \end{equation}
    is Fredholm. 
\end{proposition}
\textbf{Proof:} We just use that $D^2 = P$, together with our Fredholm estimate for $P$ to get one for $D$:
\begin{equation}
    \| u \|_{H^s_{b, \Gamma}} \lesssim \| D^2 u \|_{H^{s-2+k}_{b, \Gamma}} + \| u \|_{H^{-N}_{b, \Tilde{\Gamma}}} \lesssim \| D u \|_{H^{s-1+k}_{b, \Gamma}} + \| u \|_{H^{-N}_{b, \Tilde{\Gamma}}}
\end{equation}
We proceed similarly for the adjoint $D^{\dagger}$.
$\Box$
\vspace{0.2in}
\newline
It thus makes sense to ask what their index are on such spaces. By symmetry consideration, the index of $D$ itself is typically 0 (unless one works on weighted spaces). Indeed, with respect to the indefinite pairing $(\cdot, \cdot)$ on $V$, $D$ is formally self-adjoint:
\begin{equation}
    \int_{\Bar{M}} (Du, v) = \int_{\Bar{M}} (u, Dv) \quad u, v \in \dot{\mathcal{C}}^{\infty}
\end{equation}
On the other end, with respect to the grading $V^{+} \oplus V^{-}$ the operator splits as
\begin{equation}
D = 
\begin{pmatrix} 
    0 & D^{-} \\ 
    D^{+} & 0 
\end{pmatrix}
\end{equation}
 These may have potentially have non-trivial indexes, and this is what is typically referred as 'the' index of a Dirac operator.
\begin{definition}
\label{def:DiracOperators}
    Let $D \in \Psi^1_{b, \Gamma}$ \footnote{This holds for example if both the Clifford action $cl$ and the connection $\nabla$ are in $\Psi_{b, \Gamma}$, which means that they tend to stationary model at the boundary.} be a Dirac operator Fredholm on $D: \{ u \in H^s_{b, \Gamma}: Du \in H^{s-1+k}_{b, \Gamma} \} \rightarrow H^{s-1+k}_{b, \Gamma}$ as above. Then its index with respect to $\Gamma$ is defined by
    \begin{equation}
        ind_{\Gamma}D = \text{Fredholm Index of } D^{+} \text{on $\{ u \in H^s_{b, \Gamma}: Du \in H^{s-1+k}_{b, \Gamma} \} \rightarrow H^{s-1+k}_{b, \Gamma}$}
    \end{equation}
\end{definition}
In the local in time setting $M = [0, 1]_t \times \Sigma$, Baer-Storhmaier \cite{B_r_2019} computed the index of Dirac operators on spinor bundles \footnote{They actually proved it also for twisted spin$^c$-Dirac operators but the statement is slightly more complicated. Also, their convention for the Dirac operator is different hence why they have different signs for the Eta invariants.} under APS boundary condition as
\begin{equation}
\label{eq:IndexFormulaAPS}
    Ind[D_{APS}] = \int_{M} \widehat{A}(\nabla) + \int_{\partial M} T \widehat{A}(\nabla) - \frac{(h(A^{+}_{0}) + h(A^{+}_{1}) + \eta(A^{+}_{0}) - \eta(A^{+}_{1}))}{2}
\end{equation}
where $\widehat{A}$ is the $\widehat{A}$ form of the curvature of the Lorentzian manifold, $T\widehat{A}$ is a transgression form corresponding to a difference of $\widehat{A}$ forms for the Lorentzian metric and an auxiliary Riemannian metric, $h(A^{+}_{0}), h(A^{+}_1)$ denote the dimensions of the kernel of the induced (positive) Dirac operators at the boundaries $\{ 0 \} \times \Sigma$, $\{ 1 \} \times \Sigma$ and finally $\eta(A^{+}_{0}), \eta(A^{+}_1)$ their eta invariant. The key point is that this is the exact same formula as the original APS index theorem in the Riemannian case \cite{Atiyah1975SpectralAA} applied to a 'Wick rotated' Riemannian manifold. As mentioned in the introduction, this formula was generalized to various settings. In particular, the same authors in \cite{B_R_2025} theorem 5.5 proved a local version for general Dirac operators. Furthermore, D. Shen and M. Wrochna in \cite{Shen_Wrochna_2022} proved this exact formula (theorem 1.3) in the infinite-time setting $M = \mathbb{R} \times \Sigma$ for Dirac operators acting on spinor bundles. As it is quite similar to the setting of this paper, we give a more in-depth comparison of their result with the work in this paper. They considered globally hyperbolic spacetime whose metric was only \textit{polynomially} decaying to a static model $g_0$:
\begin{equation}
    g = g_0 + O(|t|^{-\delta}) \quad g_0 = Ndt^2 - h \quad \delta > 1
\end{equation}
and showed using time-dependent scattering theory that the Dirac operators on spinor bundles are Fredholm as operators between
\begin{equation}
    D: \{ u \in C^0_{t} L^2_{\Sigma}: Du \in \left< t \right>^{-1/2 - \epsilon} L^2 \text{ and } lim_{t \rightarrow \pm \infty} \Pi_{\pm A(t) \leq 0}u(t) = 0\} \rightarrow \left< t \right>^{-1/2 - \epsilon} L^2
\end{equation}
Here $\Pi$ denotes the spectral projector associated to the (time-dependent) induced Dirac operator $A(t)$ on the $t-$spatial slice $\{ t \} \times \Sigma$. They then showed that its Fredholm inverse, i.e. the Feynman parametrix, satisfies the appropriate wavefront set relations (theorem 7.5) and that the index of $D^{+}$ is given by equation \ref{eq:IndexFormulaAPS} (theorem 6.6) using a similar spectral flow argument as Baer-Storhmaier. Thus, the main advantage advantage of their result is that the class of metrics they can handle is much larger, since it requires only polynomial decay. From the point of view of $b-$geometry, their metrics are very irregular on the exponential compactification $\Bar{M}$ (they make this remark themselves in remark 1.2). The trade-off is that it focuses only on Dirac operators acting on spinor bundles rather than general wave-type operators, and that it requires the metric to be asymptotically static (we are going to make the same assumptions for our index computation, but the theory is well behaved without them). Furthermore, we should say that the scattering theory based proof used there is quite different than the theory presented in this paper; in terms of Fourier-Laplace transform, the polynomially weighted distributions they work with can be thought off as having Fourier transform on $\mathbb{R}$ given by the boundary value of their holomorphic Fourier-Laplace transform (here suppressing spatial variables):
\begin{equation}
    \left< t \right>^{-l} H^s \simeq \{\text{distributions whose }\widehat{u}\text{ at the past(+)/future(-) is the boundary value } \widehat{u}(\sigma \pm i0) \text{ in $\left< \sigma \right>^{-s} H^l(\mathbb{R}_{\sigma})$} \}
\end{equation}
In terms of the $H^s_{b, \Gamma}$ spaces of this paper, this can be thought of as letting the contours $\Gamma$ approach the real line and requiring $H^l$ regularity for the boundary value. It would be very interesting to incorporate the APS boundary condition in terms of these boundary values and derive a positive commutator estimate on these spaces to bring the two theories on the same footing. Anyhow, from the expansion into quasi-normal modes $\sim e^{i \sigma t}$ we derived last section \ref{prop:ModeExpansion}, it is simple to relate the asymptotic structure that our Feynman parametrix select through contours $\Gamma$ to the APS boundary condition that was used in the previous works. The next proposition says that any Feynman solution we construct will also satisfy the APS boundary condition in an exponentially decaying sens for sufficiently high Sobolev regularity $s$, while in general will satisfy a weak version of the condition.
\begin{corollary}
    Let $D \in \Psi^1_{b, \Gamma}$ be a Dirac operator on an asymptotically ultra-static spacetime whose stationary model $N(D)$ at the future and past is as follows:
    \begin{equation}
        N(D)|_{\partial \Bar{M}_i} = cl(n)(-i\partial_t - A_i) \quad A_i^{\dagger} = A_i
    \end{equation}
    where $A_i$ is the induced Dirac operator on $\Sigma$, which for our purposes here is some self-adjoint first-order elliptic operator on $\Sigma$. Suppose the contours $\Gamma_i$ correspond to the Feynman problem and cross 0 at exactly some frequency $\lambda_i$ which is not in $spec(A_i)$. More precisely, they stay uniformly bounded away from the real line $\mathbb{R}$ and satisfy
    \begin{equation}
    \begin{aligned}
        &g_0(\sigma > \lambda_0) < 0 \quad g_0(\sigma < \lambda_0) > 0 \quad g_0(\lambda_0) = 0 \\
        &g_1(\sigma > \lambda_1) < 0 \quad g_1(\sigma < \lambda_1) > 0 \quad g_1(\lambda_1) = 0 
    \end{aligned}
    \end{equation}
    Then the mode expansion of solutions to \ref{prop:ModeExpansion} for $Du = f$ compactly supported (with $N(D) - D = O(\tau^{\nu})$ satisfying a mode expansion as in proposition \ref{prop:ModeExpansion}) takes the form
    \begin{equation}
    \begin{aligned}
            &\chi_0 u = \sum_{\sigma \in spec(A_0), \sigma < \lambda_0} \sum_{\alpha, k} e^{i(\sigma - \alpha \cdot \eta) t}(-t)^{k} u^0_{\sigma,k,\alpha} + r_0 \\
            &\chi_1 u = \sum_{\sigma \in spec(A_1), \sigma >\lambda_1} \sum_{\alpha, k} e^{i(\sigma + \alpha \cdot \eta) t}t^{k} u^1_{\sigma,k,\alpha} + r_1
    \end{aligned}
    \end{equation}
    over the discrete eigenvalues of $A_0$ and $A_1$. In particular, up to a sufficiently small decaying term $O(\tau^{\delta})$ for $\delta < min(\inf Im(\Gamma_1 + i \nu), -\sup Im(\Gamma_0 - i\nu), \nu)$ (assumed to be positive by taking $\Gamma$ sufficiently close to the real line) we have
    \begin{equation}
    \label{eq:ModeExpansionSelfAdjointGenerator}
    \begin{aligned}
            &\chi_0 u = \sum_{\sigma \in spec(A_0), \sigma < \lambda_0} \sum_k e^{i\sigma t} u^0_{\sigma,k} + O(e^{-\delta|t|}) \\
            &\chi_1 u = \sum_{\sigma \in spec(A_1), \sigma > \lambda_1} \sum_{k} e^{i\sigma t} u^1_{\sigma,k} + O(e^{-\delta|t|})
    \end{aligned}
    \end{equation}
     where $u^i_{\sigma,k}$ are eigenfunctions of $A_i$. Thus $u$ satisfies the APS boundary condition in a weak Sobolev sense:
     \begin{equation}
         \Pi_{A_0 > \lambda_0}(\chi_0 u) \in \tau^{\delta}H^{-M}_{b} \quad \Pi_{A_1 < \lambda_1}(\chi_1 u) \in \tau^{\delta}H^{-M}_{b}
     \end{equation}
     where $\Pi$ are the spectral projectors of $A_0$ and $A_1$. If furthermore $s$ is large enough so that $-M > 1/2 + k$, then it holds in the limiting sense as well:
     \begin{equation}
         lim_{t \rightarrow -\infty} \Pi_{A_0 > \lambda_0}(u) = 0 \quad lim_{t \rightarrow +\infty} \Pi_{A_1 < \lambda_1}(u) = 0 
     \end{equation}
     in $C_t^0(H^k_{\Sigma})$.
\end{corollary}
\textbf{Proof:} The result follows from applying our mode expansion \ref{prop:ModeExpansion} to $D$; while it was stated for wave operators, it follows from the proof that it applies to any $P$ satisfying the assumptions. The poles of $\widehat{D}(\sigma)^{-1}$ are then identified with the spectrum of $A$, hence why we get an expansion in eigenvalues of $A$. In particular the reason there are no $t^k$ powers in \ref{eq:ModeExpansionSelfAdjointGenerator} is because the poles of $(z - A)^{-1}$ are simple. Thus, we must justify applying the mode expansion in that case. We automatically get high-energy estimates for $N(D)(\sigma)^{-1}$ in $\mathbb{C} \backslash \mathbb{R}$ by the spectral theory of $A$. The only thing to check is thus the claimed estimates on vertical line $L_{k}$. To do so, recall that the norm of the resolvent of $A$ satisfy the classical bound:
\begin{equation}
    \|(z - A)^{-1} \|_{H^s_{\Sigma} \rightarrow H^s_{\Sigma}} = \frac{1}{dist(z, spec(A))}
\end{equation}
Furthermore, Weyl's law gives a count of its eigenvalues with multiplicity (here applied to $|A|$)
\begin{equation}
    N(\Lambda) = \# \{k; \lambda_k \leq \Lambda \} = C \Lambda^{n} + o(\Lambda^{n})
\end{equation}
So for large enough $\Lambda$ we have $N(\Lambda) \leq 2C \Lambda^{n}$. By the pigeonhole principle, we can find $l_{2k}$ and $l_{2k+1}$ such that $l_{2k} \rightarrow \infty$ and $l_{2k+1} \rightarrow -\infty$ with
\begin{equation}
    dist(l_{k}, spec(A))^{-1} \lesssim \left< l_{k} \right>^{n-1}
\end{equation}
On such lines $L_i$, we thus get
\begin{equation}
    \| (\sigma - A)^{-1} \|_{H^{s}_{\partial \Bar{M}_i, \sigma} \rightarrow H^s_{\Bar{M}_i, \sigma}} \lesssim \left< l_{k} \right>^{n-1}
\end{equation}
which is the desired polynomial bound hence proves the expansion. $\Box$
\vspace{0.2in}
\newline
Thus, because in the case of Dirac operators on spinor bundles the results by Shen Wrochna in \cite{Shen_Wrochna_2022} encompasses the exponentially decaying spacetimes we've been considering here and our spaces also satisfy the APS boundary conditions, we expect the index of Dirac operators considered here to agree with formula \ref{eq:IndexFormulaAPS} as well. Thankfully, the hard part of computing the index of $b-$operators was done in the Riemannian signature by R. Melrose \cite{Melrose1993}, hence here we can use a simple Wick rotation argument to relate the two indexes. This will be done in the next section \ref{sec:WickRotation}.

\subsection{Wick Rotation}
\label{sec:WickRotation}
We now move on to explain how the Wick rotation is done and what kind of estimates one needs to get results out of it. It is important to point out that a Wick-like rotation was already a core ingredient in the previous works; Baer-Storhmaier used it in \cite{B_r_2019} to relate the Lorentzian index of $dt^2 - h^2$ to the one of $dt^2 + h$, which is how they got a transgression form $T \widehat{A}(\nabla)$ in their formula. Subsequent results based on similar spectral flow arguments also made use of it. The novelty here is that because our theory is well behaved under \textit{complex} deformations of our operator, we can do a continuous Wick rotation instead of relying on more explicit formulas to relate the Lorentzian and Riemannian indexes \footnote{The spectral flow approach could work in this setting, but doing the continuous Wick rotation was a major motivator for this theory hence why this approach is taken here.}. This will not be a 'true' Wick rotation, in the sense that we are not complexifying the time coordinate $t$; doing so would require analyticity assumptions on the coefficients of the operators, which we would like to avoid. Instead, we will work in the complexified tangent bundle without ever modifying $t$ itself. The trade-off is that properties of our operator will change as we rotate, but since we are only interested in the index, this deformation of the operator will be sufficient for our purposes. So let us begin: the plan is to turn our Lorentzian metric into a Riemannian metric in a 'continuous way', so that we can invoke the invariance of a Fredholm index under continuous deformations of the operator. The simplest way to do so is to simply rotate a timelike vector field normal to the Cauchy surfaces $\Sigma_t$ by $e^{i\theta}$. In the notation of section \ref{sec:PropagatorsWave} a general globally hyperbolic Lorentzian metric is given by
\begin{equation}
    g = N^2 dt^2 - h_{ij}(dx^i + \beta^i dt)(dx^j + \beta^j dt)    
\end{equation}
A normal vector field to the Cauchy surfaces $\Sigma_t$ is then
\begin{equation}
    n = \frac{1}{N}(\partial_t - \beta) \quad g(n, n) = 1
\end{equation}
Wick rotating it amounts to simply multiplying it by $e^{i\theta}$:
\begin{equation}
    R_{\theta} n = \frac{e^{i\theta}}{N}(\partial_t - \beta)
\end{equation}
If we denote by $R_{\theta}$ the (complex) rotation map on vectors that sends $n + v$ to $e^{i\theta} n + v$ for $v \in T\Sigma$, then the effect on the metric is given by:
\begin{equation}
    g_{\theta}(X, Y) = g(R_{\theta} X, R_{\theta} Y) = e^{2i\theta}N^2 dt^2 - h_{ij}(dx^i + \beta^i dt)(dx^j + \beta^j dt)
\end{equation}
It is now a non-degenerate complex $2-$tensor, and at the end-points $\theta = \pm \pi/2$ is negative definite, hence for all purposes Riemannian. Its inverse is given by
\begin{equation}
    g^{-1}_{\theta} = g^{-1}(R_{-\theta}\zeta, R_{-\theta} \eta) = \frac{e^{-2i\theta}}{N^2}(\partial_t - \beta^i \partial_{x^i})(\partial_t - \beta^j \partial_{x^j}) - h^{ij} \partial_{x^i} \partial_{x^j}
\end{equation}
Now, we wish to apply this to studying Dirac operators. We thus also need to rotate all the structures required to define one, as explained in definition \ref{def:DiracOperators}. This is done in the exact same fashion.
\begin{itemize}
    \item Define the rotated Clifford action by $cl_{\theta}(X) = cl(R_{\theta} X)$ on $V$ (or the in terms of covectors $cl_{\theta}(\zeta) = cl(R_{-\theta} \zeta))$. It satisfies $cl_{\theta}(X) cl_{\theta}(Y) + cl_{\theta}(Y) cl_{\theta}(X) = 2 g_{\theta}(X, Y) Id$ and is still grading reversing. \\
    \item In the case where the connection is constructed from the metric $g$, we can simply use $g_{\theta}$ to rotate the connection. This will be the case for the spinor bundle or the $d+ \delta$ construction for example. In the general case, to define the rotated connection $\nabla^{\theta}$, we will require that it is still compatible with $g_{\theta}$ and the Clifford action as we rotate, so that at the end point we are guaranteed to end up with a Riemannian Dirac operator. Other choices could be made as long as the end-point is correct. First note that $\nabla^{\theta} - \nabla$ will be an endomorphism $\Phi^{\theta}$ of the bundle, hence the requirement can be written as
    \begin{equation}
    \label{eq:RotatedConnectionCompatibility}
        [\nabla_X + \Phi^{\theta}X, cl_{\theta}(Y)] u = cl_{\theta}(\nabla^{L, \theta}_{X} Y) u
    \end{equation}
    where $\nabla^{L, \theta}$ is the Levi-Civita connection associated to $g_{\theta}$ (which is well-defined since it only requires the associated pairing to be non-degenerate, not necessarily real). To compute $\nabla^{L, \theta}$, we can use the usual formula in terms of $g_{\theta}$, and we obtain for $v, w \in T\Sigma$
    \begin{equation}
    \begin{aligned}
        &\nabla^{L, \theta}_n n = e^{2i\theta} \nabla^L_n n \\
        &\nabla^{L, \theta}_v w = \nabla^L_v w - (1 - e^{-2i\theta}) g(\nabla^L_v w, n) n \\
        &\nabla^{L, \theta}_v n = \nabla^{L}_v n \quad \nabla^{L, \theta}_n v = \nabla^{L}_n v 
    \end{aligned}
    \end{equation}
    Equipped with that, we can solve the commutation relation above. If we pick $Y = n$ and $X$ some $v \in T\Sigma$, we get $cl_{\theta}(n) = e^{i\theta} cl(n)$ hence the equation becomes 
    \begin{equation}
        e^{i\theta}[\nabla_v + \Phi^{\theta}v, cl(n)] u = cl(\nabla^{L}_{v} n) u
    \end{equation}
    Because $\nabla_X$ is already compatible with he Clifford action, we end up with
    \begin{equation}
        [\Phi^{\theta}v, cl(n)] u = (e^{-i\theta} - 1) cl(\nabla^{L}_{v} n) u \quad \implies  \quad \Phi^{\theta}v = \frac{(e^{-i\theta} - 1)}{2} cl(\nabla^{L}_{v} n) cl(n)
    \end{equation}
    where our solution $\Phi^{\theta}$ is well defined up to a term commuting with $cl(n)$. Similarly, taking $X = n$ then gives
    \begin{equation}
    e^{i\theta}[\nabla_n + \Phi^{\theta}n, cl(n)] u = e^{2i\theta} cl(\nabla^{L}_{n} n) u      
    \end{equation}
    hence using that $\nabla_n$ is compatible with the Clifford action we get
    \begin{equation}
        [\Phi^{\theta}n, cl(n)] u = (e^{i\theta} - 1) cl(\nabla^{L}_{n} n) u \quad \implies  \quad \Phi^{\theta}n = \frac{(e^{i\theta} - 1)}{2} cl(\nabla^{L}_{n} n) cl(n)
    \end{equation}
    up to a term commuting with $cl(n)$. It can be checked pretty easily that these choices also respect the full Clifford compatibility, in the sense that if we let $Y \in T\Sigma$ equation \ref{eq:RotatedConnectionCompatibility} is still satisfied. Thus, a choice of rotated connection $\nabla^{\theta}$ for $X = X_{\perp} + g(X, n) n$ is
    \begin{equation}
        \nabla^{\theta} = \nabla + \Phi^{\theta} \quad \Phi^{\theta} X =  g(X, n) \frac{(e^{i\theta} - 1)}{2} cl(\nabla^{L}_{n} n) cl(n) + \frac{(e^{-i\theta} - 1)}{2} cl(\nabla^{L}_{X_{\perp}} n) cl(n)
    \end{equation}
    Note in particular that it preserves the grading since both $\nabla$ and $cl(\nabla^L_X n)cl(n)$ do so. At the end points $\theta = \pm \pi/2$, it gives a connection compatible with the Levi-Civita connection of the Riemannian manifold $g_{\pm \pi/2}$:
    \begin{equation}
        \nabla^{\pm \pi/2}_X = \nabla_X + g(X, n) \frac{(\pm i - 1)}{2} cl(\nabla^{L}_{n} n) cl(n) - \frac{(\pm i+1)}{2} cl(\nabla^{L}_{X_{\perp}} n) cl(n)
    \end{equation}
    For us, the relevant situations will be asymptotically ultra-static spacetimes. In that case $\nabla^L n$ vanishes, hence the connection stays constant under Wick rotation.
    \item The most subtle structure to rotate is actually the indefinite pairing $(\cdot, \cdot)$ on the bundle $V$. A Riemannian Dirac operator requires a connection which preserves a positive-definite pairing and for which the Riemannian Clifford action is self-adjoint (and which makes the grading an orthogonal decomposition). We must thus manufacture one for the end-point connection $\nabla^{\pi/2}$. In the general case, this seems hard to do since our Wick rotation takes place on the complexified tangent bundle $TM^{\mathbb{C}}$ which does not need to have any relations with the bundle $V$ where our original (indefinite) pairing $(\cdot, \cdot)$ is defined. We thus restrict to a specific case where there exists such a positive definite pairing $\left< \cdot, \cdot \right>$ on $V$ preserved by $\nabla^{\pi/2}$. This is guaranteed for example if
    \begin{equation}
        \left< \cdot, \cdot \right> = (\cdot, cl(n) \cdot) \text{ is positive-definite}
    \end{equation}
    This is the case for the spinor bundle (see \cite{B_r_2019} section 1.). Our end-point $\nabla^{\pi/2}$ connection preserves this pairing. Indeed, the original one satisfies 
    \begin{equation}
        \nabla^{\pi/2}_X \left< s, t \right> = \nabla_X \left< s, t \right> = \left< \nabla_X s, t \right> + \left<s, \nabla_X t \right> + (s, [\nabla_X, cl(n)]t)
    \end{equation}
    The commutator $[\nabla_X, cl(n)]$ is $cl(\nabla^L_X n)$ by compatibility with the Clifford action. We can thus rewrite $\nabla_X$ above as
    \begin{equation}
        \nabla_X = \nabla^{\pi/2}_X + \frac{1}{2} [\nabla_X, cl(n)] cl(n) \mp \frac{i}{2} (g(X, n)cl(\nabla^L_{n}n)cl(n) - cl(\nabla^L_{X_{\perp}}n)cl(n)) = \nabla_X^{\pi/2} + \Phi
    \end{equation}
    Using the Clifford relations and that $cl(n)$ is self-adjoint with respect to the pairing $(\cdot, \cdot)$, it's pretty simple to show that the extra terms cancel by showing that conjugation by $cl(n)$ flips its sign. Finally, the Riemannian Clifford action is clearly self-adjoint with respect to this pairing (after renormalizing by $i$). If furthermore the grading was originally anti-self-adjoint with respect to $h$, i.e. $G^{\dagger} = -G$, then the grading is self-adjoint with respect to the rotated positive definite pairing $\left< \cdot, \cdot \right>$.
\end{itemize}
This Wick rotation on $\mathbb{R} \times \Sigma$ lifts to one on the exponential compactification $\Bar{M} = [0, 1] \times \Sigma$, with the only subtlety that the Wick rotation is done in 'opposite' directions at the past and future in terms of $\tau \partial_{\tau}$; indeed using $\tau_{-} = e^t$ and $\tau_{+} = e^{-t}$ we get
\begin{equation}
    n = \frac{1}{N}(\mp \tau_{\pm} \partial_{\tau_{\pm}} - \beta)
\end{equation}
hence for $\beta = 0$ say, in the future the dual coordinate of $-d\tau/\tau$ is rotated as $\sigma \rightarrow e^{i\theta}$ while in the past its changing by $-e^{i\theta}$. This is another reason why we changed convention so that $\sigma$ denotes a dual coordinate to $t$ rather than $-\tau \partial_{\tau}$. Thus, lifting everything to the exponential compactification, for each $\theta$ we obtain a first order differential operator
\begin{equation}
    D_{\theta} = -iC_{\theta} \cdot \nabla_{\theta} u \quad C_{\theta}(\zeta \otimes v) = cl_{\theta}(\zeta) v
\end{equation}
which interpolate between our original Dirac operator $D_0$ and a Riemannian geometric Dirac operator $D_{\pi/2}$. In order to use the continuity of a Fredholm index under continuous deformation, we need each of them to be Fredholm in a uniform way. For that, we consider their squares, which now have an imaginary part at the top order
\begin{equation}
    D^{2}_{\theta} = P_{\theta} + iQ_{\theta}
\end{equation}
Their principal symbols are
\begin{equation}
\begin{aligned}
    &\sigma_2(P_{\theta}) = Re(g^{-1}_{\theta}(\zeta, \zeta)) = \frac{cos(2\theta)}{N^2}(\sigma - \beta^i \xi_i)^2 - h^{-1}(\xi, \xi)\\
    &\sigma_2(Q_{\theta}) = Im(g^{-1}_{\theta}(\zeta, \zeta)) = \frac{-sin(2\theta)}{N^2}(\sigma - \beta^i \xi_i)^2
\end{aligned}
\end{equation}
So for $0 \leq \theta < \pi/4$, $P_{\theta}$ is a wave-operator associated to the Lorentzian inverse metric $Re(g^{-1}_{\theta})$. As for $Q_{\theta}$, it has a global sign corresponding to the sign of $sin(2\theta)$. Thus, the Feynman parametrix of $P + iQ_{\theta}$ will be well-behaved when $sin(2\theta) \geq 0$ while the anti-Feynman parametrix will be well-behaved when $sin(2\theta) \leq 0$. As for the causal propagators, i.e. the retarded and advanced one, they will not even be defined because in theorem \ref{theo:PropagatorConstruction} we require the imaginary part to have different signs in each halves of the characteristic set for those. This is why Wick rotation only really works for Feynman and anti-Feynman propagators. We summarize this discussion in the next proposition.
\begin{proposition}
\label{prop:UniformFredholm}
    Let $D_{\theta} = -i C_{\theta} \cdot \nabla_{\theta}$ be the Wick rotated Dirac operator as defined above. Let $\Gamma$ be choices of contours and suppose the original connection $\nabla: V \rightarrow T^{\ast}_b M \otimes V$ is in $\Psi^1_{b, \Gamma}$ and similarly that the Lorentzian Clifford action $C: T^{\ast}M \otimes V \rightarrow V$ is also in $\Psi^0_{b, \Gamma}$. Then $D_{\theta} \in \Psi^1_{b, \Gamma}$. If furthermore $\partial_t$ is timelike for the stationary models $g_0$ of $M$, and the contours $\Gamma$ are chosen to satisfy the assumptions for the Feynman (anti-Feynman) problems in theorem \ref{theo:PropagatorConstruction} for $0 \leq \theta \leq \pi/2$ ($0 \geq \theta \geq -\pi/2$), then assuming uniform high-energy estimates over $\theta$ for the indicial family $\widehat{N(P_{\theta} + iQ_{\theta})}^{-1}(\sigma)$ on $\Gamma$ yields that for $0 \leq \theta \leq \pi/2$ ($0 \geq \theta \geq -\pi/2$)
    \begin{equation}
        D_{\theta}: \{ u \in H^s_{b, \Gamma}: D_{\theta}u \in H^{s-1+k}_{b, \Gamma} \} \rightarrow H^{s-1+k}_{b, \Gamma}
    \end{equation}
    is a uniform family of Fredholm operators, in the sense that the corresponding Fredholm estimates are uniform.
\end{proposition}
\textbf{Proof:} The proof follows from our remark in theorem \ref{theo:PropagatorConstruction} that if $P_{\theta} + iQ_{\theta}$ is a family of operator which satisfy the assumptions uniformly, then the Fredholm estimates for them can be taken uniform. Indeed, the assumptions of proposition $\ref{prop:RadialPointEstimates}$ are clearly uniform for small $\theta$, while for large $\theta$ we can use that $P_{\theta} + iQ_{\theta}$ is uniformly elliptic. The assumed uniformness of the high-energy estimates then yield the missing ingredient. $\Box$
\vspace{0.2in}
\newline
Of course, the crucial thing to show is that the high-energy estimates are uniform. As we Wick rotate, what will happen is that the poles of $\widehat{N(P_{\theta} + iQ_{\theta})}^{-1}$ will also move in the complex plane, and for our estimates to be uniform, they better not cross the contour $\Gamma$ \footnote{On that note, one could also attempt to rotate $\Gamma$ so that it stays away from the pole, but that would end up heavily modifying the spaces with $\theta$ which would require careful analysis to understand, hence we do not attempt to do so here.}. In general, the roots are all spread over some finite strip $|Im(\sigma)| \leq C$ with no easy way to track them as we rotate the operator. They are given by the poles of the inverse of some quadratic family varying over $\theta$:
\begin{equation}
    \widehat{N(D_{\theta})}^2(\sigma) = \widehat{N(P_{\theta} + iQ_{\theta})}(\sigma) = \frac{e^{-i2\theta}\sigma^2}{N^2} + \sigma A_{\theta} + B_{\theta}
\end{equation}
The simplest case where we can guarantee that they never cross $\Gamma$ is if the poles at $\theta = 0$ are identified with the spectrum of some self-adjoint operator on $\Sigma$ hence all lie on the real line $\mathbb{R}$, in which case the poles get simply rotated to the imaginary axis $i \mathbb{R}$ and it is then simple to arrange $\Gamma$ so that they never touch. Thus, we now restrict to the case of $P = D^2$ having a self-adjoint generator and furthermore assume the model metric in the future and past to be ultra-static ($\beta = 0$, $N = 1$ asymptotically) \footnote{This will also be necessary to apply Melrose's $b-$index theorem as it applies to exact Riemannian $b-$metrics $\sim (d\tau/\tau)^2 + h$ which are the same as asymptotically ultra-static metrics.}. After doing the Wick rotation of the connection $\nabla^{\theta}$, we want for the Wick rotation to reduce to $\sigma \rightarrow e^{-i\theta} \sigma$ hence for the indicial family of $P_{\theta} + iQ_{\theta}$ to be
\begin{equation}
        \widehat{N(D_{\theta}^2)}(\sigma) = (e^{-i2\theta} \sigma^2 - A^2) \quad (A^2)^{\dagger} = A^2 \quad spec(A) = -spec(A)
\end{equation}
We would then have uniform control on $\widehat{N(D_{\theta}^2)}(\sigma)^{-1}$ away from some strip centered around $e^{i\theta} \mathbb{R}$ by usual spectral analysis \footnote{Here we are using that $spec(A)$ will be symmetric around $0$ for our cases using $cl(n)A + Acl(n) = 0$, hence we don't need to consider $\pm spec(A)$.}:
\begin{equation}
     \| u \|_{H^{s}_{\partial \Bar{M}_i, \sigma}} \lesssim \frac{1}{dist(e^{-i\theta}\sigma, spec(A)) | \sigma | } \| \widehat{N(D_{\theta}^2)}(\sigma) u \|_{H^{s-1}_{\partial \Bar{M}_i, \sigma}} \quad \sigma \in \mathbb{C} \backslash (e^{i\theta} spec(A) \cup \{ 0 \})
\end{equation}
hence as we rotate the poles from the real subset $spec(A)$ to the complex plane $e^{i\theta} spec(A)$, the high-energy estimates on the contours $\Gamma$ stay uniform as long as $\Gamma$ stays a positive distance away from the 'full' sweep $W = \cup_{0 \leq \theta \leq \pi/2} e^{i\theta} spec(A)$ (or $W = \cup_{-\pi/2 \leq \theta \leq 0} e^{i\theta} spec(A)$ in the anti-Feynman case where we rotate by negative angles $\theta$). To guarantee such a structure for the associated rotated wave operator and for them to fit in Melrose's framework of geometric Dirac operators, we will thus assume that the connection associated to the Dirac operator satisfies
\begin{equation} 
    \nabla_{\tau \partial_{\tau}}u|_{\partial \Bar{M}} = 0 \quad \forall u
\end{equation}
i.e. the connection reduces to $\tau \partial_{\tau}$ derivatives at the boundary. In that case, the model of the Dirac operator $D^{+}$ takes the simple form $N(D) = cl(n)(-i\partial_t - A)$, where $A$ is the induced Dirac operator on the pulled-back bundle $V|_{\partial \Bar{M}}$ to the boundary defined by the spatial Clifford action $\xi \rightarrow cl(\xi) cl(n)$. Note for further references that it preserves the grading, hence splits as
\begin{equation}
    A = A^{+} \oplus A^{-} \text{ on } V|_{\partial \Bar{M}}^{+} \oplus V|_{\partial \Bar{M}}^{-}
\end{equation}
The operator $A^{+}$ will be the one showing up in the index theorem. The Wick rotated family at the boundary can then easily be computed from our formula for $\nabla^{\theta}$ since in the ultra-static case it doesn't change with $\theta$ (by plugging $\nabla^L n|_{\partial \Bar{M}} = 0$ in our formula for $\nabla^{\theta}$):
\begin{equation}
    N(D_{\theta}) = cl(n)(-i e^{-i\theta} \partial_t - A) \implies  N(D_{\theta}^2) = (-e^{-2i\theta}\partial_t^2 - A^2) 
\end{equation}
The canonical example is the Dirac operator on the spinor bundle of an asymptotically ultra-static metric. In that case, using the induced Levi-Civita connection on spinors we get
\begin{equation}
    \widehat{N(D)}(\sigma) = cl(n)(\sigma - D_{\Sigma})
\end{equation}
where $D_{\Sigma}$ is the induced spinor Dirac operator on $\Sigma$. This formula is well-known and can be found in \cite{B_r_2019} section 1.2. As for the rotated operator, in that case we can just use the Levi-Civita connection on spinors induced by $\nabla^{L, \theta}$. The Wick rotated static model metric is obtained by $\sigma \rightarrow e^{-i\theta} \sigma$ and so is the Clifford action. As for the Levi-Civita connection, by what we computed it doesn't change in the ultra-static case. It thus follows that the model rotated Dirac operator is also obtained by such a rotation
\begin{equation}
    \widehat{N(D_{\theta})}(\sigma) = cl(n)(e^{-i\theta}\sigma - D_{\Sigma})
\end{equation}
 Our other main example is $d+\delta$, which also fits in this framework on asymptotically ultra-static spacetimes since
\begin{equation}
    \widehat{N((d+\delta))}(\sigma) = cl(n)(\sigma - A) \quad A = -i\begin{pmatrix}
  0 & d_{\Sigma} - \delta_{\Sigma} \\
  d_{\Sigma} - \delta_{\Sigma} & 0
\end{pmatrix}
\end{equation}
where we've split $k-$forms into normal and tangential component $u = u_T + \alpha dt \wedge u_N$. The formula for $d$ and $\delta$ in that splitting can be found in \cite{Hintz_2018_Cohomology} section 3, from which the above formula follows by factoring out $cl(n)$ \footnote{We would like to point out that in the work of \cite{Hintz_2018}, a Wick rotation type-argument was also used to compute the resonances; see section 3.1.}. Again, the rotated connection is the Levi-Civita connection on forms induced by $\nabla^{L, \theta}$, which is invariant in the ultra-static case, and the Clifford action is the one of $g_{\theta}$, which is again obtained by simply rotating $\sigma \rightarrow e^{-i\theta}$ hence gives:
\begin{equation}
    \widehat{N((d+\delta)_{\theta})}(\sigma) = cl(n)(e^{-i\theta}\sigma - A)
\end{equation}
In that setting, we thus have a uniform family of Fredholm operators, hence we can use that to relate the indexes of the two end-points. Pairing this with the knowledge of the index of the Riemannian Dirac operator, which is given by Melrose's theorem 9.1 \cite{Melrose1993} in the $b-$setting, we end up with a formula for the index of our Lorentzian Dirac operator $D$. We also refer to his paper (especially section 9) for the definitions of the quantities involved in the theorem.
\begin{figure}[h]
    \centering
    \includegraphics[width=0.6\linewidth]{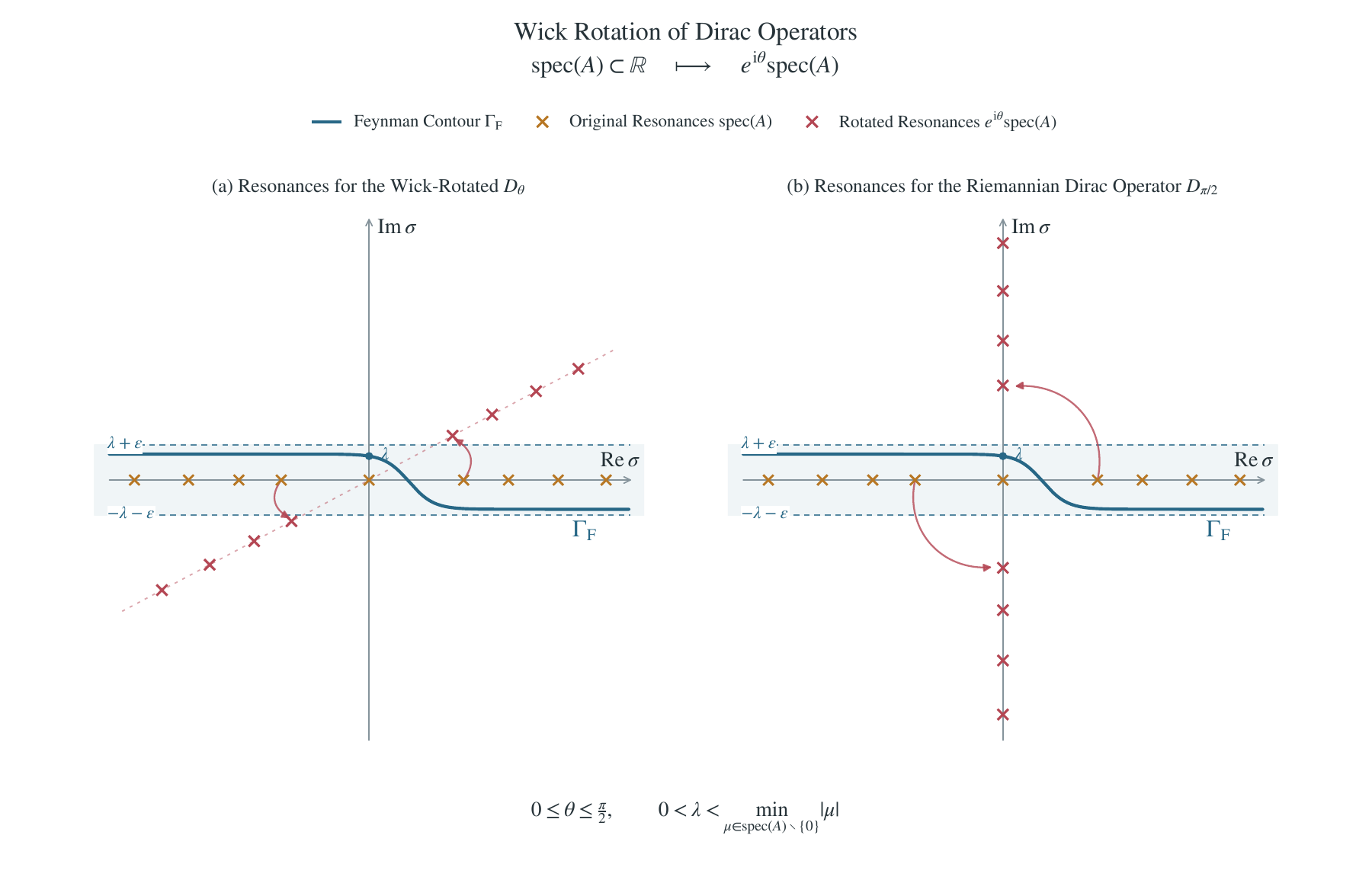}
    \caption{An illustration of the (future) contour in the Feynman case chosen for the index theorem, and how the resonances change with Wick rotation. }
    \label{fig:wick_rotation}
\end{figure}
\begin{theorem}(Indexes of Asymptotically Ultra-Static Spacetimes) 
\label{theo:IndexTheorem}
    Let $\Bar{M}$ be the exponential compactification of a globally even dimensional hyperbolic spacetime whose stationary models $g_0$ at the past and future are ultra-static. Let $D$ be a Dirac operator as defined in section \ref{sec:DiracOperators} whose associated connection satisfies $\nabla_{\tau \partial_{\tau}}|_{\partial \Bar{M}} = 0$. Thus, its Wick rotated family takes the form
    \begin{equation}
        N(D_{\theta}) = cl(n) (-i e^{-i\theta} \partial_t - A_i) \quad A_i \text{ is the induced Dirac operator on the boundary component $\partial \Bar{M}_i$}
    \end{equation}
    Let $\Gamma$ be choices of contours such that the associated connection $\nabla$ and Clifford action $C$ to $D$ are in $\Psi_{b, \Gamma}$. Suppose that the contours $\Gamma$ corresponds to the Feynman propagator of \ref{theo:PropagatorConstruction} and cross $\mathbb{R}$ exactly once at $\pm \lambda$ ($+$ in the future, $-$ in the past) for some $\lambda > 0$ smaller than the minimal nonzero eigenvalues of $|A_0|, |A_1|$, while staying a positive distance away from $e^{i\theta} spec(A)$ for $0 \leq \theta \leq \pi/2$ and having imaginary part between $\pm i(\lambda + \epsilon)$, for $\epsilon$ small. For example, the contours can be taken to be
    \begin{equation}
    \begin{aligned}
        \phi_0(\sigma) = \sigma - i l \tanh(\sigma + \lambda) \\
        \phi_1(\sigma) = \sigma - i l \tanh(\sigma - \lambda)
    \end{aligned}
    \end{equation}
    for $0 < l < 1$ small enough. If $A_0$ and $A_1$ do not have a kernel, we can take $\lambda = 0$. Then the operators $D_{\theta}$ are uniformly Fredholm with constant index on the spaces
    \begin{equation}
        D_{\theta}: \{ u \in H^s_{b, \Gamma}: D_{\theta}u \in H^{s}_{b, \Gamma} \} \rightarrow H^{s}_{b, \Gamma}
    \end{equation}
    In particular, if the Riemannian Dirac operator corresponding to the end-point connection $\nabla^{\pi/2}$ preserves some positive definite pairing with respect to which the Riemannian Clifford action is self-adjoint and the grading orthogonal (for example if its associated indefinite pairing $(\cdot, \cdot)$ has $\left< \cdot, \cdot \right> = (\cdot, cl(n) \cdot)$ positive definite and $G^{\dagger} = -G$) then the index of the Lorentzian Dirac operator at $\theta = 0$ is given by the b-index of the Riemannian one $D_{\pi/2}$ on $\tau^{\lambda} H^s_b$:
    \begin{equation}
        ind_{\lambda}(D) = {}^b\int_{\Bar{M}} AS(\nabla) - \frac{(h(A^{+}_{0}) + h(A^{+}_{1}) + \eta(A^{+}_{0}) - \eta(A^{+}_{1}))}{2}
    \end{equation}
    where $AS$ is the Atiyah-Singer integrand while $h$ and $\eta$ denote the dimension of the kernel and the eta invariant of $A^{+}_0$, $A^{+}_1$. Here, ${}^b\int_{\Bar{M}}$ denotes the regularized $b-$integral, which is basically the value of the integral with its logarithmically singular term removed (see \cite{Melrose1993}). Examples satisfying these assumptions are the Dirac operator on spinors or the $d+\delta$ Dirac operator. The index of the corresponding Dirac operator for the anti-Feynman contours $\Bar{\Gamma}$ on the dual space would be minus that. Equivalently, it would be equal to the Riemannian index of $D_{-\pi/2}$.
\end{theorem}
\textbf{Remarks:} In the case of the Dirac operator on spinors, this collapses to the previously derived formula in terms of $\widehat{A}$ forms
\begin{equation}
        ind_{\lambda}(D) = {}^b\int_{\Bar{M}} \widehat{A}(\nabla) - \frac{(h(A^{+}_{0}) + h(A^{+}_{1}) + \eta(A^{+}_{0}) - \eta(A^{+}_{1}))}{2}
\end{equation}
 We are missing the transgression form term here because in the ultra-static geometries, it vanishes at the boundary. For spinor Dirac operators, the result also holds on asymptotically \textit{static} spacetime, i.e. allowing for $N \neq 1$ asymptotically, since the spinor Dirac operator is a conformal invariant, hence we can rescale by $N^2$ the metric and have the same index. For the $d+\delta$ operator, the theorem also holds on odd-dimensional spacetimes and says that its index is the one of the Riemannian $d + \delta$ operator. An average of this index between growing and decaying spaces was computed by Melrose in \cite{Melrose1993} lemma 9.2 in terms of the relative deRham cohomology of $\Bar{M}$. For $\Bar{M} = [0, 1] \times \Sigma$, the corresponding averaged index vanishes, hence on our decaying spaces the index will be minus a half times the size of the kernel of $d+\delta$ acting on the even forms of $\Lambda^{\circ}(T^{\ast}_b \Bar{M})|_{\partial \Bar{M}}$, which is a half times minus the sum of the dimensions of the deRham complex on $\{ 0 \} \times \Sigma \cup \{1 \} \times \Sigma$. In both even and odd dimensions the index thus gives minus the total dimension of the deRham complex on $\Sigma$ (while the integral formula above only holds in even dimensions):
\begin{equation}
    ind_{\lambda}(d+\delta) = - dim(H_{dR}(\Sigma))
\end{equation}
Furthermore, during the proof we will be showing a relative index theorem, which allow us to compute the index of these Dirac operators on any weighted spaces $H^s_{b, \Gamma}$. One last remark is that given stronger Riemannian index theorems and better understanding of how the Wick rotated operators change, one could potentially remove the asymptotic ultra-static assumption of this theorem, since the Wick rotation is defined for any spacetime. 
\vspace{0.1in}
\newline
\textbf{Proof:} We split the proof into 3 steps. First we demonstrate that this satisfies the assumptions of proposition \ref{prop:UniformFredholm} to show that this is a uniform family of Fredholm operator. Secondly we show that this implies the index to be constant, and finally we relate the Riemannian index to the Lorentzian index. We will do the Feynman case, while the anti-Feynman case proceeds in the exact same fashion but doing the rotation with negative $\theta$.
\vspace{0.1in}
\newline
1) Because $N(D_{\theta})$ is $cl(n)(-ie^{-i\theta}\partial_t - A)$, the associated wave operator $P$ satisfies $N(D_{\theta}^2) = (-e^{-2i\theta}\partial_t^2 - A^2)$ hence by what we just discussed will satisfy uniform high-energy estimates as long as
\begin{equation}
    dist(\Gamma, \cup_{0 \leq \theta \leq \pi/2} e^{i\theta} spec(A)) > b > 0 \quad \text{uniformly over $\theta$}
\end{equation}
By assumptions, this holds. For our concrete examples of contours $\phi_i$ in terms of $\tanh$ imaginary parts, see figure \ref{fig:wick_rotation}. Thus, we can apply proposition \ref{prop:UniformFredholm} to obtain that the $D_{\theta}$ forms a uniform family of Fredholm operators. The index concerns the Fredholm index of its positive part $D^{+}_{\theta}$. These are also uniformly Fredholm on the same spaces, as they are given by $ D^{+}_{\theta} = \Pi_{-} D_{\theta} \Pi_{+}$ where $\Pi_{\pm}$ are the smooth projections onto $V^{\pm}$.
\vspace{0.1in}
\newline
2) Next we show that this implies the index to be constant. An argument of this form was done previously in \cite{Gell_Redman_2016} theorem 4.3, hence we simply recall it here. We begin with assuming that the $D^{+}$ is invertible. In that case, suppose $D_{\theta}^{+}$ was not invertible for small $\theta$, and pick a sequence $u_{\theta_j}$ of elements of the kernels. Then by the uniformity of our Fredholm estimate in $\theta$, $\| u_{\theta_j} \| \lesssim \| K u_{\theta_j} \|$ where $K$ is a compact map. After rescaling, we can pick a weakly-converging subsequence $u_{\theta_{k}}$, which by the estimate above is also strongly converging to some nonzero $u$. By continuity of the family, this would be an element of the kernel of $D^{+}_{\theta}$, which is a contradiction. Doing the same argument with the adjoint gives that the cokernel must be empty too for small $\theta$. Thus the family is invertible for small $\theta$, hence in particular has index $0$. The index is then $0$ for all $\theta$, since for $\theta > 0$ the family is elliptic hence Fredholm on the same spaces $H^s_{b, \Gamma} \rightarrow H^{s-1}_{b, \Gamma}$. To cover the general case of $D^{+}$ Fredholm, we then append finite dimensional vector spaces $W_1$ and $W_2$ so that
\begin{equation}
    \Tilde{D}_{\theta}^{+}: \{ u \in H^s_{b, \Gamma}: D_{\theta}^{+}u \in H^{s}_{b, \Gamma} \} \oplus W_1 \rightarrow H^{s}_{b, \Gamma} \oplus W_2
\end{equation}
is invertible at $\theta = 0$. By what we just showed the index of $\Tilde{D}^{+}_{\theta}$ is then $0$ for all $\theta$, hence the index of $D^{+}_{\theta}$ is also constant, equal to $dim(W_2) - dim(W_1) = ind(D^{+})$ as desired.
\vspace{0.1in}
\newline
3) We thus obtain the index of our Lorentzian operator to be the index of a Riemannian Dirac operator $D^{+}_{\pi/2}$ on 
\begin{equation}
    D^{+}_{\pi/2}: H^s_{b, \Gamma} \rightarrow H^{s-1}_{b, \Gamma} \quad N(D_{\pi/2}) = cl(n)(-\partial_t - A_i)
\end{equation}
The contour $\Gamma_0$ and $\Gamma_1$ can be deformed to the horizontal lines $Im(\sigma) = \mp \lambda$ without crossing any resonances of the indicial family, hence we expect the indexes to agree with the index of 
\begin{equation}
    D^{+}_{\pi/2}: \tau^{\lambda} H^s_b \rightarrow \tau^{\lambda} H^{s-1}_b
\end{equation}
which by Melrose's APS index theorem 9.1 \cite{Melrose1993} is given by 
\begin{equation}
        ind_{\lambda}(D_{\pi/2}) = {}^b\int_{\Bar{M}} AS(\nabla^{\pi/2}) - \frac{h(A^{+}) + \eta(A^{+})}{2}
\end{equation}
where $A^{+}$ denotes the induced Dirac operator on the full boundary $\{ 0 \} \times \Sigma$ and $\{ 1 \} \times \Sigma$. In particular, the dimension of its kernel is $h(A^{+}_0) + h(A^{+}_1)$ while its Eta invariant $\eta(A^{+})$ is $\eta(A^{+}_{0}) - \eta(A^{+}_{1})$, where the future contribution comes with a minus sign due to our convention for the induced Dirac operator being minus the one of Melrose there (who uses $-i cl(d\tau/\tau)(\tau \partial_{\tau} + A)$ everywhere hence the opposite at the future). We then write the integrant in terms of the original connection $\nabla^{\pi/2}$ using the transgression form $AS(\nabla^{\pi/2}) = AS(\nabla) + dTAS(\nabla, \nabla^{\pi/2})$, which when integrating gives a boundary term
\begin{equation}
    \int_{\partial \Bar{M}} TAS(\nabla, \nabla^{\pi/2})
\end{equation}
This vanishes in our ultra-static case because at the boundary, $\nabla^{\theta}|_{\partial \Bar{M}}$ is invariant under our Wick rotation due to $\nabla^L n$ being $0$ there, hence the transgression form vanishes there. This gives the claimed formula. Thus, we only have to show a 'relative index theorem' relating the index on $H^s_{b, \Gamma}$ and $\tau^{\lambda} H^s_{b}$ to conclude the proof. Because there are no resonances between $\Gamma_0$, $\Gamma_1$ and the horizontal lines $Im(\sigma) = \mp \lambda$ we expect the index to agree. The following relative index theorem will show that, by applying it to $D^{+}_{\pi/2}$. This will also allow computing the index of non-elliptic operators on different contours $\Gamma'$, by counting the number of resonances one crosses between $\Gamma$ and $\Gamma'$.
\begin{lemma}(Relative Index Theorem on $H^s_{b, \Gamma})$ Let $P$ be a differential operator in $\Psi^m_{b, \Gamma} \cap \Psi^m_{b, \Tilde{\Gamma}}$ which is Fredholm on
\begin{equation}
    P: \{ u \in H^s_{b, \Gamma}: Pu \in H^{s-m+1} \} \rightarrow H^{s-m+1}_{b, \Gamma}
\end{equation}
in the sense of theorem \ref{theo:NonEllipticOperatorsAreFredholm} (or $P$ elliptic), for the two different choices of contours $\Gamma$,$\Tilde{\Gamma}$. Assume that $\Gamma_0$ is below $\Tilde{\Gamma}_0$ and that $\Gamma_1$ is above $\Tilde{\Gamma}_1$, with a lower bounded gap between them. In other words, $\Gamma$ implies more decay than $\Tilde{\Gamma}$ hence  $H^s_{b, \Gamma} \hookrightarrow H^s_{b, \Tilde{\Gamma}}$. Suppose furthermore that we have polynomial estimates along vertical lines between $\Gamma$ and $\Tilde{\Gamma}$ for the indicial family $\widehat{N(P)}(\sigma)^{-1}$, as in proposition \ref{prop:ModeExpansion}, and that there is only finitely many resonances of $\widehat{N(P)}(\sigma)^{-1}$ between $\Gamma_i$ and $\Tilde{\Gamma}_i$ for $i=0, 1$ outside of which we have high-energy estimates on (curved) horizontal strips \footnote{In other words, we can justify contour deformation between them.}. Then the difference in index between $P$ on $H^s_{b, \Gamma}$ and $P$ on $H^s_{b, \Tilde{\Gamma}}$ is given by the total number of resonances (counted with multiplicity) between the two contours:
\begin{equation}
    ind_{\Tilde{\Gamma}}(P) - ind_{\Gamma}(P) = \sum_{i = 0}^{1} \sum_{\sigma_l \text{ resonance between }\Gamma_i, \Tilde{\Gamma}_i} m(\sigma_l)
\end{equation}
In particular the indexes are the same if there is no resonances between them.
\end{lemma}
\textbf{Proof:} This is effectively the same proof as Melrose's relative index theorem \cite{Melrose1993} (theorem 6.5) but in the non-elliptic setting, hence we proceed in the exact same way. By subdividing the region between $\Tilde{\Gamma}$ and $\Gamma$ into thin (curved) horizontal strips and summing the contributions of each index change at the end, we can assume that the gap between the two contours is small enough so that $P - N(P) = O(\tau^{\nu})$ has $\nu$ bigger than this gap, which implies that $P - N(P)$ maps $H^s_{b, \Tilde{\Gamma}}$ to $H^{s-m}_{b, \Gamma}$. The key step will be then to use that for some cutoffs $\chi$, $\psi$ near the boundary $\partial \Bar{M}$ we have the statement
\begin{equation}
    u \in H^s_{b, \Tilde{\Gamma}}, Pu \in H^{s-m+1}_{b, \Gamma} \implies \chi u = \sum_{\sigma \text{ resonances between } \Gamma, \Tilde{\Gamma}} \sum_{k} \tau^{-i\sigma}(-\log(\tau))^{k} \psi(\tau) u^0_{\sigma,k} + r \quad r \in H^s_{b, \Gamma}
\end{equation}
This was shown in our mode expansion proposition \ref{prop:ModeExpansion}, except that $r$ was in some weak space $H^{-M}_{b, \Gamma}$. To improve its regularity, we write $r = \chi u - S$ where $S$ is the sum of modes above. Then for $Pu = f \in H^{s-m+1}_{b, \Gamma}$, $S = \psi(\tau) S_0$ and using that the quasinormal modes are smooth we have
\begin{equation}
    Pr = \chi f + [P, \chi]u - (P - N(P))S -[N(P), \psi] S_0 \in H^{s-m+1}_{b, \Gamma}
\end{equation}
We can then apply the strong version of our radial set estimate \ref{prop:RadialPointEstimates} plus propagation of singularity to conclude that $r \in H^s_{b, \Gamma}$. As a consequence of this statement, elements of the kernel of $P$ on $H^s_{b, \Gamma}$ and $H^s_{b, \Tilde{\Gamma}}$ agree up to the $m(\sigma)$ resonance terms, and $u \in H^s_{b, \Gamma}$ iff $S$ vanishes. Consider the sequence
\begin{equation}
    0 \rightarrow kerP|_{H^s_{b, \Gamma}} \rightarrow kerP|_{H^s_{b, \Tilde{\Gamma}}} \rightarrow \mathcal{S} \rightarrow cokerP|_{H^s_{b, \Gamma}} \rightarrow cokerP|_{H^s_{b, \Tilde{\Gamma}}} \rightarrow 0
\end{equation}
where $\mathcal{S}$ is the space of sums of resonant states $S$ between $\Gamma$ and $\Tilde{\Gamma}$ as above, and $kerP|_{H^s_{b, \Tilde{\Gamma}}}$ maps to it by mapping $u$ to the corresponding sum $S$ in $\chi u$, while the map from $\mathcal{S}$ to the cokernel is given by applying $P$. This sequence is exact (the first half by what we discussed, and the second half follows from $P(S) = 0$ in $cokerP|_{H^s_{b, \Gamma}}$ implying that $P(S) = Pv$ for some $v \in H^s_{b, \Gamma}$ hence $P(S-v) = 0$, i.e. $S$ comes from the data of some $S - v \in ker P|_{H^s_{b, \Tilde{\Gamma}}}$) hence the alternating sum of dimension vanishes, which says exactly that the indexes on $H^s_{b, \Gamma}$ and $H^s_{b, \Tilde{\Gamma}}$ differ by the dimension of $\mathcal{S}$, which is the total number of resonances between the two contours.
$\Box$

\vspace{0.2in}
\newpage
\section{Acknowledgments}

I would like to acknowledge my advisor's help (Andras Vasy) for the many discussions on the project. They were crucial in making progress and properly structuring the theory. I would also like to thank Stanford's math department and the FRQ (`Fond de Recherche du Québec'; Doctoral Research Scholarship \#332587) for funding my PhD and allowing me to pursue my research. Finally, I am grateful to the Erwin Schrödinger International Institute for Mathematics and Physics (ESI) for its hospitality in the summer of 2023 for the 'Spectral Theory and Mathematical Relativity' program, where I first heard about index theorems for Lorentzian Dirac operators in a talk by Christian Bär. 

\section{Declarations}
\textbf{Use of AI}
\vspace{0.1in}
\newline
\noindent The figures in this paper were done on Python with ChatGPT Pro Edu (using GPT-6), which was provided for me by Stanford University's math department following a donation. The mathematical work was mostly done in the winter and spring quarters of 2026, during which a personal (free, unsponsored) account running on GPT-5 was used for gathering references to various results in the literature and helping with recalling formulas. Furthermore, in the later stage of writing this paper, the AI overview feature of Google powered by Gemini was also used for a similar purpose. Finally, before submission, this paper was reviewed by ChatGPT for catching typos and inconsistencies. This paper contains no AI written content.
\vspace{0.2in}
\newline
\noindent \textbf{Funding}
\vspace{0.1in}
\newline
\noindent The research leading to these results received funding from the FRQ (' Fond de Recherche du Québec') under grant number \#332587 and from the NSF ('National Science Foundation') under grant numbers DMS-2247004 and DMS-2553664.
\vspace{0.2in}
\newline
\noindent \textbf{Conflict of Interest}
\vspace{0.1in}
\newline
\noindent The author has no relevant financial or non-financial interests to disclose.
\vspace{0.2in}
\newline
\noindent \textbf{Ethical Statement}
\vspace{0.1in}
\newline
\noindent This research did not involve human or animal participants and therefore was not subjected to compliance with ethical standards concerning them.
\vspace{0.2in}
\newline
\noindent \textbf{Data Availability}
\vspace{0.1in}
\newline
\noindent The author has no relevant data to share.

\bibliographystyle{acm}
\bibliography{main}

\label{AppendixHighEnergy}
\end{document}